\documentclass[aos]{imsart}

\RequirePackage{amsthm,amsmath,amsfonts,amssymb}
\RequirePackage[authoryear]{natbib}

\startlocaldefs

\usepackage{booktabs}
\usepackage{threeparttable}
\usepackage{amsfonts}
\usepackage{amsmath}
\usepackage{mathrsfs}
\usepackage{fancyhdr}
\usepackage{epsfig,morefloats}

\usepackage{amsthm}
\usepackage{amssymb}
\usepackage{color}
\usepackage{lscape}
\usepackage{rotating}
\usepackage{caption}
\usepackage{subfigure}
\usepackage{graphicx}
\usepackage[compact]{titlesec}
\usepackage{lineno}
\usepackage{cancel}
\usepackage{multirow}

\newcommand{\p}{{\rm p}}

\definecolor{mygray}{gray}{0.6}

\newcommand{\Date}[1]{\def\@Date{#1}}
\def\today{\number\day~\ifcase\month\or
	January\or February\or March\or April\or May\or June\or
	July\or August\or September\or October\or November\or December\fi~\number\year}

\def\be{\begin{equation}}
	\def\ee{\end{equation}}
\def\bea{\begin{eqnarray}}
	\def\eea{\end{eqnarray}}
\def\bd{\begin{displaymath}}
	\def\ed{\end{displaymath}}
\def\bda{\begin{eqnarray*}}
	\def\eda{\end{eqnarray*}}
\def\bsm{\begin{small}}
	\def\esm{\end{small}}

\def\t0{\theta_0}

\def\ha1{\hat \beta_1}

\def\bnt{\begin{enumerate}}
	\def\ent{\end{enumerate}}
\def\T{{\scriptscriptstyle \top}}
\def\p{{ \mathrm{p} }}

\def\bsc{\begin{scriptsize}}
	\def\esc{\end{scriptsize}}

\newtheorem{theorem}{Theorem}

\newtheorem{lemma}{Lemma}

\newtheorem{proposition}{Proposition}

\theoremstyle{definition}

\newtheorem{remark}{Remark}

\makeatletter
\newcommand{\figcaption}{\def\@captype{figure}\caption}
\newcommand{\tabcaption}{\def\@captype{table}\caption}
\makeatother

\newcommand{\calI}{{\mathcal I}}

\def\ga{\gamma}
\def\de{\delta}

\newcommand{\ve}{{\varepsilon}}

\newcommand{\bB}{{\mathbf B}}

\newcommand{\bI}{{\mathbf I}}

\newcommand{\bU}{{\mathbf U}}
\newcommand{\bV}{{\mathbf V}}

\newcommand{\bX}{{\mathbf X}}
\newcommand{\bY}{{\mathbf Y}}
\newcommand{\bZ}{{\mathbf Z}}
\newcommand{\ba}{{\mathbf a}}

\newcommand{\bc}{{\mathbf c}}
\newcommand{\bfe}{{\mathbf e}}

\newcommand{\bh}{{\mathbf h}}
\newcommand{\bg}{{\mathbf g}}

\newcommand{\br}{{\mathbf r}}

\newcommand{\bu}{{\mathbf u}}
\newcommand{\bv}{{\mathbf v}}
\newcommand{\bw}{{\mathbf w}}

\newcommand{\by}{{\mathbf y}}

\newcommand{\bdelta} {\boldsymbol{\delta}}

\newcommand{\bTheta} {\boldsymbol{\Theta}}

\newcommand{\btheta} {\boldsymbol{\theta}}
\newcommand{\bxi} {\boldsymbol{\xi}}

\newcommand{\bzeta} {\boldsymbol{\zeta}}

\newcommand{\bzero}{{\mathbf 0}}

\renewcommand{\theequation}{\thesection.\arabic{equation}}

\allowdisplaybreaks

\numberwithin{equation}{section}
\usepackage{amsmath, amssymb, amsthm, latexsym,color, natbib}
\usepackage[mathscr]{euscript}
\usepackage{graphicx}
\usepackage{fancyhdr}
\usepackage{algpseudocode}
\usepackage{algorithm}
\usepackage{algpseudocode}
\usepackage{enumitem}

\theoremstyle{definition}
\newtheorem{condition}{Condition}
\newtheorem{definition}{Definition}
\def\ga{\gamma}
\def\de{\delta}

\newcommand{\calX}{{\mathcal X}}

\def\6bullets{\bullet\bullet\bullet\bullet\bullet\bullet}

\DeclareMathAlphabet\EuScriptBF{U}{eus}{b}{n}

\endlocaldefs

\begin{document}

\begin{frontmatter}
\title{Edge Differentially Private Estimation in the $\beta$-model via Jittering and Method of Moments}
\runtitle{Edge Differentially Private Estimation}

\begin{aug}
\author[A,B]{\fnms{Jinyuan} \snm{Chang}\ead[label=e1,mark]{changjinyuan@swufe.edu.cn}},
\author[A]{\fnms{Qiao} \snm{Hu}\ead[label=e2,mark]{huqiao@smail.swufe.edu.cn}},
\author[C]{\fnms{Eric D.} \snm{Kolaczyk}\ead[label=e3,mark]{eric.kolaczyk@mcgill.ca}},
\author[D]{\fnms{Qiwei} \snm{Yao}\ead[label=e4,mark]{q.yao@lse.ac.uk}}
\and
\author[E]{\fnms{Fengting} \snm{Yi}\ead[label=e5,mark]{yifengting@ynu.edu.cn}}
\address[A]{Joint Laboratory
				of Data Science and Business Intelligence, Southwestern University of Finance and Economics, Chengdu, 611130, China, \printead{e1,e2}}

    \address[B]{Academy of
Mathematics and Systems Science, Chinese Academy of Sciences, Beijing, 100190, China}

\address[C]{Department of Mathematics and Statistics, McGill University, Montreal, QC, H3A 0B8, Canada, \printead{e3}}

\address[D]{Department of Statistics, London School of Economics and Political Science, London, WC2A 2AE, UK, \printead{e4}}

\address[E]{Yunnan Key Laboratory of Statistical Modeling and Data Analysis, Yunnan University, Kunming, 650500, China, \printead{e5}}
\end{aug}

\begin{abstract}
A standing challenge in data privacy is the trade-off between the level of privacy and the efficiency of statistical inference. Here we conduct an in-depth study of this trade-off for parameter estimation in the $\beta$-model \citep{ChatterjeeDiaconisSly_2011} for edge differentially private network data released via jittering \citep{KarwaKrivitskySlavkovic_2017}.  Unlike most previous approaches based on maximum likelihood estimation for this network model, we proceed via method-of-moments.  This choice facilitates our exploration of a substantially broader range of privacy levels -- corresponding to stricter privacy -- than has been to date.  Over this new range we discover our proposed estimator for the parameters exhibits an interesting phase transition, with both its convergence rate and asymptotic variance following one of three different regimes of behavior depending on the level of privacy.  Because identification of the operable regime is difficult if not impossible in practice, we devise a novel adaptive bootstrap procedure to construct uniform inference across different phases.  In fact, leveraging this bootstrap we are able to provide for simultaneous inference of all parameters in the $\beta$-model (i.e., equal to the number of nodes), which, to our best knowledge, is the first result of its kind.  Numerical experiments confirm the competitive and reliable finite sample performance of the proposed inference methods, next to a comparable maximum likelihood method, as well as significant advantages in terms of computational speed and memory.
\end{abstract}



\end{frontmatter}

\section{Introduction.}
	
	In this information age, data is one of the most important
	assets. With ever-advancing machine learning technology,
	collecting, sharing and using data yield great societal and economic
	benefits, while the abundance and granularity of personal data bring new
	risks of potential exposure of
	sensitive personal or financial information which
	may lead to adverse consequences. Therefore, continuous and
	conscientious effort has been made to formulate concepts of sensitivity of the
	data and privacy guarantee in data usage, and those concepts evolve along with
	the technological advancement. At present, one of most commonly used
	formulations of data privacy is the so-called differential privacy \citep{Dwork_2006,
		WassermanZhou_2010}. This paper is devoted to studying statistical estimation in the context of edge differential privacy for network data.

	In network data, individuals (e.g., persons or firms) are
	typically represented by nodes and their
	inter-relationships are represented by edges. Therefore, network data often contain
	sensitive individual information. On the other hand, for analysis purposes the
	information of interest in the data should be sufficiently preserved. Hence, the
	primary concern for data privacy 
	is two-folded: (a) to release
	only a sanitized version of the original network data to protect privacy,
	and (b) the sanitized data should preserve the information of interest
	such that analysis based on the sanitized data is still effective.		
	
	To protect privacy, the conventional approach is to release some noised
	version of summary statistics of interest. Normally the summary
	statistics used are of (much) lower dimension
	than the original data. In the context of network data, the chosen
	summary statistics can be the node degree sequence
	\citep{KarwaSlavkovic_2016} or subgraph counts \citep{Blockietal_2013}.
	To achieve differential privacy, only a noised version of the summary
	statistics is released. The noised version of the statistics is generated based on
	some appropriate release mechanism, which depends critically on the
	so-called sensitivity of the adopted statistics.
	One of the most frequently used data release schemes is the Laplace mechanism
	of \cite{Dworketal_2006}. See also
	Section 2 of \cite{WassermanZhou_2010}, and
	Section 3 of \cite{KarwaSlavkovic_2016}.
	\cite{KarwaSlavkovic_2016} consider edge differential privacy for the
	$\beta$-model \citep{ChatterjeeDiaconisSly_2011}, where only the
	node degree sequence, which is a sufficient statistic,  is released with
	added noise generated from a discrete Laplace mechanism. However, a noisy
	degree sequence may no longer be a legitimate degree sequence.
	Even for a legitimate degree sequence, the maximum likelihood estimator (MLE)
	may not exist.
	\cite{KarwaSlavkovic_2016} propose a two-step procedure that entails `de-noising' the
	noisy sequence first and then estimating the parameters using the
	de-noised data by MLE.
	
	A radically different approach is to release a noisy version of an entire network.	\cite{KarwaKrivitskySlavkovic_2017} offer what they call a generalized random response mechanism for doing so and present empirical results of its use with maximum likelihood estimation in exponential random graph models. The structure of this release mechanism is same as the noisy network setting of \cite{ChangKolaczykYao_2018}, where the edge status of each pair of nodes is known only up to some binary noise and method-of-moments was used to estimate certain network summary statistics.  As noted by \cite{ChangKolaczykYao_2020}, this noisy network setting in turn is essentially analogous to the idea of jittering in the analysis of classical Euclidean data, where each original data point is released with added noise. In this paper we study this jittering release mechanism for network data, and we do so in the specific context of parameter estimation for the $\beta$-model.  However, importantly, we note that unlike approaches based on releasing noised versions of specific, pre-determined summary statistics, jittering allows for the possibility of multiple statistics to be calculated and/or quantities to be estimated from the same released network.

	Specifically, we conduct an in-depth study on the statistical inference for the $\beta$-model based on the edge $\pi$-differentially private data generated via jittering, where $\pi>0$ reflects the privacy level; the smaller $\pi$, the greater the level of privacy. Unlike most previous approaches to inference under this model, based on maximum likelihood estimation, we proceed via method of moments. This choice facilitates our exploration of a substantially broader range of privacy levels $\pi$ than has been to date. Let $p$ be the number of nodes in the network. Our major contributions are as follows.
 \begin{itemize}
\item First, we develop the asymptotic theory when $p\to \infty$ and $\pi \to 0$, and find that (i) in order to achieve consistency of the newly proposed moment-based estimator, $\pi$ should decay to zero slower than $p^{-1/3}\log^{1/6}p$, while (ii) both the convergence rate and the asymptotic variance of our proposed estimator depend intimately on the interplay between $p$ and $\pi$. In particular, the asymptotic behavior of these quantities exhibits an interesting phase transition phenomenon, as $\pi$ decays to zero as a function of $p$, following one of three different regimes of behavior: $\pi\gg p^{-1/4}$, $\pi\asymp p^{-1/4}$, and $p^{-1/4}\gg \pi\gg p^{-1/3}\log^{1/6}p$.

\item Second, because identification of the operable regime is difficult if not impossible in practice, we devise a novel adaptive bootstrap procedure to construct uniform inference across different phases.

\item Third, leveraging this bootstrap we are able to provide for simultaneous inference of all parameters in the $\beta$-model (i.e., equal to the number of nodes). This, to our best knowledge, is  the first result of its kind, which requires a substantially different and more nuanced technical investigation than those for finite-dimensional results.

\item Lastly, numerical experiments confirm the competitive and reliable finite sample performance of the proposed inference methods, next to a comparable maximum likelihood method, as well as significant advantages in terms of computational speed and memory.
 \end{itemize}
	
	The dichotomy of `dense' versus `sparse' networks is an important one in network science, as sparsity of edges is a property encountered widely in practice with real-world networks.  In recent years, theoretical properties of sparse $\beta$-models have been successfully considered, extending the original developments for dense $\beta$-models (such as cited above).  See, for example, \cite{Mukherjeeetal_2018}, \cite{ChenKatoLeng_2021}, \cite{SteinLeng_2021}, and~\cite{Zhangetal_2021}, which in turn build on earlier work of~\cite{Rinaldoetal_2013}.  \cite{FanZhangYan_2020} have addressed estimation in an edge-weighted version of the sparse $\beta$-model (as well as in the dense case) under the differential privacy mechanism of~\cite{KarwaSlavkovic_2016}.
	Here in this paper we conduct the majority of our development in the dense case, after which we then extend our results to the sparse case.
	
	The rest of the paper is organized as follows. Section \ref{se2}
	introduces the concept of edge $\pi$-differential privacy for networks, and the data release mechanism by jittering \citep{KarwaKrivitskySlavkovic_2017}. Section \ref{sec4} addresses
	inference for the $\beta$-model based on edge differentially
	private data, introducing the method-of-moments estimator and characterizing its asymptotic behavior. Section \ref{sec44} develops the adaptive bootstrap inference that makes inference feasible in practice (i.e., despite the phase transition), and presents the accompanying results on simultaneous inference. Some numerical results are reported in Section \ref{sec6}. Section \ref{sec:sparsemodel} illustrates
 how to extend the proposed moment-based method and the associated theory to  sparse $\beta$-models. We
	relegate all the technical proofs to the supplementary material. 

	{\it Notation}. For any integer $d\geq1$, we write $[d]=\{1,\ldots,d\}$, and denote by $\bI_d$ the $d\times d$ identity matrix. We denote by $I(\cdot)$ the indicator function. For a
	vector $\bh=(h_1,\ldots,h_d)^\T$, we write $|\bh|_0=\sum_{j=1}^dI(h_j\neq0)$
	and $|\bh|_\infty=\max_{j\in[d]}|h_j|$ for its $L_0$-norm and
	$L_\infty$-norm, respectively. For a countable set $\mathcal{S}$, we use
	$\#\mathcal{S}$ or $|\mathcal{S}|$ to denote its cardinality. For two
	sequences of positive numbers $\{a_p\}_{p\geq1}$ and $\{c_p\}_{p\geq1}$,
	we write $a_p\lesssim c_p$ or $c_p\gtrsim a_p$ if
	$\lim\sup_{p\rightarrow\infty}a_p/c_p<\infty$, and write $a_p\asymp c_p$
	if and only if $a_p\lesssim c_p$ and $c_p\lesssim a_p$ hold simultaneously.
	We also write $a_p\ll c_p$ or $c_p\gg a_p$ if
	$\lim\sup_{p\rightarrow\infty}a_p/c_p=0$.

	\section{Edge differential privacy.}\label{se2}
	
	\subsection{Definition.} \label{sec21}
	We consider simple networks in the sense that there are
	no self-loops and there exists at most one edge from one node to another
	for a directed network, and at most one edge between two nodes for an undirected network.
	Such a network with $p$ nodes can be represented by an adjacency matrix
	$\bX = (X_{i,j})_{p\times p}$, where $X_{i,i}\equiv 0$, and $X_{i,j}=1$ indicating
	an edge from the $i$-th node to the $j$-th node, and 0 otherwise.
	For undirected networks,  $X_{i,j}=X_{j,i}$. In this paper,
	we always assume that the $p$ nodes are fixed and are labeled as $1, \ldots, p$.
	Then a simple network can be represented entirely by its adjacency matrix.
	To simplify statements, we often refer to an adjacency matrix $\bX$ as a network.
	
	Let $\calX$ be the set consisting of the adjacency matrices of all the simple
	and directed (or undirected) networks with $p$ nodes. For any $\bX= (X_{i,j})_{p\times p}\in \calX$ and $\bY= (Y_{i,j})_{p\times p} \in \calX$,
	the Hamming distance between $\bX$ and $\bY$ is defined as
	\begin{equation} \label{a0}
		\de(\bX, \bY) = \# \{ (i,j) \in \calI : X_{i,j} \ne Y_{i,j} \}\,,
	\end{equation}
	where $\calI=\{ (i,j): 1 \le i \ne j \le p\}$ for directed networks, and $ \calI=
	\{ (i,j): 1 \le i < j \le p\}$ for undirected networks.
	To protect privacy, the original network $\bX$ is not released directly.
	Instead we release a sanitized version $\bZ=(Z_{i,j})_{p\times p} \in \calX$ of the network, where $\bZ$ is generated according to some conditional distribution $Q(\cdot\, |\, \bX)$. Here $Q$ is also called a release mechanism \citep{WassermanZhou_2010}.
	\begin{definition}
		[Edge differential privacy] For any $\pi >0$, a release mechanism
		(i.e. a conditional probability distribution) $Q$ satisfies $\pi$-edge differential
		privacy if
		\begin{equation} \label{a4}
			\sup_{\bX, \bY \in \calX:\,\delta(\bX, \bY)=1} \,\sup_{\bZ \in \calX:\,Q(\bZ\,|\,\bX)>0}\frac{ Q(\bZ\,|\,\bY)}{Q(\bZ\,|\,\bX) }  \le  e^{\pi}\,.
		\end{equation}
	\end{definition}
	\noindent The definition above equates privacy with the inability to distinguish	two close networks. The privacy parameter $\pi$ controls the amount of randomness added to released data; the smaller $\pi$ is the more protection on privacy. Notice that (\ref{a4}) is much more stringent than requiring $|Q(\bZ\,|\,\bY) - Q(\bZ\,|\,\bX)|$ to be small for any $\bX, \bY\in \calX$ with $\delta(\bX, \bY)=1$. In practice $\pi$ is often chosen to
	be small. Then it follows from (\ref{a4}) that
	\[
	\sup_{\bX, \bY \in \calX:\, \delta(\bX, \bY)=1}\,
	\sup_{\bZ \in \calX:\,
		Q(\bZ\,|\,\bX)>0}
	\frac{|Q(\bZ\,|\,\bY) - Q(\bZ\,|\,\bX)|}{ Q(\bZ\,|\,\bX)} \leq e^\pi -1 \approx \pi\,.
	\]
	Note that multiple notions of privacy have been introduced for networks; see \cite{JiangEtal_2020} for a recent survey.  In this paper we focus on the notion of edge differential privacy (e.g., \cite{NissimRaskhodnikovaSmith_2007}).
	At the same time, there is a connection between differential privacy and hypothesis testing.
	\begin{proposition}
		[Wasserman and Zhou, 2010]\label{pn:1} Let the released network
		$\bZ \sim Q( \cdot\,|\,\bX)$ and
		$Q$ satisfy $\pi$-edge differential privacy for some $\pi > 0$. For
		any given $i \ne j$, consider
		hypotheses $H_0: X_{i,j} = 1$ versus $H_1: X_{i,j} =0$.
		Then the power of any test at the significance level
		$\ga$ and based on $\bZ$, $Q$ and the distribution of $\bX$ is
		bounded from above by $\gamma e^\pi$, provided that $X_{i,j}$ is
		independent of $\{ X_{k,\ell}: ( k,  \ell) \in \calI  \mbox{ and } (k,
		\ell) \ne (i,j) \}$.
	\end{proposition}

 Proposition \ref{pn:1} implies that if $\bZ$ is released through $Q$ which
	satisfies $\pi$-edge differential privacy and $\pi$ is sufficiently
	small, it is virtually impossible to identify whether an edge exists (i.e. $X_{i,j} =1$) or not (i.e. $X_{i,j} =0$) in the original network	through statistical tests, as the power of any test is bounded by its significance level multiplied by $e^\pi$.
	The independence condition in Proposition 1 is satisfied by
	the Erd\"os-R\'enyi class of models for which all edges are independent, including the $\beta$-model and the well-known
	stochastic block model. Proposition \ref{pn:1} follows almost immediately from the Neyman-Pearson lemma for the optimality of likelihood ratio tests for simple null and simple alternative hypotheses.	It was first proved by \cite{WassermanZhou_2010} with independent observations. Since their proof can be adapted to our setting in a straightforward manner, we omit the details.
	
	For further discussion on
	differential privacy under more general settings,
	we refer to \cite{Dworketal_2006} and \cite{WassermanZhou_2010}.

	\subsection{Edge privacy via jittering.} \label{sec22}
	
	Now we introduce the data release mechanism of \cite{KarwaKrivitskySlavkovic_2017}, which is formally the same as the noisy network structure adopted in \cite{ChangKolaczykYao_2018}. This approach releases a jittered version of the entire network. The word `jittering' means that a small amount of noise is added to every single data point \citep{Hennig_2007}.

	For $\mathcal{I}$ specified just after \eqref{a0} above, we define a data release mechanism as follows:
	\begin{equation} \label{a1}
		Z_{i,j} = X_{i,j} I(\ve_{i,j} =0) + I(\ve_{i,j}=1)
	\end{equation}
	for each $(i,j) \in \calI$. In the above expression, $\{\ve_{i,j}\}_{(i,j) \in \calI}$ are independent random variables only taking three possible values $-1, 0$ and $1$ with
	\begin{equation} \label{a2}
		\mathbb{P}(\ve_{i,j}=1) = \alpha\,, ~
		\mathbb{P}(\ve_{i,j}=0) = 1 - \alpha - \beta~~ \textrm{and}~~ \mathbb{P}(\ve_{i,j}=-1) = \beta\,,
	\end{equation}
 where $\alpha,\beta\in[0,0.5]$.
	For an undirected network, $Z_{i,j}=Z_{j,i}$ for $j>i$.
	Then it follows from (\ref{a1}) and (\ref{a2}) that
	\begin{equation} \label{a3}
		\mathbb{P}(Z_{i,j}=1\,|\,X_{i,j}=0) = \alpha~~ \textrm{and} ~~
		\mathbb{P}(Z_{i,j}=0\,|\,X_{i,j}=1) = \beta\,.
	\end{equation}
	Furthermore the proposition below follows from (\ref{a4}) and (\ref{a3}) immediately.
	See also Proposition 1 of \cite{KarwaKrivitskySlavkovic_2017}.
	
	\begin{proposition}\label{pn:2} The data release mechanism {\rm(\ref{a1})} satisfies $\pi$-edge
		differential privacy with
		\[
		\pi = \log \bigg\{ \max\bigg(
		\frac{\alpha}{1- \beta}\,, \;
		\frac{\beta}{1-\alpha}\,, \;
		\frac{1-\alpha}{\beta}\,, \;
		\frac{1- \beta}{\alpha} \bigg) \bigg\}\,.
		\]
	\end{proposition}
	
	\begin{remark}\label{rk1}
		Notice that
  \[
  \begin{split}
&\frac{\alpha}{1-\beta}=1-\frac{1-\alpha-\beta}{1-\beta}\,,~~~~~\frac{\beta}{1-\alpha}=1-\frac{1-\alpha-\beta}{1-\alpha}\,,\\
&\frac{1-\alpha}{\beta}=1+\frac{1-\alpha-\beta}{\beta}\,,~~~~~\frac{1-\beta}{\alpha}=1+\frac{1-\alpha-\beta}{\alpha}\,,
  \end{split}
  \]
  where $1-\alpha-\beta\geq0$.
  Then the differential privacy parameter $\pi$ given in Proposition \ref{pn:2} can be reformulated as
\begin{align*}
\pi&=\log\bigg\{1+(1-\alpha-\beta)\max\bigg(-\frac{1}{1-\beta}\,,\;-\frac{1}{1-\alpha}\,,\;\frac{1}{\beta}\,,\;\frac{1}{\alpha}\bigg)\bigg\}\\
&=\log\bigg\{1+(1-\alpha-\beta)\max\bigg(\frac{1}{\beta}\,,\;\frac{1}{\alpha}\bigg)\bigg\}=\log\bigg\{1+\frac{1-\alpha-\beta}{\min(\alpha,\beta)}\bigg\}\,.
\end{align*}
Recall $\alpha,\beta\in[0, 0.5]$. The maximum privacy
		is achieved by setting
		$\alpha = \beta=0.5$, as then
		$\pi =0$. By (\ref{a1})
		and (\ref{a2}), $Z_{i,j} = I(\ve_{i,j}=1)$ then, i.e. $\bZ$ carries no information
		about $\bX$. In order to achieve high privacy, we need to use
		large $\alpha$ and $\beta$.
Due to $\alpha,\beta\in[0, 0.5]$, when $\pi\rightarrow0$, $\min(\alpha,\beta)$ cannot converge to zero, which means  there exists a constant $\epsilon\in(0, 0.5)$ such that $\min(\alpha,\beta)>\epsilon$ when $\pi\rightarrow0$. Hence, when $\pi\rightarrow0$, we have
		$\pi\asymp1-\alpha-\beta$.
		%
		In Section \ref{sec4} below we will develop statistical inference approaches for the original network $\bX$ based on the released data $\bZ$
		with $\pi\rightarrow0$.
		
	\end{remark}
	
	\section{Differentially private inference for the $\beta$-model.}\label{sec4}
	
	In this section we introduce a new method-of-moments estimator for the parameters of the network $\beta$-model and characterize the asymptotic behavior of this estimator, through which we discover an interesting phase transition.
	
	\subsection{The $\beta$-model.}\label{secbeta}
	
	The so-called $\beta$-model \citep{ChatterjeeDiaconisSly_2011} for undirected networks is characterized by
	$p$ parameters
	$\btheta = (\theta_1, \ldots, \theta_p)^\T \in \mathbb{R}^p$
	which define the probability function
	\begin{equation} \label{d2}
		\mathbb{P}(X_{i,j}=1) = \frac{\exp ( \theta_i + \theta_j)}{ 1 + \exp (
			\theta_i + \theta_j)}\,, \quad i\ne j\,.
	\end{equation}
	The parameter $\theta_i$ in this model has a natural interpretation as it measures the propensity of node $i$ to have connections with other nodes. Namely, the larger $\theta_i$ is, the more likely node $i$ is connected to other nodes. The likelihood function for $\beta$-model is given by
	\[
	f(\bX;\btheta)=\prod_{i,j:\,i< j}  {\exp\{( \theta_i + \theta_j)X_{i,j} \}\over  1
		+\exp(\theta_i + \theta_j)} \propto \exp ( U_1 \theta_1 + \cdots + U_p \theta_p)\,,
	\]
	where $U_i = \sum_{j:\,j\ne i}X_{i,j}$ is the degree of the $i$-th node.
	Hence the degree sequence $\bU= (U_1, \ldots, U_p)^\T$ is a sufficient
	statistic.
	
	Denote by $\tilde{\btheta}(\bU)=\{\tilde{\theta}_1(\bU),\ldots,\tilde{\theta}_p(\bU)\}^\T$ the MLE for $\btheta$ based on $\bU$. For given degree sequence $\bU$, $\tilde{\btheta}(\bU)$ must satisfy the following moment equations:
	\begin{align*}
		U_i=\sum_{j:\,j\neq i}\frac{\exp\{\tilde{\theta}_i(\bU)+\tilde{\theta}_j(\bU)\}}{1+\exp\{\tilde{\theta}_i(\bU)+\tilde{\theta}_j(\bU)\}}\,,\quad i\in[p]\,.
	\end{align*}
	Unfortunately $\tilde{\btheta}(\bU)$ may not exist; see Theorem
	1 of \cite{KarwaSlavkovic_2016} for necessary and sufficient conditions for the existence of $\tilde{\btheta}(\bU)$. When $\tilde{\btheta}(\bU)$ exists,
	\cite{ChatterjeeDiaconisSly_2011} show that
	\begin{equation} \label{d1}
		| \tilde{\btheta}(\bU) - \btheta|_\infty \le C_*\sqrt{\frac{\log p}{p}}
	\end{equation}
	with probability at least $1- C_*p^{-2}$,
	where $C_* >0$ is a constant depending only on
	$|\btheta|_\infty$. For any fixed integer $s\geq1$ and distinct
	$\ell_1,\ldots,\ell_s\in[p]$, \cite{yx13} establish the asymptotic
	normality of
	$\{\tilde{\theta}_{\ell_1}(\bU),\ldots,\tilde{\theta}_{\ell_s}(\bU)\}^\T$
	as $p\rightarrow\infty$, which can be used to construct joint
	confidence regions for $(\theta_{\ell_1},\ldots,\theta_{\ell_s})^\T$.
	However, to our best knowledge, simultaneous inference for all $p$ parameters in the $\beta$-model remains unresolved in the literature. 	
	
	\cite{KarwaSlavkovic_2016} consider differentially private MLE for $\btheta$ based on a noisy version of the degree sequence. More specifically,
	the noisy degree sequence in their setting is defined as
	$\bU+\bV$, where the components of $\bV=(V_1,\ldots,V_p)^\T$ are
	drawn independently from a discrete Laplace distribution with the
	probability mass function $$\mathbb{P}(V=v)=\frac{(1-\kappa)\kappa^{|v|}}{1+\kappa}$$ for any
	integer $v$ with $\kappa=\exp(-\pi/2)$.
	\cite{KarwaSlavkovic_2016} propose a two-step procedure:
	(a) find the MLE $\bU^*$ for $\bU$ based on $\bU+\bV$,
	and (b) estimate $\btheta$ by $\tilde{\btheta}(\bU^*)$.
	For any fixed integer $s\geq1$ and distinct
	$\ell_1,\ldots,\ell_s\in[p]$, Theorem 4 of
	\cite{KarwaSlavkovic_2016} shows that
	$\{\tilde{\theta}_{\ell_1}(\bU^*),\ldots,\tilde{\theta}_{\ell_s}(\bU^*)\}^\T$
	shares the same asymptotic normality as
	$\{\tilde{\theta}_{\ell_1}(\bU),\ldots,\tilde{\theta}_{\ell_s}(\bU)\}^\T$ when $\pi\asymp (\log p)^{-1/2}$. To appreciate this `free privacy' result, let us assume first that
	$|\btheta|_\infty\leq C$ for some universal constant $C>0$. Then
	there exists a universal constant $\tilde{C}>1$ such
	that  $\tilde{C}^{-1}p\leq \min_{i\in[p]}U_i\leq\max_{i\in[p]}U_i\leq
	\tilde{C}p$ holds almost surely as $p\rightarrow\infty$. On the other
	hand, when $\pi\asymp(\log p)^{-1/2}$, Lemma C in the supplementary material of \cite{KarwaSlavkovic_2016}
	indicates that $|\bU^*-\bU|_\infty\leq \sqrt{6}p^{1/2}\log^{1/2}p$ holds
	almost surely as $p\rightarrow\infty$, which implies that $\bU^*$ is
	dominated by $\bU$. Based on this result, Theorem 3 of \cite{KarwaSlavkovic_2016} shows that $\tilde{\btheta}(\bU^*)$ exists and is unique and can be used to estimate $\btheta$ with uniform accuracy in all coordinates when $\pi\asymp(\log p)^{-1/2}$. 
However, when $\pi\ll (\log p)^{-1/2}$, the asymptotic behavior of $\tilde{\btheta}(\bU^*)$
	is unknown.
	
	Our interest in this paper is on differentially private estimation based on released data  $\bZ=(Z_{i,j})_{p\times p}$ generated by the more general jittering mechanism \eqref{a1}.  Remark \ref{rk1} in Section \ref{sec22} shows that $\bZ$ is $\pi$-differentially private with $\pi\asymp 1-\alpha-\beta$. To gain more appreciation of the impact of the privacy level $\pi$ on the efficiency
	of inference, we introduce a new moment-based estimation for
	$\btheta$ based on $\bZ$.  We then establish
	the asymptotic theory under the setting that $p\to \infty$ and
	$\pi$ may vary with respect to $p$.
	Of particular interest is the findings when $\pi \to 0$ together
	with $p\to \infty$. It turns out the asymptotic distribution of the new proposed
	estimator depends intimately on the interplay between
	$\pi$ and $p$,
	exhibiting interesting phase transition in the convergence rate and the asymptotic variance as $\pi$ decays to zero as a
	function of $p$.  See Theorem \ref{tm:asymnormal} and Remark \ref{rk3}(a) in Section \ref{sec43}. To overcome
	the complexity in inference due to the phase transition,
	a novel bootstrap method is proposed, which provides a uniform
	inference regardless different phases.
	In addition, it also facilitates
	the simultaneous inference for all the $p$ components of $\btheta$ as
	$p \to \infty$.

	\subsection{A new moment-based estimator.} \label{sec42}
	
	Under the $\beta$-model (\ref{d2}), it holds that
	$$
	\frac{\mathbb{P}(X_{i,j}=1)}{\mathbb{P}(X_{i,j}=0)}=\exp(\theta_i+\theta_j)$$ for
	any $i\neq j$, which implies
	\begin{equation}\label{eq:est1}
		\frac{\mathbb{P}(X_{i,\ell}=1)\mathbb{P}(X_{i,j}=0)\mathbb{P}(X_{\ell,j}=1)}{\mathbb{P}(X_{i,\ell}=0)\mathbb{P}(X_{i,j}=1)\mathbb{P}(X_{\ell,j}=0)}=\exp(2\theta_{\ell})
		\,, \quad i\neq j\neq \ell\,.
	\end{equation}
	Since only the
	sanitized network $\bZ= (Z_{i,j})_{p\times p}$, defined as in
	\eqref{a1}--\eqref{a3}, is
	available, we represent (\ref{eq:est1})
	in terms of the probabilities of $Z_{i,j}$.
	For $\tau\in\{0,1\}$, put
	\begin{equation*}
		\varphi_{\tau}(x)=(x-\alpha)^{\tau}(1-\beta-x)^{1-\tau}
	\end{equation*}
 with $x\in\{0, 1\}$.
	Then for any $i\neq j$,
	\begin{align}\label{eq:exvarphi}
	    	\mathbb{P}(X_{i,j}=0)=\frac{\mathbb{E}\{\varphi_{0}(Z_{i,j})\}}{1-\alpha-\beta}~~\textrm{and}~~
	    	 \mathbb{P}(X_{i,j}=1)=\frac{\mathbb{E}\{\varphi_{1}(Z_{i,j})\}}{1-\alpha-\beta}\,.
	\end{align}

	To simplify the notation, 
	we write $\varphi_{\tau}(Z_{i,j})$ as $\varphi_{(i,j),\tau}$ for any $i\neq j$ and $\tau\in\{0,1\}$. 
 Since $\{Z_{i,j}:i<j\}$ is a sequence of independent random variables and $Z_{i,j}=Z_{j,i}$ for any $i\neq j$,
	it follows from (\ref{eq:est1}) that
	\begin{equation}\label{eq:estimating1}
		\frac{\mathbb{E}\{\varphi_{(i,\ell),1}\varphi_{(i,j),0}\varphi_{(\ell,j),1}\}}{\mathbb{E}\{\varphi_{(i,\ell),0}\varphi_{(i,j),1}\varphi_{(\ell,j),0}\}}=\exp(2\theta_{\ell})\,, \quad i\neq j\neq \ell\,.
	\end{equation}
	For each $\ell\in[p]$, let
	\begin{align}
	\mu_{\ell,1}=&\,\frac{1}{|\mathcal{H}_{\ell}|}\sum_{(i,j)\in\mathcal{H}_{\ell}}\mathbb{E}\{\varphi_{(i,\ell),1}\varphi_{(i,j),0}\varphi_{(\ell,j),1}\}\,,\label{eq:muell1}\\ \mu_{\ell,2}=&\,\frac{1}{|\mathcal{H}_{\ell}|}\sum_{(i,j)\in\mathcal{H}_{\ell}}\mathbb{E}\{\varphi_{(i,\ell),0}\varphi_{(i,j),1}\varphi_{(\ell,j),0}\}\,,\label{eq:muell2}
	\end{align}
	where $\mathcal{H}_{\ell}=\{(i,j):i,j\neq\ell~\textrm{such that}~i< j\}$. By \eqref{eq:estimating1}, we have
	\[
	\theta_\ell=\frac{1}{2}\log\bigg(\frac{\mu_{\ell,1}}{\mu_{\ell,2}}\bigg)\,.
	\]
	Hence a moment-based estimator for $\theta_\ell$ can
	be defined as
	\begin{equation} \label{eq:esttheta}
		\hat \theta_\ell = {1\over 2} \log\bigg(\frac{\hat{\mu}_{\ell,1}}{\hat{\mu}_{\ell,2}}\bigg)\,,
	\end{equation}
	where
	\begin{align} \hat{\mu}_{\ell,1}=&\,\frac{1}{|\mathcal{H}_\ell|}\sum_{(i,j)\in\mathcal{H}_{\ell}}\varphi_{(i,\ell),1}\varphi_{(i,j),0}\varphi_{(\ell,j),1}\,,\label{eq:hatmuell1}\\ \hat{\mu}_{\ell,2}=&\,\frac{1}{|\mathcal{H}_\ell|}\sum_{(i,j)\in\mathcal{H}_{\ell}}\varphi_{(i,\ell),0}\varphi_{(i,j),1}\varphi_{(\ell,j),0}\,.\label{eq:hatmuell2}
	\end{align}
	
	\subsection{Asymptotic properties and phase transition.} \label{sec43}
	
	We always confine
	$(\alpha,\beta)
	\in\mathcal{M}(\gamma;C_1)$ with
	\begin{equation*}\label{eq:Mgamma}
			\mathcal{M}(\gamma;C_1)=\big\{(\alpha,\beta) : C_1< \alpha, \beta<0.5\,,
			1-\alpha-\beta=\gamma\big\}
	\end{equation*}
	for some $\gamma\in (0,1]$ and $C_1\in(0,0.5)$.
	Our theoretical analysis allows $\gamma$ to be a
	constant,
	or to vary with respect to $p$. Of particular interest are the cases
	when $\gamma \to 0 $ (at different rates) together with $p\to \infty$. When $(\alpha,\beta)
	\in\mathcal{M}(\gamma;C_1)$ for some fixed constants $C_1\in(0,0.5)$, it follows from Remark \ref{rk1} in Section \ref{sec22}
	that the privacy level $\pi \asymp\gamma$.
	
	\subsubsection{Consistency.} \label{sec431}
	Proposition \ref{tm:1} below presents the consistency for the moment-based estimator
	$\hat{\theta}_\ell$ defined in (\ref{eq:esttheta}), which indicates
	that $\theta_\ell$ can be estimated consistently under the edge
	$\pi$-differential privacy with $\pi \to 0$,
	as long as $\pi\gg p^{-1/3}\log^{1/6}p$.
	
	
	\begin{condition}\label{cond1}
		There exists a universal constant $C_3>0$ such that $|\btheta|_\infty\leq C_3$.
	\end{condition}
	
	
	%
	%

	\begin{proposition}\label{tm:1}
		Let Condition {\rm\ref{cond1}} hold and
		$(\alpha,\beta)\in \mathcal{M}
		(\gamma;C_1)$ for some fixed constant $C_1\in(0,0.5)$. 
		If $\gamma\gg p^{-1/3}\log^{1/6}p$, it then holds that
			\begin{align*}
	\max_{\ell\in[p]}|\hat{\theta}_\ell-\theta_\ell|=O_{\p}\bigg(\frac{\log^{1/2} p}{\gamma^3p}\bigg)+O_{\p}\bigg(\frac{\log^{1/2}p}{\gamma p^{1/2}}\bigg)\,.
	\end{align*}
	\end{proposition}
	
	\begin{remark}
	(a) By Condition \ref{cond1} and \eqref{eq:exvarphi}, we know
 $$
 \min_{\tau\in\{0, 1\}}\min_{i,j:\,i\neq j}\mathbb{E}\{\varphi_{(i,j),\tau}\}\asymp\gamma\asymp \max_{\tau\in\{0,  1\}}\max_{i,j:\,i\neq j}\mathbb{E}\{\varphi_{(i,j),\tau}\}\,,
 $$
 which implies
 $$\min_{k\in\{1, 2\}}\min_{\ell\in[p]}\mu_{\ell,k}\asymp \gamma^3\asymp \max_{k\in\{1, 2\}}\max_{\ell\in[p]}\mu_{\ell,k}\,.$$ Lemma \ref{la:3} in the supplementary material shows that $$\max_{k\in\{1, 2\}}\max_{\ell\in[p]}|\hat{\mu}_{\ell,k}-\mu_{\ell,k}|=O_{\p}\bigg(\frac{\log^{1/2}p}{p}\bigg)+O_{\p}\bigg(\frac{\gamma^2\log^{1/2}p}{p^{1/2}}\bigg)+O_{\p}\bigg(\frac{\gamma \log p}{p}\bigg)\,.$$ To make $(\hat{\mu}_{\ell,1},\hat{\mu}_{\ell,2})$ be a valid estimate of $(\mu_{\ell,1},\mu_{\ell,2})$, we need to require $p^{-1}\log^{1/2}p=o(\gamma^3)$, $\gamma^2p^{-1/2}\log^{1/2}p=o(\gamma^3)$ and $\gamma p^{-1}\log p=o(\gamma^3)$. Hence, we need the restriction $\gamma\gg p^{-1/3}\log^{1/6}p$.
	Notice that the privacy level $\pi\asymp\gamma$. In order to
	ensure the consistency of $\hat
	\theta_\ell$, the edge differential privacy level $\pi$ must
	satisfy condition $\pi \gg  p^{-1/3}\log^{1/6}p$.
	
	(b) Recall $\varepsilon_{i,j}$ involved in the data release mechanism \eqref{a1} for $Z_{i,j}$ is a discrete random variable that only takes three possible values $-1, 0$ and $1$. When $\alpha=\beta=0$, $\varepsilon_{i,j}\equiv0$, and our moment-based estimator \eqref{eq:esttheta} is then constructed
	based on the original network $\bX$. By setting $\gamma=1$ in our proof of Proposition \ref{tm:1},
	we can establish the following convergence rate for our moment-based estimator based on the original network $\bX$:
 \begin{align*}
	\max_{\ell\in[p]}|\hat{\theta}_\ell-\theta_\ell|=O_{\p}\bigg(\frac{\log^{1/2}p}{ p^{1/2}}\bigg)\,,
	\end{align*}
 which shares the same convergence rate of the MLE
	of \cite{ChatterjeeDiaconisSly_2011}; see (\ref{d1}) in Section \ref{secbeta}.
	\end{remark}

	
	\subsubsection{Asymptotic normality.} \label{sec432}
	
	Put $N=(p-1)(p-2)$. Proposition \ref{pn:asyexp} gives the asymptotic expansion of $\hat{\theta}_\ell-\theta_\ell$, which can be obtained from the proof of Theorem \ref{tm:asymnormal} in Section \ref{sec:pfthm1} of the supplementary material. For any $i\neq \ell$, let
	\begin{align}\label{eq:lambdaiell}
			\lambda_{i,\ell}=\frac{1}{p-2}\sum_{j:\,j\neq \ell,i}\bigg[\frac{1}{\mu_{\ell,1}}\mathbb{E}\{\varphi_{(\ell,j),1}\}\mathbb{E}\{\varphi_{(i,j),0}\}+\frac{1}{\mu_{\ell,2}}\mathbb{E}\{\varphi_{(\ell,j),0}\}\mathbb{E}\{\varphi_{(i,j),1}\}\bigg]\,.
	\end{align}
	\begin{proposition}\label{pn:asyexp}
		For any $i\neq j$, write $\mathring{Z}_{i,j}=Z_{i,j}-\mathbb{E}(Z_{i,j})$. Let Condition {\rm\ref{cond1}} hold and
		$(\alpha,\beta)\in \mathcal{M}
		(\gamma;C_1)$ for some fixed constant $C_1\in(0,0.5)$. If $\gamma\gg p^{-1/3}\log^{1/6}p$, it then holds that
		$$
		\hat{\theta}_\ell-\theta_\ell=\tilde{T}_{\ell,1}+\tilde{T}_{\ell,2}+\tilde{R}_\ell\,,$$
		where
			\begin{align*}
	\tilde{T}_{\ell,1}=-\frac{1}{N}\sum_{i,j:\,i\neq j,\,i,j\neq \ell}\bigg(\frac{\mu_{\ell,1}+\mu_{\ell,2}}{2\mu_{\ell,1}\mu_{\ell,2}}\bigg)\mathring{Z}_{i,\ell}\mathring{Z}_{\ell,j}\mathring{Z}_{i,j}~~and~~ \tilde{T}_{\ell,2}=\frac{1}{p-1}\sum_{i:\,i\neq \ell}\lambda_{i,\ell}\mathring{Z}_{i,\ell}
	\end{align*}
		satisfy $\tilde{T}_{\ell,1}=O_{\p}(\gamma^{-3}p^{-1})$ and $\tilde{T}_{\ell,2}=O_{\p}(\gamma^{-1}p^{-1/2})$,
		and the remainder term $\tilde{R}_\ell$ satisfies $\tilde{R}_\ell=O_{\p}(\gamma^{-6}p^{-2})+O_{\p}(\gamma^{-2}p^{-1}\log p)$.
	\end{proposition}
	
	The leading term in the asymptotic expansion of $\hat{\theta}_\ell-\theta_\ell$ will be different for different scenarios of $\gamma$: $\tilde{T}_{\ell,2}$, a partial sum of independent random variables, serves as the leading term if $\gamma\gg p^{-1/4}$, $\tilde{T}_{\ell,1}+\tilde{T}_{\ell,2}$ is the leading term if $\gamma\asymp p^{-1/4}$, and $\tilde{T}_{\ell,1}$, a generalized $U$-statistic, is the leading term if $p^{-1/4}\gg\gamma\gg p^{-1/3}\log^{1/6}p$. Such characteristic will lead to a phase transition phenomenon in the limiting distribution of the proposed moment-based estimator. Put
	\begin{align}\label{eq:bell}
		&~~~~~~~~~~~~~~~~~~~~~~~~~~~b_\ell=\frac{1}{p-1}\sum_{i:\,i\neq \ell}\lambda_{i,\ell}^2{\rm Var}(Z_{i,\ell})\,,\\
		\label{eq:bellti}
		&\tilde{b}_\ell=\frac{1}{2N}\bigg(\frac{\mu_{\ell,1}+\mu_{\ell,2}}{\mu_{\ell,1}\mu_{\ell,2}}\bigg)^2\sum_{i,j:\,i\neq j,\,i,j\neq\ell}{\rm Var}(Z_{i,\ell}){\rm Var}(Z_{\ell,j}){\rm Var}(Z_{i,j})\,.
	\end{align}

	
	\begin{theorem}\label{tm:asymnormal}
		Let  Condition {\rm\ref{cond1}} hold
		and $(\alpha,\beta)\in\mathcal{I}\in\mathcal{M}(\gamma;C_1)$ for some fixed constant $C_1\in(0,0.5)$.
		Let $1\le \ell_1 < \cdots < \ell_s \le p$ be any $s$ given indices for some fixed integer $s\ge 1$.
		As $p\rightarrow\infty$, the following three assertions hold.
		
		{\rm (a)} If $\gamma\gg p^{-1/4}$, then
		$$
		(p-1)^{1/2}{\rm diag}(b_{\ell_1}^{-1/2},\ldots,b_{\ell_s}^{-1/2})(\hat{\theta}_{\ell_1}-\theta_{\ell_1},\ldots,\hat{\theta}_{\ell_s}-\theta_{\ell_s})^\T\rightarrow\mathcal{N}(\bzero,\bI_s)$$ in distribution.
		
		{\rm (b)} If $p^{-1/4}\gg \gamma\gg p^{-1/3}\log^{1/6}p$, then $$
		N^{1/2}\,{\rm diag}(\tilde{b}_{\ell_1}^{-1/2},\ldots,\tilde{b}_{\ell_s}^{-1/2})(\hat{\theta}_{\ell_1}-\theta_{\ell_1},\ldots,\hat{\theta}_{\ell_s}-\theta_{\ell_s})^\T\rightarrow\mathcal{N}(\bzero,\bI_s)$$ in distribution.
		
		{\rm (c)} If $\gamma\asymp p^{-1/4}$, then \begin{align*}
		&N^{1/2}\,{\rm diag}[\{(p-2)b_{\ell_1}+\tilde{b}_{\ell_1}\}^{-1/2},\ldots,\{(p-2)b_{\ell_s}+\tilde{b}_{\ell_s}\}^{-1/2}]\\
&~~~~~~~~~~~~~~~~~~\times(\hat{\theta}_{\ell_1}-\theta_{\ell_1},\ldots,\hat{\theta}_{\ell_s}-\theta_{\ell_s})^\T\rightarrow\mathcal{N}(\bzero,\bI_s)
  \end{align*}in distribution.
	\end{theorem}
	
	\begin{remark} \label{rk3}
		(a) Theorem \ref{tm:asymnormal} presents
		the asymptotic normality of the proposed estimator
		when $p\to \infty$ and also, possibly, $\pi \asymp \gamma \to 0$. It can be shown that $b_\ell\asymp \gamma^{-2}$ and
		$\tilde{b}_\ell\asymp\gamma^{-6}$ under Condition {\rm\ref{cond1}}.
		The limiting distribution depends on the relative
		rates of $p$ and $\gamma$ intimately; yielding an interesting
		phase transition phenomenon in the convergence rate.
		More precisely, when $\gamma\gg
		p^{-1/4}$ (including the case  $\gamma$ is a fixed constant),
		we have
		$|\hat{\theta}_\ell-\theta_\ell|=O_{\p}(p^{-1/2}\gamma^{-1})$.
		On the other hand, $|\hat{\theta}_\ell-\theta_\ell|=O_{\p}(p^{-1/4})$
		when $\gamma\asymp p^{-1/4}$,
		and $O_{\p}(p^{-1}\gamma^{-3})$
		when $p^{-1/4}\gg
		\gamma\gg p^{-1/3}\log^{1/6}p$.
		
		
		(b) The asymptotic normality of the proposed moment-based estimator with the original network $\bX$ can
		be also established. By setting $\gamma =1$ (i.e. $\alpha=\beta=0$) in our technical proof of Theorem \ref{tm:asymnormal}(a), we can show
		%
		$
		p^{1/2}b_\ell^{-1/2}(\hat{\theta}_\ell-\theta_\ell)\rightarrow\mathcal{N}(0,1)$ in distribution.
		
		
		(c) Theorem \ref{tm:asymnormal} cannot be used to construct confidence intervals for $\theta_\ell$ directly since
		we would have to overcome two
		obstacles:
		(i) to identify the most appropriate phase in terms of relative
		sizes  between $\gamma $ and $p$, and (ii) to estimate $b_\ell$
		and $\tilde b_\ell$ which determine the asymptotic variances.
		For (ii), we give their estimates in the Appendix.
		Unfortunately (i) is extremely difficult if not impossible, as
		in practice we only have one $\gamma$ and one $p$.
		Proposition \ref{tm:varc} in the Appendix shows that (ii) is only
		partially attainable,
		as, for example, $b_\ell$ cannot be estimated consistently when
		$p^{-1/4} \lesssim \gamma \lesssim p^{-1/4} \log^{1/4}p$. In practice, we always need  $\pi\rightarrow0$ for retaining the privacy. With $\pi \rightarrow 0$, (i) can be overcome from a new perspective. More specifically, let $
	\nu_\ell=(p-2)b_\ell+\tilde{b}_\ell
	$ for any $\ell\in[p]$. Note that $b_\ell\asymp \gamma^{-2}$ and $\tilde{b}_\ell\asymp
	\gamma^{-6}$ under Condition \ref{cond1}, and $\gamma\asymp\pi$ when $\pi\rightarrow0$. Then  $(p-2)b_\ell/\nu_\ell\rightarrow1$ when $1\gg\gamma\gg
	p^{-1/4}$, and $\tilde{b}_\ell/\nu_\ell\rightarrow1$ when $\gamma\ll p^{-1/4}$. Recall $N=(p-1)(p-2)$. Hence, as $\gamma \asymp \pi \to 0$, the three asymptotic assertions in
	Theorem \ref{tm:asymnormal} admit a uniform representation:
	\begin{equation*}
		N^{1/2}\,{\rm diag}(\nu_{\ell_1}^{-1/2},\ldots,\nu_{\ell_s}^{-1/2})(\hat{\theta}_{\ell_1}-\theta_{\ell_1},\ldots,\hat{\theta}_{\ell_s}-\theta_{\ell_s})^\T\rightarrow\mathcal{N}(\bzero,\bI_s)
	\end{equation*}
	in distribution. However, even with the additional requirement $\pi\rightarrow0$, we still cannot obtain a consistent estimate for $\nu_\ell$ by the plug-in method with estimating $b_\ell$ and $\tilde{b}_\ell$ separately for all $\gamma\gg p^{-1/3}\log^{1/6}p$. A novel adaptive bootstrap procedure will be developed in
		Section \ref{sec44}, which provides a unified estimation procedure for $\nu_\ell$ when
		$\gamma \asymp \pi \to 0$
		across the three different phases. On the other hand, the inference with $\gamma$ being a fixed constant
		can be obtained based on Theorem \ref{tm:asymnormal} with the estimated
		$\hat{b}_\ell$ specified in the Appendix.
	\end{remark}

	\section{Adaptive bootstrap inference.} \label{sec44}
	
	The goal of this section is primarily two-fold. First, we construct
	a novel bootstrap confidence interval for $\theta_\ell$ which is
	automatically adaptive to the three phases identified in Theorem \ref{tm:asymnormal}. Second, we leverage the new bootstrap procedure with Gaussian approximation to provide simultaneous
	inference for all $p$ components of $\btheta$ as $p\rightarrow\infty$. Additionally, we provide an algorithm for data-adaptive selection of a working parameter in our approach.  In the sequel, we always assume that the privacy level $\pi \to 0$ together with the number of nodes $p\to \infty$.
	
	\subsection{Bootstrap algorithm and simultaneous inference.}\label{sec:bootstrap}
	
As we have discussed in Remark \ref{rk3}(c), it holds that 
	\begin{equation}\label{eq:unified}
		N^{1/2}\,{\rm diag}(\nu_{\ell_1}^{-1/2},\ldots,\nu_{\ell_s}^{-1/2})(\hat{\theta}_{\ell_1}-\theta_{\ell_1},\ldots,\hat{\theta}_{\ell_s}-\theta_{\ell_s})^\T\rightarrow\mathcal{N}(\bzero,\bI_s)
	\end{equation}
	in distribution as $\gamma \asymp \pi \to 0$, where
	$
	\nu_\ell=(p-2)b_\ell+\tilde{b}_\ell
	$. Now we reproduce this structure in a bootstrap
	world based on the available network $\bZ$. The goal is to
	estimate $\nu_\ell$ adaptively regardless of the decay rate of $\gamma$.
	
	Recall $\mathcal{I}=\{(i,j): 1\le i<j \le p\}$. For a
	given constant $\delta\in(0,0.5)$, we draw bootstrap samples $\bZ^\dag=(Z_{i,j}^\dag)_{p\times p}$ according to
	\begin{equation}\label{eq:dag1}
		Z_{i,j}^\dag \equiv  Z_{j,i}^\dag = Z_{i,j}I(\eta_{i,j}=0)+I(\eta_{i,j}=1)\,, \quad (i,j) \in \calI\,,
	\end{equation}
	where $\{\eta_{i,j}\}_{(i,j)\in\mathcal{I}}$ are independent and identically distributed random variables only taking three possible values $-1, 0$ and $1$ with
		\begin{equation*}
		\mathbb{P}(\eta_{i,j}=0)=1-2\delta\,, ~~
		\mathbb{P}(\eta_{i,j}=1)=\delta ~~ \textrm{and}~~ \mathbb{P}(\eta_{i,j}=-1)=\delta\,.
	\end{equation*}
	For $i\neq j$ and $\tau\in\{0,1\}$, put
	\begin{align*}
		\varphi_{\tau}^\dag(x)=\{x-\delta-\alpha(1-2\delta)\}^{\tau}\{1-\delta-\beta(1-2\delta)-x\}^{1-\tau}
	\end{align*}
 with $x\in\{0, 1\}$.
	To simplify the notation, we write $\varphi_{\tau}^\dag(Z_{i,j}^\dag)$ as $\varphi_{(i,j),\tau}^\dag$ for any $i\neq j$ and $\tau\in\{0,1\}$. Note that
	\[
\mathbb{P}(X_{i,j}=0)=\frac{\mathbb{E}\{\varphi_{(i,j),0}^\dag\}}{(1-2\delta)(1-\alpha-\beta)}~~\textrm{and}~~
\mathbb{P}(X_{i,j}=1)=\frac{\mathbb{E}\{\varphi_{(i,j),1}^\dag\}}{(1-2\delta)(1-\alpha-\beta)}\,.
\]
For any given $(i,j)$ such that $i\neq j$, we know $Z_{i,j}^\dag$ is independent of
	$\{Z_{\tilde{i},\tilde{j}}^\dag:|\{\tilde{i},\tilde{j}\}\cap\{i,j\}|\leq1\}$. Hence, it follows from (\ref{eq:est1}) that
	\begin{equation}\label{eq:estimating2} \frac{\mathbb{E}\{\varphi_{(i,\ell),1}^\dag\varphi_{(i,j),0}^\dag\varphi_{(\ell,j),1}^\dag\}}{\mathbb{E}\{\varphi_{(i,\ell),0}^\dag\varphi_{(i,j),1}^\dag\varphi_{(\ell,j),0}^\dag\}}=\exp(2\theta_{\ell})\,, \quad i\neq j\neq \ell\,,
	\end{equation}
	which is a bootstrap analogue
	of \eqref{eq:estimating1}.
	Similarly, we define a bootstrap estimator for $\theta_\ell$ as:
	\begin{equation}\label{eq:estthetadag}
		\hat{\theta}_{\ell}^\dag={1\over 2} \log\bigg( \frac{\hat{\mu}_{\ell,1}^\dag}{
			\hat{\mu}_{\ell,2}^\dag}\bigg)\,,
	\end{equation}
	where
	\begin{align*}
	  \hat{\mu}_{\ell,1}^\dag=&\,\frac{1}{|\mathcal{H}_\ell|}\sum_{(i,j)\in\mathcal{H}_{\ell}}\varphi_{(i,\ell),1}^\dag\varphi_{(i,j),0}^\dag\varphi_{(\ell,j),1}^\dag\,,\\
\hat{\mu}_{\ell,2}^\dag=&\,\frac{1}{|\mathcal{H}_\ell|}\sum_{(i,j)\in\mathcal{H}_{\ell}}\varphi_{(i,\ell),0}^\dag\varphi_{(i,j),1}^\dag\varphi_{(\ell,j),0}^\dag\,.	    		\end{align*}
Such defined $\hat{\mu}_{\ell,1}^\dag$ and $\hat{\mu}_{\ell,2}^\dag$ are, respectively, the bootstrap analogues of $\hat{\mu}_{\ell,1}$ and $\hat{\mu}_{\ell,2}$ defined as \eqref{eq:hatmuell1} and \eqref{eq:hatmuell2}. For $\mu_{\ell,1}$, $\mu_{\ell,2}$ and $\lambda_{i,\ell}$ defined as \eqref{eq:muell1}, \eqref{eq:muell2} and \eqref{eq:lambdaiell}, we define their bootstrap analogues, respectively, as
	\begin{align*}
&~~~~~~~~~~~~~~~~~~~~~~~~~~~~\mu_{\ell,1}^\dag=\frac{1}{|\mathcal{H}_{\ell}|}\sum_{(i,j)\in\mathcal{H}_{\ell}}\mathbb{E}\big\{\varphi_{(i,\ell),1}^\dag\varphi_{(i,j),0}^\dag\varphi_{(\ell,j),1}^\dag\big\}\,,\\
&~~~~~~~~~~~~~~~~~~~~~~~~~~~~\mu_{\ell,2}^\dag=\frac{1}{|\mathcal{H}_{\ell}|}\sum_{(i,j)\in\mathcal{H}_{\ell}}\mathbb{E}\big\{\varphi_{(i,\ell),0}^\dag\varphi_{(i,j),1}^\dag\varphi_{(\ell,j),0}^\dag\big\}\,,\\
&\lambda_{i,\ell}^\dag=\frac{1}{p-2}\sum_{j:\,j\neq \ell,i}\bigg[\frac{1}{\mu_{\ell,1}^\dag}\mathbb{E}\big\{\varphi_{(\ell,j),1}^\dag\big\}\mathbb{E}\big\{\varphi_{(i,j),0}^\dag\big\}+\frac{1}{\mu_{\ell,2}^\dag}\mathbb{E}\big\{\varphi_{(\ell,j),0}^\dag\big\}\mathbb{E}\big\{\varphi_{(i,j),1}^\dag\big\}\bigg]\,.
\end{align*}
	Then $\hat{\theta}_{\ell}^\dag$ admits a similar asymptotic property
	as (\ref{eq:unified}). To present it explicitly, we
	let
	\begin{align}\label{eq:nulldag}
		\nu_{\ell}^\dag=(p-2)b_{\ell}^\dag+\tilde{b}_{\ell}^\dag\,, \quad \ell \in[p]\,,
	\end{align}
	where
	\begin{align*}
		&~~~~~~~~~~~~~~~~~~~~~~~~~~~~~b_\ell^\dag=\frac{1}{p-1}\sum_{i:\,i\neq \ell}\lambda_{i,\ell}^{\dag,2}{\rm Var}(Z_{i,\ell}^\dag)\,,\\
		&\tilde{b}_\ell^\dag=\frac{1}{2N}\bigg(\frac{\mu_{\ell,1}^\dag+\mu_{\ell,2}^\dag}{\mu_{\ell,1}^\dag\mu_{\ell,2}^\dag}\bigg)^2\sum_{i,j:\,i\neq j,\,i,j\neq\ell}{\rm Var}(Z_{i,\ell}^\dag){\rm Var}(Z_{\ell,j}^\dag){\rm Var}(Z_{i,j}^\dag)\,.
	\end{align*}

	\begin{theorem}\label{tm:2}
		Let the conditions of Theorem {\rm\ref{tm:asymnormal}} hold,
		and $\delta\in(0,c]$ for some positive constant $c<0.5$.
		As $p\to \infty$, if $1\gg\gamma \gg p^{-1/3}\log^{1/6} p $,
		the following two assertions hold.

{\rm (a)} Let $1\le \ell_1 < \cdots < \ell_s \le p$ be any $s$ given indices for some fixed integer $s\ge 1$. Then
		$$
		    N^{1/2}\,{\rm diag}(\nu_{\ell_1}^{\dag,-1/2},\ldots,\nu_{\ell_s}^{\dag,-1/2})
		(\hat{\theta}_{\ell_1}^\dag-\theta_{\ell_1},\ldots,\hat{\theta}_{\ell_s}^\dag-\theta_{\ell_s})^\T\rightarrow\mathcal{N}(\bzero,\bI_s)$$
		in distribution.

{\rm(b)} $\max_{\ell\in[p]}|\nu_\ell^\dag\nu_\ell^{-1}-1|=O(\delta)$, where $\nu_\ell$ is specified in {\rm(\ref{eq:unified})}.
	\end{theorem}

	Theorem \ref{tm:2} indicates that
	$\nu_\ell^\dag/\nu_\ell\rightarrow1$ for any $1\gg\gamma\gg
	p^{-1/3}\log^{1/6}p$ provided that we set $\delta=o(1)$. For fixed $s\geq1$ and given $1\leq \ell_1<\cdots<\ell_s\leq p$, we can draw bootstrap samples $\bZ^{\dag}$ as in \eqref{eq:dag1} with some $\delta=o(1)$, and compute the bootstrap estimate $(\hat{\theta}_{\ell_1}^\dag,\ldots,\hat{\theta}_{\ell_s}^\dag)^\T$ defined in \eqref{eq:estthetadag} based on $\bZ^\dag$. We repeat this procedure $M$ times for some large integer $M$ and compute $$
\hat{\nu}_{\ell_k}^\dag=\frac{N}{M}\sum_{m=1}^M\{\hat{\theta}_{\ell_k}^{\dag,(m)}-\bar{\hat{\theta}}_{\ell_k}^\dag\}^2\,,~~~k\in[s]\,,
$$
with $\bar{\hat{\theta}}_{\ell_k}^\dag=M^{-1}\sum_{m=1}^M\hat{\theta}_{\ell_k}^{\dag,(m)}$, where $\{\hat{\theta}_{\ell_1}^{\dag,(m)},\ldots,\hat{\theta}_{\ell_s}^{\dag,(m)}\}^\T$ is the associated bootstrap estimate in the $m$-th repetition. Then a confidence region for $(\theta_{\ell_1},\ldots,\theta_{\ell_s})^\T$ can be constructed based on the asymptotic approximation
$$N^{1/2}\,{\rm
		diag}(\hat{\nu}_{\ell_1}^{\dag,-1/2},\ldots,\hat{\nu}_{\ell_s}^{\dag,-1/2})(\hat{\theta}_{\ell_1}-\theta_{\ell_1},\ldots,\hat{\theta}_{\ell_s}-\theta_{\ell_s})^\T \overset{{\rm d}}{\approx} \mathcal{N}(\bzero,\bI_s)\,.$$
	
	Importantly, we note that in both Theorems \ref{tm:asymnormal} and \ref{tm:2}, $s$ is a fixed integer when $p \to \infty$.
	Hence the inference methods presented so far are not applicable to
	all $p$ components of $\btheta$ simultaneously.  However,
	a breakthrough can be had via the Gaussian approximation in
	Theorem \ref{tm:ga} below. To our best knowledge, this is the first
	method for simultaneous inference for all the $p$ components of
	$\btheta$ in the $\beta$-model.
	Write $\hat \btheta =(\hat \theta_1, \ldots, \hat\theta_p)^\T$ where $\hat{\theta}_\ell$ is the proposed moment-based
	estimator given in \eqref{eq:esttheta} based on the sanitized data $\bZ$. As shown in Proposition \ref{pn:asyexp}, the leading term of $\hat{\btheta}-\btheta$ cannot be formulated as
	a partial sum of independent (or weakly dependent) random vectors, which is different from the standard framework of Gaussian approximation
	\citep{CCK_2013,ChangChenWu_2021}. Hence the existing results of Gaussian approximation cannot be applied
	directly to obtain Theorem \ref{tm:ga}, which requires significant technical challenge to be overcome in our theoretical analysis. 
	
	\begin{theorem}\label{tm:ga}
		Let Condition {\rm\ref{cond1}} hold and $(\alpha,\beta)\in\mathcal{M}(\gamma;C_1)$ for some fixed constant $C_1\in(0,0.5)$. As $p\to\infty$, if $0< \delta\ll (p\log p)^{-1}$ and $1\gg\gamma \gg
		p^{-1/3}\log^{1/2} p $, then
			\[
	\sup_{\bu\in\mathbb{R}^p}\big|\mathbb{P}\big\{N^{1/2}({\bV}^{\dag})^{-1/2}(\widehat{\btheta}-\btheta)\leq \bu\big\}-\mathbb{P}(\bxi\leq \bu)\big|\rightarrow0\,,
	\]
		where $\bV^{\dag}={\rm diag}({\nu}_1^\dag,\ldots,{\nu}_p^\dag)$,
		and $\bxi\sim\mathcal{N}(\bzero,\bI_p)$.
	\end{theorem}
	
	Let $\bxi=(\xi_1,\ldots,\xi_p)^\T\sim\mathcal{N}(\bzero,\bI_p)$. For any $\mathcal{J}=\{\ell_1,\ldots,\ell_s\}\subset[p]$, write  \begin{align*}
	&\bV_{\mathcal{J}}^{\dag}={\rm diag}(\nu_{\ell_1}^\dag,\ldots,
	\nu_{\ell_s}^\dag)\,,~~ \hat{\btheta}_{\mathcal{J}}=(\hat{\theta}_{\ell_1},\ldots,\hat{\theta}_{\ell_s})^\T\,,\\ &~~~~\btheta_{\mathcal{J}}=(\theta_{\ell_1},\ldots,\theta_{\ell_s})^\T\,,~~\bxi_{\mathcal{J}}=(\xi_{\ell_1},\ldots,\xi_{\ell_s})^\T\,.
	\end{align*}
	Following the same arguments in the proof of Proposition 1 in the supplementary material of \cite{CZZZ_2017}, we can obtain from Theorem \ref{tm:ga} that
	\[
	\sup_{\mathcal{J}}\sup_{u\in\mathbb{R}}\big|\mathbb{P}\big\{N^{1/2}|({\bV}_{\mathcal{J}}^{\dag})^{-1/2}(\hat{\btheta}_{\mathcal{J}}-\btheta_{\mathcal{J}})|_\infty\leq u \big\}-\mathbb{P}(|\bxi_{\mathcal{J}}|_\infty\leq u)\big|\rightarrow0
	\]
	as $p\rightarrow\infty$. Given $\alpha\in(0,1)$ and $\mathcal{J}\subset[p]$,
	\begin{align}\label{eq:confreg}
		\bTheta_{\mathcal{J},\alpha}:=\bigg\{\ba\in\mathbb{R}^{|\mathcal{J}|}:N^{1/2}|({\bV}_{\mathcal{J}}^{\dag})^{-1/2}(\hat{\btheta}_{\mathcal{J}}-\ba)|_\infty\leq \Phi^{-1}\bigg(\frac{1+\alpha^{1/|\mathcal{J}|}}{2}\bigg) \bigg\}
	\end{align}
	is a $100\cdot\alpha\%$ confidence region for
	$\btheta_{\mathcal{J}}$, where $\Phi(\cdot)$ is the cumulative distribution function of the standard normal distribution. We refer to Section 4 of \cite{CQYZ_2018} for applications of this type of confidence region in simultaneous
	inference. If $\gamma$ is a fixed constant, Theorem \ref{tm:ga} still holds with replacing $\bV^\dag$ by $(p-2)\cdot{\rm diag}(\hat{b}_1,\ldots,\hat{b}_p)$ where $\hat{b}_\ell$ is given in \eqref{eq:hatbell} in the Appendix. If we set $\alpha=\beta=0$ in the jittering mechanism \eqref{a1}--\eqref{a3}, then $\gamma=1$ in this case and the released data $\bZ$ is identical to the original data $\bX$. Our simultaneous inference procedure still also works in this case.	
	
	\subsection{Adaptive selection of $\delta$.}
	
	The tuning parameter $\delta$ plays a key role in our simultaneous
	inference procedure. We propose a data-driven method in Algorithm \ref{alg:1} to select $\delta$.
	To illustrate the basic idea,
	we denote by $\nu_\ell^\dag(\delta)$ the associated
	$\nu_\ell^\dag$ defined in \eqref{eq:nulldag} with $\delta$ used in generating
	the bootstrap samples $\bZ^\dag$ in \eqref{eq:dag1}. If
	$\{\nu_\ell\}_{\ell\in\mathcal{J}}$ are known, the ideal selection for
	the tuning parameter $\delta$ should be $$\delta_{{\rm
			opt}}=\arg\min_{\delta>0}\max_{\ell\in\mathcal{J}}|\nu_{\ell}^\dag(\delta)-\nu_\ell|\,.$$
	Unfortunately,  $\{\nu_\ell\}_{\ell\in\mathcal{J}}$ are unknown in
	practice, as they depend on the unknown parameters
	$\theta_1,\ldots,\theta_p$. A natural idea is to replace $\nu_\ell$'s by	their estimates. Recall
	$\hat{\theta}_\ell=2^{-1}\log(\hat{\mu}_{\ell,1}\hat{\mu}_{\ell,2}^{-1})$
	with
	\begin{align*}
 \hat{\mu}_{\ell,1}=&~\frac{1}{|\mathcal{H}_\ell|}\sum_{(i,j)\in\mathcal{H}_{\ell}}
	\varphi_{(i,\ell),1}\varphi_{(i,j),0}\varphi_{(\ell,j),1}\,,\\
	\hat{\mu}_{\ell,2}=&~\frac{1}{|\mathcal{H}_\ell|}\sum_{(i,j)\in\mathcal{H}_{\ell}}
	\varphi_{(i,\ell),0}\varphi_{(i,j),1}\varphi_{(\ell,j),0}\,.
 \end{align*}
 Due
 to the
	nonlinear function $\log(\cdot)$ and the ratio between
	$\hat{\mu}_{\ell,1}$ and $\hat{\mu}_{\ell,2}$, $\hat{\theta}_\ell$
	usually includes some high-order bias term. More specifically,
	\[
	\hat{\theta}_{\ell}-\theta_{\ell}=\frac{\hat{\mu}_{\ell,1}-{\mu}_{\ell,1}}{2\mu_{\ell,1}}-\frac{\hat{\mu}_{\ell,2}-{\mu}_{\ell,2}}{2\mu_{\ell,2}}+\underbrace{\frac{(\hat{\mu}_{\ell,2}-{\mu}_{\ell,2})^2}{4\mu_{\ell,2}^2}-\frac{(\hat{\mu}_{\ell,1}-{\mu}_{\ell,1})^2}{4\mu_{\ell,1}^2}}_{\text{high-order bias}}+\hat{R}_{\ell}\,,
	\]
	where $\hat{R}_{\ell}$ is a negligible term in comparison to the high-order bias.
	Although the high-order bias
	has little impact on the estimation of $\theta_\ell$,
	it may lead to a bad estimate of $\nu_\ell$ if we just plug-in
	$\hat{\theta}_1,\ldots,\hat{\theta}_p$ in the nonlinear function
	$\nu_\ell$ which depends on $\theta_1,\ldots,\theta_p$. Hence, when we replace $\{\nu_\ell\}_{\ell\in\mathcal{J}}$ in Algorithm \ref{alg:1}, we use their associated estimates with bias-corrected $\hat{\theta}_1^{{\rm bc}},\ldots,\hat{\theta}_p^{{\rm bc}}$.
	Based
	on the optimal $\hat{\delta}_{\rm opt}$ selected in Algorithm \ref{alg:1}, we can replace the values $\{\nu_{\ell}^\dag\}_{\ell\in\mathcal{J}}$ in \eqref{eq:confreg} by $\{\hat{\nu}_{\ell}^\dag(\hat{\delta}_{\rm opt})\}_{\ell\in\mathcal{J}}$ specified in Algorithm \ref{alg:1} to construct a $100\cdot \alpha\%$ simultaneous confidence region for $\btheta_{\mathcal{J}}$ in practice.

	\begin{algorithm}[!htb]
		\caption{Selecting tuning parameter $\delta$}\label{alg:1}
		\begin{algorithmic}[1]
			\State obtain $\{\hat{\theta}_\ell\}_{\ell=1}^p$, $\{\hat{\mu}_{\ell,1}\}_{\ell=1}^p$ and $\{\hat{\mu}_{\ell,2}\}_{\ell=1}^p$ based on \eqref{eq:esttheta}, \eqref{eq:hatmuell1}	and \eqref{eq:hatmuell2}, respectively.
			\State calculate
			\begin{align*}
				\hat{\varphi}_{(i,j,\ell),1}
				&= \frac{(1-\alpha-\beta)\exp(\hat{\theta}_i+\hat{\theta}_{\ell})}{1+\exp(\hat{\theta}_i+\hat{\theta}_{\ell})}\frac{1-\alpha-\beta}{1+\exp(\hat{\theta}_i+\hat{\theta}_j)}\frac{(1-\alpha-\beta)\exp(\hat{\theta}_{\ell}+\hat{\theta}_j)}{1+\exp(\hat{\theta}_{\ell}+\hat{\theta}_j)}\, ,\\
				\hat{\varphi}_{(i,j,\ell),2}
				&=  \frac{1-\alpha-\beta}{1+\exp(\hat{\theta}_i+\hat{\theta}_{\ell})}\frac{(1-\alpha-\beta)\exp(\hat{\theta}_i+\hat{\theta}_j)}{1+\exp(\hat{\theta}_i+\hat{\theta}_j)}\frac{1-\alpha-\beta}{1+\exp(\hat{\theta}_{\ell}+\hat{\theta}_j)}\, .
			\end{align*}	
			\Repeat
			\State  leave out one $(i,j)\in\mathcal{H}_\ell$ randomly and denote by $\mathcal{H}_\ell^{-}$ the set including the rest elements in $\mathcal{H}_\ell$.
			\State  calculate $$\tilde{\mu}_{\ell,1}=\frac{1}{|\mathcal{H}_\ell^{-}|}\sum_{(i,j)\in\mathcal{H}_{\ell}^{-}}\hat{\varphi}_{(i,j,\ell),1}~~\textrm{and}~~\tilde{\mu}_{\ell,2}=\frac{1}{|\mathcal{H}_\ell^{-}|}\sum_{(i,j)\in\mathcal{H}_{\ell}^{-}}\hat{\varphi}_{(i,j,\ell),2}\,,$$ which provide the estimates of $\mu_{\ell,1}$ and $\mu_{\ell,2}$, respectively.
			\State calculate $${\rm bias}_\ell=4^{-1}\tilde{\mu}_{\ell,2}^{-2}(\hat{\mu}_{\ell,2}-\tilde{\mu}_{\ell,2})^2-4^{-1}\tilde{\mu}_{\ell,1}^{-2}(\hat{\mu}_{\ell,1}-\tilde{\mu}_{\ell,1})^2\,.$$
			\Until $M$ replicates obtained, for a large integer $M$, and get ${\rm bias}_\ell^{(1)},\ldots,{\rm bias}_\ell^{(M)}$.
			\State approximate the high-order bias in $\hat{\theta}_\ell$ by $$\widehat{{\rm bias}}_\ell=\frac{1}{M}\sum_{m=1}^M{\rm bias}_\ell^{(m)}\,,$$ and obtain $\hat{\theta}_\ell^{\rm bc}=\hat{\theta}_\ell-\widehat{{\rm bias}}_\ell$, the bias-correction for $\hat{\theta}_\ell$.
			\State calculate $$\tilde{\mu}_{\ell,1}^{\rm bc}=\frac{1}{|\mathcal{H}_\ell|}\sum_{(i,j)\in\mathcal{H}_\ell}\tilde{\varphi}_{(i,j,\ell),1}~~\textrm{and}~~ \tilde{\mu}_{\ell,2}^{\rm bc}=\frac{1}{|\mathcal{H}_\ell|}\sum_{(i,j)\in\mathcal{H}_\ell}\tilde{\varphi}_{(i,j,\ell),2}\,,$$ where $\tilde{\varphi}_{(i,j,\ell),1}$ and $\tilde{\varphi}_{(i,j,\ell),2}$ are defined in the same manner as $\hat{\varphi}_{(i,j,\ell),1}$ and $\hat{\varphi}_{(i,j,\ell),2}$, respectively, with replacing $\{\hat{\theta}_\ell\}_{\ell=1}^p$ by $\{\hat{\theta}_\ell^{\rm bc}\}_{\ell=1}^p$.
			\State calculate $\hat{\nu}_\ell^{\rm bc}=(p-2)\hat{b}_\ell^{\rm bc}+\hat{\tilde{b}}_\ell^{\rm bc}$, where $\hat{b}_\ell^{\rm bc}$ and $\hat{\tilde{b}}_\ell^{\rm bc}$ are defined in the same manner of $\hat{b}_\ell$ and $\hat{\tilde{b}}_\ell$ specified as \eqref{eq:hatbell} in the Appendix with replacing $(\hat{\mu}_{\ell,1},\hat{\mu}_{\ell,2},\{\hat{\theta}_k\}_{k=1}^p)$ by $(\tilde{\mu}_{\ell,1}^{\rm bc},\tilde{\mu}_{\ell,2}^{\rm bc},\{\hat{\theta}_k^{\rm bc}\}_{k=1}^p)$.
			\Repeat
			\State  given $\delta>0$ and draw bootstrap samples $\bZ^\dag=(Z_{i,j}^\dag)_{p\times p}$ as in (\ref{eq:dag1}),
			calculate the bootstrap estimate
			$\hat{\theta}_{\ell}^\dag$ defined in \eqref{eq:estthetadag} based on the bootstrap samples $\bZ^\dag$.
			\Until $M$ replicates obtained, for a large integer $M$, and get $\hat{\theta}_{\ell}^{\dag,(1)},\ldots,\hat{\theta}_{\ell}^{\dag,(M)}$.
			
			\State calculate $$\hat{\nu}_{\ell}^\dag(\delta)=\frac{p^2}{M}\sum_{m=1}^M\{\hat{\theta}_{\ell}^{\dag,(m)}-\bar{\hat{\theta}}_{\ell}^\dag\}^2$$ with $\bar{\hat{\theta}}_{\ell}^\dag=M^{-1}\sum_{m=1}^M\hat{\theta}_{\ell}^{\dag,(m)}$.
			
			\State select $$\hat{\delta}_{\rm opt}=\arg\min_{\delta>0}\max_{\ell\in\mathcal{J}}|\hat{\nu}_\ell^\dag(\delta)-\hat{\nu}_\ell^{\rm bc}|\,.$$

		\end{algorithmic}
	\end{algorithm}

	\section{Numerical study.}\label{sec6}

\subsection{Simulation.}

	In this section we illustrate the finite sample
	properties of our proposed method of estimation and inference for the unknown parameters in the
	$\beta$-model by simulation. For $p\in\{1000,2000\}$, we draw
$\theta_1,\ldots,\theta_p$ independently from $ \mathcal{N}(0,0.2)$, and then
generate the adjacency matrix $\bX$ according to the $\beta$-model \eqref{d2}. For a given original network $\bX$, we set $\alpha = \beta
	\in\{0,0.1,0.2,0.3\}$
	 in the data release mechanism \eqref{a1} and \eqref{a2} to
	generate $\bZ$. Note that $\bZ = \bX$
	when $\alpha = \beta = 0$.

Based on the released data $\bZ$, we applied the moment-based
method \eqref{eq:esttheta} to estimate
$\btheta=(\theta_1,\ldots,\theta_p)^\T$, and then calculated the
estimation error
$L(\hat{\btheta})=p^{-1}|\hat{\btheta}-\btheta|_2^2$. For comparison, we also considered to apply the MLE of \cite{KarwaSlavkovic_2016} to the degree sequence of the released data $\bZ$.  Table \ref{tab1} reports the
	averages, medians and  standard deviations of the estimation
errors 
over 500 replications.
The proposed moment-based estimation performed competitively in relation to the MLE, though
the MLE is slightly more accurate overall.
However the MLE method
is memory-demanding
when $p$ is large.
For example with $p=1000$ and $\alpha=\beta=0.1$, the step generating a graph with given degree
sequence (i.e. Algorithm 2 of \cite{KarwaSlavkovic_2016})
occupied 3.91 GB memory. In contrast, the newly proposed moment-based estimation
only used 38.19 MB memory.
Furthermore, the MLE is excessively time-consuming computationally when $p$ is large.
See Table \ref{tab1} for the recorded average CPU times for each realization on
an Intel(R) Xeon(R) Platinum 8160 processor
(2.10GHz). With $p=1000$,  the average required CPU time for computing
the MLE once is over 471 minutes with the original data $\bX$ (i.e.
$\alpha=\beta=0$) and is almost
double with the sanitized data $\bZ$ (i.e. $\alpha, \beta >0$).
It is practically infeasible to conduct the simulation (with replications) for
all scenarios with $p=2000$, for which we only report the results with $\alpha=\beta=0$
with the average CPU time 5095 minutes per estimation.

We note that Algorithm 2 of \cite{KarwaSlavkovic_2016} might be made more efficient if it is modified to directly estimate the node degree sequence without actually producing the intermediate graph, the latter step which requires MCMC. Additionally, such an approach might also help with convergence issues.  In particular, and as an important caveat to the above results, we note that in order to achieve MLE estimates for $500$ trials in our simulations it was necessary to discard a nontrivial fraction of trials for which the MCMC algorithm failed to converge.  Specifically, when $\alpha=\beta=0.1$, $0.2$, and $0.3$, the proportion  of trials that needed to be discarded were, respectively, $3\%$, $10\%$ and $21\%$.  That is, MLE convergence was increasingly problematic with increasing noise level and hence with increasing privacy.  No trials were discarded for our proposed moment-based approach.

Based on our moment-based estimator $\hat{\btheta}$, we also constructed
the simultaneous confidence regions
\eqref{eq:confreg} for all the $p$ components
$\theta_1,\ldots,\theta_p$. To determine the tuning parameter $\delta$,
we applied the data-driven Algorithm \ref{alg:1} with $M=500$.
Table \ref{tab3} lists the relative frequencies, in 500 replications for
each settings, of the occurrence
of the event that
the constructed confidence region contains the true value
of $\btheta$.
At each of the three nominal levels,
those relative frequencies are always close to the corresponding nominal level.
	

\begin{table}[thp]
		\centering
		\caption{Estimation errors of the proposed moment-based
estimation and the maximum likelihood estimation for $\btheta$ in the $\beta$-model \eqref{d2}.
Also reported are the average CPU times (in minutes) for completing the estimation once for each
of the two methods. }\label{tab1}
\resizebox{\textwidth}{!}{
		\begin{threeparttable}
			\begin{tabular}{ccccccccccccc}
				\toprule
				\midrule
				&& \multicolumn{4}{c}{Proposed method} & \multicolumn{4}{c}{Maximum likelihood estimation}\\
				\midrule
				$p$& Summary statistics	& $\alpha=\beta=0$ & $\alpha=\beta=0.1$ &$\alpha=\beta= 0.2$& $\alpha=\beta=0.3$ &$\alpha=\beta=0$ & $\alpha=\beta=0.1$ &$\alpha=\beta= 0.2$& $\alpha=\beta=0.3$\\
				\midrule
				1000&	Average & 0.0041  & 0.0065  & 0.0117&0.0274   & 0.0062  & 0.0057&0.0107& 0.0239 \\
				&	Median & 0.0041  & 0.0065  & 0.0117& 0.0274  & 0.0041  &0.0057 &0.0107&0.0239   \\
				&	Standard deviation &0.0002   & 0.0003  & 0.0006&0.0012  & 0.0085  & 0.0002&0.0005&0.0015     \\
				&	Time (min) & 1.0340  & 1.0439  & 0.9191& 0.8540& 471.6290  & 850.2369&754.2811& 780.9615 \\
				\midrule
				2000&	Average & 0.0020  & 0.0032  & 0.0058&0.0133  & 0.0058  &NA&NA& NA \\
				&	Median & 0.0020  &  0.0032 & 0.0058& 0.0133 &0.0043   &NA&NA& NA  \\
				&	Standard deviation &  0.0001 &  0.0001 &0.0002 & 0.0004 &0.0019   &NA&NA&NA    \\
				&	Time (min) & 4.2333  & 4.8707  & 3.7540& 3.7256 & 5095.0520  &NA&NA&NA   \\	
				\midrule
				\toprule
			\end{tabular}%
		\end{threeparttable}
}
\end{table}
	
	\begin{table}[htbp]
		\centering
		\caption{Empirical frequencies of the constructed simultaneous confidence regions for $\btheta$ covering the truth in the $\beta$-model \eqref{d2}.}\label{tab3}
		\begin{threeparttable}
			\setlength{\tabcolsep}{4mm}
			\begin{tabular}{ccccccccccccc}
				\toprule
				\midrule
				\multicolumn{1}{c}{$p$} & \multicolumn{1}{c}{Level} &     	 \multicolumn{1}{c}{$\alpha=\beta=0$} &       \multicolumn{1}{c}{$\alpha=\beta=0.1$}       &\multicolumn{1}{c}{$\alpha=\beta=0.2$}
				&\multicolumn{1}{c}{$\alpha=\beta=0.3$}
				  \\
				\midrule
				1000&90\%  & 0.876   &0.868    & 0.910& 0.888  \\			
				&95\%  & 0.932    & 0.928  & 0.958& 0.948  \\
				&99\%  &   0.984  & 0.982   & 0.982 & 0.992 \\
				\midrule
				2000&90\%  & 0.900   & 0.876   & 0.898& 0.896   \\			
				&95\%  &   0.950  & 0.956  &0.946 & 0.952  \\
				&99\%  & 0.988    &  0.990  & 0.996 & 0.992\\
				\midrule
				\toprule
			\end{tabular}%
		\end{threeparttable}
	\end{table}%

\subsection{Real data analysis.}

Facebook, a social networking site launched
	in February 2004, now overwhelms numerous aspects of everyday life, and has become an immensely popular societal obsession. The Facebook friendships define a network of undirected edges that connect individual users. In this section, we analyze a small Facebook friendship network dataset
	available at {\tt http://wwwlovre.appspot.com/support.jsp}. 
	The network consists of 334 nodes and 2218 edges. 

	We fit the $\beta$-model to this network. As an illustration
	on the impact of the `jittering', we identify the nodes with the associated parameters equal to 0 based on both the original network and some sanitized versions.
	More specifically, we first consider the multiple hypothesis tests:
		\begin{align*}
		H_{0,\ell} : \theta_{\ell} = 0 ~~~\mbox{versus}~~~ H_{1,\ell}: \theta_{\ell} \neq 0 \,
	\end{align*}
	for $1\le \ell \le 334.$
	The moment-based estimate
	$\hat{\btheta}=(\hat{\theta}_1,\ldots,\hat{\theta}_{334})^{\T}$ based on the original data $\bX$ is calculated according
	to \eqref{eq:esttheta}. Theorem \ref{tm:asymnormal} indicates that the p-value for the $\ell$-th test is given by $2\{1-\Phi(\sqrt{333}\hat{b}^{-1/2}_\ell|\hat{\theta}_\ell|)\}$ with $\hat{b}_\ell$ defined as \eqref{eq:hatbell} in the Appendix. Note that $\hat{\theta}_{\ell_1}$ and $\hat{\theta}_{\ell_2}$ are asymptotically independent for any $\ell_1\neq \ell_2$.  The BH procedure \citep{Benjamini_1995} at the rate $1\%$ for the 334 multiple tests identifies the 10 nodal parameters ($\theta_2, \theta_{21}, \theta_{33}, \theta_{51}, \theta_{78}, \theta_{186}, \theta_{202}, \theta_{211}, \theta_{263}, \theta_{272}$) being not significantly different from 0.  Put $\mathcal{J}=\{2,21,33,51,78,186,202,211,263,272\}$. We consider now the testing problem for the single hypothesis setting
	\begin{align}\label{eq:globaltest}
	H_{0} : \btheta_{\mathcal{J}} ={\bf 0} ~~~\mbox{versus}~~~ H_{1}: \btheta_{\mathcal{J}} \neq {\bf 0} \,
	\end{align}
	based on both the original network $\bX$ and
	its sanitized versions $\bZ$ via  jittering mechanism \eqref{a1} with $\alpha=\beta=0.1, 0.2$ and $0.3$.
	Let $\bzeta_1,\dots, \bzeta_{1000}$ be independent and $\mathcal{N}(\bzero,\bI_{10})$. By Theorem \ref{tm:ga}, the p-value of the test for \eqref{eq:globaltest} based on $\bZ$ is approximately
		$$
	\frac{1}{1000}\sum_{m=1}^{1000} I\{|\bzeta_m|_{\infty}\geq \sqrt{333\times332}|\widehat{\bV}_{\mathcal{J}}^{-1/2}\hat{\btheta}_{\mathcal{J}}^{(\bZ)}|_{\infty}\}\,,$$
	where $\hat{\btheta}_{\mathcal{J}}^{(\bZ)}$ is the estimate of $\btheta_{\mathcal{J}}$ based on $\bZ$ by the moment-based method \eqref{eq:esttheta}, and $\widehat{\bV}_{\mathcal{J}}$ is the estimate of the asymptotic covariance of $\sqrt{333\times 332}\{\hat{\btheta}_{\mathcal{J}}^{(\bZ)}-\btheta_{\mathcal{J}}\}$. When $\alpha=\beta=0$, $\bZ=\bX$, the p-value for testing \eqref{eq:globaltest} based on $\bX$ is then 0.1019.
	As the test based on $\bZ$ depends on a particular realization when $\alpha=\beta=0.1, \, 0.2$ and $0.3$, we repeat the test 500 times for each setting.
	The average p-values of those 500 tests (based on $\bZ$)
	with $\alpha=\beta=0.1, \, 0.2$ and $0.3$ are, respectively,
	$0.1276$, $0.1522$ and $0.1874$, which are reasonably close to the p-value based on $\bX$. The standard errors of the 500 p-values are $0.0795$, $0.1281$ and $0.1408$, respectively, for $\alpha=\beta=0.1, 0.2$ and $0.3$.
	
This small illustration suggests that, with increasing edge noise (and hence increasing privacy), the resulting p-value is increasingly over-estimated with increasing standard error.  Both trends are to be expected -- since with increasing edge-noise the signal will be weakened -- and merit future study.
	

\section{Extension to the sparse $\beta$-model.} \label{sec:sparsemodel}
Under Condition \ref{cond1} imposed on the $\beta$-model \eqref{d2}, we have
	\[
	\min_{i,j:\,i<j}\mathbb{P}(X_{i,j}=1)\geq \frac{\exp(-c)}{1+\exp(-c)}\asymp 1
	\]
	for some positive constant $c$, which implies the expected number of edges of the network should be of order at least $p^2$ and thus the network will be dense. In this last section we illustrate how our results may be extended to the case of sparse networks, through several additional results. A full generalization of our results for the dense case, inclusive of the bootstrap-based inferential procedure, is beyond the present scope.

 To model the sparse networks, \cite{ChenKatoLeng_2021} consider the sparse $\beta$-model defined as
\begin{align}\label{eq:sparsebeta}
\mathbb{P}(X_{i,j}=1)=\frac{\exp(\xi+\check{\theta}_i+\check{\theta}_j)}{1+\exp(\xi+\check{\theta}_i+\check{\theta}_j)}\,,
\end{align}
where $\xi\in\mathbb{R}$ and $\check{\btheta}=(\check{\theta}_1,\ldots,\check{\theta}_p)^\T\in\mathbb{R}_+^p$ are both unknown parameters with $|\check{\btheta}|_0\ll p$ and $\min_{\ell\in[p]}\check{\theta}_\ell=0$. Denote by $\mathcal{S}$ the support of $\check{\btheta}$, that is $\mathcal{S}=\{\ell\in[p]:\check{\theta}_\ell\neq 0\}$. Write  $|\mathcal{S}|=s$. Given some constants $\omega_1\in[0,2)$ and $\omega_2\in[0,1)$ such that $0\leq \omega_1-\omega_2<1$, \cite{ChenKatoLeng_2021} consider the reparametrization
\[
\xi=-\omega_1\log p+\xi^
+~~\textrm{and}~~\check{\theta}_\ell=\omega_2\log p+\check{\theta}_\ell^+~\textrm{for all}~\ell\in\mathcal{S}\,,
\]
where $|\xi^+|=o(\log p)$ and $\max_{\ell\in\mathcal{S}}|\check{\theta}_\ell^+|=o(\log p)$.
Let
\begin{equation}\label{eq:thetaellxi}
\theta_\ell=\frac{\xi}{2}+\check{\theta}_\ell\,,~~~~\ell\in[p]\,.
\end{equation}
The sparse $\beta$-model \eqref{eq:sparsebeta} can be reformulated as the standard $\beta$-model \eqref{d2} with
\[
|\btheta|_\infty\left\{\begin{aligned}
   \sim \bigg|\frac{\omega_1}{2}-\omega_2\bigg|\log p\,,~~&\textrm{if}~\omega_1\neq 2\omega_2\,,  \\
   =o(\log p)\,,~~~~~~~&\textrm{if}~\omega_1=2\omega_2\,.
\end{aligned}\right.
\]

Applying the estimation procedure given in Section \ref{sec42} to the sanitized network $\bZ=(Z_{i,j})_{p\times p}$ defined as in \eqref{a1}--\eqref{a3}, we can also obtain the moment-based estimator  $\hat{\theta}_\ell$ defined as \eqref{eq:esttheta} for the unknown parameter $\theta_\ell$ given in \eqref{eq:thetaellxi}. For the positive stochastic sequence $\{a_p\}$ and the positive sequence $\{c_p\}$, we write $a_p=\tilde{O}_{\p}(c_p)$ if $a_p=O_{\p}(p^{\epsilon}c_p)$ for some sufficiently small fixed constant $\epsilon>0$. Proposition \ref{tm:4} gives the convergence rate of  $\max_{\ell\in[p]}|\hat{\theta}_{\ell}-\theta_{\ell}|$  under the sparse $\beta$-model.

\begin{proposition}\label{tm:4}
Let $(\alpha,\beta)\in\mathcal{M}(\gamma, C_1)$ for some fixed constant $C_1\in(0,0.5)$. Write $\chi_{p}=\exp(-|\xi^+|\vee\max_{\ell\in\mathcal{S}}|\check{\theta}_\ell^+|)$. If $0\leq \omega_2\leq \omega_1<1/2$, then
\begin{align*}
    \max_{\ell \in [p]}|\hat{\theta}_{\ell}-\theta_{\ell}|=
    &~\tilde{O}_{\p}\bigg(\frac{\log^{1/2}p}{\gamma p^{1/2-\omega_1}}\bigg)+\tilde{O}_{\p}\bigg(\frac{s\log^{1/2}p}{\gamma p^{3/2-\omega_1-\omega_2}}\bigg)
    +\tilde{O}_{\p}\bigg(\frac{\log^{1/2}p}{\gamma^{3}p^{1-2\omega_1}}\bigg)\,.
\end{align*}
provided that
	$
 \gamma\gg\chi_p^{-8}(sp^{-3/2+\omega_1+\omega_2}\log^{1/2}p+p^{-1/3+2\omega_1/3}\log^{1/6}p)
$.
	\end{proposition}

	\begin{remark}\label{rek4}
Under the assumption $|\xi^+|\vee\max_{\ell\in\mathcal{S}}|\check{\theta}_\ell^+|=o(\log p)$, we know $\chi_p^{-1}=\exp\{o(\log p)\}$. As shown in Section \ref{sec:pftm4} of the supplementary material, there exists some universal positive constant $c$ such that
\begin{align*}
    \max_{\ell \in [p]}|\hat{\theta}_{\ell}-\theta_{\ell}|=&~\chi_p^{-c}\cdot\bigg\{{O}_{\p}\bigg(\frac{\log^{1/2}p}{\gamma p^{1/2-\omega_1}}\bigg)+{O}_{\p}\bigg(\frac{s\log^{1/2}p}{\gamma p^{3/2-\omega_1-\omega_2}}\bigg)+{O}_{\p}\bigg(\frac{\log^{1/2}p}{\gamma^{3}p^{1-2\omega_1}}\bigg)\bigg\}
\end{align*}
provided that $\gamma\gg \chi_p^{-8}(sp^{-3/2+\omega_1+\omega_2}\log^{1/2}p+p^{-1/3+2\omega_1/3}\log^{1/6}p)$.
If the network $\bX$ is dense with $\omega_1=0$, $|\xi^+|\leq C$ and $ \max_{\ell\in\mathcal{S}}|\check{\theta}_\ell^+|\leq C$ for some universal positive constant $C$, it follows from Proposition \ref{tm:4} that
\begin{align*}
\max_{\ell\in[p]}|\hat{\theta}_\ell-\theta_\ell|=O_{\rm p}\bigg(\frac{\log^{1/2}p}{\gamma p^{1/2}}\bigg)+O_{\rm p}\bigg(\frac{\log^{1/2}p}{\gamma^{3}p}\bigg)
\end{align*}
provided that $\gamma\gg p^{-1/3}\log^{1/6}p$, which is identical to the result in Proposition \ref{tm:1}.
	\end{remark}
	

By \eqref{eq:thetaellxi} and $s\ll p$ in the sparse $\beta$-model, we can estimate $\xi$ and $\check{\theta}_\ell$ as follows:
\begin{equation}\label{eq:hatxihatthetaell}
\hat{\xi}=\frac{2}{p}\sum_{\ell\in[p]}\hat{\theta}_\ell~~~\textrm{and}~~~\hat{\check{\theta}}_\ell=\hat{\theta}_\ell-\frac{\hat{\xi}}{2}\,.
\end{equation}
Due to $|\hat{\check{\theta}}_\ell-\check{\theta}_\ell|\leq|\hat{\theta}_\ell-\theta_\ell|+|\hat{\xi}-\xi|/2$ and
\begin{align*}
|\hat{\xi}-\xi|=\bigg|\frac{2}{p}\sum_{\ell\in[p]}(\hat{\theta}_\ell-\theta_\ell+\check{\theta}_\ell)\bigg|\leq2\max_{\ell\in[p]}|\hat{\theta}_\ell-\theta_\ell|+O\bigg(\frac{s\log p}{p}\bigg)\,,
\end{align*}
by Proposition \ref{tm:4}, we have the following theorem.

\begin{theorem}\label{tm:5}
Let $(\alpha,\beta)\in\mathcal{M}(\gamma, C_1)$ for some fixed constant $C_1\in(0,0.5)$. Write $\chi_{p}=\exp(-|\xi^+|\vee\max_{\ell\in\mathcal{S}}|\check{\theta}_\ell^+|)$. If $0\leq \omega_2\leq \omega_1<1/2$, then
\begin{align*}
    |\hat{\xi}-\xi|=
    &~\tilde{O}_{\p}\bigg(\frac{\log^{1/2}p}{\gamma p^{1/2-\omega_1}}\bigg)+\tilde{O}_{\p}\bigg(\frac{s\log^{1/2}p}{\gamma p^{3/2-\omega_1-\omega_2}}\bigg)\\
    &
    +\tilde{O}_{\p}\bigg(\frac{\log^{1/2}p}{\gamma^{3}p^{1-2\omega_1}}\bigg)+O\bigg(\frac{s\log p}{p}\bigg)=\max_{\ell\in[p]}|\hat{\check{\theta}}_\ell-\check{\theta}_\ell|
\end{align*}
provided that
	$\gamma\gg \chi_p^{-8}(sp^{-3/2+\omega_1+\omega_2}\log^{1/2}p+p^{-1/3+2\omega_1/3}\log^{1/6}p)$.
\end{theorem}	

\begin{remark}
For known $(\omega_1,\omega_2)$ and $\mathcal{S}$, Theorem 1 of \cite{ChenKatoLeng_2021} specifies the convergence rates of the MLE for $\xi^+$ and $\{\check{\theta}_\ell^+\}_{\ell\in\mathcal{S}}$ based on the true network $\bX$ rather than the sanitized network $\bZ$ (i.e., $\alpha=\beta=0$ in our setting). Denote by $\tilde{\xi}^+$ and $\tilde{\check{\theta}}_\ell^+$, respectively, the MLE of $\xi^+$ and $\check{\theta}_\ell^+$ proposed in \cite{ChenKatoLeng_2021}. To simplify our comparison, we assume $|\xi^+|\vee \max_{\ell\in\mathcal{S}}|\check{\theta}_\ell^+|=O(1)$. Under the restriction $s=O\{ p^{(1-\omega_2)/2-c}\}$ for some sufficiently small constant $c>0$, Theorem 1(ii) of \cite{ChenKatoLeng_2021} implies $|\tilde{\xi}^+-\xi^+|=O_{\rm p}(p^{-1+\omega_1/2})$ and $|\tilde{\check{\theta}}_\ell^+-\check{\theta}_\ell^+|=O_{\rm p}\{p^{-1/2+(\omega_1-\omega_2)/2}\}$ for any $\ell\in\mathcal{S}$. With known $(\omega_1,\omega_2)$, we can obtain the following estimators for $\xi^+$ and $\check{\theta}_\ell^+$ based on $\hat{\xi}$ and $\hat{\check{\theta}}_\ell$ given in \eqref{eq:hatxihatthetaell}:
\[
\hat{\xi}^+=\hat{\xi}+\omega_1\log p~~~\textrm{and}~~~\hat{\check{\theta}}_\ell^+=\hat{\check{\theta}}_\ell-\omega_2\log p\,.
\]
Recall $\gamma=1-\alpha-\beta$. By Theorem \ref{tm:5} with $\gamma=1$ and $s=O\{ p^{(1-\omega_2)/2-c}\}$ for some sufficiently small constant $c>0$, it holds that
$
 |\hat{\xi}^+-\xi^+|=\tilde{O}_{\p}(p^{-1/2+\omega_1}\log ^{1/2}p)$ and $\max_{\ell\in[p]}|\hat{\check{\theta}}_\ell^+-\check{\theta}_\ell^+|=\tilde{O}_{\p}(p^{-1/2+\omega_1}\log ^{1/2}p)$, which are slower than the convergence rates of the MLE considered in \cite{ChenKatoLeng_2021}. Their method cannot be implemented directly with unknown $(\omega_1,\omega_2)$ while our moment-based method can still work.
\end{remark}

\begin{appendix}
\section*{}
A brief discussion of the fundamental issue of estimating asymptotic variances in Theorem \ref{tm:asymnormal} is provided here. 
	If we know the decay rate of $\gamma$ falls into which region, we may consider to construct the confidence region of $(\theta_{\ell_1},\ldots,\theta_{\ell_s})^\T$ based on Theorem \ref{tm:asymnormal} with the plug-in method. To do this, we need to estimate $b_{\ell_k}$'s and $\tilde{b}_{\ell_k}$'s first. By (\ref{eq:lambdaiell}), we can estimate $\lambda_{i,\ell}$ by
	\begin{equation}\label{eq:hatlambdaiell}
		\hat{\lambda}_{i,\ell}=\frac{1}{p-2}\sum_{j:\,j\neq \ell,i}\bigg\{\frac{1}{\hat{\mu}_{\ell,1}}\varphi_{(\ell,j),1}\varphi_{(i,j),0}+\frac{1}{\hat{\mu}_{\ell,2}}\varphi_{(\ell,j),0}\varphi_{(i,j),1}\bigg\}
	\end{equation}
	with $\hat{\mu}_{\ell,1}$ and $\hat{\mu}_{\ell,2}$ specified in \eqref{eq:hatmuell1} and \eqref{eq:hatmuell2}, respectively.
	By the definition of $Z_{i,j}$, we have
	\[
{\rm Var}(Z_{i,j})=\frac{\alpha+(1-\beta)\exp(\theta_i+\theta_j)}{1+\exp(\theta_i+\theta_j)}\cdot\frac{1-\alpha+\beta\exp(\theta_i+\theta_j)}{1+\exp(\theta_i+\theta_j)}
\]
	for any $i\neq j$. We can estimate ${\rm Var}(Z_{i,j})$ by
	\begin{equation}\label{eq:hatvar}
		\widehat{{\rm Var}}(Z_{i,j})=\frac{\alpha+(1-\beta)\exp(\hat{\theta}_i+\hat{\theta}_j)}{1+\exp(\hat{\theta}_i+\hat{\theta}_j)}\cdot\frac{1-\alpha+\beta\exp(\hat{\theta}_i+\hat{\theta}_j)}{1+\exp(\hat{\theta}_i+\hat{\theta}_j)}\,.
	\end{equation}
	Based on (\ref{eq:hatlambdaiell}) and (\ref{eq:hatvar}), we can estimate $b_\ell$ and $\tilde{b}_\ell$, respectively, by
	\begin{align}
		&~~~~~~~~~~~~~~~~~~~~~~~~~~~\hat{b}_\ell=\frac{1}{p-1}\sum_{i:\,i\neq \ell}\hat{\lambda}_{i,\ell}^2\widehat{{\rm Var}}(Z_{i,\ell})\,,\label{eq:hatbell}\\
		&\hat{\tilde{b}}_\ell=\frac{1}{2N}\bigg(\frac{\hat{\mu}_{\ell,1}+\hat{\mu}_{\ell,2}}{\hat{\mu}_{\ell,1}\hat{\mu}_{\ell,2}}\bigg)^2\sum_{i,j:\,i\neq j,\,i,j\neq\ell}\widehat{{\rm Var}}(Z_{i,\ell})\widehat{{\rm Var}}(Z_{\ell,j})\widehat{{\rm Var}}(Z_{i,j})\,.\notag
	\end{align}
	The convergence rates of such estimates are presented in Proposition \ref{tm:varc}. The proof of Proposition \ref{tm:varc} is given in Section \ref{sec:pfpn5} of the supplementary material.
	
	\begin{proposition}\label{tm:varc}
		Let Condition {\rm\ref{cond1}} hold and $(\alpha,\beta)\in\mathcal{M}(\gamma;C_1)$ for some fixed constant $C_1\in(0,0.5)$. If $\gamma\gg p^{-1/3}\log^{1/6}p$, for any given $\ell\in[p]$, it holds that
			\begin{align*}
	&~~\bigg|\frac{\hat{b}_\ell}{b_\ell}-1\bigg|=O_{\p}\bigg(\frac{\log p}{\gamma^4p}\bigg)+O_{\p}\bigg(\frac{\log^{1/2}p}{\gamma^2p^{1/2}}\bigg)\,,\\
&\bigg|\frac{\hat{\tilde{b}}_\ell}{\tilde{b}_\ell}-1\bigg|=O_{\p}\bigg(\frac{\log^{1/2}p}{\gamma^3p}\bigg)+O_{\p}\bigg(\frac{\log^{1/2}p}{\gamma p^{1/2}}\bigg)\,.
	\end{align*}
	
	\end{proposition}
	
	For any fixed integer $s\geq1$, Theorem \ref{tm:asymnormal} and Proposition \ref{tm:varc} imply that
$$p^{1/2}{\rm diag}(\hat{b}_{\ell_1}^{-1/2},\ldots,\hat{b}_{\ell_s}^{-1/2})
	(\hat{\theta}_{\ell_1}-\theta_{\ell_1},\ldots,\hat{\theta}_{\ell_s}-\theta_{\ell_s})^\T\rightarrow\mathcal{N}(\bzero,\bI_s)$$ in distribution
	if $\gamma\gg p^{-1/4}\log^{1/4}p$,
and
$$
    	p\,{\rm diag}(\hat{\tilde{b}}_{\ell_1}^{-1/2},\ldots, \hat{\tilde{b}}_{\ell_s}^{-1/2})(\hat{\theta}_{\ell_1}-\theta_{\ell_1},\ldots,\hat{\theta}_{\ell_s}-\theta_{\ell_s})^\T\rightarrow\mathcal{N}(\bzero,\bI_s)$$ in distribution
if $p^{-1/3}\log^{1/6}p\ll \gamma\ll p^{-1/4}$.
Unfortunately, such plug-in method does not work in the scenario $p^{-1/4}\lesssim \gamma\lesssim p^{-1/4}\log^{1/4}p$ since $\hat{b}_{\ell}$ is no longer a valid estimate for $b_\ell$. On the other hand, it is difficult to judge which regime the decay rate of $\gamma$ falls into in practice with finite samples. Hence, the plug-in method is powerless practically.

\end{appendix}

\begin{acks}[Acknowledgments]
The authors are grateful to the editor, an associate editor and three referees for their helpful suggestions. The authors also thank Vishesh Karwa for sharing code for differentially maximum likelihood estimation in the $\beta$-model.
\end{acks}
\begin{funding}
Chang and Hu were supported in part
				by the National Natural Science Foundation of China (grant nos.~71991472, 72125008
				and 12326360). Chang was also supported by  the Center of Statistical Research at Southwestern University of
				Finance and Economics.  Kolaczyk was supported in part by U.S. National Science Foundation award SES-2120115 and Canadian NSERC RGPIN-2023-03566 and NSERC DGDND-2023-03566 awards. Yao was	supported in part by the U.K. Engineering and Physical Sciences Research
				Council (grant nos. EP/V007556/1 and EP/X002195/1). Yi was supported in part by the National Key R\&D Program of China (grant no. 2022YFA1003701), the National Natural Science Foundation of China (grant nos. 12071416, 12271472 and 12301363),
the Fundamental Research Funds for the Central Universities (grant no. JBK2101013), and the China Postdoctoral Science Foundation (grant no. 2021M692663).
\end{funding}

\begin{supplement}
\stitle{Supplement to ``Edge Differentially Private Estimation in the $\beta$-model via Jittering and Method of Moments''}
\sdescription{This supplement contains all the technical proofs of the theoretical results in this paper.}
\end{supplement}


\clearpage
\setcounter{equation}{0}
\setcounter{table}{0}
\setcounter{lemma}{0}
\newtheorem{props}{Proposition}
\newtheorem{inequ}{Inequality}
\setcounter{props}{0}
\setcounter{page}{1}
\setcounter{section}{0}
\renewcommand{\theequation}{S.\thesection.\arabic{equation}}
\renewcommand{\theprops}{P\arabic{props}}
\renewcommand{\thepage}{S\arabic{page}}
\renewcommand{\thetable}{S\arabic{table}}

	\begin{frontmatter}
		\title{Supplement to ``Edge Differentially Private Estimation in the $\beta$-model via Jittering and Method of Moments''}
		\runtitle{Supplement to ``Edge Differentially Private Estimation''}
		
		\begin{aug}
			\author[A,B]{\fnms{Jinyuan} \snm{Chang}},
\author[A]{\fnms{Qiao} \snm{Hu}\ead[label=e2,mark]{huqiao@smail.swufe.edu.cn}},
\author[C]{\fnms{Eric D.} \snm{Kolaczyk}},
\author[D]{\fnms{Qiwei} \snm{Yao}}
\and
\author[E]{\fnms{Fengting} \snm{Yi}}
\address[A]{Joint Laboratory of Data Science and Business Intelligence, Southwestern University of Finance and Economics, Chengdu, 611130, China, \printead{e1,e2}}

\address[B]{Academy of
Mathematics and Systems Science, Chinese Academy of Sciences, Beijing, 100190, China}

\address[C]{Department of Mathematics and Statistics, McGill University, Montreal, QC, H3A 0B8, Canada, \printead{e3}}

\address[D]{Department of Statistics, London School of Economics and Political Science, London, WC2A 2AE, UK, \printead{e4}}

\address[E]{Yunnan Key Laboratory of Statistical Modeling and Data Analysis, Yunnan University, Kunming, 650500, China, \printead{e5}}
	\end{aug}

\end{frontmatter}

	Throughout the supplementary material, we use $C$ and $\tilde{C}$ to denote generic positive finite universal constants that may be different in different uses. The following two inequalities will be used in our proofs.

\begin{inequ}[Decoupling inequality, Theorem 1 of \cite{DM_1995}] Let $\{U_i\}$ be a sequence of independent random variables on a measurable space $(\mathcal{J},S)$ and let $\{U_{i}^{(j)}\}$, $j\in[k]$, be $k$ independent copies of $\{U_i\}$. Let $f_{i_1,\ldots,i_k}$ be families of functions of $k$ variables taking $(S\times \cdots \times S)$ into a Banach space $(\mathcal{B},\|\cdot\|)$. Then, for all $n\geq k\geq 2$ and $t>0$, there exists a numerical constant $C_k$ depending on $k$ only so that
\begin{align*}
&\mathbb{P}\bigg\{\bigg\|\sum_{1\leq i_1\neq\cdots\neq i_k\leq n}f_{i_1,\ldots,i_k}(U_{i_1}^{(1)},U_{i_2}^{(1)},\ldots,U_{i_k}^{(1)})\bigg\|\geq t\bigg\}\\
&~~~~~~~~~~~~~\leq C_k\mathbb{P}\bigg\{C_k\bigg\|\sum_{1\leq i_1\neq\cdots\neq i_k\leq n}f_{i_1,\ldots,i_k}(U_{i_1}^{(1)},U_{i_2}^{(2)},\ldots,U_{i_k}^{(k)})\bigg\|\geq t\bigg\}\,.
\end{align*}
\end{inequ}

\begin{inequ}[Theorem 3.3 of \cite{Gineetal_2000}]
There exists a universal constant $L>0$ such that, if $h_{i,j}$ are bounded canonical kernels of two variables for independent random variables $U_i^{(1)}$ and $U_j^{(2)}$, $i,j\in[n]$, then
\begin{align*}
&\mathbb{P}\bigg\{\bigg|\sum_{i,j=1}^nh_{i,j}(U_i^{(1)},U_j^{(2)})\bigg|\geq t\bigg\}\leq L\exp\bigg\{-\frac{1}{L}\min\bigg(\frac{t^2}{E^2},\frac{t}{D},\frac{t^{2/3}}{B^{2/3}},\frac{t^{1/2}}{A^{1/2}}\bigg)\bigg\}
\end{align*}
for all $t>0$, where
\begin{align*}
&~~~~~~~~~~~~~~~~~~~~A=\max_{i,j\in[n]}\|h_{i,j}\|_\infty\,,~~~E^2=\sum_{i,j=1}^n\mathbb{E}\{h_{i,j}^2(U_i^{(1)},U_j^{(2)})\}\,,\\
&~~~~B^2=\max_{i,j\in[n]}\bigg[\bigg\|\sum_{i=1}^n\mathbb{E}_1\{h_{i,j}^2(U_i^{(1)},y)\}\bigg\|_\infty,~\bigg\|\sum_{j=1}^n\mathbb{E}_2\{h_{i,j}^2(x,U_j^{(2)})\}\bigg\|_\infty\bigg]\,,\\
&D=\sup\bigg[\mathbb{E}\bigg\{\sum_{i,j=1}^nh_{i,j}(U_i^{(1)},U_j^{(2)})f_i(U_i^{(1)})g_j(U_i^{(2)})\bigg\}:\\
&~~~~~~~~~~~~~~~~~~~~~~~~~~~~~~~~~~~~~~~~~~~\mathbb{E}\bigg\{\sum_{i=1}^nf_i^2(U_i^{(1)})\bigg\}\leq 1,~ \mathbb{E}\bigg\{\sum_{j=1}^ng_j^2(U_j^{(2)})\bigg\}\leq 1\bigg]\,.
\end{align*}
\end{inequ}

	\section{Proof of Proposition \ref{tm:1}.}\label{appB}


	 Define
	\begin{equation}\label{eq:hatzetaj}
		\hat{\zeta}_{\ell}=\frac{\hat{\mu}_{\ell,1}}{\hat{\mu}_{\ell,2}}~~\textrm{and}~~
		\zeta_{\ell}=\frac{\mu_{\ell,1}}{\mu_{\ell,2}}\,.
	\end{equation}
	To prove Proposition \ref{tm:1}, we need the following lemma whose proof is given in Section \ref{sec:pflem1}.
	
	\begin{lemma}\label{la:3}
		Let Condition {\rm\ref{cond1}} hold and $(\alpha,\beta)\in\mathcal{M}(\gamma;C_1)$ for some fixed constant $C_1\in(0,0.5)$. Then
	\[
	\begin{split}
	\max_{k\in\{1,2\}}\max_{\ell\in[p]}|\hat{\mu}_{\ell,k}-\mu_{\ell,k}|=&\,O_{\p}\bigg(\frac{\log^{1/2}p}{p}\bigg)+O_{\p}\bigg(\frac{\gamma^2\log^{1/2}p}{p^{1/2}}\bigg)+O_{\p}\bigg(\frac{\gamma\log p}{p}\bigg)\,.
	\end{split}
	\]
	\end{lemma}

 By Condition \ref{cond1},
 \[
 \min_{\ell \in[p]}\mu_{\ell,1}\asymp\gamma^3\asymp\max_{\ell\in[p]}\mu_{\ell,2}\,.
 \]
If  $\gamma\gg p^{-1/3}\log^{1/6}p$, Lemma \ref{la:3} implies
	$$
	\max_{\ell\in[p]}|\hat{\mu}_{\ell,1}-\mu_{\ell,1}|=o_{\p}(\gamma^3)=\max_{\ell\in[p]}|\hat{\mu}_{\ell,2}-\mu_{\ell,2}|\,.$$
	By (\ref{eq:hatzetaj}), it holds that
	\[
	\begin{split}
\hat{\zeta}_{\ell}-\zeta_{\ell}=&\,\frac{\hat{\mu}_{\ell,1}-\mu_{\ell,1}}{\hat{\mu}_{\ell,2}}-\frac{\mu_{\ell,1}(\hat{\mu}_{\ell,2}-\mu_{\ell,2})}{\hat{\mu}_{\ell,2}\mu_{\ell,2}}\\
=&\,\frac{\hat{\mu}_{\ell,1}-\mu_{\ell,1}}{\mu_{\ell,2}}-\frac{\mu_{\ell,1}}{\mu_{\ell,2}^2}(\hat{\mu}_{\ell,2}-\mu_{\ell,2})+R_{\ell,1}\,,
\end{split}
\]
	where
	$
	\max_{\ell\in[p]}|R_{\ell,1}|=O_{\p}(\gamma^{-6}p^{-2}\log p)+O_{\p}(\gamma^{-2}p^{-1}\log p)$.
	Thus,
	\[
\max_{\ell\in[p]}|\hat{\zeta}_\ell-\zeta_\ell|=O_{\p}\bigg(\frac{\log^{1/2}p}{\gamma^3p}\bigg)+O_{\p}\bigg(\frac{\log^{1/2}p}{\gamma p^{1/2}}\bigg)=o_{\p}(1)\,.
\]
Since $\theta_{\ell}=\log(\zeta_{\ell})/2$ and $\hat{\theta}_{\ell}=\log(\hat{\zeta}_{\ell})/2$, by Taylor expansion, we have that
	\begin{align}\label{eq:thetaexp}
		\hat{\theta}_{\ell}-\theta_{\ell}=\frac{1}{2\zeta_{\ell}}(\hat{\zeta}_{\ell}-\zeta_{\ell})+R_{\ell,2}=\frac{\hat{\mu}_{\ell,1}-\mu_{\ell,1}}{2\mu_{\ell,1}}-\frac{\hat{\mu}_{\ell,2}-\mu_{\ell,2}}{2\mu_{\ell,2}}+R_{\ell,3}\,,
	\end{align}
	where
	\[
\max_{\ell\in[p]}|R_{\ell,2}|=O_{\p}\bigg(\frac{\log p}{\gamma^6p^2}\bigg)+O_{\p}\bigg(\frac{\log p}{\gamma^2p}\bigg)=\max_{\ell\in[p]}|R_{\ell,3}|\,.
\]
Therefore,
\[
\max_{\ell\in[p]}|\hat{\theta}_\ell-\theta_\ell|=O_{\p}\bigg(\frac{\log^{1/2}p}{\gamma^3p}\bigg)+O_{\p}\bigg(\frac{\log^{1/2}p}{\gamma p^{1/2}}\bigg)=o_{\p}(1)\,.
\]
	We complete the proof of Proposition \ref{tm:1}. $\hfill\Box$

	\section{Proof of Theorem \ref{tm:asymnormal}.}\label{sec:pfthm1}
For any $i,j,\ell\in[p]$, let  $\psi_1(i,j;\ell)=\varphi_{(i,\ell),1}\varphi_{(i,j),0}\varphi_{(\ell,j),1}$ and $\psi_2(i,j;\ell)=\varphi_{(i,\ell),0}\varphi_{(i,j),1}\varphi_{(\ell,j),0}$. Note that $\gamma\gg p^{-1/3}\log^{1/6}p$. Write $N=(p-1)(p-2)$ and $\mathring{\varphi}_{(i,j),\tau}=\varphi_{(i,j),\tau}-\mathbb{E}\{\varphi_{(i,j),\tau}\}$. As shown in Section \ref{sec:pflem1} for the proof of Lemma \ref{la:3}, 
	\begin{align*}
		\hat{\mu}_{\ell,1}-\mu_{\ell,1}=&~\underbrace{\frac{1}{|\mathcal{H}_\ell|}\sum_{(i,j)\in\mathcal{H}_\ell}[\psi_1(i,j;\ell)-\mathbb{E}\{\psi_1(i,j;\ell)\,|\,\mathscr{F}_\ell\}]}_{I_{\ell,1,1}}\\
		&+\underbrace{\frac{2}{N}\sum_{i,j:\,i\neq j,\,i,j\neq \ell}\mathring{\varphi}_{(i,\ell),1}\mathbb{E}\{\varphi_{(\ell,j),1}\}\mathbb{E}\{\varphi_{(i,j),0}\}}_{N^{-1}I_{\ell,1,2}(1)}\\
		&+\underbrace{\frac{1}{N}\sum_{i,j:\,i\neq j,\,i,j\neq \ell}\mathring{\varphi}_{(i,\ell),1}\mathring{\varphi}_{(\ell,j),1}\mathbb{E}\{\varphi_{(i,j),0}\}}_{N^{-1}I_{\ell,1,2}(2)}\,,\\
		\hat{\mu}_{\ell,2}-\mu_{\ell,2}=&~\underbrace{\frac{1}{|\mathcal{H}_\ell|}\sum_{(i,j)\in\mathcal{H}_\ell}[\psi_2(i,j;\ell)-\mathbb{E}\{\psi_2(i,j;\ell)\,|\,\mathscr{F}_\ell\}]}_{I_{\ell,2,1}}\\
		&+\underbrace{\frac{2}{N}\sum_{i,j:\,i\neq j,\,i,j\neq \ell}\mathring{\varphi}_{(i,\ell),0}\mathbb{E}\{\varphi_{(\ell,j),0}\}\mathbb{E}\{\varphi_{(i,j),1}\}}_{N^{-1}I_{\ell,2,2}(1)}\\
		&+\underbrace{\frac{1}{N}\sum_{i,j:\,i\neq j,\,i,j\neq \ell}\mathring{\varphi}_{(i,\ell),0}\mathring{\varphi}_{(\ell,j),0}\mathbb{E}\{\varphi_{(i,j),1}\}}_{N^{-1}I_{\ell,2,2}(2)}\,,
	\end{align*}
where $\mathscr{F}_\ell=\{Z_{i,\ell},Z_{\ell,j}:(i,j)\in\mathcal{H}_\ell\}$. For each given $\ell\in[p]$, following the proof of Lemma \ref{la:3} in Section \ref{sec:pflem1}, we know that
$I_{\ell,1,1}=O_{\p}(p^{-1})=I_{\ell,2,1}$, $I_{\ell,1,2}(1)=O_{\p}(\gamma^2p^{3/2})=I_{\ell,2,2}(1)$ and $I_{\ell,1,2}(2)=O_{\p}(\gamma p)=I_{\ell,2,2}(2)$.
Also, the remainder term $R_{\ell,3}$ in \eqref{eq:thetaexp} satisfies $R_{\ell,3}=O_{\p}(\gamma^{-6}p^{-2})+O_{\p}(\gamma^{-2}p^{-1})$ for each given $\ell\in[p]$. By \eqref{eq:thetaexp}, we have
	\begin{align}\label{eq:asymexp}	\hat{\theta}_\ell-\theta_\ell=\underbrace{\frac{I_{\ell,1,1}}{2\mu_{\ell,1}}-\frac{I_{\ell,2,1}}{2\mu_{\ell,2}}}_{T_{\ell,1}}+\underbrace{\frac{I_{\ell,1,2}(1)}{2\mu_{\ell,1}N}-\frac{I_{\ell,2,2}(1)}{2\mu_{\ell,2}N}}_{T_{\ell,2}}+R_{\ell,4}\,,
	\end{align}
 where $R_{\ell,4}=O_{\p}(\gamma^{-6}p^{-2})+O_{\p}(\gamma^{-2}p^{-1})$. Write $\mathring{Z}_{i,j}=Z_{i,j}-\mathbb{E}(Z_{i,j})$. Then $\mathring{\varphi}_{(i,j),1}=\mathring{Z}_{i,j}$ and $\mathring{\varphi}_{(i,j),0}=-\mathring{Z}_{i,j}$. It holds that
 \begin{align}
 &T_{\ell,1}=-\frac{1}{N}\sum_{i,j:\,i\neq j,\,i,j\neq\ell}\bigg\{\frac{\varphi_{(i,\ell),1}\varphi_{(\ell,j),1}}{2\mu_{\ell,1}}+\frac{\varphi_{(i,\ell),0}\varphi_{(\ell,j),0}}{2\mu_{\ell,2}}\bigg\}\mathring{Z}_{i,j}\,,\label{eq:Tell1}\\
 &~~~~~~~~~~~~~~~~~~~~~~~~~~~T_{\ell,2}=\frac{1}{p-1}\sum_{i:\,i\neq\ell}\lambda_{i,\ell}\mathring{Z}_{i,\ell}\,,\label{eq:Tell2}
 \end{align}
 where $\lambda_{i,\ell}$ is defined as \eqref{eq:lambdaiell}. In Sections \ref{sec:case1}--\ref{sec:case3}, we will prove Theorem \ref{tm:asymnormal}(a)--\ref{tm:asymnormal}(c) based on \eqref{eq:asymexp}--\eqref{eq:Tell2}, respectively.

	\subsection{Case 1: $\gamma\gg p^{-1/4}$.}\label{sec:case1}
	
	Due to $T_{\ell,1}=O_{\p}(\gamma^{-3}p^{-1})$, by \eqref{eq:asymexp} and \eqref{eq:Tell2},
	\begin{align*}
\hat{\theta}_\ell-\theta_\ell=&~T_{\ell,2}+R_{\ell,5}=\frac{1}{p-1}\sum_{i:\,i\neq\ell}\lambda_{i,\ell}\mathring{Z}_{i,\ell}+R_{\ell,5}\,,
\end{align*}
	where $R_{\ell,5}=O_{\p}(\gamma^{-3}p^{-1})$. Write $\lambda_{i,\ell}^*=\gamma\lambda_{i,\ell}$.
	Under Condition \ref{cond1}, $\min_{\ell\in[p]}\min_{i:\,i\neq\ell}\lambda_{i,\ell}^*\asymp1\asymp\max_{\ell\in[p]}\max_{i:\,i\neq\ell}\lambda_{i,\ell}^*$. Then
	$$
	\gamma(\hat{\theta}_\ell-\theta_\ell)=\frac{1}{p-1}\sum_{i:\,i\neq \ell}\lambda_{i,\ell}^*\mathring{Z}_{i,\ell}+R_{\ell,6}$$
	with $R_{\ell,6}=O_{\p}(\gamma^{-2}p^{-1})$.
	Given $s$ different $\ell_1,\ldots,\ell_s\in[p]$, we define an $s$-dimensional vector $\bw_i=(W_{i,1},\ldots,W_{i,s})^\T$ with $W_{i,j}=\lambda_{i,\ell_j}^*\mathring{Z}_{i,\ell_j}$. Then it holds that
	\begin{align}\label{eq:wi}
	    \gamma(\hat{\theta}_{\ell_1}-\theta_{\ell_1},\ldots,\hat{\theta}_{\ell_s}-\theta_{\ell_s})^\T=\frac{1}{p-1}\sum_{i:\,i\neq \ell_1,\ldots,\ell_s}\bw_i+\br
	\end{align}
	 with $|\br|_\infty=O_{\p}(\gamma^{-2}p^{-1})$. Let $
	b_{\ell,*}=(p-1)^{-1}\sum_{i:\,i\neq \ell}\lambda_{i,\ell}^{*,2}{\rm Var}(Z_{i,\ell})$. Notice that
	 \begin{align*}
	   \frac{1}{p-1
}\sum_{i:\,i\neq\ell_1,\ldots,\ell_s}{\rm Var}(\bw_i)={\rm diag}(b_{\ell_1,*},\ldots,b_{\ell_s,*})+O(p^{-1})\,,
	 \end{align*}
	  and $\{\bw_i\}_{i\neq \ell_1,\ldots,\ell_s}$ is an independent sequence. By the Central Limit Theorem,
	 \begin{align*}
	     	{(p-1)}^{1/2}\gamma\,{\rm diag}(b_{\ell_1,*}^{-1/2},\ldots,b_{\ell_s,*}^{-1/2})(\hat{\theta}_{\ell_s}-\theta_{\ell_1},\ldots,\hat{\theta}_{\ell_s}-\theta_{\ell_s})^\T\xrightarrow{d}\mathcal{N}(\bzero,\bI_s)\,.
	 \end{align*}
	Due to $b_{\ell,*}=\gamma^2b_\ell$, we complete the proof of Theorem \ref{tm:asymnormal}(a). $\hfill\Box$
	
	\subsection{Case 2: $p^{-1/4}\gg \gamma\gg p^{-1/3}\log^{1/6}p$.}\label{se:case2}
	Due to $T_{\ell,2}=O_{\p}(\gamma^{-1}p^{-1/2})$, by \eqref{eq:asymexp},  
	\begin{align}\label{eq:thetaexpansion}
		\hat{\theta}_\ell-\theta_\ell&=T_{\ell,1}+R_{\ell,7}
	\end{align}
	with $R_{\ell,7}=O_{\p}(\gamma^{-6}p^{-2})+O_{\p}(\gamma^{-1}p^{-1/2})$. Recall $\mathring{\varphi}_{(i,j),1}=\mathring{Z}_{i,j}$ and $\mathring{\varphi}_{(i,j),0}=-\mathring{Z}_{i,j}$. By \eqref{eq:Tell1}, we have
	\begin{align}\label{eq:jell}
		T_{\ell,1}=&-\frac{1}{N}\sum_{i,j:\,i\neq j,\,i,j\neq\ell}\bigg(\frac{1}{2\mu_{\ell,1}}+\frac{1}{2\mu_{\ell,2}}\bigg)\mathring{Z}_{i,\ell}\mathring{Z}_{\ell,j}\mathring{Z}_{i,j}\notag\\
		&-\underbrace{\frac{2}{N}\sum_{i,j:\,i\neq j,\,i,j\neq\ell}\bigg[\frac{\mathbb{E}\{\varphi_{(\ell,j),1}\}}{2\mu_{\ell,1}}-\frac{\mathbb{E}\{\varphi_{(\ell,j),0}\}}{2\mu_{\ell,2}}\bigg]\mathring{Z}_{i,\ell}\mathring{Z}_{i,j}}_{J_{\ell,1}}\\
		&-\underbrace{\frac{1}{N}\sum_{i,j:\,i\neq j,\,i,j\neq\ell}\bigg[\frac{\mathbb{E}\{\varphi_{(i,\ell),1}\}\mathbb{E}\{\varphi_{(\ell,j),1}\}}{2\mu_{\ell,1}}+\frac{\mathbb{E}\{\varphi_{(i,\ell),0}\}\mathbb{E}\{\varphi_{(\ell,j),0}\}}{2\mu_{\ell,2}}\bigg]\mathring{Z}_{i,j}}_{J_{\ell,2}}\notag\,.
	\end{align}
	Due to $\mathring{Z}_{i,j}=\mathring{Z}_{j,i}$ and $\mathbb{E}\{\varphi_{(i,j),\tau}\}=\mathbb{E}\{\varphi_{(j,i),\tau}\}$, we have
\[
J_{\ell,2}=\frac{1}{|\mathcal{H}_\ell|}\sum_{(i,j)\in\mathcal{H}_\ell}\bigg[\frac{\mathbb{E}\{\varphi_{(i,\ell),1}\}\mathbb{E}\{\varphi_{(\ell,j),1}\}}{2\mu_{\ell,1}}+\frac{\mathbb{E}\{\varphi_{(i,\ell),0}\}\mathbb{E}\{\varphi_{(\ell,j),0}\}}{2\mu_{\ell,2}}\bigg]\mathring{Z}_{i,j}\,.
\]
Under Condition \ref{cond1},
	\begin{align*}
	&\min_{\ell \in [p]}\min_{i,j:\,i\neq j,\,i,j\neq \ell}\bigg[\frac{\mathbb{E}\{\varphi_{(i,\ell),1}\}\mathbb{E}\{\varphi_{(\ell,j),1}\}}{2\mu_{\ell,1}}+\frac{\mathbb{E}\{\varphi_{(i,\ell),0}\}\mathbb{E}\{\varphi_{(\ell,j),0}\}}{2\mu_{\ell,2}}\bigg]\\
	&~~~~~\asymp \gamma^{-1}\asymp \max_{\ell\in[p]}\max_{i,j:\,i\neq j,\,i,j\neq \ell}\bigg[\frac{\mathbb{E}\{\varphi_{(i,\ell),1}\}\mathbb{E}\{\varphi_{(\ell,j),1}\}}{2\mu_{\ell,1}}+\frac{\mathbb{E}\{\varphi_{(i,\ell),0}\}\mathbb{E}\{\varphi_{(\ell,j),0}\}}{2\mu_{\ell,2}}\bigg]\,.
	\end{align*}
Notice that $\{\mathring{Z}_{i,j}\}_{(i,j)\in\mathcal{H}_\ell}$ is an independent sequence with $|\mathcal{H}_\ell|=(p-1)(p-2)/2$. By Bernstein inequality, we have
	\begin{equation}\label{eq:convJell2}
	|J_{\ell,2}|=O_{\p}(\gamma^{-1}p^{-1})\,.
\end{equation}
For $J_{\ell,1}$, we can reformulate it as follows:
	\begin{align}\label{eq:jell1}
		J_{\ell,1}=&\underbrace{\frac{1}{p-1}\sum_{i:\,i\neq\ell}\mathring{Z}_{i,\ell}\bigg[\frac{1}{p-2}\sum_{j:\,j\neq i,\ell}\frac{\mathbb{E}\{\varphi_{(\ell,j),1}\}}{\mu_{\ell,1}}\mathring{Z}_{i,j}\bigg]}_{J_{\ell,1}(1)}\notag\\
		&-\underbrace{\frac{1}{p-1}\sum_{i:\,i\neq\ell}\mathring{Z}_{i,\ell}\bigg[\frac{1}{p-2}\sum_{j:\,j\neq i,\ell}\frac{\mathbb{E}\{\varphi_{(\ell,j),0}\}}{\mu_{\ell,2}}\mathring{Z}_{i,j}\bigg]}_{J_{\ell,1}(2)}\,.
	\end{align}
	Let $\tilde{A}_{i,\ell} = (p-2)^{-1}\sum_{j:\,j\neq i,\ell}\mathbb{E}\{\varphi_{(\ell,j),1}\}\mu_{\ell,1}^{-1}\mathring{Z}_{i,j}$.
	Due to
	\[
\min_{\ell \in [p]}\min_{j:\,j\neq\ell}\frac{\mathbb{E}\{\varphi_{(\ell,j),1}\}}{\mu_{\ell,1}}\asymp\gamma^{-2}\asymp \max_{\ell\in[p]}\max_{j:\,j\neq\ell}\frac{\mathbb{E}\{\varphi_{(\ell,j),1}\}}{\mu_{\ell,1}}
\]
under Condition \ref{cond1}, by Bernstein inequality,
	\begin{equation}\label{eq:bdAiell}
\mathbb{P}(|\tilde{A}_{i,\ell}|>u)\lesssim\exp(-{C\gamma^4pu^2})
\end{equation}
	for any $u=o(1)$. Given a sufficiently large constant $C_*>0$, define
	\[
\mathcal{E}_\ell(C_*)=\bigg\{\max_{i:\,i\neq\ell}|\tilde{A}_{i,\ell}|\leq\frac{C_*\log^{1/2}p}{\gamma^{2}p^{1/2}}\bigg\}\,.
\]
	By \eqref{eq:bdAiell}, we have
	\begin{align}\label{eq:bj1}
		\mathbb{P}\{|J_{\ell,1}(1)|>u\}\leq&~ \mathbb{P}\bigg\{\bigg|\frac{1}{p-1}\sum_{i:\,i\neq\ell}\mathring{Z}_{i,\ell}\tilde{A}_{i,\ell}\bigg|>u\,,\mathcal{E}_\ell(C_*)\bigg\}+p^{-\bar{C}}\,.
	\end{align}
	The constant $\bar{C}>0$ in \eqref{eq:bj1} can be sufficiently large if we select a sufficiently large constant $C_*$. Write $\mathcal{F}_{-\ell}=\{Z_{i,j}:i,j\neq \ell\}$.
	Conditional on $\mathcal{F}_{-\ell}$, by Bernstein inequality,
	\begin{align*}
&\mathbb{P}\bigg\{\bigg|\frac{1}{p-1}\sum_{i:\,i\neq\ell}\mathring{Z}_{i,\ell}\tilde{A}_{i,\ell}\bigg|>u\,,\mathcal{E}_\ell(C_*)\,\bigg|\,\mathcal{F}_{-\ell}\bigg\}\\
&~~~~~~~~~~~~~~\lesssim \exp\bigg(-\frac{Cp^2u^2}{\tilde{C}\sum_{i:\,i\neq\ell}\tilde{A}_{i,\ell}^2+pu\max_{i:\,i\neq\ell}|\tilde{A}_{i,\ell}|}\bigg)I\bigg(\max_{i:\,i\neq\ell}|\tilde{A}_{i,\ell}|\leq \frac{C_*\log^{1/2}p}{\gamma^{2}p^{1/2}}\bigg)
\end{align*}
	for any $u>0$. Selecting $u=C_{**}\gamma^{-2}p^{-1}\log p$, we have
	\begin{align*}
		&\mathbb{P}\bigg\{\bigg|\frac{1}{p-1}\sum_{i:\,i\neq\ell}\mathring{Z}_{i,\ell}\tilde{A}_{i,\ell}\bigg|>\frac{C_{**}\log p}{\gamma^{2}p}\,,\mathcal{E}_\ell(C_*)\,\bigg|\,\mathcal{F}_{-\ell}\bigg\}\\
		&~~~~~~~~~~~~~~\lesssim p^{-\check{C}}\cdot I\bigg(\max_{i:\,i\neq\ell}|\tilde{A}_{i,\ell}|\leq \frac{C_*\log^{1/2}p}{\gamma^{2}p^{1/2}}\bigg)\leq p^{-\check{C}}\,,
	\end{align*}
	which implies
	\begin{align*}
	    	\mathbb{P}\bigg\{\bigg|\frac{1}{p-1}\sum_{i:\,i\neq\ell}\mathring{Z}_{i,\ell}\tilde{A}_{i,\ell}\bigg|>\frac{C_{**}\log p}{\gamma^{2}p}\,,\mathcal{E}_\ell(C_*)\bigg\} \lesssim p^{-\check{C}}\rightarrow0\,.
	\end{align*}
	Here the constant $\check{C}>0$ can be sufficiently large if we select a sufficiently large constant $C_{**}$. Thus, \eqref{eq:bj1} implies $\max_{\ell\in[p]}|J_{\ell,1}(1)|=O_{\p}(\gamma^{-2}p^{-1}\log p)$. Analogously, we also have $\max_{\ell\in[p]}|J_{\ell,1}(2)|=O_{\p}(\gamma^{-2}p^{-1}\log p)$. By \eqref{eq:jell1},
\begin{align}\label{eq:bdJell1}
\max_{\ell\in[p]}|J_{\ell,1}|=O_{\p}(\gamma^{-2}p^{-1}\log p)\,.
 \end{align}
Together with \eqref{eq:convJell2}, by \eqref{eq:thetaexpansion} and \eqref{eq:jell}, we have
	\begin{align*}
	    	\hat{\theta}_{\ell}-\theta_\ell=-\frac{1}{N}\sum_{i,j:\,i\neq j,\,i,j\neq\ell}\bigg(\frac{1}{2\mu_{\ell,1}}+\frac{1}{2\mu_{\ell,2}}
	    	\bigg)\mathring{Z}_{i,\ell}\mathring{Z}_{\ell,j}\mathring{Z}_{i,j}+R_{\ell,8}
	\end{align*}
	with $R_{\ell,8}=O_{\p}(\gamma^{-6}p^{-2})+O_{\p}(\gamma^{-1}p^{-1/2})$, which implies
	\begin{align}\label{eq:expcase2}
	    -\frac{2\mu_{\ell,1}\mu_{\ell,2}}{{\mu_{\ell,1}+\mu_{\ell,2}}}(\hat{\theta}_\ell-\theta_\ell) = \underbrace{\frac{1}{N}\sum_{i,j:\,i\neq j,\,i,j\neq\ell}\mathring{Z}_{i,\ell}\mathring{Z}_{\ell,j}\mathring{Z}_{i,j}}_{\Delta_\ell}+R_{\ell,9}
	\end{align}
	with $R_{\ell,9}=O_{\p}(\gamma^{-3}p^{-2})+O_{\p}(\gamma^2 p^{-1/2})$.
	
In the sequel, we will specify the limiting distribution of $\sqrt{N}(\Delta_{\ell_1},\ldots,\Delta_{\ell_s})$ for given $s$ different $\ell_1,\ldots,\ell_s\in[p]$. For given $k\in[s]$, we have
	\begin{equation}\label{eq:theasymp2}
		\begin{split}
			\Delta_{\ell_k}=&\underbrace{\frac{1}{N}\sum_{i,j:\,i\neq j,\atop i,j\neq\ell_1,\ldots,\ell_s}\mathring{Z}_{i,\ell_k}\mathring{Z}_{\ell_k,j}\mathring{Z}_{i,j}}_{M_{\ell_k,1}}+\underbrace{\frac{1}{N}\sum_{i,j:\,i\neq j,\,i,j\neq \ell_k,\atop \{i,j\}\cap\{\ell_1,\ldots,\ell_{k-1},\ell_{k+1},\ldots,\ell_s\}\neq\emptyset}\mathring{Z}_{i,\ell_k}\mathring{Z}_{\ell_k,j}\mathring{Z}_{i,j}}_{M_{\ell_k,2}}\,.
		\end{split}
	\end{equation}
	Notice that
	\begin{align*}
	    	 M_{\ell_k,2}
	    	=&\,\frac{2}{N}\sum_{k':\,k'\neq k}\sum_{i:\,i\neq\ell_1,\ldots,\ell_s}\mathring{Z}_{i,\ell_{k}}\mathring{Z}_{\ell_k,\ell_{k'}}\mathring{Z}_{i,\ell_{k'}}\\
&+\frac{1}{N}\sum_{k',k'':\,k'\neq k'',\,k',k''\neq k}\mathring{Z}_{\ell_{k'},\ell_k}\mathring{Z}_{\ell_k,\ell_{k''}}\mathring{Z}_{\ell_{k'},\ell_{k''}}\,.
	\end{align*}
	Since
	$
	\max_{k\in[s]}\max_{k',k'':\,k'\neq k'', \, k',k''\neq k}|\mathring{Z}_{\ell_{k'},\ell_k}\mathring{Z}_{\ell_k,\ell_{k''}}\mathring{Z}_{\ell_{k'},\ell_{k''}}|\lesssim 1
	$,
	it holds that
	\begin{align*}
	    M_{\ell_k,2}=2\sum_{k':\,k'\neq k}\mathring{Z}_{\ell_k,\ell_{k'}}\bigg(\frac{1}{N}\sum_{i:\,i\neq\ell_1,\ldots,\ell_s}\mathring{Z}_{i,\ell_{k}}\mathring{Z}_{i,\ell_{k'}}\bigg)+R_{\ell_k,10}
	\end{align*}
	with $R_{\ell_k,10}=O(p^{-2
})$. For given $k'$ such that $k'\neq k$, since $\{\mathring{Z}_{i,\ell_{k}}\mathring{Z}_{i,\ell_{k'}}\}_{i\neq\ell_1,\ldots,\ell_s}$ is an independent sequence, by Bernstein inequality, we have
	\begin{align*}
	    \mathbb{P}\bigg(\frac{1}{p-s}\bigg|\sum_{i:\,i\neq\ell_1,\ldots,\ell_s}\mathring{Z}_{i,\ell_{k}}\mathring{Z}_{i,\ell_{k'}}\bigg|>u \bigg)\lesssim \exp(-{Cpu^2})
	\end{align*}
	for any $u=o(1)$. Therefore, it holds that
	\[
\max_{k',k:\,k'\neq k}\bigg|\frac{1}{N}\sum_{i:\,i\neq\ell_1,\ldots,\ell_s}\mathring{Z}_{i,\ell_{k}}\mathring{Z}_{i,\ell_{k'}}\bigg|=O_{\p}(p^{-3/2})\,,
\]
	which implies
	$
	M_{\ell_{k},2}=O_{\p}(p^{-3/2})$.
	By \eqref{eq:theasymp2},
	\begin{align}\label{eq:expcase21}
	\Delta_{\ell_k}=M_{\ell_k,1}+R_{\ell_k,11}
\end{align}
	with $R_{\ell_k,11}=O_{\p}(p^{-3/2})$. Given $\ell_1,\ldots,\ell_s$, write
	\begin{align}\label{eq:bellk**}
b_{\ell_k,**}=&\,\frac{1}{N}\sum_{i,j:\,i\neq j,\atop i,j\neq \ell_1,\ldots,\ell_s}{\rm Var}(Z_{i,\ell_k}){\rm Var}(Z_{\ell_k,j}){\rm Var}(Z_{i,j})\,,\\
f(t_1,\ldots,t_s)=&\,\mathbb{E}\bigg\{\exp\bigg(\iota \sqrt{N}\sum_{k=1}^st_kb_{\ell_k,**}^{-1/2}M_{\ell_k,1}\bigg)\bigg\}~~\textrm{with}~~\iota^2=-1\,.\notag
 \end{align}
 Let $\mathcal{F}_{\ell_1,\ldots,\ell_s}^*=\cup_{k=1}^s\{Z_{i,\ell_k},Z_{\ell_k,j}:i,j\neq\ell_k\}$. Then
	\begin{align}
		&\mathbb{E}\bigg\{\exp\bigg(\iota \sqrt{N}\sum_{k=1}^st_kb_{\ell_k,**}^{-1/2}M_{\ell_k,1}\bigg)\,\bigg|\,\mathcal{F}_{\ell_1,\ldots,\ell_s}^*\bigg\}\\
		&~~~~~~~~~~~~=\mathbb{E}\bigg[\exp\bigg\{\sum_{i,j:\,i<j,\atop i,j\neq \ell_1,\ldots,\ell_s}\iota\bigg(\sum_{k=1}^s\frac{2t_kb_{\ell_k,**}^{-1/2}}{\sqrt{N}}\mathring{Z}_{i,\ell_k}\mathring{Z}_{\ell_k,j}\bigg)\mathring{Z}_{i,j}\bigg\}\,\bigg|\,\mathcal{F}_{\ell_1,\ldots,\ell_s}^*\bigg]\\
		&~~~~~~~~~~~~=\prod_{i,j:\,i<j,\atop i,j\neq \ell_1,\ldots,\ell_s}\mathbb{E}\bigg[\exp\bigg\{\iota\bigg(\sum_{k=1}^s\frac{2t_kb_{\ell_k,**}^{-1/2}}{\sqrt{N}}\mathring{Z}_{i,\ell_k}\mathring{Z}_{\ell_k,j}\bigg)\mathring{Z}_{i,j}\bigg\}\,\bigg|\,\mathcal{F}_{\ell_1,\ldots,\ell_s}^*\bigg]\,.
	\end{align}
Write $\bar{A}_{i,j} = (\sum_{k=1}^s{2t_kb_{\ell_k,**}^{-1/2}}{{N}^{-1/2}}\mathring{Z}_{i,\ell_k}\mathring{Z}_{\ell_k,j})\mathring{Z}_{i,j}$. By Taylor expansion,
	\begin{align*}
	    \exp (\iota\bar{A}_{i,j})=1+\iota\bar{A}_{i,j}
	    -\frac{1}{2}\bar{A}_{i,j}^2+\tilde{R}_{i,j}
	\end{align*}
	with $|\tilde{R}_{i,j}|\leq CN^{-3/2}(|t_1|+\cdots+|t_s|)^3$, which implies
	\begin{align*}
	  	\mathbb{E}\{\exp (\iota\bar{A}_{i,j})|\,\mathcal{F}_{\ell_1,\ldots,\ell_s}^*\}
	=&1-\frac{1}{2}\bigg(\sum_{k=1}^s\frac{2t_kb_{\ell_k,**}^{-1/2}}{\sqrt{N}}\mathring{Z}_{i,\ell_k}\mathring{Z}_{\ell_k,j}\bigg)^2{\rm Var}({Z}_{i,j})+\tilde{R}_{i,j}^*
	\end{align*}
with $|\tilde{R}_{i,j}^*|\leq CN^{-3/2}(|t_1|+\cdots+|t_s|)^3$ for any $i,j\neq\ell_1,\ldots,\ell_k$.
	Due to the fact $|\prod_{k=1}^mz_k-\prod_{k=1}^mw_k|\leq \sum_{k=1}^m|z_k-w_k|$ for any $z_k,w_k\in\mathbb{C}$ with $|z_k|\leq1$ and $|w_k|\leq 1$, we have
	\begin{align*}
&\bigg|\mathbb{E}\bigg\{\exp\bigg(\iota \sqrt{N}\sum_{k=1}^st_kb_{\ell_k,**}^{-1/2}M_{\ell_k,1}\bigg)\,\bigg|\,\mathcal{F}_{\ell_1,\ldots,\ell_s}^*\bigg\}\\
&~~~~~~-\prod_{i,j:\,i<j,\atop i,j\neq \ell_1,\ldots,\ell_s}\bigg\{1-\frac{1}{2}\bigg(\sum_{k=1}^s\frac{2t_kb_{\ell_k,**}^{-1/2}}{\sqrt{N}}\mathring{Z}_{i,\ell_k}\mathring{Z}_{\ell_k,j}\bigg)^2{\rm Var}({Z}_{i,j})\bigg\}\bigg|\lesssim
\frac{(|t_1|+\cdots+|t_s|)^3}{\sqrt{N}}\,.
\end{align*}
	It also holds that
\begin{align*}
&\bigg|\prod_{i,j:\,i<j,\atop i,j\neq \ell_1,\ldots,\ell_s}\exp\bigg\{-\frac{1}{2}\bigg(\sum_{k=1}^s\frac{2t_kb_{\ell_k,**}^{-1/2}}{\sqrt{N}}\mathring{Z}_{i,\ell_k}\mathring{Z}_{\ell_k,j}\bigg)^2{\rm Var}({Z}_{i,j})\bigg\}\\
&~~~~~~-\prod_{i,j:\,i<j,\atop i,j\neq \ell_1,\ldots,\ell_s}\bigg\{1-\frac{1}{2}\bigg(\sum_{k=1}^s\frac{2t_kb_{\ell_k,**}^{-1/2}}{\sqrt{N}}\mathring{Z}_{i,\ell_k}\mathring{Z}_{\ell_k,j}\bigg)^2{\rm Var}({Z}_{i,j})\bigg\}\bigg|\lesssim \frac{(|t_1|+\cdots+|t_s|)^4}{N}\,.
\end{align*}
By Triangle inequality,
	\begin{align}
		&\bigg|\exp\bigg\{-\frac{1}{2}\sum_{i,j:\,i<j,\atop i,j\neq \ell_1,\ldots,\ell_s}\bigg(\sum_{k=1}^s
\frac{2t_kb_{\ell_k,**}^{-1/2}}{\sqrt{N}}\mathring{Z}_{i,\ell_k}\mathring{Z}_{\ell_k,j}\bigg)^2{\rm Var}({Z}_{i,j})\bigg\}\notag\\
		&~~~~~~~~-\mathbb{E}\bigg\{\exp\bigg(\iota \sqrt{N}\sum_{k=1}^st_kb_{\ell_k,**}^{-1/2}M_{\ell_k,1}\bigg)\,\bigg|\,\mathcal{F}_{\ell_1,\ldots,\ell_s}^*\bigg\}\bigg|\notag\\
&~~~~~~~\lesssim \frac{(|t_1|+\cdots+|t_s|)^3}{\sqrt{N}}+\frac{(|t_1|+\cdots+|t_s|)^4}{N}\,.\label{eq:chardiff1}
	\end{align}
	Define
	\[
\begin{split}
Q=&\sum_{i,j:\,i\neq j,\atop i,j\neq \ell_1,\ldots,\ell_s}\bigg(\sum_{k=1}^s\frac{t_kb_{\ell_k,**}^{-1/2}}{\sqrt{N}}\mathring{Z}_{i,\ell_k}\mathring{Z}_{\ell_k,j}\bigg)^2{\rm Var}({Z}_{i,j})-\sum_{k=1}^st_k^2\,.
\end{split}
\]
By \eqref{eq:chardiff1}, we have
	\begin{align}
	  f(t_1,\ldots,t_s)=&\,\mathbb{E}\bigg[\mathbb{E}\bigg\{\exp\bigg(\iota \sqrt{N}\sum_{k=1}^st_kb_{\ell_k,**}^{-1/2}M_{\ell_k,1}\bigg)\,\bigg|\,\mathcal{F}_{\ell_1,\ldots,\ell_s}^*\bigg\}\bigg]\notag\\
=&\,\exp\bigg(-\sum_{k=1}^st_k^2\bigg)\mathbb{E}\{\exp(-Q)\}+R(t_1,\ldots,t_s)\label{eq:ft1ts}
	\end{align}
	with $|R(t_1,\ldots,t_s)|\lesssim N^{-1/2}(|t_1|+\cdots+|t_s|)^3+N^{-1}(|t_1|+\cdots+|t_s|)^4$. 
Note that
	\begin{align*}
Q
=&~\sum_{k=1}^s\frac{2t_k^2}{Nb_{\ell_k,**}}\sum_{i:\,i\neq \ell_1,\ldots,\ell_s}\{\mathring{Z}_{i,\ell_k}^2-\mathbb{E}(\mathring{Z}_{i,\ell_k}^2)\}\bigg\{\sum_{j:\,j\neq i, \ell_1,\ldots,\ell_s}\mathbb{E}(\mathring{Z}_{\ell_k,j}^2){\rm Var}(Z_{i,j})\bigg\}\\
&+\sum_{k=1}^s\frac{t_k^2}{Nb_{\ell_k,**}}\sum_{i,j:\,i\neq j,\atop i,j\neq \ell_1,\ldots,\ell_s}\{\mathring{Z}_{i,\ell_k}^2-\mathbb{E}(\mathring{Z}_{i,\ell_k}^2)\}\{\mathring{Z}_{\ell_k,j}^2-\mathbb{E}(\mathring{Z}_{\ell_k,j}^2)\}{\rm Var}(Z_{i,j})\\
&+\sum_{k,k':\,k\neq k'}\frac{t_{k}t_{k'}}{Nb_{\ell_k,**}^{1/2}b_{\ell_{k'},**}^{1/2}}\sum_{i,j:\,i\neq j,\atop i,j\neq \ell_1,\ldots,\ell_s}\mathring{Z}_{i,\ell_k}\mathring{Z}_{i,\ell_{k'}}\mathring{Z}_{\ell_k,j}\mathring{Z}_{\ell_{k'},j}{\rm Var}(Z_{i,j})\,.
\end{align*}
Recall
\[
\min_{i:\,i\neq \ell_1,\ldots,\ell_s}\sum_{j:\,j\neq i, \ell_1,\ldots,\ell_s}\mathbb{E}(\mathring{Z}_{\ell_k,j}^2){\rm Var}(Z_{i,j})\asymp p\asymp \max_{i:\,i\neq \ell_1,\ldots,\ell_s}\sum_{j:\,j\neq i, \ell_1,\ldots,\ell_s}\mathbb{E}(\mathring{Z}_{\ell_k,j}^2){\rm Var}(Z_{i,j})\,.
\]
By Bernstein inequality,
\[
\mathbb{P}\bigg[\bigg|\frac{1}{N}\sum_{i:\,i\neq \ell_1,\ldots,\ell_s}\{\mathring{Z}_{i,\ell_k}^2-\mathbb{E}(\mathring{Z}_{i,\ell_k}^2)\}\bigg\{\sum_{j:\,j\neq i, \ell_1,\ldots,\ell_s}\mathbb{E}(\mathring{Z}_{\ell_k,j}^2){\rm Var}(Z_{i,j})\bigg\}\bigg|>u\bigg]\lesssim \exp(-Cpu^2)
\]
for any $u\rightarrow0$. By the decoupling inequality of \cite{DM_1995} and Theorem 3.3 of \cite{Gineetal_2000}, for any $u\rightarrow0$ but $pu\rightarrow\infty$,
	\begin{align*}
&\mathbb{P}\bigg[\bigg|\frac{1}{N}\sum_{i,j:\,i\neq j,\atop i,j\neq \ell_1,\ldots,\ell_s}\{\mathring{Z}_{i,\ell_k}^2-\mathbb{E}(\mathring{Z}_{i,\ell_k}^2)\}\{\mathring{Z}_{\ell_k,j}^2-\mathbb{E}(\mathring{Z}_{\ell_k,j}^2)\}{\rm Var}(Z_{i,j})\bigg|>u\bigg]\lesssim \exp(-Cpu)\,,\\
&~~~~~~~~~~~~~	    	\mathbb{P}\bigg\{\bigg|\frac{1}{N}\sum_{i,j:\,i\neq j,\atop i,j\neq \ell_1,\ldots,\ell_s}\mathring{Z}_{i,\ell_k}\mathring{Z}_{i,\ell_{k'}}\mathring{Z}_{\ell_k,j}\mathring{Z}_{\ell_{k'},j}{\rm Var}(Z_{i,j})\bigg|>u\bigg\}\lesssim \exp(-Cpu)
	\end{align*}
with $k\neq k'$, which implies
	\begin{align}\label{eq:bdQ}
	\mathbb{P}\{|Q|>(|t_1|+\cdots+|t_s|)^2u\}\lesssim\exp(-Cpu^2)
	\end{align}
	for any $u\rightarrow0$ but $pu\rightarrow\infty$. For a sufficiently large constant $C_*>0$, define
	\begin{align*}
	    \mathcal{E}=\bigg\{|Q|\leq \frac{C_*\log^{1/2}p}{p^{1/2}}\bigg\}\,.
	\end{align*}
For given $(t_1,\ldots,t_s)$, we have $\mathbb{P}(\mathcal{E}^c)\rightarrow0$ as $p\rightarrow\infty$. By \eqref{eq:ft1ts},
	\begin{align*}
	    f(t_1,\ldots,t_s)=&\,\exp\bigg(-\sum_{k=1}^st_k^2\bigg)\big[\mathbb{E}\{\exp(-Q)I(\mathcal{E})\}
	+\mathbb{E}\{\exp(-Q)I(\mathcal{E}^c)\}\big]\\
&+R(t_1,\ldots,t_s)\,.
	\end{align*}
Notice that $Q+\sum_{k=1}^st_k^2\geq0$. For given $(t_1,\ldots,t_s)$, as $p\rightarrow\infty$,
	due to the facts
	\begin{align*}
	    	&~~~~~R(t_1,\ldots,t_s)\rightarrow0\,,~~~~~0\leq \mathbb{E}\bigg\{\exp\bigg(-Q-\sum_{k=1}^st_k^2\bigg)I(\mathcal{E}^c)\bigg\}\leq \mathbb{P}(\mathcal{E}^c)\rightarrow0\,,\\
	    	&1\leftarrow\exp\bigg(-\frac{C_*\log^{1/2} p}{p^{1/2}}\bigg)\mathbb{P}(\mathcal{E})\leq\mathbb{E}\{\exp(-Q)I(\mathcal{E})\}\leq\exp\bigg(\frac{C_*\log^{1/2} p}{p^{1/2                                                                      }}\bigg)\mathbb{P}(\mathcal{E})\rightarrow1\,,
	\end{align*}
	 we have
	$
	f(t_1,\ldots,t_s)\rightarrow\exp(-\sum_{k=1}^st_k^2)
	$, which implies
\[
\sqrt{N}(b_{\ell_1,**}^{-1/2}M_{\ell_1,1},\ldots,b_{\ell_s,**}^{-1/2}M_{\ell_s,1})^\T\xrightarrow{d}\mathcal{N}(\bzero,2\bI_s)\,.
\]
 For $\tilde{b}_{\ell_k}$ specified in \eqref{eq:bellti}, since $
	\sqrt{2}{\mu_{\ell_k,1}\mu_{\ell_k,2}b_{\ell_k,**}^{-1/2}}{(\mu_{\ell_k,1}+\mu_{\ell_k,2})^{-1}\tilde{b}_{\ell_k}^{1/2}}\rightarrow1
	$, then
\begin{align*}
&\sqrt{N}{\rm diag}(\tilde{b}_{\ell_1}^{-1/2},\ldots,\tilde{b}_{\ell_s}^{-1/2})\bigg\{\frac{(\mu_{\ell_1,1}+\mu_{\ell_1,2})M_{\ell_1,1}}{\sqrt{2}\mu_{\ell_1,1}\mu_{\ell_1,2}},\ldots,\frac{(\mu_{\ell_s,1}+\mu_{\ell_s,2})M_{\ell_s,1}}{\sqrt{2}\mu_{\ell_s,1}\mu_{\ell_s,2}}\bigg\}^\T\\
&~~~~~~~~~~~~~~~~~~~~~~~\xrightarrow{d}\mathcal{N}(\bzero,2\bI_s)\,.
\end{align*}
Together with \eqref{eq:expcase2} and \eqref{eq:expcase21},	we complete the proof of Theorem \ref{tm:asymnormal}(b). $\hfill\Box$

	\subsection{Case 3: $\gamma\asymp p^{-1/4}$.}\label{sec:case3}
	
	For $T_{\ell,1}$ defined as \eqref{eq:Tell1}, as shown in Section \ref{se:case2},
	\[
T_{\ell,1}=-\frac{1}{N}\sum_{i,j:\,i\neq j,\,i,j\neq\ell}\bigg(\frac{1}{2\mu_{\ell,1}}+\frac{1}{2\mu_{\ell,2}}\bigg)\mathring{Z}_{i,\ell}\mathring{Z}_{\ell,j}\mathring{Z}_{i,j}+O_{\p}\bigg(
\frac{\log p}{p^{1/2}}\bigg)
\]
when $\gamma\asymp p^{-1/4}$.
By  \eqref{eq:asymexp} and \eqref{eq:Tell2}, we have
	\[
\begin{split}
\hat{\theta}_\ell-\theta_\ell=&~\frac{1}{p-1}\sum_{i:\,i\neq \ell}\lambda_{i,\ell}\mathring{Z}_{i,\ell}-\frac{1}{N}\sum_{i,j:\,i\neq j,\,i,j\neq\ell}\bigg(\frac{1}{2\mu_{\ell,1}}+\frac{1}{2\mu_{\ell,2}}\bigg)\mathring{Z}_{i,\ell}\mathring{Z}_{\ell,j}\mathring{Z}_{i,j}+R_{\ell,12}
\end{split}
\]
	with $R_{\ell,12}=O_{\p}(p^{-1/2}\log p)$,  where $\lambda_{i,\ell}$ is defined as \eqref{eq:lambdaiell}. Define
	\[
\tilde{\lambda}_{i,\ell}=\frac{2\lambda_{i,\ell}\mu_{\ell,1}\mu_{\ell,2}}{\mu_{\ell,1}+\mu_{\ell,2}}\,.
\]
	Then
	\begin{align*}
	   -\frac{2\mu_{\ell,1}\mu_{\ell,2}}{{\mu_{\ell,1}+\mu_{\ell,2}}}(\hat{\theta}_\ell-\theta_\ell)=-\frac{1}{p-1}\sum_{i:\,i\neq \ell}\tilde{\lambda}_{i,\ell}\mathring{Z}_{i,\ell}+\frac{1}{N}\sum_{i,j:\,i\neq j,\,i,j\neq\ell}\mathring{Z}_{i,\ell}\mathring{Z}_{\ell,j}\mathring{Z}_{i,j}+R_{\ell,13}
	\end{align*}
	with $R_{\ell,13}=O_{\p}(p^{-5/4}\log p)$.
	Given $s$ different $\ell_1,\ldots,\ell_s$, as shown in Section \ref{se:case2},
	\[
\frac{1}{N}\sum_{i,j:\,i\neq j,\,i,j\neq\ell_k}\mathring{Z}_{i,\ell_k}\mathring{Z}_{\ell_k,j}\mathring{Z}_{i,j}=M_{\ell_k,1}+O_{\p}\bigg(\frac{1}{p^{3/2}}\bigg)\,,
\]
	where $M_{\ell_k,1}$ is specified in \eqref{eq:theasymp2}. We also have
	\begin{align*}
	    	-\frac{1}{p-1}\sum_{i:\,i\neq \ell_k}\tilde{\lambda}_{i,\ell_k}\mathring{Z}_{i,\ell_k}&=-\frac{1}{p-1}\sum_{i:\,i\neq \ell_1,\ldots,\ell_s}\tilde{\lambda}_{i,\ell_k}\mathring{Z}_{i,\ell_k}+O\bigg(\frac{1}{p^{3/2}}\bigg)\\
	    	&=:{M}_{\ell_k,1}^*+O\bigg(\frac{1}{p^{3/2}}\bigg)\,.
	\end{align*}
	Then
	\begin{align*}
  	-\frac{2\mu_{\ell_k,1}\mu_{\ell_k,2}}{\mu_{\ell_k,1}+\mu_{\ell_k,2}}(\hat{\theta}_{\ell_k}-\theta_{\ell_k})=M_{\ell_k,1}+M_{\ell_k,1}^*+O_{\p}\bigg(\frac{\log p}{p^{5/4}}\bigg)\,.
	\end{align*}
	For any $k\in[s]$, write
	\begin{align}\label{eq:bellk***}
b_{\ell_k,***}=2\bigg(\frac{\mu_{\ell_k,1}\mu_{\ell_k,2}}{\mu_{\ell_k,1}+\mu_{\ell_k,2}}\bigg)^2(p-2)b_{\ell_k}\,,
\end{align}
	 where $b_{\ell_k}$ is defined in \eqref{eq:bell}. Write $\check{b}_{\ell_k}=b_{\ell_k,**}+b_{\ell_k,***}$ with $b_{\ell_k,**}$ specified in \eqref{eq:bellk**}. Let
	\begin{align*}
	   \check{f}(t_1,\ldots,t_s)=\mathbb{E}\bigg[\exp\bigg\{\iota \sqrt{N}\sum_{k=1}^{s}t_k\check{b}_{\ell_k}^{-1/2}(M_{\ell_k,1}+M_{\ell_k,1}^*)\bigg\}\bigg]~~\textrm{with}~~\iota^2=-1\,.
	\end{align*}
	 Recall $\mathcal{F}_{\ell_1,\ldots,\ell_s}^*=\cup_{k=1}^s\{Z_{i,\ell_k},Z_{\ell_k,j}:i,j\neq\ell_k\}$. Same as \eqref{eq:chardiff1},
\begin{align*}
&\bigg|\exp\bigg\{-\frac{1}{2}\sum_{i,j:\,i<j,\atop i,j\neq \ell_1,\ldots,\ell_s}\bigg(\sum_{k=1}^s\frac{2t_k\check{b}_{\ell_k}^{-1/2}}{\sqrt{N}}\mathring{Z}_{i,\ell_k}\mathring{Z}_{\ell_k,j}\Bigg)^2{\rm Var}({Z}_{i,j})\bigg\}-\\
&~~~~~~~-\mathbb{E}\bigg\{\exp\bigg(\iota \sqrt{N}\sum_{k=1}^st_k\check{b}_{\ell_k}^{-1/2}M_{\ell_k,1}\bigg)\,\bigg|\,\mathcal{F}_{\ell_1,\ldots,\ell_s}^*\bigg\}\bigg|\\
&~~~~~~
\lesssim \frac{(|t_1|+\cdots+|t_s|)^3}{\sqrt{N}}+\frac{(|t_1|+\cdots+|t_s|)^4}{N}\,,
\end{align*}
	which implies
	\begin{align}
\check{f}(t_1,\ldots,t_s)=&~\mathbb{E}\bigg[\exp\bigg\{-\sum_{i,j:\,i\neq j,\atop i,j\neq \ell_1,\ldots,\ell_s}\bigg(\sum_{k=1}^s\frac{t_k\check{b}_{\ell_k}^{-1/2}}{\sqrt{N}}\mathring{Z}_{i,\ell_k}\mathring{Z}_{\ell_k,j}\bigg)^2{\rm Var}({Z}_{i,j})\bigg\}\notag\\
&~~~~~~~\times\exp\bigg(\iota \sqrt{N}\sum_{k=1}^st_k\check{b}_{\ell_k}^{-1/2}M_{\ell_k,1}^*\bigg)\bigg]+O(N^{-1/2})\label{eq:hatf}
\end{align}
	for any given $(t_1,\ldots,t_s)$.
	Define
	\[
\begin{split}
\tilde{Q}=&\sum_{i,j:\,i\neq j,\atop i,j\neq \ell_1,\ldots,\ell_s}\bigg(\sum_{k=1}^s\frac{t_k\check{b}_{\ell_k}^{-1/2}}{\sqrt{N}}\mathring{Z}_{i,\ell_k}\mathring{Z}_{\ell_k,j}\bigg)^2{\rm Var}({Z}_{i,j})-\sum_{k=1}^s\frac{t_k^2b_{\ell_k,**}}{\check{b}_{\ell_k}}\,.
\end{split}
\]
Same as \eqref{eq:bdQ}, we have
 \begin{align*}
     	\mathbb{P}\{|\tilde{Q}|>(|t_1|+\cdots+|t_s|)^2u^2\}\lesssim\exp(-Cpu^2)
 \end{align*}
	for any $u\rightarrow0$ but $pu\rightarrow\infty$. For a sufficiently large constant $C_*>0$, define $\tilde{\mathcal{E}}=\{|\tilde{Q}|\leq C_*p^{-1/2}\log^{1/2} p\}$. Then $\mathbb{P}(\tilde{\mathcal{E}}^c)\rightarrow0$ as $p\rightarrow\infty$. 
By \eqref{eq:hatf},
\[
\begin{split}
\check{f}(t_1,\ldots,t_s)=&~\mathbb{E}\bigg\{\exp\bigg(\iota \sqrt{N}\sum_{k=1}^st_k\check{b}_{\ell_k}^{-1/2}M_{\ell_k,1}^*\bigg)\exp(-\tilde{Q})I(\tilde{\mathcal{E}})\bigg\}\\
&~~~~~~~~~~\times \exp\bigg(-\sum_{k=1}^s\frac{t_k^2b_{\ell_k,**}}{\check{b}_{\ell_k}}\bigg)+o(1)\,.
\end{split}
\]
	Note that
	\[
(M_{\ell_1,1}^*,\ldots,M_{\ell_s,1}^*)^\T=-\frac{1}{\gamma}\cdot{\rm diag}\bigg(\frac{2\mu_{\ell_1,1}\mu_{\ell_1,2}}{\mu_{\ell_1,1}+\mu_{\ell_1,2}},\ldots,\frac{2\mu_{\ell_s,1}\mu_{\ell_s,2}}{\mu_{\ell_s,1}+\mu_{\ell_s,2}}\bigg)\cdot\frac{1}{p-1}\sum_{i:\,i\neq\ell_1,\ldots,\ell_s}\bw_i
\]
for $\bw_i$ specified in \eqref{eq:wi}. As shown in Section \ref{sec:case1},
	\begin{align*}
	    	(p-1)^{-1/2}{\rm diag}(b_{\ell_1,*}^{-1/2},\ldots,b_{\ell_s,*}^{-1/2})\sum_{i:\,i\neq\ell_1,\ldots,\ell_s}\bw_i\xrightarrow{d}\mathcal{N}(\bzero,\bI_s)\,.
	\end{align*}
For $b_{\ell_k,***}$ defined as \eqref{eq:bellk***}, due to $b_{\ell_k,*}=\gamma^2 b_{\ell_k}$, then
	\begin{align*}
&\sqrt{N}{\rm diag}\bigg(\frac{b_{\ell_1,***}^{-1/2}}{\sqrt{2}},\ldots,\frac{b_{\ell_k,***}^{-1/2}}{\sqrt{2}}\bigg)(M_{\ell_1,1}^*,\ldots,M_{\ell_s,1}^*)^{\T}\xrightarrow{d}\mathcal{N}(\bzero,\bI_s)\,.
\end{align*}
Notice that $\check{b}_{\ell_k}=b_{\ell_k,**}+b_{\ell_k,***}$.	By the Dominated Convergence Theorem,
	$
	\check{f}(t_1,\ldots,t_s)\rightarrow\exp(-\sum_{k=1}^st_k^2)
	$ for any given $(t_1,\ldots,t_s)$, which implies
\[
\begin{split}
&\sqrt{N}{\rm diag}\bigg\{\frac{\sqrt{2}\mu_{\ell_1,1}\mu_{\ell_1,2}}{(\mu_{\ell_1,1}+\mu_{\ell_1,2})\check{b}_{\ell_1}^{1/2}},\ldots,\frac{\sqrt{2}\mu_{\ell_s,1}\mu_{\ell_s,2}}{(\mu_{\ell_s,1}+\mu_{\ell_s,2})\check{b}_{\ell_s}^{1/2}}\bigg\}(\hat{\theta}_{\ell_1}-\theta_{\ell_1},\ldots,\hat{\theta}_{\ell_s}-\theta_{\ell_s})^\T\\
&~~~~~~~~~~~~~~~~~~~~~~\xrightarrow{d}\mathcal{N}(\bzero,\bI_s)\,.
\end{split}
\]
Due to $b_{\ell_k}\asymp \sqrt{p}$, $\tilde{b}_{\ell_k}\asymp p^{3/2}$ and
\begin{align*}
&\frac{(\mu_{\ell_k,1}+\mu_{\ell_k,2})^2\check{b}_{\ell_k}}{2\mu_{\ell_k,1}^2\mu_{\ell_k,2}^2}=\frac{(\mu_{\ell_k,1}+\mu_{\ell_k,2})^2{b}_{\ell_k,***}}{2\mu_{\ell_k,1}^2\mu_{\ell_k,2}^2}+\frac{(\mu_{\ell_k,1}+\mu_{\ell_k,2})^2{b}_{\ell_k,**
}}{2\mu_{\ell_k,1}^2\mu_{\ell_k,2}^2}\\
&~~~~~~~~~~~~=(p-2)b_{\ell_k}+\frac{(\mu_{\ell_k,1}+\mu_{\ell_k,2})^2}{2N\mu_{\ell_k,1}^2\mu_{\ell_k,2}^2}\sum_{i,j:\,i\neq j,\atop i,j\neq \ell_1,\ldots,\ell_s}{\rm Var}(Z_{i,\ell_k}){\rm Var}(Z_{\ell_k,j}){\rm Var}(Z_{i,j})\\
&~~~~~~~~~~~~=(p-2)b_{\ell_k}+\tilde{b}_{\ell_k}+O(\sqrt{p})\,,
\end{align*}
then
\[
\frac{\sqrt{2}\mu_{\ell_k,1}\mu_{\ell_k,2}}{(\mu_{\ell_k,1}+\mu_{\ell_k,2})\check{b}_{\ell_k}^{1/2}}=\{(p-2)b_{\ell_k}+\tilde{b}_{\ell_k}\}^{-1/2}+O(p^{-7/4})\,.
\]
	Hence, we complete the proof of Theorem \ref{tm:asymnormal}(c). $\hfill\Box$

		\section{Proof of Proposition \ref{tm:varc}.}\label{sec:pfpn5}
	
	To construct Proposition \ref{tm:varc}, we need the following lemmas whose proofs are given in Sections \ref{sec:pflemma2} and \ref{sec:pflemma3}, respectively.
	\begin{lemma}\label{la:4}
		Let Condition {\rm\ref{cond1}} hold and $(\alpha,\beta)\in\mathcal{M}(\gamma;C_1)$ for some fixed constant $C_1\in(0,0.5)$. If $\gamma\gg p^{-1/3}\log^{1/6}p$, it holds that
		$$
		\max_{i:\,i\neq\ell}|\hat{\lambda}_{i,\ell}-\lambda_{i,\ell}|=O_{\p}\bigg(\frac{\log^{1/2}p}{\gamma^{3}p^{1/2}}\bigg)$$
		for any given $\ell\in[p]$, where $\hat{\lambda}_{i,\ell}$ is defined in \eqref{eq:hatlambdaiell}.
	\end{lemma}

	\begin{lemma}\label{la:5}
		Let Condition {\rm\ref{cond1}} hold and $(\alpha,\beta)\in\mathcal{M}(\gamma;C_1)$ for some fixed constant $C_1\in(0,0.5)$. If $\gamma\gg p^{-1/3}\log^{1/6}p$, it holds that
			\[
	\max_{i:\,i\neq \ell}|\widehat{{\rm Var}}(Z_{i,\ell})-{\rm Var}(Z_{i,\ell})|=(3\gamma+\alpha-\beta)\bigg\{O_{\p}\bigg(\frac{\log^{1/2}p}{\gamma^2p}\bigg)+O_{\p}\bigg(\frac{\log^{1/2}p}{ p^{1/2}}\bigg)\bigg\}
	\]
		for any given $\ell\in[p]$, where $\widehat{{\rm Var}}(Z_{i,\ell})$ is defined in {\rm(\ref{eq:hatvar})}.
	\end{lemma}
	
 Recall that
	\[
b_\ell=\frac{1}{p-1}\sum_{i:\,i\neq \ell}\lambda_{i,\ell}^2{\rm Var}(Z_{i,\ell})~~\textrm{and}~ ~\hat{b}_\ell=\frac{1}{p-1}\sum_{i:\,i\neq \ell}\hat{\lambda}_{i,\ell}^2\widehat{{\rm Var}}(Z_{i,\ell})\,.
\]
	Then
	\begin{align*}
\hat{b}_\ell-b_\ell=&~\frac{1}{p-1}\sum_{i:\,i\neq \ell}(\hat{\lambda}_{i,\ell}^2-\lambda_{i,\ell}^2){\rm Var}(Z_{i,\ell})+\frac{1}{p-1}\sum_{i:\,i\neq \ell}\lambda_{i,\ell}^2\{\widehat{{\rm Var}}(Z_{i,\ell})-{\rm Var}(Z_{i,\ell})\}\\
&+\frac{1}{p-1}\sum_{i:\,i\neq \ell}(\hat{\lambda}_{i,\ell}^2-\lambda_{i,\ell}^2)\{\widehat{{\rm Var}}(Z_{i,\ell})-{\rm Var}(Z_{i,\ell})\}\,.
\end{align*}
	Note that $\lambda_{i,\ell}\asymp \gamma^{-1}$ and ${\rm Var}(Z_{i,\ell})\asymp 1$.
	By Lemmas \ref{la:4} and \ref{la:5},
	\begin{align*}
|\hat{b}_\ell-b_\ell|\lesssim&\,\max_{i:\,i\neq\ell}|\hat{\lambda}_{i,\ell}^2-\lambda_{i,\ell}^2|+\gamma^{-2}\max_{i:\,i\neq\ell}|\widehat{{\rm Var}}(Z_{i,\ell})-{\rm Var}(Z_{i,\ell})|\\
=&~O_{\p}\bigg(\frac{\log p}{\gamma^6p}\bigg)+O_{\p}\bigg(\frac{\log^{1/2}p}{\gamma^4p^{1/2}}\bigg)\,.
\end{align*}
 Since $b_\ell\asymp\gamma^{-2}$, we have
\[
\bigg|\frac{\hat{b}_\ell}{b_\ell}-1\bigg|=O_{\p}\bigg(\frac{\log p}{\gamma^4p}\bigg)+O_{\p}\bigg(\frac{\log^{1/2}p}{\gamma^2p^{1/2}}\bigg)\,.
\]
Analogously, we have
\[
|\hat{\tilde{b}}_\ell-\tilde{b}_\ell|=O_{\p}\bigg(\frac{\log^{1/2}p}{\gamma^9p}\bigg)+O_{\p}\bigg(\frac{\log^{1/2}p}{\gamma^7p^{1/2}}\bigg)\,.
\]
Recall that $\tilde{b}_\ell\asymp \gamma^{-6}$. It holds that
\[
\bigg|\frac{\hat{\tilde{b}}_\ell}{\tilde{b}_\ell}-1\bigg|=O_{\p}\bigg(\frac{\log^{1/2}p}{\gamma^3p}\bigg)+O_{\p}\bigg(\frac{\log^{1/2}p}{\gamma p^{1/2}}\bigg)\,.
\]
	We complete the proof of Proposition \ref{tm:varc}. $\hfill\Box$

	\section{Proof of Theorem \ref{tm:2}.}
	The proof of Theorem \ref{tm:2}(a) is almost identical to that of Theorem \ref{tm:asymnormal} given in Section \ref{sec:pfthm1} with replacing $(\hat{\mu}_{\ell,1},\hat{\mu}_{\ell,2},\mu_{\ell,1},\mu_{\ell,2})$ there by the bootstrap analogues
$(\hat{\mu}_{\ell,1}^\dag,\hat{\mu}_{\ell,2}^\dag,\mu_{\ell,1}^\dag,\mu_{\ell,2}^\dag)$ and is thus omitted here.  We only prove Theorem \ref{tm:2}(b).
	Due to $\nu_{\ell}^{\dag}=(p-2)b_{\ell}^{\dag}+\tilde{b}_{\ell}^{\dag}$ and $\nu_{\ell}=(p-2)b_{\ell}+\tilde{b}_{\ell}$, then
	$
	    	\nu_{\ell}^{\dag}-\nu_{\ell}=(p-2)(b_{\ell}^{\dag}-b_{\ell})+\tilde{b}_{\ell}^{\dag}-\tilde{b}_{\ell}$.
Note that
$
		\lambda_{i,\ell}^{\dag}=(1-2\delta)^{-1}\lambda_{i,\ell}$, $\mu_{\ell,1}^{\dag}=(1-2\delta)^{3}\mu_{\ell,1}$ and
$\mu_{\ell,2}^{\dag}=(1-2\delta)^{3}\mu_{\ell,2}$.
 Then
	\begin{align*}
	    	&~~~~~~~~~~~~~~~~~~~~~b_{\ell}^{\dag}-b_\ell=\frac{1}{p-1}\sum_{i:\,i\neq \ell}\lambda_{i,\ell}^2\bigg\{\frac{{\rm Var}(Z_{i,\ell}^{\dag})}{(1-2\delta)^{2}}-{\rm Var}(Z_{i,\ell})\bigg\}\,,\\
	    	&\tilde{b}_{\ell}^{\dag}-\tilde{b}_\ell=\frac{1}{2N}\bigg(\frac{\mu_{\ell,1}+\mu_{\ell,2}}{\mu_{\ell,1}\mu_{\ell,2}}\bigg)^2\sum_{i,j:\, i\neq j,i,j\neq \ell}\bigg\{\frac{{\rm Var}(Z_{i,\ell}^{\dag})}{(1-2\delta)^2}\frac{{\rm Var}(Z_{\ell,j}^{\dag})}{(1-2\delta)^2}\frac{{\rm Var}(Z_{i,j}^{\dag})}{(1-2\delta)^2}\\
&~~~~~~~~~~~~~~~~~~~~~~~~~~~~~~~~~~~~~~~~~~~~~~~~~~~~~~~~~~~~~~~~~~~~~~~~~~~~-{\rm Var}(Z_{i,\ell}){\rm Var}(Z_{\ell,j}){\rm Var}(Z_{i,j})\bigg\}\,.
	\end{align*}
	Recall $\delta\in (0,c]$ with $c<0.5$.	For any $i\neq \ell$, we have
	\begin{align*}
	    \bigg|\frac{{\rm Var}(Z_{i,\ell}^{\dag})}{(1-2\delta)^{2}}-{\rm Var}(Z_{i,\ell})\bigg|=\frac{\delta(1-\delta)}{(1-2\delta)^{2}}\lesssim\delta\,.
	\end{align*}
	Under  Condition \ref{cond1}, we have
	\begin{align*}
	    \min_{\ell\in[p]}\min_{i:\,i\neq\ell}\lambda_{i,\ell}\asymp\gamma^{-1}\asymp\max_{\ell\in[p]}\max_{i:\,i\neq\ell}\lambda_{i,\ell}~~\textrm{and}~~
	  \min_{\ell\in[p]}\mu_{\ell,1}\asymp\gamma^3\asymp\max_{\ell\in[p]}\mu_{\ell,2}\,.
	\end{align*}
	 Then
	 	$
		(p-2)\max_{\ell\in[p]}|b_{\ell}^{\dag}-b_{\ell}|\lesssim \gamma^{-2}p\delta$ and $\max_{\ell\in[p]}|\tilde{b}_{\ell}^{\dag}-\tilde{b}_{\ell}|\lesssim \gamma^{-6}\delta$,
 which implies
	\begin{align*}
	    	\max_{\ell \in[p]}|\nu_{\ell}^{\dag}-\nu_{\ell}|\lesssim \frac{p\delta}{\gamma^{2}}+\frac{\delta}{\gamma^{6}}\,.
	\end{align*}
Note that $\min_{\ell\in[p]}\nu_{\ell} \asymp p\ga^{-2}+\gamma^{-6}\asymp\max_{\ell\in[p]}\nu_\ell$. Then
	 \begin{align*}
	    \max_{\ell \in[p]}\bigg|\frac{\nu_\ell^\dag}{\nu_\ell}-1\bigg|=O(\delta)\,.
	 \end{align*}
We complete the proof of Theorem \ref{tm:2}(b).	$\hfill\Box$

		\section{Proof of Theorem \ref{tm:ga}.}\label{appF}
	
	Recall $N=(p-1)(p-2)$. As shown in \eqref{eq:asymexp}, it holds that
	\begin{align}\label{eq:asymexpga1}	\hat{\theta}_\ell-\theta_\ell=\underbrace{\frac{I_{\ell,1,1}}{2\mu_{\ell,1}}-\frac{I_{\ell,2,1}}{2\mu_{\ell,2}}}_{T_{\ell,1}}+\underbrace{\frac{I_{\ell,1,2}(1)}{2\mu_{\ell,1}N}-\frac{I_{\ell,2,2}(1)}{2\mu_{\ell,2}N}}_{T_{\ell,2}}+R_{\ell,4}\,,
	\end{align}
	where
	$\max_{\ell\in[p]}|R_{\ell,4}|=O_{\p}(\gamma^{-6}p^{-2}\log p)+O_{\p}(\gamma^{-2}p^{-1}\log p)$. By \eqref{eq:jell},
\begin{align*}
    T_{\ell,1}=-\frac{1}{N}\sum_{i,j:\,i\neq j,\,i,j\neq\ell}\bigg(\frac{\mu_{\ell,1}+\mu_{\ell,2}}{2\mu_{\ell,1}\mu_{\ell,2}}\bigg)\mathring{Z}_{i,\ell}\mathring{Z}_{\ell,j}\mathring{Z}_{i,j}-J_{\ell,1}-J_{\ell,2}\,.
   \end{align*}
	Following the arguments for deriving \eqref{eq:convJell2}, we have $\max_{\ell\in[p]}|J_{\ell,2}|=O_{\p}(\gamma^{-1}p^{-1}\log^{1/2}p)$. Together with \eqref{eq:bdJell1}, we have $\max_{\ell\in[p]}|J_{\ell,1}+J_{\ell,2}|=O_{\p}(\gamma^{-2}p^{-1}\log p)$. Recall
	$
	T_{\ell,2}=(p-1)^{-1}\sum_{i:\,i\neq\ell}\lambda_{i,\ell}\mathring{Z}_{i,\ell}
	$
	with $\lambda_{i,\ell}$ specified in \eqref{eq:lambdaiell}. By \eqref{eq:asymexpga1}, we have
	\[
\hat{\theta}_\ell-\theta_\ell=-\frac{1}{N}\sum_{i,j:\,i\neq j,\,i,j\neq\ell}\bigg(\frac{\mu_{\ell,1}+\mu_{\ell,2}}{2\mu_{\ell,1}\mu_{\ell,2}} \bigg)\mathring{Z}_{i,\ell}\mathring{Z}_{\ell,j}\mathring{Z}_{i,j}+\frac{1}{p-1}\sum_{i:\,i\neq\ell}\lambda_{i,\ell}\mathring{Z}_{i,\ell}+R_{\ell,14}\,,
\]
	where $\max_{\ell\in[p]}|R_{\ell,14}|=O_{\p}(\gamma^{-6}p^{-2}\log p)+O_{\p}(\gamma^{-2}p^{-1}\log p)$. Recall $\nu_\ell=(p-2)b_\ell+\tilde{b}_\ell\asymp p\gamma^{-2}+\gamma^{-6}$. Then
	\[
\begin{split}
\sqrt{N}\nu_\ell^{-1/2}(\hat{\theta}_\ell-\theta_\ell)=&-\frac{1}{\sqrt{N}}\sum_{i,j:\,i\neq j,\,i,j\neq\ell}\nu_\ell^{-1/2}\bigg\{\bigg(\frac{\mu_{\ell,1}+\mu_{\ell,2}}{2\mu_{\ell,1}\mu_{\ell,2}}\bigg)\mathring{Z}_{i,\ell}\mathring{Z}_{\ell,j}\mathring{Z}_{i,j}-\lambda_{i,\ell}\mathring{Z}_{i,\ell}\bigg\}\\
&+R_{\ell,15}\,,
\end{split}
\]
	where $\max_{\ell\in[p]}|R_{\ell,15}|=O_{\p}(\gamma^{-1}p^{-1/2}\log p)+O_{\p}(\gamma^{-3}p^{-1}\log p)$.
Notice that
\begin{align*}
\sqrt{N}(\nu_\ell^\dag)^{-1/2}(\hat{\theta}_\ell-\theta_\ell)=&\,\sqrt{N}\nu_\ell^{-1/2}(\hat{\theta}_\ell-\theta_\ell)+\sqrt{N}\nu_\ell^{-1/2}\bigg\{\bigg(\frac{\nu_\ell}{\nu_\ell^\dag}\bigg)^{1/2}-1\bigg\}(\hat{\theta}_\ell-\theta_\ell)\\
=&\,\sqrt{N}\nu_\ell^{-1/2}(\hat{\theta}_\ell-\theta_\ell)+R_{\ell,16}\,.
\end{align*}
By Theorem \ref{tm:2}(b), $\max_{\ell\in[p]}|(\nu_\ell/\nu_\ell^\dag)^{1/2}-1|=O(\delta)$. Together with Proposition \ref{tm:1}, we have
$
\max_{\ell\in[p]}|R_{\ell,16}|=O_{\p}(\delta\log^{1/2}p)$.
Given $(i,j)$ such that $i\neq j$, we define
	\begin{align}\label{eq:Yijell}
Y_{(i,j),\ell}=-\nu_\ell^{-1/2}\bigg\{\bigg(\frac{\mu_{\ell,1}+\mu_{\ell,2}}{2\mu_{\ell,1}\mu_{\ell,2}}
\bigg)\mathring{Z}_{i,\ell}\mathring{Z}_{\ell,j}\mathring{Z}_{i,j}-\lambda_{i,\ell}\mathring{Z}_{i,\ell}\bigg\}
\end{align}
	for any $\ell\neq i,j$, and $Y_{(i,j),\ell}=0$ for $\ell=i$ or $j$. Write $R_{\ell,17}=R_{\ell,15}+R_{\ell,16}$. Hence,
	\[
\sqrt{N}(\nu_\ell^\dag)^{-1/2}(\hat{\theta}_\ell-\theta_\ell)=\frac{1}{\sqrt{N}}\sum_{i,j:\,i\neq j}Y_{(i,j),\ell}+R_{\ell,17}=:{h_\ell}+R_{\ell,17}\,,
\]
where $\max_{\ell\in[p]}|R_{\ell,17}|=O_{\p}(\gamma^{-1}p^{-1/2}\log p)+O_{\p}(\gamma^{-3}p^{-1}\log p)+O_{\p}(\delta\log^{1/2}p)$.

	Write $\bV^{\dag}={\rm diag}(\nu_1^\dag,\ldots,\nu_p^\dag)$, $\bh=(h_1,\ldots,h_p)^\T$ and $\tilde{\br}=(R_{1,17},\ldots,R_{p,17})^\T$. Then
	\begin{equation}\label{eq:asymexpleading}
		\sqrt{N}(\bV^\dag)^{-1/2}(\hat{\btheta}-\btheta)=\bh+\tilde{\br}\,.
	\end{equation}
	Lemma \ref{la:vnorm} states the property of $\bB={\rm Cov}(\bh)$, whose proof is given in Section \ref{se:proofla4}.
	
	\begin{lemma}\label{la:vnorm}
		Write
		$\bB=(B_{\ell_1,\ell_2})_{p\times p}
		$. Then $\max_{1\leq\ell_1\neq\ell_2\leq p}|B_{\ell_1,\ell_2}|\lesssim p^{-1}$ and $B_{\ell,\ell}=1$ for any $\ell\in [p]$.
	\end{lemma}

Define
	\[
\varrho=\sup_{\bu\in\mathbb{R}^{p}}\big|\mathbb{P}(\bh\leq\bu)-\mathbb{P}(\bxi\leq \bu)\big|
\]
	with $\bxi\sim\mathcal{N}(\bzero,\bI_p)$. By \eqref{eq:asymexpleading}, we have 
	\[
\begin{split}
\mathbb{P}\{\sqrt{N}(\bV^\dag)^{-1/2}(\hat{\btheta}-\btheta)\leq \bu\}=\mathbb{P}(\bh+\tilde{\br}\leq \bu, |\tilde{\br}|_\infty\leq\epsilon)+\mathbb{P}(\bh+\tilde{\br}\leq \bu, |\tilde{\br}|_\infty>\epsilon)
\end{split}
\]
	for any $\epsilon>0$, which implies 
	\begin{align*}
&~~~~~~~~~~~~~~~~~~~~~\mathbb{P}\{\sqrt{N}(\bV^\dag)^{-1/2}(\hat{\btheta}-\btheta)\leq \bu\}\leq \mathbb{P}(\bh\leq \bu+\epsilon)+\mathbb{P}(|\tilde{\br}|_\infty>\epsilon)\,,\\
&\mathbb{P}\{\sqrt{N}(\bV^\dag)^{-1/2}(\hat{\btheta}-\btheta)\leq \bu\}\geq\mathbb{P}(\bh\leq\bu-\epsilon,|\tilde{\br}|_\infty\leq\epsilon)\geq\mathbb{P}(\bh\leq\bu-\epsilon)-\mathbb{P}(|\tilde{\br}|_\infty>\epsilon)\,.
\end{align*}
	Therefore,
	\begin{align}
&\mathbb{P}\{\sqrt{N}(\bV^\dag)^{-1/2}(\hat{\btheta}-\btheta)\leq \bu\}-\mathbb{P}(\bxi\leq \bu)\notag\\
&~~~~~~~~~\leq \mathbb{P}(\bh\leq \bu+\epsilon)-\mathbb{P}(\bxi\leq \bu+\epsilon)+\mathbb{P}(\bxi\leq \bu+\epsilon)-\mathbb{P}(\bxi\leq \bu)+\mathbb{P}(|\tilde{\br}|_\infty>\epsilon)\label{eq:bb}\\
&~~~~~~~~~\leq \varrho+\mathbb{P}(\bxi\leq \bu+\epsilon)-\mathbb{P}(\bxi\leq \bu)+\mathbb{P}(|\tilde{\br}|_\infty>\epsilon)\,,\notag\\
&\mathbb{P}\{\sqrt{N}(\bV^\dag)^{-1/2}(\hat{\btheta}-\btheta)\leq \bu\}-\mathbb{P}(\bxi\leq \bu)\notag\\
&~~~~~~~~~\geq \mathbb{P}(\bh\leq \bu-\epsilon)-\mathbb{P}(\bxi\leq \bu-\epsilon)+\mathbb{P}(\bxi\leq \bu-\epsilon)-\mathbb{P}(\bxi\leq \bu)-\mathbb{P}(|\tilde{\br}|_\infty>\epsilon)\notag\\
&~~~~~~~~~\geq -\varrho+\mathbb{P}(\bxi\leq \bu-\epsilon)-\mathbb{P}(\bxi\leq \bu)-\mathbb{P}(|\tilde{\br}|_\infty>\epsilon)\,,\notag
\end{align}
	which implies that
	\begin{align}\label{eq:bb}
		&\sup_{\bu\in\mathbb{R}^p}\big|\mathbb{P}\{\sqrt{N}(\bV^\dag)^{-1/2}(\hat{\btheta}-\btheta)\leq \bu\}-\mathbb{P}(\bxi\leq \bu)\big|\notag\\
		&~~~~~~~~~~~~~~~\leq \varrho+\sup_{\bu\in\mathbb{R}^p}\big|\mathbb{P}(\bxi\leq \bu+\epsilon)-\mathbb{P}(\bxi\leq \bu)\big|+\mathbb{P}(|\tilde{\br}|_\infty>\epsilon)\\
		&~~~~~~~~~~~~~~~\leq \varrho+C\epsilon\log^{1/2}p+\mathbb{P}(|\tilde{\br}|_\infty>\epsilon)\notag\,.
	\end{align}
	The last step in above inequality is due to Nazarov's inequality \citep[Lemma A.1]{CCK2017}. As we will show in Section \ref{sec:pfrho0}, $\varrho\rightarrow0$ as $p\rightarrow\infty$. Recall $|\tilde{\br}|_\infty=O_{\p}(p^{-1/2}\gamma^{-1}\log p)+O_{\p}(p^{-1}\gamma^{-3}\log p)+O_{\p}(\delta\log^{1/2}p)$. Since $\gamma\gg p^{-1/3}\log^{1/2}p$ and $\delta\ll (\log p)^{-1}$, we know $p^{-1/2}\gamma^{-1}\log^{3/2}p+p^{-1}\gamma^{-3}\log^{3/2}p+\delta\log p\rightarrow0$ as $p\rightarrow\infty$. Therefore,
 there exists $\epsilon>0$ satisfying  $p^{-1/2}\gamma^{-1}\log p+p^{-1}\gamma^{-3}\log p+\delta\log^{1/2}p\ll \epsilon\ll (\log p)^{-1/2}$. For such selected $\epsilon$, we have $C\epsilon\log^{1/2}p+\mathbb{P}(|\tilde{\br}|_\infty>\epsilon)\rightarrow0$ as $p\rightarrow\infty$. Then we complete the proof of Theorem \ref{tm:ga} by \eqref{eq:bb}. $\hfill\Box$

\subsection{To show $\varrho\rightarrow0$ as $p\rightarrow\infty$.}\label{sec:pfrho0}
	Let $\bg\sim\mathcal{N}(\bzero,\bB)$ with $\bB$ specified in Lemma \ref{la:vnorm}, and define
	\begin{align*}
	    	\bar{\varrho}=\sup_{\bu\in\mathbb{R}^{p}}|\mathbb{P}(\bh\leq\bu)-\mathbb{P}(\bg\leq \bu)|\,.
	\end{align*}
	Then
	\begin{align*}
	    	\varrho\leq \bar{\varrho}+\sup_{\bu\in\mathbb{R}^{p}}|\mathbb{P}(\bxi\leq\bu)-\mathbb{P}(\bg\leq \bu)|\,.
	\end{align*}
	Recall $\bxi\sim\mathcal{N}(\bzero,\bI_p)$ and $\bg\sim\mathcal{N}(\bzero,\bB)$ with $|\bI_p-\bB|_\infty\lesssim p^{-1}$. By Lemma 1 of \cite{ChangChenWu_2021}, we have
	\begin{align*}
	    \sup_{\bu\in\mathbb{R}^{p}}|\mathbb{P}(\bxi\leq\bu)-\mathbb{P}(\bg\leq \bu)|\lesssim |\bI_p-\bB|_\infty^{1/3}\log^{2/3}p\lesssim p^{-1/3}\log^{2/3}p\rightarrow0\,.
	\end{align*}
	To show $\varrho\rightarrow0$ as $p\rightarrow\infty$, it suffices to show $\bar{\varrho}\rightarrow0$ as $p\rightarrow\infty$. Define
	\begin{align*}
	    \varrho_*=\sup_{\bu\in\mathbb{R}^{p},v\in[0,1]}|\mathbb{P}(\sqrt{v}\bh+\sqrt{1-v}\bg\leq\bu)-\mathbb{P}(\bg\leq \bu)|
	\end{align*}
	with $\bg\sim\mathcal{N}(\bzero,\bB)$.
	It is obvious that $\bar{\varrho}\leq \varrho_*$. In the sequel, we will show $\varrho_*\rightarrow0$ as $p\rightarrow\infty$. Let $\beta:=\phi\log p$. For a given $\bu=(u_1,\ldots,u_{p})^\T\in\mathbb{R}^{p}$, we define
	\begin{equation}\label{eq:Fbeta}
		F_\beta(\by):=\beta^{-1}\log\bigg[\sum_{\ell=1}^{p}\exp\{\beta(y_\ell-u_\ell)\}\bigg]
	\end{equation}
	for any $\by=(y_1,\ldots,y_{p})^\T\in\mathbb{R}^{p}$. Such defined function $F_\beta(\by)$ satisfies the property
	\begin{align*}
	    0\leq F_\beta(\by)-\max_{\ell\in[p]}(y_\ell-u_\ell)\leq \frac{\log p}{\beta}=\frac{1}{\phi}
	\end{align*}
	 for any $\by\in\mathbb{R}^{p}$. Select a thrice continuously differentiable function $f_0:\mathbb{R}\rightarrow[0,1]$ whose derivatives up to the third order are all bounded such that $f_0(t)=1$ for $t\leq 0$ and $f_0(t)=0$ for $t\geq1$. Define $f(t):=f_0(\phi t)$ for any $t\in\mathbb{R}$, and $q(\by):=f\{F_\beta(\by)\}$ for any $\by\in\mathbb{R}^{p}$. To simplify the notation, we write $q_\ell(\by)=\partial q(\by)/\partial y_\ell$, $q_{\ell k}(\by)=\partial^2q(\by)/\partial y_\ell\partial y_k$ and $q_{\ell kl}(\by)=\partial^3q(\by)/\partial y_\ell\partial y_k \partial y_l$. Let $\tilde{\bg}$ be an independent copy of $\bg$. Define
	 \begin{align*}
	     \mathcal{T}:=q(\sqrt{v}\bh+\sqrt{1-v}\bg)-q(\tilde{\bg})\,.
	 \end{align*}
Lemma \ref{lemma5} gives an upper bound for $\sup_{v\in[0,1]}|\mathbb{E}(\mathcal{T})|$, whose proof is given in Section \ref{sec:pflemma5}.

\begin{lemma}\label{lemma5}
If $\phi\ll p^{1/2}(\log p)^{-3/2}$, then $\sup_{v\in[0,1]}|\mathbb{E}(\mathcal{T})|\lesssim p^{-1/2}\phi^3\log^{7/2}p$.
\end{lemma}

	Write $\bdelta=\sqrt{v}\bh+\sqrt{1-v}\bg$. Notice that
\begin{align*}
\mathbb{P}(\bdelta\leq \bu-\phi^{-1})\leq&~\mathbb{P}\{F_\beta(\bdelta)\leq 0\}\leq \mathbb{E}\{q(\bdelta)\}\leq\mathbb{P}\{F_\beta(\tilde{\bg})\leq \phi^{-1}\}+\mathbb{E}(\mathcal{T})\\
\leq&~\mathbb{P}(\tilde{\bg}\leq \bu+\phi^{-1})+|\mathbb{E}(\mathcal{T})|\\
\leq&~\mathbb{P}(\tilde{\bg}\leq \bu-\phi^{-1})+C\phi^{-1}\log^{1/2}p+|\mathbb{E}(\mathcal{T})|
\end{align*}
	and
	$
	  \mathbb{P}(\bdelta\leq \bu-\phi^{-1})\geq\mathbb{P}(\tilde{\bg}\leq \bu-\phi^{-1})-C\phi^{-1}\log^{1/2}p-|\mathbb{E}(\mathcal{T})|$.
Together with Lemma \ref{lemma5}, we have
	\begin{align*}
		\varrho_*=&\,\sup_{\bu\in\mathbb{R}^{p},v\in[0,1]}|\mathbb{P}(\bdelta\leq\bu)-\mathbb{P}(\bg\leq \bu)|\\
		\leq&\, C\phi^{-1}\log^{1/2}p+\sup_{v\in[0,1]}|\mathbb{E}(\mathcal{T})|\lesssim \phi^{-1}\log^{1/2}p+p^{-1/2}\phi^3\log^{7/2}p\,.
	\end{align*}
Selecting $\phi=p^{1/8}\log^{-3/4}p$, we have
	$
	\varrho_*
	\lesssim{p^{-1/8}}{\log^{5/4}p}$ as $p\rightarrow\infty$. Hence, it holds that $\varrho\rightarrow0$ as $p\rightarrow\infty$. $\hfill\Box$

\section{Proof of Proposition \ref{tm:4}.} \label{sec:pftm4}

Recall $s=|\mathcal{S}|$ and $\chi_{p}=\exp(-|\xi^+|\vee\max_{\ell\in\mathcal{S}}|\check{\theta}_\ell^+|)$. Define $A_{i,\ell}=2\sum_{j:\,j\neq i,\ell}\mathbb{E}\{\varphi_{(\ell,j),1}\}\mathbb{E}\{\varphi_{(i,j),0}\}$. To prove Proposition \ref{tm:4}, we need the following lemmas whose proofs are given in Sections \ref{sec:pfla:mu1mul2} and \ref{sec:pfl7}, respectively.

\begin{lemma}\label{la:mu1mul2}
Let $(\alpha,\beta)\in\mathcal{M}(\gamma, C_1)$ for some fixed constant $C_1\in(0,0.5)$.
It holds that
	\begin{align*}
  &\min_{\ell \in\mathcal{S}}\mu_{\ell,1}\gtrsim p^{2\omega_2-2\omega_1}\chi_{p}^5\gamma^3\,, ~~~~\min_{\ell \in\mathcal{S}^c}\mu_{\ell,1}\gtrsim p^{-2\omega_1}\chi_{p}^3\gamma^3\,,
		\\
		&~~~~\min_{\ell \in\mathcal{S}}\mu_{\ell,2}\gtrsim p^{-\omega_1}\chi_{p}^5\gamma^3\,,~~~~
		\min_{\ell \in\mathcal{S}^c}\mu_{\ell,2}\gtrsim p^{-\omega_1}\chi_{p}^3\gamma^3\,.
\end{align*}
	\end{lemma}


\begin{lemma}\label{la:9}
Let $(\alpha,\beta)\in\mathcal{M}(\gamma, C_1)$ for some fixed constant $C_1\in(0,0.5)$. It holds that
\begin{align*}
	\max_{\ell\in\mathcal{S}}|\hat{\mu}_{\ell,1}-\mu_{\ell,1}|=&~O_{\p}\big[\gamma^2p^{-3/2}\{sp^{\min(2\omega_2-\omega_1,\,0)} \chi_{p}^{-1}+p^{1+\omega_2-\omega_1}\}\chi_{p}^{-2}\log^{1/2}s\big]\\
 &+O_{\p}(\gamma p^{-1}\log s)+O_{\p}(p^{-1}\log^{1/2}s)\,,\\
	\max_{\ell\in\mathcal{S}^c}|\hat{\mu}_{\ell,1}-\mu_{\ell,1}|=&~O_{\p}\big[\gamma^2s^{1/2}p^{-2}\{sp^{\min(-\omega_2,\,\omega_2-\omega_1)} \chi_{p}^{-4}+p^{1-\omega_1}\}\chi_{p}^{-1}\log^{1/2}p\big]\\
 &+O_{\p}\big\{\gamma^2p^{-3/2}(sp^{\omega_2-\omega_1} \chi_{p}^{-1}+p^{1-\omega_1})\chi_{p}^{-1}\log^{1/2}p\big\}\\
 &+O_{\p}(\gamma p^{-1}\log p)+O_{\p}(p^{-1}\log^{1/2}p)\,,\\
\max_{\ell\in\mathcal{S}}|\hat{\mu}_{\ell,2}-\mu_{\ell,2}|=&~O_{\p}\big[\gamma^2s^{1/2}p^{-2}\{sp^{\min(2\omega_2-\omega_1,\,\omega_1-2\omega_2)}\chi_{p}^{-1}+p^{1+\omega_2-\omega_1}\}\chi_{p}^{-2}\log^{1/2}s\big]\\
&+O_{\p}\big[\gamma^2p^{-3/2}\{sp^{\min(-\omega_2,\,\omega_2-\omega_1)} \chi_{p}^{-4}+p^{1-\omega_1}\}\chi_{p}^{-1}\log^{1/2}s\big]\\
&+O_{\p}\big[\gamma s p^{-2}\{p^{2\omega_2-\omega_1}\chi_{p}^{-3}  I(\omega_1>2\omega_2)+I(\omega_1\leq2\omega_2)\}\log s\big]\\
&+O_{\p}\big[\gamma s^{1/2} p^{-3/2}\{p^{\omega_2-\omega_1}\chi_{p}^{-2}I(\omega_1>\omega_2)+I(\omega_1=\omega_2)\}\log s\big]\\
&+O_{\p}\big[\gamma p^{-1}\{p^{-\omega_1}\chi_{p}^{-1}I(\omega_1>0)+I(\omega_1=0)\}\log s\big]\\
&+O_{\p}(p^{-1}\log^{1/2}s)\,,\\
	\max_{\ell\in\mathcal{S}^c}|\hat{\mu}_{\ell,2}-\mu_{\ell,2}|=&~O_{\p}\big[\gamma^2s^{1/2}p^{-2}\{sp^{\min(2\omega_2-\omega_1,\,0)} \chi_{p}^{-1}+p^{1+\omega_2-\omega_1}\}\chi_{p}^{-2}\log^{1/2}p\big]\\
 &+O_{\p}\big\{\gamma^2p^{-3/2}(sp^{\omega_2-\omega_1} \chi_{p}^{-1}+p^{1-\omega_1})\chi_{p}^{-1}\log^{1/2}p\big\}\\
 &+O_{\p}\big[\gamma sp^{-2}\{p^{2\omega_2-\omega_1}\chi_{p}^{-3}I(\omega_1>2\omega_2)+ I(\omega_1\leq2\omega_2)\}\log p\big]\\
 &+O_{\p}\big[\gamma s^{1/2} p^{-3/2}\{p^{\omega_2-\omega_1}\chi_{p}^{-2} I(\omega_1>\omega_2)+ I(\omega_1=\omega_2)\}\log p\big]\\
 &+O_{\p}\big[\gamma p^{-1}\{p^{-\omega_1}\chi_{p}^{-1} I(\omega_1>0)+I(\omega_1=0)\}\log p\big]\\
 &+O_{\p}(p^{-1}\log^{1/2}p)\,.
	\end{align*}
\end{lemma}

Now we begin to prove Proposition \ref{tm:4}.
	Recall
	$
		\hat{\zeta}_{\ell}=\hat{\mu}_{\ell,1}/\hat{\mu}_{\ell,2}$ and $
		\zeta_{\ell}=\mu_{\ell,1}/\mu_{\ell,2}$.
	Then
	\[
	\begin{split}
		\hat{\zeta}_{\ell}-\zeta_{\ell}
		=&~\frac{\hat{\mu}_{\ell,1}-\mu_{\ell,1}}{\mu_{\ell,2}}-\frac{\mu_{\ell,1}}{\mu_{\ell,2}^2}(\hat{\mu}_{\ell,2}-\mu_{\ell,2})\\
  &-\frac{(\hat{\mu}_{\ell,1}-\mu_{\ell,1})(\hat{\mu}_{\ell,2}-\mu_{\ell,2})}{\mu_{\ell,2}\hat{\mu}_{\ell,2}}+\frac{\mu_{\ell,1}(\hat{\mu}_{\ell,2}-\mu_{\ell,2})^2}{\hat{\mu}_{\ell,2}{\mu}_{\ell,2}^2}\,,
	\end{split}
	\]
	which implies
	\[
	\begin{split}
		\max_{\ell \in \mathcal{S}}\bigg|\frac{1}{\zeta_{\ell}}(\hat{\zeta}_{\ell}-\zeta_{\ell})\bigg|
		\leq\bigg(\frac{\max_{\ell \in \mathcal{S}}|\hat{\mu}_{\ell,1}-\mu_{\ell,1}|}{\min_{\ell\in\mathcal{S}}\mu_{\ell,1}}+\frac{\max_{\ell \in \mathcal{S}}|\hat{\mu}_{\ell,2}-\mu_{\ell,2}|}{\min_{\ell\in\mathcal{S}}\mu_{\ell,2}}\bigg)\cdot O_{\p}(1)=o_{\p}(1)\,
	\end{split}
	\]
	provided that \begin{align}\label{eq:maxmulS}
		\max_{\ell\in\mathcal{S}}|\hat{\mu}_{\ell,1}-\mu_{\ell,1}|=o_{\p}\Big(\min_{\ell\in\mathcal{S}}\mu_{\ell,1}\Big)~~\textrm{and}~~\max_{\ell\in\mathcal{S}}|\hat{\mu}_{\ell,2}-\mu_{\ell,2}|=o_{\p}\Big(\min_{\ell\in\mathcal{S}}\mu_{\ell,2}\Big)\,.
	\end{align}
	Since $\theta_{\ell}=\log(\zeta_{\ell})/2$ and $\hat{\theta}_{\ell}=\log(\hat{\zeta}_{\ell})/2$, under \eqref{eq:maxmulS}, it follows from the Taylor expansion that
	\begin{align}\label{eq:thetaS}
		\notag\max_{\ell \in \mathcal{S}}|\hat{\theta}_{\ell}-\theta_{\ell}|\leq&~ O_{\p}(1)\cdot\max_{\ell \in \mathcal{S}}\bigg|\frac{1}{\zeta_{\ell}}(\hat{\zeta}_{\ell}-\zeta_{\ell})\bigg|\\
  \leq&~\bigg(\frac{\max_{\ell \in \mathcal{S}}|\hat{\mu}_{\ell,1}-\mu_{\ell,1}|}{\min_{\ell\in\mathcal{S}}\mu_{\ell,1}}+\frac{\max_{\ell \in \mathcal{S}}|\hat{\mu}_{\ell,2}-\mu_{\ell,2}|}{\min_{\ell\in\mathcal{S}}\mu_{\ell,2}}\bigg)\cdot O_{\p}(1)\,.
	\end{align}
	Analogously, we also have
	\begin{align}\label{eq:thetaSc}
		\max_{\ell \in \mathcal{S}^c}|\hat{\theta}_{\ell}-\theta_{\ell}|\leq \bigg(\frac{\max_{\ell \in \mathcal{S}^c}|\hat{\mu}_{\ell,1}-\mu_{\ell,1}|}{\min_{\ell\in\mathcal{S}^c}\mu_{\ell,1}}+\frac{\max_{\ell \in \mathcal{S}^c}|\hat{\mu}_{\ell,2}-\mu_{\ell,2}|}{\min_{\ell\in\mathcal{S}^c}\mu_{\ell,2}}\bigg)\cdot O_{\p}(1)
	\end{align}
	provided that \begin{align}\label{eq:maxmulSc}
	\max_{\ell\in\mathcal{S}^c}|\hat{\mu}_{\ell,1}-\mu_{\ell,1}|=o_{\p}\Big(\min_{\ell\in\mathcal{S}^c}\mu_{\ell,1}\Big)~~\textrm{and}~~\max_{\ell\in\mathcal{S}^c}|\hat{\mu}_{\ell,2}-\mu_{\ell,2}|=o_{\p}\Big(\min_{\ell\in\mathcal{S}^c}\mu_{\ell,2}\Big)\,.
\end{align}
As shown in Section \ref{sec:S}, we need to require $(\omega_1,\omega_2)$ to satisfy $0\leq \omega_1-\omega_2<1/2$ and $0\leq \omega_2\leq \omega_1<1$, and the restrictions in \eqref{eq:maxmulS} can be simplified as
\begin{equation}\label{eq:gamms11}
 \gamma\gg\left\{\begin{aligned}
&sp^{-3/2+\min(2\omega_1-2\omega_2,\,\omega_1)}\chi_{p}^{-8}\log^{1/2}s\,,\\
&p^{-1/3+\max(2\omega_1-2\omega_2, \,\omega_1)/3}\chi_{p}^{-5/3}\log ^{1/6}s\,,  \\
 &s^{1/2}p^{-1+\omega_2}\chi_{p}^{-7}\log ^{1/2}s \,.
\end{aligned}\right.
\end{equation}
As shown in Section \ref{sec:Sc}, we need to require $(\omega_1,\omega_2)$ to satisfy  $0\leq \omega_2\leq \omega_1<1/2$, and the restrictions in \eqref{eq:maxmulSc} can be simplified as
\begin{equation}\label{eq:gammasc2}
\gamma\gg \left\{\begin{aligned}
&s^{3/2}p^{-2+\min(2\omega_1-\omega_2,\,\omega_1+\omega_2)}\chi_{p}^{-8}\log^{1/2}p\,,\\
 &sp^{-3/2+\omega_1+\omega_2}\chi_{p}^{-5}\log^{1/2}p\,,\\
& p^{-1/3+2\omega_1/3}\chi_{p}^{-1}\log^{1/6}p\,.
\end{aligned}\right.
\end{equation}
Due to $s\ll p$ and $\chi_p^{-1}=\exp\{o(\log p)\}$, if $0\leq \omega_2\leq \omega_1<1/2$, we have
 \begin{align*}
&sp^{-3/2+\min(2\omega_1-2\omega_2,\,\omega_1)}\chi_{p}^{-8}\log^{1/2}s
\lesssim
p^{-1/3+2\omega_1/3}\chi_{p}^{-1}\log^{1/6}p\,,
\\
&~~~~~~~~~s^{1/2}p^{-1+\omega_2}\chi_{p}^{-7}\log ^{1/2}s\lesssim p^{-1/3+2\omega_1/3}\chi_{p}^{-1}\log^{1/6}p\,.
\end{align*}
Therefore, combining \eqref{eq:gamms11} and \eqref{eq:gammasc2}, if $0\leq \omega_2\leq \omega_1<1/2$, \eqref{eq:maxmulS} and \eqref{eq:maxmulSc} hold provided that
\begin{equation}\label{eq:gammaf}
 \gamma\gg\left\{\begin{aligned}
&p^{-1/3+\max(2\omega_1-2\omega_2, \,\omega_1)/3}\chi_{p}^{-5/3}\log ^{1/6}s\,,  \\
 &sp^{-3/2+\omega_1+\omega_2}\chi_{p}^{-5}\log^{1/2}p\,,\\
 &s^{3/2}p^{-2+\min(2\omega_1-\omega_2,\,\omega_1+\omega_2)}\chi_{p}^{-8}\log^{1/2}p\,,\\
& p^{-1/3+2\omega_1/3}\chi_{p}^{-1}\log^{1/6}p\,.
\end{aligned}\right.
\end{equation}
Under \eqref{eq:gammaf} with $0\leq \omega_2\leq \omega_1<1/2$, by \eqref{eq:thetaS} and Lemmas \ref{la:mu1mul2} and \ref{la:9},  it holds that
\begin{align*}
     \max_{\ell \in \mathcal{S}}|\hat{\theta}_{\ell}-\theta_{\ell}|\notag=&~O_{\p}\big(\gamma^{-1}p^{-1/2+\omega_1-\omega_2}\chi_{p}^{-7}\log ^{1/2}s\big)+O_{\p}\big(\gamma^{-2}p^{-1+2\omega_1-2\omega_2}\chi_{p}^{-5}\log s\big)
     \\
     &+O_{\p}\big\{\gamma^{-1}sp^{-3/2+\min(2\omega_1-2\omega_2,\,\omega_1)}\chi_{p}^{-8}\log^{1/2}s\big\}
      \\&+O_{\p}\big\{\gamma^{-3}p^{-1+\max(2\omega_1-2\omega_2,\, \omega_1)}\chi_{p}^{-5}\log ^{1/2}s\big\}
      \\&+O_{\p}\big(\gamma^{-1}s^{1/2}p^{-1+\omega_2}\chi_{p}^{-7}\log^{1/2}s\big)+O_{\p}\big(\gamma^{-1}p^{-1/2}\chi_{p}^{-6}\log^{1/2} s\big)
      \\&+O_{\p}\big\{\gamma^{-1}s^{3/2}p^{-2+\min(2\omega_2,\,2\omega_1-2\omega_2)}\chi_{p}^{-8}\log^{1/2}s\big\}
       \\&+O_{\p}\big\{\gamma^{-1}sp^{-3/2+\min(\omega_1-\omega_2,\,\omega_2)}\chi_{p}^{-10}\log^{1/2}s\big\}
         \\&+O_{\p}\big[\gamma^{-2} s p^{-2+\omega_1}\{p^{2\omega_2-\omega_1}\chi_{p}^{-3}I(\omega_1>2\omega_2)+ I(\omega_1\leq2\omega_2)\}\chi_{p}^{-5}\log s\big]
         \\&+O_{\p}\big[\gamma^{-2} s^{1/2} p^{-3/2+\omega_1}\{p^{\omega_2-\omega_1}\chi_{p}^{-2}I(\omega_1>\omega_2)+I(\omega_1=\omega_2)\}\chi_{p}^{-5}\log s\big]
         \\&+O_{\p}\big[\gamma^{-2} p^{-1+\omega_1}\{p^{-\omega_1}\chi_{p}^{-1}I(\omega_1>0)+ I(\omega_1=0)\}\chi_{p}^{-5}\log s\big]\,,\\
    \max_{\ell \in \mathcal{S}^c}|\hat{\theta}_{\ell}-\theta_{\ell}|=
    &~O_{\p}\big(\gamma^{-1}s^{1/2}p^{-1+\omega_1}\chi_{p}^{-4}\log ^{1/2}p\big)+O_{\p}\big(\gamma^{-1}sp^{-3/2+\omega_1+\omega_2}\chi_{p}^{-5}\log^{1/2}p\big) \\
    &+O_{\p}\big\{\gamma^{-1}s^{3/2}p^{-2+\min(2\omega_1-\omega_2,\,\omega_1+\omega_2)}\chi_{p}^{-8}\log^{1/2}p\big\}\\
    &+O_{\p}\big(\gamma^{-1}p^{-1/2+\omega_1}\chi_{p}^{-4}\log ^{1/2}p\big)+O_{\p}\big(\gamma^{-2}p^{-1+2\omega_1}\chi_{p}^{-3}\log p\big) \\
    &+O_{\p}\big(\gamma^{-1}s^{1/2}p^{-1+\omega_2}\chi_{p}^{-5}\log^{1/2}p\big)+O_{\p}\big(\gamma^{-3}p^{-1+2\omega_1}\chi_{p}^{-3}\log^{1/2}p\big)
    \\
    &+O_{\p}\big\{\gamma^{-1}s^{3/2}p^{-2+\min(2\omega_2,\,\omega_1)} \chi_{p}^{-6}\log^{1/2}p\big\}\,.
\end{align*}
Due to $\chi_{p}\in (0,1]$ and $s\ll p$, under \eqref{eq:gammaf}  with $0\leq \omega_2\leq \omega_1<1/2$, we have
\begin{align*}
&\gamma^{-2}p^{-1+2\omega_1-2\omega_2}\chi_{p}^{-5}\log s +\gamma^{-1}p^{-1/2}\chi_{p}^{-6}\log^{1/2} s\\
&~~~~~~~~~~~~~~~~~~~~\lesssim \gamma^{-1}p^{-1/2+\omega_1-\omega_2}\chi_{p}^{-7}\log ^{1/2}s\,,\\
&\gamma^{-1}s^{3/2}p^{-2+\min(2\omega_2,\,2\omega_1-2\omega_2)}\chi_{p}^{-8}\log^{1/2}s\\
&~~~~~~~~~~~~~~~~~~~~\lesssim \gamma^{-1}sp^{-3/2+\min(2\omega_1-2\omega_2,\,\omega_1)}\chi_{p}^{-8}\log^{1/2}s\,,\\
&\gamma^{-2} s p^{-2+\omega_1}\{p^{2\omega_2-\omega_1}\chi_{p}^{-3}I(\omega_1>2\omega_2)+ I(\omega_1\leq2\omega_2)\}\chi_{p}^{-5}\log s \\&~~~~~~~~~~~~~~~~~~~~\lesssim \gamma^{-1}s^{1/2}p^{-1+\omega_2}\chi_{p}^{-7}\log^{1/2}s\,,\\
&\gamma^{-2} s^{1/2} p^{-3/2+\omega_1}\{p^{\omega_2-\omega_1}\chi_{p}^{-2}I(\omega_1>\omega_2)+I(\omega_1=\omega_2)\}\chi_{p}^{-5}\log s\\&~~~~~~~~~~~~~~~~~~~~\lesssim \gamma^{-1}s^{1/2}p^{-1+\omega_2}\chi_{p}^{-7}\log^{1/2}s\,,\\
&\gamma^{-2} p^{-1+\omega_1}\{p^{-\omega_1}\chi_{p}^{-1}I(\omega_1>0)+ I(\omega_1=0)\}\chi_{p}^{-5}\log s\\
&~~~~~~~~~~~~~~~~~~~~\lesssim \gamma^{-1}p^{-1/2+\omega_1-\omega_2}\chi_{p}^{-7}\log ^{1/2}s\,,\\
&\gamma^{-1}s^{1/2}p^{-1+\omega_1}\chi_{p}^{-4}\log ^{1/2}p +\gamma^{-2}p^{-1+2\omega_1}\chi_{p}^{-3}\log p \\
&~~~~~~~~~~~~~~~~~~~~\lesssim \gamma^{-1}p^{-1/2+\omega_1}\chi_{p}^{-4}\log ^{1/2}p\,,\\
&\gamma^{-1}s^{3/2}p^{-2+\min(2\omega_2,\,\omega_1)} \chi_{p}^{-6}\log^{1/2}p \\
&~~~~~~~~~~~~~~~~~~~~\lesssim \gamma^{-1}s^{3/2}p^{-2+\min(2\omega_1-\omega_2,\,\omega_1+\omega_2)}\chi_{p}^{-8}\log^{1/2}p\,,
\end{align*}
which implies
\begin{align*}
    \max_{\ell \in \mathcal{S}}|\hat{\theta}_{\ell}-\theta_{\ell}|\notag=&~O_{\p}\big(\gamma^{-1}p^{-1/2+\omega_1-\omega_2}\chi_{p}^{-7}\log ^{1/2}s\big)+O_{\p}\big(\gamma^{-1}s^{1/2}p^{-1+\omega_2}\chi_{p}^{-7}\log ^{1/2}s\big)\\
    &+O_{\p}\big\{\gamma^{-1}sp^{-3/2+\min(2\omega_1-2\omega_2,\,\omega_1)}\chi_{p}^{-8}\log^{1/2}s\big\}
    \\
    &+O_{\p}\big\{\gamma^{-1}sp^{-3/2+\min(\omega_1-\omega_2,\,\omega_2)}\chi_{p}^{-10}\log^{1/2}s\big\}\\
    &+O_{\p}\big\{\gamma^{-3}p^{-1+\max(2\omega_1-2\omega_2,\, \omega_1)}\chi_{p}^{-5}\log ^{1/2}s\big\}\,,\\
    \max_{\ell \in \mathcal{S}^c}|\hat{\theta}_{\ell}-\theta_{\ell}|=&~O_{\p}\big\{\gamma^{-1}s^{3/2}p^{-2+\min(2\omega_1-\omega_2,\,\omega_1+\omega_2)}\chi_{p}^{-8}\log^{1/2}p\big\}\\
    &+O_{\p}\big(\gamma^{-1}p^{-1/2+\omega_1}\chi_{p}^{-4}\log ^{1/2}p\big)+O_{\p}\big(\gamma^{-1}sp^{-3/2+\omega_1+\omega_2}\chi_{p}^{-5}\log^{1/2}p\big)\\
    &+O_{\p}\big(\gamma^{-1}s^{1/2}p^{-1+\omega_2}\chi_{p}^{-5}\log^{1/2}p\big)+O_{\p}\big(\gamma^{-3}p^{-1+2\omega_1}\chi_{p}^{-3}\log^{1/2}p\big)\,.
\end{align*}
Hence, it holds that
\begin{align*}
    \max_{\ell \in [p]}|\hat{\theta}_{\ell}-\theta_{\ell}|\notag=
    &~\tilde{O}_{\p}\big(\gamma^{-1}p^{-1/2+\omega_1}\log ^{1/2}p\big)+\tilde{O}_{\p}\big(\gamma^{-1}sp^{-3/2+\omega_1+\omega_2}\log^{1/2}p\big)\\
    &
    +\tilde{O}_{\p}\big(\gamma^{-3}p^{-1+2\omega_1}\log^{1/2}p\big)
\end{align*}
provided that
\begin{align*}
 \gamma\gg\chi_p^{-8}\cdot\left\{\begin{aligned}
 &sp^{-3/2+\omega_1+\omega_2}\log^{1/2}p\,,\\
& p^{-1/3+2\omega_1/3}\log^{1/6}p\,.
\end{aligned}\right.
\end{align*}
We complete the proof of Proposition \ref{tm:4}. $\hfill\Box$

\subsection{Sufficient conditions for \eqref{eq:maxmulS}. }\label{sec:S}

By Lemmas \ref{la:mu1mul2} and \ref{la:9}, to make \eqref{eq:maxmulS} hold, it suffices to require
\begin{equation*}
\left\{\begin{aligned}
 &\gamma^2p^{-1/2+\omega_2-\omega_1}\chi_{p}^{-2}\log^{1/2}s\ll p^{2\omega_2-2\omega_1}\chi_{p}^{5} \gamma^3\,,\\
&\gamma^2sp^{-3/2+\min(2\omega_2-\omega_1,\,0)} \chi_{p}^{-3}\log^{1/2}s\ll p^{2\omega_2-2\omega_1}\chi_{p}^{5} \gamma^3\,,  \\
 &\gamma p^{-1}\log s \ll p^{2\omega_2-2\omega_1}\chi_{p}^{5} \gamma^3 \,,  \\
 &p^{-1}\log^{1/2}s \ll p^{2\omega_2-2\omega_1}\chi_{p}^{5} \gamma^3\,,\\
 &\gamma^2s^{1/2}p^{-1+\omega_2-\omega_1}\chi_{p}^{-2}\log^{1/2}s \ll p^{-\omega_1}\chi_{p}^{5} \gamma^3\,,\\
&\gamma^2s^{3/2}p^{-2+\min(2\omega_2-\omega_1,\,\omega_1-2\omega_2)}\chi_{p}^{-3}\log^{1/2}s \ll p^{-\omega_1}\chi_{p}^{5} \gamma^3  \,,  \\
 &\gamma^2p^{-1/2-\omega_1}\chi_{p}^{-1}\log^{1/2}s \ll p^{-\omega_1}\chi_{p}^{5} \gamma^3\,,
 \\&\gamma^2sp^{-3/2+\min(-\omega_2,\,\omega_2-\omega_1)} \chi_{p}^{-5}\log^{1/2}s\ll p^{-\omega_1}\chi_{p}^{5} \gamma^3\,,
 \\&\gamma s p^{-2}\{p^{2\omega_2-\omega_1}\chi_{p}^{-3}I(\omega_1>2\omega_2)+ I(\omega_1\leq2\omega_2)\}\log s\ll p^{-\omega_1}\chi_{p}^{5} \gamma^3\,,
 \\&\gamma s^{1/2} p^{-3/2}\{p^{\omega_2-\omega_1}\chi_{p}^{-2}I(\omega_1>\omega_2)+I(\omega_1=\omega_2)\}\log s \ll p^{-\omega_1}\chi_{p}^{5} \gamma^3\,,
 \\&\gamma p^{-1}\{p^{-\omega_1}\chi_{p}^{-1}I(\omega_1>0)+ I(\omega_1=0)\}\log s\, \ll p^{-\omega_1}\chi_{p}^{5} \gamma^3 \,,
 \\&p^{-1}\log^{1/2}s\ll p^{-\omega_1}\chi_{p}^{5} \gamma^3\,,
\end{aligned}\right.
\end{equation*}
%
which is equivalent to
\begin{equation}\label{eq:gamms1}
 \gamma\gg\left\{\begin{aligned}
&p^{-1/2+\omega_1-\omega_2}\chi_{p}^{-7} \log^{1/2}s\,,\\
&sp^{-3/2+\min(2\omega_1-2\omega_2,\,\omega_1)}\chi_{p}^{-8}\log^{1/2}s\,,\\
&p^{-1/3+\max(2\omega_1-2\omega_2, \,\omega_1)/3}\chi_{p}^{-5/3}\log ^{1/6}s\,,  \\
 &s^{1/2}p^{-1+\omega_2}\chi_{p}^{-7}\log ^{1/2}s \,,  \\
 &s^{3/2}p^{-2+\min(2\omega_2,\,2\omega_1-2\omega_2)}\chi_{p}^{-8}\log ^{1/2}s\,,\\
 &sp^{-3/2+\min(\omega_1-\omega_2,\,\omega_2)}\chi_{p}^{-10}\log^{1/2}s\,,\\
 &s^{1/4}p^{-3/4}\chi_{p}^{-5/2}\{p^{\omega_2/2}\chi_p^{-1}I(\omega_1>\omega_2)+p^{\omega_1/2}I(\omega_1=\omega_2)\}\log ^{1/2}s\,.
\end{aligned}\right.
\end{equation}
Due to $s\leq p$, $\omega_1\geq\omega_2$, $\omega_2\in[0,1)$ and $\chi_p^{-1}=\exp\{o(\log p)\}$, then
\begin{align*}
&s^{1/4}p^{-3/4}\chi_{p}^{-5/2}\{p^{\omega_2/2}\chi_p^{-1}I(\omega_1>\omega_2)+p^{\omega_1/2}I(\omega_1=\omega_2)\}\log ^{1/2}s\\
&~~~~~~~~~~~~~\leq s^{1/4}p^{-3/4+\omega_2/2}\chi_{p}^{-7/2}\log ^{1/2}s\\
&~~~~~~~~~~~~~\lesssim p^{-1/3+\omega_2/3}\chi_p^{-5/3}\log^{1/6}s\leq p^{-1/3+\max(2\omega_1-2\omega_2, \,\omega_1)/3}\chi_{p}^{-5/3}\log ^{1/6}s\,,\\
	&sp^{-3/2+\min(\omega_1-\omega_2,\,\omega_2)}\chi_{p}^{-10}\log^{1/2}s\\
 &~~~~~~~~~~~~~\lesssim 	\left\{\begin{aligned}
	 sp^{-3/2+\min(2\omega_1-2\omega_2,\,\omega_1)}\chi_{p}^{-8}\log^{1/2}s\,,~~~~&\textrm{if}~\omega_1>\omega_2\,,
\\
  p^{-1/3+\max(2\omega_1-2\omega_2, \,\omega_1)/3}\chi_{p}^{-5/3}\log ^{1/6}s\,,~&\textrm{if}~\omega_1=\omega_2\,,
		\end{aligned}\right.
 \\& s^{3/2}p^{-2+\min(2\omega_2,\,2\omega_1-2\omega_2)}\chi_{p}^{-8}\log ^{1/2}s
 \\
 &~~~~~~~~~~~~~\lesssim \left\{\begin{aligned}
sp^{-3/2+\omega_1}\chi_{p}^{-8}\log^{1/2}s\,,~~~~~&\textrm{if}~\omega_1\geq 2\omega_2\,,
\\
  sp^{-3/2+2\omega_1-2\omega_2}\chi_{p}^{-8}\log^{1/2}s\,,~&\textrm{if}~\omega_1<2\omega_2\,,
		\end{aligned}\right.
  \\
 &~~~~~~~~~~~~~=sp^{-3/2+\min(2\omega_1-2\omega_2,\,\omega_1)}\chi_{p}^{-8}\log^{1/2}s\,.
\end{align*}
Notice that $\gamma=O(1)$ and $0\leq\omega_1-\omega_2<1$ with $\omega_2\in[0,1)$. We need to require $(\omega_1,\omega_2)$ to satisfy $0\leq \omega_1-\omega_2<1/2$ and $0\leq \omega_2\leq \omega_1<1$. Under such restrictions, due to $\chi_p^{-1}=\exp\{o(\log p)\}$, we have
\begin{align*}
p^{-1/2+\omega_1-\omega_2}\chi_{p}^{-7} \log^{1/2}s\lesssim p^{-1/3+\max(2\omega_1-2\omega_2, \,\omega_1)/3}\chi_{p}^{-5/3}\log ^{1/6}s\,.
\end{align*}
Then \eqref{eq:gamms1} can be simplified as
\begin{equation*}
 \gamma\gg\left\{\begin{aligned}
&sp^{-3/2+\min(2\omega_1-2\omega_2,\,\omega_1)}\chi_{p}^{-8}\log^{1/2}s\,,\\
&p^{-1/3+\max(2\omega_1-2\omega_2, \,\omega_1)/3}\chi_{p}^{-5/3}\log ^{1/6}s\,,  \\
 &s^{1/2}p^{-1+\omega_2}\chi_{p}^{-7}\log ^{1/2}s \,,
\end{aligned}\right.
\end{equation*}
which gives the sufficient conditions for \eqref{eq:maxmulS}. $\hfill\Box$

\subsection{Sufficient conditions for \eqref{eq:maxmulSc}. }\label{sec:Sc}
By Lemmas \ref{la:mu1mul2} and \ref{la:9}, to make \eqref{eq:maxmulSc} hold, it suffices to require
\begin{equation*}\label{eq:gammasc1}
\left\{\begin{aligned}
&\gamma^2s^{1/2}p^{-1-\omega_1}\chi_{p}^{-1}\log^{1/2}p\ll p^{-2\omega_1}\chi_{p}^3\gamma^3\,,\\
&\gamma^2s^{3/2}p^{-2+\min(-\omega_2,\,\omega_2-\omega_1)} \chi_{p}^{-5}\log^{1/2}p\ll p^{-2\omega_1}\chi_{p}^3\gamma^3\,,\\
 &\gamma^2p^{-1/2-\omega_1}\chi_{p}^{-1}\log^{1/2}p\ll p^{-2\omega_1}\chi_{p}^3\gamma^3\,,\\
&\gamma^2sp^{-3/2+\omega_2-\omega_1} \chi_{p}^{-2}\log^{1/2}p\ll p^{-2\omega_1}\chi_{p}^3\gamma^3 \,,  \\
 & \gamma p^{-1}\log p \ll p^{-2\omega_1}\chi_{p}^3\gamma^3\,,  \\
 &p^{-1}\log^{1/2}p \ll p^{-2\omega_1}\chi_{p}^3\gamma^3\,,\\
  &\gamma^2s^{1/2}p^{-1+\omega_2-\omega_1}\chi_{p}^{-2}\log^{1/2}p \ll p^{-\omega_1}\chi_{p}^3\gamma^3\,,\\
& \gamma^2s^{3/2}p^{-2+\min(2\omega_2-\omega_1,\,0)} \chi_{p}^{-3}\log^{1/2}p \ll p^{-\omega_1}\chi_{p}^3\gamma^3 \,,
\end{aligned}\right.
\end{equation*}
which is equivalent to
\begin{equation}\label{eq:gammasc1}
\gamma\gg \left\{\begin{aligned}
&s^{3/2}p^{-2+\min(2\omega_1-\omega_2,\,\omega_1+\omega_2)}\chi_{p}^{-8}\log^{1/2}p\,,\\
 &p^{-1/2+\omega_1}\chi_{p}^{-4}\log ^{1/2}p\,,\\
 &sp^{-3/2+\omega_1+\omega_2}\chi_{p}^{-5}\log^{1/2}p\,,\\
& p^{-1/3+2\omega_1/3}\chi_{p}^{-1}\log^{1/6}p\,,  \\
&s^{1/2}p^{-1+\omega_2}\chi_{p}^{-5}\log^{1/2}p\,,\\
&s^{3/2}p^{-2+\min(2\omega_2,\,\omega_1)}\chi_{p}^{-6}\log^{1/2}p\,.
\end{aligned}\right.
\end{equation}
Due to $s\leq p$, $\omega_1\geq\omega_2$, $\omega_2\in[0,1)$ and $\chi_p^{-1}=\exp\{o(\log p)\}$, then
\begin{align*}
&s^{3/2}p^{-2+\min(2\omega_2,\,\omega_1)}\chi_{p}^{-6}\log^{1/2}p
 \\
 &~~~~~~~~~~~~~\lesssim \left\{\begin{aligned}
s^{3/2}p^{-2+\omega_1+\omega_2}\chi_{p}^{-8}\log^{1/2}p\,,~~&\textrm{if}~\omega_1\geq 2\omega_2\,,
\\
  s^{3/2}p^{-2+2\omega_1-\omega_2}\chi_{p}^{-8}\log^{1/2}p\,,~&\textrm{if}~\omega_1<2\omega_2\,,
		\end{aligned}\right.
  \\
 &~~~~~~~~~~~~~=s^{3/2}p^{-2+\min(2\omega_1-\omega_2,\,\omega_1+\omega_2)}\chi_{p}^{-8}\log^{1/2}p\,.
\end{align*}
Notice that $\gamma=O(1)$ and $0\leq\omega_1-\omega_2<1$ with $\omega_2\in[0,1)$. We need to require $(\omega_1,\omega_2)$ to satisfy  $0\leq \omega_2\leq \omega_1<1/2$. Under such restrictions, due to $\chi_p^{-1}=\exp\{o(\log p)\}$, we have
\begin{align*}
p^{-1/2+\omega_1}\chi_{p}^{-4}\log ^{1/2}p\lesssim&~ p^{-1/3+2\omega_1/3}\chi_{p}^{-1}\log^{1/6}p\,,\\
s^{1/2}p^{-1+\omega_2}\chi_{p}^{-5}\log^{1/2}p\lesssim&~ p^{-1/3+2\omega_1/3}\chi_{p}^{-1}\log^{1/6}p\,.
\end{align*}
Then \eqref{eq:gammasc1} can be simplified as
\begin{equation*}
\gamma\gg \left\{\begin{aligned}
&s^{3/2}p^{-2+\min(2\omega_1-\omega_2,\,\omega_1+\omega_2)}\chi_{p}^{-8}\log^{1/2}p\,,\\
 &sp^{-3/2+\omega_1+\omega_2}\chi_{p}^{-5}\log^{1/2}p\,,\\
& p^{-1/3+2\omega_1/3}\chi_{p}^{-1}\log^{1/6}p\,,
\end{aligned}\right.
\end{equation*}
which gives the sufficient conditions for \eqref{eq:maxmulSc}. $\hfill\Box$

\section{Proofs of Lemmas \ref{la:3}--\ref{la:9}.}
\subsection{Proof of Lemma {\rm\ref{la:3}}.}\label{sec:pflem1} 	 For any $i,j,\ell\in[p]$, let $\mathring{\psi}_1(i,j;\ell)=\psi_1(i,j;\ell)-\mathbb{E}\{\psi_1(i,j;\ell)\}$ and  $\mathring{\psi}_2(i,j;\ell)=\psi_2(i,j;\ell)-\mathbb{E}\{\psi_2(i,j;\ell)\}$, where  $\psi_1(i,j;\ell)=\varphi_{(i,\ell),1}\varphi_{(i,j),0}\varphi_{(\ell,j),1}$ and $\psi_2(i,j;\ell)=\varphi_{(i,\ell),0}\varphi_{(i,j),1}\varphi_{(\ell,j),0}$. Then
	 \[
	 \hat{\mu}_{\ell,1}-\mu_{\ell,1}=\frac{1}{|\mathcal{H}_\ell|}\sum_{(i,j)\in\mathcal{H}_\ell}\mathring{\psi}_1(i,j;\ell)~~\textrm{and}~~\hat{\mu}_{\ell,2}-\mu_{\ell,2}=\frac{1}{|\mathcal{H}_\ell|}\sum_{(i,j)\in\mathcal{H}_\ell}\mathring{\psi}_2(i,j;\ell)\,.
	 \]Write  $\mathscr{F}_{\ell}=\{Z_{i,\ell},Z_{\ell,j}:(i,j)\in\mathcal{H}_{\ell}\}$. For any $(i,j)\in\mathcal{H}_\ell$, we have $\mathbb{E}\{\psi_1(i,j;\ell)\,|\,\mathscr{F}_\ell\}=\varphi_{(i,\ell),1}\varphi_{(\ell,j),1}\mathbb{E}\{\varphi_{(i,j),0}\}$ and
	\[
	\begin{split}
		\frac{1}{|\mathcal{H}_\ell|}\sum_{(i,j)\in\mathcal{H}_\ell}\mathring{\psi}_1(i,j;\ell)=&~\underbrace{\frac{1}{|\mathcal{H}_\ell|}\sum_{(i,j)\in\mathcal{H}_\ell}[\psi_1(i,j;\ell)-\mathbb{E}\{\psi_1(i,j;\ell)\,|\,\mathscr{F}_\ell\}]}_{I_{\ell,1,1}}\\
		&+\underbrace{\frac{1}{|\mathcal{H}_\ell|}\sum_{(i,j)\in\mathcal{H}_\ell}[\mathbb{E}\{\psi_1(i,j;\ell)\,|\,\mathscr{F}_\ell\}-\mathbb{E}\{\psi_1(i,j;\ell)\}]}_{I_{\ell,1,2}}\,.
	\end{split}
	\]
	As we will show in Sections \ref{sec:pflem11} and \ref{sec:pflem12}, respectively, that
		\begin{align}
	    &~~~~~~~~~~~~~~~\max_{\ell\in[p]}|I_{\ell,1,1}|=O_{\p}\bigg(\frac{\log^{1/2}p}{p}\bigg)\,,\label{eq:lem11}\\
	   &\max_{\ell\in[p]}|I_{\ell,1,2}|=O_{\p}\bigg(\frac{\gamma^2\log^{1/2}p}{p^{1/2}}\bigg)+O_{\p}\bigg(\frac{\gamma\log p}{p}\bigg)\,, \label{eq:lem12}
	\end{align}
	we then have
	\[		\max_{\ell\in[p]}|\hat{\mu}_{\ell,1}-\mu_{\ell,1}|=O_{\p}\bigg(\frac{\log^{1/2}p}{p}\bigg)+O_{\p}\bigg(\frac{\gamma^2\log^{1/2}p}{p^{1/2}}\bigg)+O_{\p}\bigg(\frac{\gamma\log p}{p}\bigg)\,.
	\]
Similarly, we can also prove another result. $\hfill\Box$	
	

	\subsubsection{Convergence rate of $\max_{\ell\in[p]}|I_{\ell,1,1}|$.}\label{sec:pflem11}
	
	Conditional on $\mathscr{F}_\ell$,  $\{\psi_1(i,j;\ell)\}_{(i,j)\in\mathcal{H}_\ell}$ is an independent sequence. For any $(i,j)\in\mathcal{H}_{\ell}$, write
	$$
	\sigma_{(i,j),\ell,1}^2:=\textrm{Var}\{{\psi}_1(i,j;\ell)\,|\,\mathscr{F}_{\ell}\}=\varphi_{(i,\ell),1}^2\varphi_{(\ell,j),1}^2{\rm Var}(Z_{i,j})\,.$$
	Due to $(\alpha,\beta)
	\in\mathcal{M}(\gamma;C_1)$, we have
	\[
	\begin{split}
	&|{\rm Var}(Z_{i,j})-\alpha(1-\alpha)|\\
	&~~~~~~~\leq (1-\alpha-\beta)(1-2\alpha)\mathbb{E}(X_{i,j})+(1-\alpha-\beta)^2\mathbb{E}^2(X_{i,j})\leq 2
	\end{split}
	\]
	holds uniformly over $(i,j)\in\mathcal{H}_\ell$. Notice that $|\varphi_{(i,j),1}|\leq 1-\alpha$ and $|\varphi_{(i,j),0}|\leq 1-\beta$. Then $\max_{(i,j)\in\mathcal{H}_\ell}|\psi_1(i,j;\ell)|\leq C$ and $\max_{(i,j)\in\mathcal{H}_\ell}\sigma_{(i,j),\ell,1}^2\leq C$. Due to $|\mathcal{H}_\ell|=(p-1)(p-2)/2$, by Bernstein inequality, we have
	\begin{align*}
	    \mathbb{P}(|I_{\ell,1,1}|>u\,|\,\mathscr{F}_\ell)\lesssim \exp(-Cp^2u^2)
	\end{align*}
	for any $0<u\leq o(1)$,	which implies that
	\begin{align}\label{eq:bdIell11}
	    \mathbb{P}(|I_{\ell,1,1}|>u)=\mathbb{E}\{\mathbb{P}(|I_{\ell,1,1}|>u\,|\,\mathscr{F}_\ell)\}\lesssim \exp(-Cp^2u^2)
	\end{align}
	for any $0<u\leq o(1)$. Therefore, we have \eqref{eq:lem11}. $\hfill\Box$

	\subsubsection{Convergence rate of $\max_{\ell\in[p]}|I_{\ell,1,2}|$.} \label{sec:pflem12}
	
	Define $\mathring{\varphi}_{(i,j),\tau}=\varphi_{(i,j),\tau}-\mathbb{E}\{\varphi_{(i,j),\tau}\}$. Due to $\psi_1(i,j;\ell)=\psi_1(j,i;\ell)$ for any $i\neq j$, it then holds that
	\begin{align}
		(p-1)(p-2)I_{\ell,1,2}=&~\sum_{i,j:\,i\neq j,\,i,j\neq \ell}\big[\varphi_{(i,\ell),1}\varphi_{(\ell,j),1}-\mathbb{E}\{\varphi_{(i,\ell),1}\varphi_{(\ell,j),1}\}\big]\mathbb{E}\{\varphi_{(i,j),0}\}\notag\\
		=&~\underbrace{2\sum_{i,j:\,i\neq j,\,i,j\neq \ell}\mathring{\varphi}_{(i,\ell),1}\mathbb{E}\{\varphi_{(\ell,j),1}\}\mathbb{E}\{\varphi_{(i,j),0}\}}_{I_{\ell,1,2}(1)}\label{eq:Iell12}\\
		&+\underbrace{\sum_{i,j:\,i\neq j,\,i,j\neq \ell}\mathring{\varphi}_{(i,\ell),1}\mathring{\varphi}_{(\ell,j),1}\mathbb{E}\{\varphi_{(i,j),0}\}}_{I_{\ell,1,2}(2)}\,.\notag
	\end{align}

	For $I_{\ell,1,2}(1)$, we have
	\begin{align*}
	I_{\ell,1,2}(1)&=\sum_{i:\,i\neq \ell}\mathring{\varphi}_{(i,\ell),1}\bigg[2\sum_{j:\,j\neq i,\ell}\mathbb{E}\{\varphi_{(\ell,j),1}\}\mathbb{E}\{\varphi_{(i,j),0}\}\bigg]=:\sum_{i:\,i\neq \ell}\mathring{\varphi}_{(i,\ell),1}A_{i,\ell}\,.
	\end{align*}
Recall
	\begin{align*}
	    	\mathbb{E}\{\varphi_{(i,j),0}\}=&~\gamma\mathbb{P}(X_{i,j}=0)=\frac{\gamma}{1+\exp(\theta_i+\theta_j)}\,,\\
	    	 \mathbb{E}\{\varphi_{(i,j),1}\}=&~\gamma\mathbb{P}(X_{i,j}=1)=\frac{\gamma\exp(\theta_i+\theta_j)}{1+\exp(\theta_i+\theta_j)}\,.
	\end{align*}	By Condition \ref{cond1}, we have
 \begin{align*}
 \min_{i,j:\,i\neq j}\mathbb{E}\{\varphi_{(i,j),0}\}\asymp\gamma\asymp\max_{i,j:\,i\neq j}\mathbb{E}\{\varphi_{(i,j),0}\}\,,\\
 \min_{i,j:\,i\neq j}\mathbb{E}\{\varphi_{(i,j),1}\}\asymp\gamma\asymp\max_{i,j:\,i\neq j}\mathbb{E}\{\varphi_{(i,j),1}\}\,,
 \end{align*}
 which implies
 $$\min_{\ell\in [p]}\min_{i:\,i\neq \ell}A_{i,\ell}\asymp p\gamma^2\asymp\max_{\ell\in[p]}\max_{i:\,i\neq\ell}A_{i,\ell}\,.$$ Note that $\max_{\ell\in[p]}\max_{i:\,i\neq\ell}{\rm Var}\{\mathring{\varphi}_{(i,\ell),1}\}=\max_{\ell\in[p]}\max_{i:\,i\neq\ell}{\rm Var}(Z_{i,\ell})\leq C$. Given $\ell\in[p]$, since $\{\mathring{\varphi}_{(i,\ell),1}\}_{i:\,i\neq\ell}$ is an independent sequence, by Bernstein inequality,
	\begin{align}\label{eq:bdIell121}
	    \mathbb{P}\{|I_{\ell,1,2}(1)|>u\}\lesssim \exp\bigg(-\frac{Cu^2}{p^{3}\gamma^{4}}\bigg)
	\end{align}
	for any $0<u\leq o(p^2\gamma^2)$. Thus,
	\begin{align}\label{eq:Iell121}
	    \max_{\ell\in[p]}|I_{\ell,1,2}(1)|=O_{\p}(\gamma^2{p^{3/2}}\log^{1/2}p)\,.
	\end{align}

	For $I_{\ell,1,2}(2)$, letting $B_{(i,j)}=\gamma^{-1}\mathbb{E}\{\varphi_{(i,j),0}\}$, then
	$$
	\gamma^{-1}I_{\ell,1,2}(2)=\sum_{i,j:\,i\neq j,\,i,j\neq \ell}\mathring{\varphi}_{(i,\ell),1}\mathring{\varphi}_{(\ell,j),1}B_{(i,j)}\,.$$
	Under Condition \ref{cond1},
	\begin{align*}
	    	\min_{\ell\in[p]}\min_{i,j:\,i\neq j,\,i,j\neq \ell}B_{(i,j)}\asymp1\asymp \max_{\ell\in[p]}\max_{i,j:\,i\neq j,\,i,j\neq \ell}B_{(i,j)}\,.
	\end{align*}
 By the decoupling inequality of \cite{DM_1995} and Theorem 3.3 of \cite{Gineetal_2000}, we have
	\begin{align}\label{eq:bdIell122}
	    \max_{\ell\in[p]}\mathbb{P}\bigg\{\bigg|\sum_{i,j:\,i\neq j,\,i,j\neq \ell}\mathring{\varphi}_{(i,\ell),1}\mathring{\varphi}_{(\ell,j),1}B_{(i,j)}\bigg|>u\bigg\}\lesssim\exp\bigg(-\frac{Cu}{p}\bigg)
	\end{align}
	 for any $p\ll u\ll p^2$, 	
	which implies
	$
	\max_{\ell\in[p]}|I_{\ell,1,2}(2)|=O_{\p}({\gamma p\log p})$. Together with \eqref{eq:Iell121}, we can obtain \eqref{eq:lem12} by \eqref{eq:Iell12}. $\hfill\Box$


\subsection{Proof of Lemma {\rm\ref{la:4}}.}\label{sec:pflemma2}
As shown in Section \ref{sec:pflem1},
	$
	    \hat{\mu}_{\ell,1}-\mu_{\ell,1}=I_{\ell,1,1}+I_{\ell,1,2}$
for $I_{\ell,1,1}$ and $I_{\ell,1,2}$ specified there. By \eqref{eq:bdIell11}, $I_{\ell,1,1}=O_{\p}(p^{-1})$. Let $N=(p-1)(p-2)$. Recall $I_{\ell,1,2}=N^{-1}\{I_{\ell,1,2}(1)+I_{\ell,1,2}(2)\}$ with $I_{\ell,1,2}(1)$ and $I_{\ell,1,2}(2)$ defined in \eqref{eq:Iell12}. By \eqref{eq:bdIell121} and \eqref{eq:bdIell122}, $I_{\ell,1,2}(1)=O_{\p}(\gamma^2p^{3/2})$ and $I_{\ell,1,2}(2)=O_{\p}(\gamma p)$. Hence, $|\hat{\mu}_{\ell,1}-\mu_{\ell,1}|=O_{\p}(p^{-1})+O_{\p}(\gamma^2p^{-1/2})$. Analogously, we also have
$
	    |\hat{\mu}_{\ell,2}-\mu_{\ell,2}|=O_{\p}(p^{-1})+O_{\p}(\gamma^2p^{-1/2})$. For given $(\ell,i)$ such that $i\neq \ell$, we know $\{\varphi_{(\ell,j),1}\varphi_{(i,j),0}\}_{j:\,j\neq \ell,i}$ and $\{\varphi_{(\ell,j),0}\varphi_{(i,j),1}\}_{j:\,j\neq \ell,i}$ are two independent and bounded sequences. By Bernstein inequality,
	\begin{align*}
&\max_{i:\,i\neq\ell}\bigg|\frac{1}{p-2}\sum_{j:\,j\neq \ell,i}\big[\varphi_{(\ell,j),1}\varphi_{(i,j),0}-\mathbb{E}\{\varphi_{(\ell,j),1}\}\mathbb{E}\{\varphi_{(i,j),0}\}\big]\bigg|\\
&~~~~=O_{\p}\bigg(\frac{\log^{1/2}p}{p^{1/2}}\bigg)=\max_{i:\,i\neq\ell}\bigg|\frac{1}{p-2}\sum_{j:\,j\neq \ell,i}\big[\varphi_{(\ell,j),0}\varphi_{(i,j),1}-\mathbb{E}\{\varphi_{(\ell,j),0}\}\mathbb{E}\{\varphi_{(i,j),1}\}\big
 ]\bigg|\,.
\end{align*}
	Notice that $\mu_{\ell,1}\asymp\gamma^3\asymp\mu_{\ell,2}$. Based on the definition of $\lambda_{i,\ell}$ and $\hat{\lambda}_{i,\ell}$ given, respectively, in (\ref{eq:lambdaiell}) and (\ref{eq:hatlambdaiell}), we complete the proof of Lemma \ref{la:4}. $\hfill\Box$

\subsection{Proof of Lemma {\rm\ref{la:5}}.}\label{sec:pflemma3}
 Due to
\begin{align*}
{\rm Var}(Z_{i,\ell})=&\,(1-\beta)\beta+\frac{\gamma(\gamma+\alpha-\beta)}{1+\exp(\theta_i+\theta_\ell)}-\frac{\gamma^2}{\{1+\exp(\theta_i+\theta_\ell)\}^2}\,,\\
\widehat{{\rm Var}}(Z_{i,\ell})=&\,(1-\beta)\beta+\frac{\gamma(\gamma+\alpha-\beta)}{1+\exp(\hat{\theta}_i+\hat{\theta}_\ell)}-\frac{\gamma^2}{\{1+\exp(\hat{\theta}_i+\hat{\theta}_\ell)\}^2}\,,
\end{align*}
it holds that
\begin{align*}
|\widehat{{\rm Var}}(Z_{i,\ell})-{\rm Var}(Z_{i,\ell})|\leq&\, \gamma(\gamma+\alpha-\beta)\bigg|\frac{1}{1+\exp(\hat{\theta}_i+\hat{\theta}_\ell)}-\frac{1}{1+\exp(\theta_i+\theta_\ell)}\bigg|\\
&+2\gamma^2\bigg|\frac{1}{1+\exp(\hat{\theta}_i+\hat{\theta}_\ell)}-\frac{1}{1+\exp(\theta_i+\theta_\ell)}\bigg|\\
\leq&\,\gamma(3\gamma+\alpha-\beta)\bigg|\frac{1}{1+\exp(\hat{\theta}_i+\hat{\theta}_\ell)}-\frac{1}{1+\exp(\theta_i+\theta_\ell)}\bigg|\,.
\end{align*}
Define $f(x)=(1+e^x)^{-1}$ for $x\in\mathbb{R}$. Then $\sup_{x\in\mathbb{R}}|f'(x)|\leq 1$. By Proposition \ref{tm:1}, we have
	\begin{align*}
&\max_{i:\,i\neq \ell}\bigg|\frac{1}{1+\exp(\hat{\theta}_i+\hat{\theta}_\ell)}-\frac{1}{\exp(\theta_i+\theta_\ell)}\bigg|\\
&~~~~~~~~~~\leq 2\max_{\ell\in[p]}|\hat{\theta}_\ell-\theta_\ell|=O_{\p}\bigg(\frac{\log^{1/2}p}{\gamma^3p}\bigg)+O_{\p}\bigg(\frac{\log^{1/2}p}{\gamma p^{1/2}}\bigg)\,,
\end{align*}
which implies
\[
\max_{i:\,i\neq\ell}|\widehat{{\rm Var}}(Z_{i,\ell})-{\rm Var}(Z_{i,\ell})|=(3\gamma+\alpha-\beta)\bigg\{O_{\p}\bigg(\frac{\log^{1/2}p}{\gamma^2p}\bigg)+O_{\p}\bigg(\frac{\log^{1/2}p}{ p^{1/2}}\bigg)\bigg\}\,.
\]
We complete the proof of Lemma \ref{la:5}. $\hfill\Box$

	\subsection{Proof of Lemma {\rm\ref{la:vnorm}}.}\label{se:proofla4}
	
	Recall $\mathbb{E}\{Y_{(i,j),\ell}\}=0$ for any $(i,j,\ell)$ such that $i\neq j$ and $\ell\neq i,j$. For any $\ell_1,\ell_2\in[p]$, we have
	\begin{align*}
		B_{\ell_1,\ell_2}=&~\underbrace{\frac{1}{N}\sum_{i,j:\,i\neq j,\,i,j\neq\ell_1,\ell_2}\mathbb{E}\{Y_{(i,j),\ell_1}Y_{(i,j),\ell_2}\}}_{B_{\ell_1,\ell_2}^{(1)}}\\
&+\underbrace{\frac{1}{N}\sum_{i,j_1,j_2:\,i\neq j_1,j_2,\,j_1\neq j_2,\atop \,i,j_1\neq\ell_1,\,i,j_2\neq\ell_2}\mathbb{E}\{Y_{(i,j_1),\ell_1}Y_{(i,j_2),\ell_2}\}}_{B_{\ell_1,\ell_2}^{(2)}}\\
&+\underbrace{\frac{1}{N}\sum_{i_1,i_2,j:\,j\neq i_1,i_2,\,i_1\neq i_2,\atop \,j,i_1\neq\ell_1,\,j,i_2\neq\ell_2}\mathbb{E}\{Y_{(i_1,j),\ell_1}Y_{(i_2,j),\ell_2}\}}_{B_{\ell_1,\ell_2}^{(3)}}\\
&+\underbrace{\frac{1}{N}\sum_{i_1,j_1,i_2,j_2:\,i_1\neq j_1,\,i_2\neq j_2,\atop i_1\neq i_2,\, j_1\neq j_2,\,i_1,j_1\neq\ell_1,\,i_2,j_2\neq\ell_2}\mathbb{E}\{Y_{(i_1,j_1),\ell_1}Y_{(i_2,j_2),\ell_2}\}}_{B_{\ell_1,\ell_2}^{(4)}}\,.
	\end{align*}
By \eqref{eq:Yijell}, it holds that
	\begin{align}\label{eq:explictform}		
&\nu_{\ell_1}^{1/2}\nu_{\ell_2}^{1/2}\mathbb{E}\{Y_{(i_1,j_1),\ell_1}Y_{(i_2,j_2),\ell_2}\}\notag\\	&~~~~~~~~~=\bigg(\frac{\mu_{\ell_1,1}+\mu_{\ell_1,2}}{2\mu_{\ell_1,1}\mu_{\ell_1,2}}\bigg)\bigg(\frac{\mu_{\ell_2,1}+\mu_{\ell_2,2}}{2\mu_{\ell_2,1}\mu_{\ell_2,2}}\bigg)\mathbb{E}(\mathring{Z}_{i_1,\ell_1}\mathring{Z}_{\ell_1,j_1}\mathring{Z}_{i_1,j_1}\mathring{Z}_{i_2,\ell_2}\mathring{Z}_{\ell_2,j_2}\mathring{Z}_{i_2,j_2})\notag\\	&~~~~~~~~~~~~~-\bigg(\frac{\mu_{\ell_2,1}+\mu_{\ell_2,2}}{2\mu_{\ell_2,1}\mu_{\ell_2,2}}\bigg)\lambda_{i_1,\ell_1}\mathbb{E}(\mathring{Z}_{i_1,\ell_1}\mathring{Z}_{i_2,\ell_2}\mathring{Z}_{\ell_2,j_2}\mathring{Z}_{i_2,j_2})\\	&~~~~~~~~~~~~~-\bigg(\frac{\mu_{\ell_1,1}+\mu_{\ell_1,2}}{2\mu_{\ell_1,1}\mu_{\ell_1,2}}\bigg)\lambda_{i_2,\ell_2}\mathbb{E}(\mathring{Z}_{i_2,\ell_2}\mathring{Z}_{i_1,\ell_1}\mathring{Z}_{\ell_1,j_1}\mathring{Z}_{i_1,j_1})\notag\\
&~~~~~~~~~~~~~+\lambda_{i_1,\ell_1}\lambda_{i_2,\ell_2}\mathbb{E}(\mathring{Z}_{i_1,\ell_1}\mathring{Z}_{i_2,\ell_2})\notag
	\end{align}
	for any $(i_1,j_1,\ell_1,i_2,j_2,\ell_2)$ such that $i_1\neq j_1$, $i_2\neq j_2$, $i_1,j_1\neq \ell_1$ and $i_2,j_2\neq \ell_2$. As we will show in Sections \ref{sec:le4case1}--\ref{sec:pfle4case4},
\begin{align}
&~~~~~~~~~~~~~~~~~~~~~~~\nu_{\ell_1}^{1/2}\nu_{\ell_2}^{1/2}B_{\ell_1,\ell_2}^{(1)}=\bigg(\frac{\tilde{b}_{\ell_1}}{2}+b_{\ell_1}\bigg)I(\ell_1=\ell_2)\,,\label{eq:Bell1ell21}\\
&~~~~~\nu_{\ell_1}^{1/2}\nu_{\ell_2}^{1/2}B_{\ell_1,\ell_2}^{(2)}=O(\gamma^{-6}p^{-1})I(\ell_1\neq \ell_2)+(p-3)b_{\ell_1} I(\ell_1=\ell_2)\,,\label{eq:Bell1ell22}\\
&~~~~~~~~~~~~~~~~~~~~~~~~\nu_{\ell_1}^{1/2}\nu_{\ell_2}^{1/2}B_{\ell_1,\ell_2}^{(3)}=O(\gamma^{-6}p^{-1})I(\ell_1\neq\ell_2)\,,\label{eq:Bell1ell23}\\
&\nu_{\ell_1}^{1/2}\nu_{\ell_2}^{1/2}B_{\ell_1,\ell_2}^{(4)}=\{O(\gamma^{-6}p^{-1})+O(\gamma^{-2})\}I(\ell_1\neq \ell_2)+\frac{\tilde{b}_{\ell_1}}{2}I(\ell_1=\ell_2)\,.\label{eq:Bell1ell24}
\end{align}
By \eqref{eq:Bell1ell21}--\eqref{eq:Bell1ell24}, we have
$
\nu_{\ell_1}^{1/2}\nu_{\ell_2}^{1/2}B_{\ell_1,\ell_2}=\{\tilde{b}_{\ell_1}+(p-2)b_{\ell_1}\}I(\ell_1=\ell_2)+\{O(\gamma^{-6}p^{-1})+O(\gamma^{-2})\}I(\ell_1\neq \ell_2)$.
Recall $\nu_\ell=(p-2)b_\ell+\tilde{b}_\ell\asymp p\gamma^{-2}+\gamma^{-6}$. Hence, $B_{\ell,\ell}=1$ for any $\ell\in[p]$, and
$
\max_{1\leq\ell_1\neq\ell_2\leq p}|B_{\ell_1,\ell_2}|\lesssim p^{-1}$.
 We complete the proof of Lemma \ref{la:vnorm}. $\hfill\Box$

	\subsubsection{Proof of \eqref{eq:Bell1ell21}.}\label{sec:le4case1} 
For any $(i,j,\ell_1,\ell_2)$ such that $i\neq j$, $\ell_1\neq\ell_2$ and $i,j\neq\ell_1,\ell_2$, by \eqref{eq:explictform}, it holds that  $\nu_{\ell_1}^{1/2}\nu_{\ell_2}^{1/2}\mathbb{E}\{Y_{(i,j),\ell_1}Y_{(i,j),\ell_2}\}=0$, which implies
$
		\nu_{\ell_1}^{1/2}\nu_{\ell_2}^{1/2}B_{\ell_1,\ell_2}^{(1)}=0$
	for any $\ell_1\neq\ell_2$. For any $(i,j,\ell)$ such that $i\neq j$ and $i,j\neq\ell$, by \eqref{eq:explictform}, it holds that
\begin{align*}
\nu_{\ell}\mathbb{E}\{Y_{(i,j),\ell}Y_{(i,j),\ell}\}
=\bigg(\frac{\mu_{\ell,1}+\mu_{\ell,2}}{2\mu_{\ell,1}\mu_{\ell,2}}\bigg)^2{\rm Var}(Z_{i,\ell}){\rm Var}(Z_{\ell,j}){\rm Var}(Z_{i,j})+\lambda_{i,\ell}^2{\rm Var}({Z}_{i,\ell})\,,
\end{align*}
which implies
	\begin{align*}
\nu_\ell B_{\ell,\ell}^{(1)}=&~\frac{1}{N}\bigg(\frac{\mu_{\ell,1}+\mu_{\ell,2}}{2\mu_{\ell,1}\mu_{\ell,2}}\bigg)^2\sum_{i,j:\,i\neq j,\,i,j\neq\ell}{\rm Var}(Z_{i,\ell}){\rm Var}(Z_{\ell,j}){\rm Var}(Z_{i,j})\\
&+\frac{1}{N}\sum_{i,j:\,i\neq j,\,i,j\neq\ell}\lambda_{i,\ell}^2{\rm Var}({Z}_{i,\ell})\\
=&~\frac{\tilde{b}_\ell}{2}+ b_\ell\,.
\end{align*}
We complete the proof of \eqref{eq:Bell1ell21}. $\hfill\Box$
	
	\subsubsection{Proof of \eqref{eq:Bell1ell22}.}\label{sec:pfBell1ell22}
	 For any $(i,j_1,j_2,\ell_1,\ell_2)$ such that $j_1,j_2\neq i$, $j_1\neq j_2$, $\ell_1\neq\ell_2$, $i,j_1\neq\ell_1$ and $i,j_2\neq\ell_2$, by \eqref{eq:explictform}, it holds that
\begin{align}
&\nu_{\ell_1}^{1/2}\nu_{\ell_2}^{1/2}\mathbb{E}\{Y_{(i,j_1),\ell_1}Y_{(i,j_2),\ell_2}\}=\bigg(\frac{\mu_{\ell_1,1}+\mu_{\ell_1,2}}{2\mu_{\ell_1,1}\mu_{\ell_1,2}}\bigg)\bigg(\frac{\mu_{\ell_2,1}+\mu_{\ell_2,2}}{2\mu_{\ell_2,1}\mu_{\ell_2,2}}\bigg){\rm Var}(Z_{i,\ell_1})\notag\\
&~~~~~~~~~~~~~~~~~~~~~~~~~~~~~~~~~~~~~~~~~~~~~~~~~~~~~~~~~~~~~~\times{\rm Var}(Z_{\ell_1,\ell_2}){\rm Var}(Z_{i,\ell_2})I(j_2=\ell_1,j_1=\ell_2)\,,\label{eq:rebd1}
\end{align}
	  which implies
		\begin{align*}
			\nu_{\ell_1}^{1/2}\nu_{\ell_2}^{1/2}B_{\ell_1,\ell_2}^{(2)} =&\,\bigg(\frac{\mu_{\ell_1,1}+\mu_{\ell_1,2}}{2\mu_{\ell_1,1}\mu_{\ell_1,2}}\bigg)\bigg(\frac{\mu_{\ell_2,1}+\mu_{\ell_2,2}}{2\mu_{\ell_2,1}\mu_{\ell_2,2}}\bigg)\\
&~~~~~~~~~~~~~\times\frac{1}{N}\sum_{i:\,i\neq \ell_1,\ell_2}{\rm Var}(Z_{i,\ell_1}){\rm Var}(Z_{\ell_1,\ell_2}){\rm Var}(Z_{i,\ell_2})\\
			\asymp&\,\gamma^{-6}p^{-1}
		\end{align*}
	for any $\ell_1\neq\ell_2$. For any $(i,j_1,j_2,\ell)$ such that $j_1,j_2\neq i$, $j_1\neq j_2$ and $i,j_1,j_2\neq \ell$, by \eqref{eq:explictform}, it holds that
$
\nu_{\ell}\mathbb{E}\{Y_{(i,j_1),\ell}Y_{(i,j_2),\ell}\}=\lambda_{i,\ell}^2{\rm Var}(Z_{i,\ell})$, which implies
		\begin{align*}
			&\nu_\ell B_{\ell,\ell}^{(2)}=\frac{1}{N}\sum_{i,j_1,j_2:\,i\neq j_1,j_2,\,j_1\neq j_2,\atop \,i,j_1,j_2\neq\ell}\lambda_{i,\ell}^2{\rm Var}(Z_{i,\ell})=\frac{p-3}{p-1}\sum_{i:\,i\neq\ell}\lambda_{i,\ell}^2{\rm Var}(Z_{i,\ell})=(p-3)b_\ell\,.
		\end{align*}
We complete the proof of \eqref{eq:Bell1ell22}. $\hfill\Box$	
	
	\subsubsection{Proof of \eqref{eq:Bell1ell23}.}\label{sec:pfBell1ell23} For any $(i_1,i_2,j,\ell_1,\ell_2)$ such that $i_1,i_2\neq j$, $i_1\neq i_2$, $\ell_1\neq\ell_2$, $j,i_1\neq\ell_1$ and $j,i_2\neq\ell_2$, by \eqref{eq:explictform}, it holds that
\begin{align}
\nu_{\ell_1}^{1/2}\nu_{\ell_2}^{1/2}\mathbb{E}\{Y_{(i_1,j),\ell_1}Y_{(i_2,j),\ell_2}\}=&\,\lambda_{\ell_1,\ell_2}\lambda_{\ell_2,\ell_1}{\rm Var}(Z_{\ell_1,\ell_2})I(i_1=\ell_2,i_2=\ell_1)\notag\\
&
+\bigg(\frac{\mu_{\ell_1,1}+\mu_{\ell_1,2}}{2\mu_{\ell_1,1}\mu_{\ell_1,2}}\bigg)\bigg(\frac{\mu_{\ell_2,1}+\mu_{\ell_2,2}}{2\mu_{\ell_2,1}\mu_{\ell_2,2}}\bigg){\rm Var}(Z_{j,\ell_1})\label{eq:rebd2}\\
&~~~~~~~~~~~~~~\times{\rm Var}(Z_{\ell_1,\ell_2}){\rm Var}(Z_{j,\ell_2})I(i_1=\ell_2,i_2=\ell_1)\,,\notag
\end{align}
which implies
	\begin{align*}	\nu_{\ell_1}^{1/2}\nu_{\ell_2}^{1/2}B_{\ell_1,\ell_2}^{(3)}=&~\frac{\lambda_{\ell_1,\ell_2}\lambda_{\ell_2,\ell_1}}{N}\sum_{j:\,j\neq \ell_1,\ell_2}{\rm Var}(Z_{\ell_1,\ell_2})\\
&+\bigg(\frac{\mu_{\ell_1,1}+\mu_{\ell_1,2}}{2\mu_{\ell_1,1}\mu_{\ell_1,2}}\bigg)\bigg(\frac{\mu_{\ell_2,1}+\mu_{\ell_2,2}}{2\mu_{\ell_2,1}\mu_{\ell_2,2}}\bigg)\\
&~~~~~~~~~~~~~~~~~~~~~~\times\frac{1}{N}\sum_{j:\,j\neq \ell_1,\ell_2}{\rm Var}(Z_{j,\ell_1}){\rm Var}(Z_{\ell_1,\ell_2}){\rm Var}(Z_{j,\ell_2})\\
\asymp&\,\gamma^{-2}p^{-1}+\gamma^{-6}p^{-1}\asymp \gamma^{-6}p^{-1}\notag
	\end{align*}
	for any $\ell_1\neq\ell_2$. For any $(i_1,i_2,j,\ell)$ such that $i_1,i_2\neq j$, $i_1\neq i_2$ and $j,i_1,i_2\neq \ell$, by \eqref{eq:explictform}, it holds that $
\nu_{\ell}\mathbb{E}\{Y_{(i_1,j),\ell}Y_{(i_2,j),\ell}\}=0$, which implies $\nu_{\ell}B_{\ell,\ell}^{(3)}=0$. We complete the proof of \eqref{eq:Bell1ell23}. $\hfill\Box$

	\subsubsection{Proof of \eqref{eq:Bell1ell24}.}\label{sec:pfle4case4} For any $(i_1,j_1,i_2,j_2,\ell_1,\ell_2)$ such that $i_1\neq j_1\neq \ell_1$, $i_2\neq j_2\neq \ell_2$, $i_1\neq i_2$, $j_1\neq j_2$ and $\ell_1\neq\ell_2$, by \eqref{eq:explictform}, it holds that
\begin{align}
\nu_{\ell_1}^{1/2}\nu_{\ell_2}^{1/2}\mathbb{E}\{Y_{(i_1,j_1),\ell_1}Y_{(i_2,j_2),\ell_2}\}=&\,\bigg(\frac{\mu_{\ell_1,1}+\mu_{\ell_1,2}}{2\mu_{\ell_1,1}\mu_{\ell_1,2}}\bigg)\bigg(\frac{\mu_{\ell_2,1}+\mu_{\ell_2,2}}{2\mu_{\ell_2,1}\mu_{\ell_2,2}}\bigg){\rm Var}(Z_{i_1,\ell_1})\notag\\
&~~~~\times{\rm Var}(Z_{\ell_1,\ell_2})
{\rm Var}(Z_{i_1,\ell_2})I(j_2=i_1,i_2=\ell_1,j_1=\ell_2)\notag\\
 &+\bigg(\frac{\mu_{\ell_1,1}+\mu_{\ell_1,2}}{2\mu_{\ell_1,1}\mu_{\ell_1,2}}\bigg)\bigg(\frac{\mu_{\ell_2,1}+\mu_{\ell_2,2}}{2\mu_{\ell_2,1}\mu_{\ell_2,2}}\bigg){\rm Var}(Z_{\ell_1,\ell_2})\label{eq:rebd3}\\
 &~~~~\times{\rm Var}(Z_{\ell_1,i_2}){\rm Var}(Z_{\ell_2,i_2})
 I(i_1=\ell_2,j_2=\ell_1,i_2=j_1)\notag\\
 &+\lambda_{\ell_1,\ell_2}\lambda_{\ell_2,\ell_1}{\rm Var}(Z_{\ell_1,\ell_2})I(i_1=\ell_2,i_2=\ell_1)\,,\notag
\end{align}
which implies
\begin{align*}
\nu_{\ell_1}^{1/2}\nu_{\ell_2}^{1/2}B_{\ell_1,\ell_2}^{(4)}=&\,\bigg(\frac{\mu_{\ell_1,1}+\mu_{\ell_1,2}}{2\mu_{\ell_1,1}\mu_{\ell_1,2}}\bigg)\bigg(\frac{\mu_{\ell_2,1}+\mu_{\ell_2,2}}{2\mu_{\ell_2,1}\mu_{\ell_2,2}}\bigg)\\
&~~~~~~~~~~~~~~~~~~~~\times\frac{2}{N}\sum_{i:\,i\neq \ell_1,\ell_2}{\rm Var}(Z_{i,\ell_1}){\rm Var}(Z_{\ell_1,\ell_2}){\rm Var}(Z_{i,\ell_2})\\
&+\frac{\lambda_{\ell_1,\ell_2}\lambda_{\ell_2,\ell_1}}{N}\sum_{j_1,j_2:\,j_1\neq j_2,\,j_1,j_2\neq\ell_1,\ell_2}{\rm Var}(Z_{\ell_1,\ell_2})\notag\\
\asymp&\,\gamma^{-6}p^{-1}+\gamma^{-2}\notag
\end{align*}
for any $\ell_1\neq\ell_2$. For any $(i_1,j_1,i_2,j_2,\ell)$ such that $i_1\neq j_1$, $i_2\neq j_2$, $i_1\neq i_2$, $j_1\neq j_2$ and $i_1,i_2,j_1,j_2\neq\ell$, by \eqref{eq:explictform}, it holds that
\begin{align*}
\nu_{\ell}\mathbb{E}\{Y_{(i_1,j_1),\ell}Y_{(i_2,j_2),\ell}\}=&\,\bigg(\frac{\mu_{\ell,1}+\mu_{\ell,2}}{2\mu_{\ell,1}\mu_{\ell,2}}\bigg)^2{\rm Var}(Z_{i_1,\ell})\\
&~~~~~~~~~~~\times{\rm Var}(Z_{\ell,i_2}){\rm Var}(Z_{i_1,i_2})I(i_1=j_2,i_2=j_1)\,,
\end{align*}
which implies
\begin{align*}
\nu_\ell B_{\ell,\ell}^{(4)}=\bigg(\frac{\mu_{\ell,1}+\mu_{\ell,2}}{2\mu_{\ell,1}\mu_{\ell,2}}\bigg)^2\frac{1}{N}\sum_{i_1,i_2:\,i_1\neq i_2,\,i_1,i_2\neq\ell}{\rm Var}(Z_{i_1,\ell}){\rm Var}(Z_{\ell,i_2}){\rm Var}(Z_{i_1,i_2})=\frac{\tilde{b}_\ell}{2}\,.
\end{align*}
We complete the proof of \eqref{eq:Bell1ell24}. $\hfill\Box$
	
\subsection{Proof of Lemma {\rm\ref{lemma5}}.}\label{sec:pflemma5}

For $Y_{(i,j),\ell}$ defined as \eqref{eq:Yijell}, let $\{V_{(i,j),\ell}\}_{i,j,\ell:\,i\neq j\neq\ell}$ be mean zero normal distributed random variables  that independent of the sequence $\{Y_{(i,j),\ell}\}_{i,j,\ell:\,i\neq j\neq\ell}$ and
	$
	{\rm Cov}\{V_{(i_1,j_1),\ell_1},V_{(i_2,j_2),\ell_2}\}={\rm Cov}\{Y_{(i_1,j_1),\ell_1},Y_{(i_2,j_2),\ell_2}\}$
	for any $(i_1,j_1,\ell_1)$ and $(i_2,j_2,\ell_2)$ such that $i_1\neq j_1\neq\ell_1$ and $i_2\neq j_2\neq\ell_2$. We set $V_{(i,j),\ell}=0$ if $\ell=i$ or $j$. Let $\{W_{(i,j),\ell}\}_{i,j,\ell:\,i\neq j\neq\ell}$ be an independent copy of $\{V_{(i,j),\ell}\}_{i,j,\ell:\,i\neq j\neq\ell}$. We also set $W_{(i,j),\ell}=0$ if $\ell=i$ or $j$.
 For each pair $(i,j)$ such that $i\neq j$, define three $p$-dimensional vectors $\by_{(i,j)}=\{Y_{(i,j),1},\ldots,Y_{(i,j),p}\}^\T$, $\bv_{(i,j)}=\{V_{(i,j),1},\ldots,V_{(i,j),p}\}^\T$ and $\bw_{(i,j)}=\{W_{(i,j),1},\ldots,W_{(i,j),p}\}^\T$. For $\bh$ specified in \eqref{eq:asymexpleading}, we know $\bh=N^{-1/2}\sum_{i,j:\,i\neq j}\by_{(i,j)}$. Recall ${\rm Cov}(\bh)=\bB$. Then	
 \begin{align*}
		\bg:=\frac{1}{\sqrt{N}}\sum_{i,j:\,i\neq j}\bv_{(i,j)}\sim\mathcal{N}(\bzero,\bB)~~\textrm{and}~~\tilde{\bg}:=\frac{1}{\sqrt{N}}\sum_{i,j:\,i\neq j}\bw_{(i,j)}\sim\mathcal{N}(\bzero,\bB)\,.
	\end{align*}
	
	Define $\bc(t)=\sum_{i,j:\,i\neq j}\bc_{(i,j)}(t)$ for any $t\in[0,1]$, where
	\begin{align}\label{eq:cijt1}
	    	\bc_{(i,j)}(t):=\frac{1}{\sqrt{N}}\big[\sqrt{t}\{\sqrt{v}\,\by_{(i,j)}+\sqrt{1-v}\,\bv_{(i,j)}\}+\sqrt{1-t}\,\bw_{(i,j)}\big]\,.
	\end{align}
	Write $\bc_{(i,j)}(t)=\{c_{(i,j),1}(t),\ldots,c_{(i,j),p}(t)\}^\T$. Then $\bc(1)=\sqrt{v}\bh+\sqrt{1-v}\bg$ and $\bc(0)=\tilde{\bg}$. Define
\begin{align}\label{eq:dotcijt1}
\dot{\bc}_{(i,j)}(t):=\frac{1}{\sqrt{N}}\bigg[\frac{1}{\sqrt{t}}\{\sqrt{v}\,\by_{(i,j)}+\sqrt{1-v}\,\bv_{(i,j)}\}-\frac{1}{\sqrt{1-t}}\,\bw_{(i,j)}\bigg]
\end{align}
	and write $\dot{\bc}_{(i,j)}(t)=\{\dot{c}_{(i,j),1}(t),\ldots,\dot{c}_{(i,j),p}(t)\}^\T$. Then
	\[
	\mathcal{T}=q\{\bc(1)\}-q\{\bc(0)\}=\int_0^1\frac{{\rm d}q\{\bc(t)\}}{{\rm d}t}{\rm d}t=\frac{1}{2}\sum_{i,j:\,i\neq j}\sum_{\ell=1}^{p}\int_0^1q_\ell\{\bc(t)\}\dot{c}_{(i,j),\ell}(t)\,{\rm d}t\,,
	\]
	which implies
	\begin{equation*}
\mathbb{E}(\mathcal{T})=\frac{1}{2}\sum_{i,j:\,i\neq j}\sum_{\ell=1}^{p}\int_0^1\mathbb{E}\big[q_\ell\{\bc(t)\}\dot{c}_{(i,j),\ell}(t)\big]\,{\rm d}t\,.
\end{equation*}
	Recall
	\[
\dot{c}_{(i,j),\ell}(t)=\frac{1}{\sqrt{N}}\bigg[\frac{1}{\sqrt{t}}\{\sqrt{v}{Y}_{(i,j),\ell}+\sqrt{1-v}V_{(i,j),\ell}\}-\frac{1}{\sqrt{1-t}}W_{(i,j),\ell}\bigg]
\]
and $Y_{(i,j),\ell}=V_{(i,j),\ell}=W_{(i,j),\ell}=0$ for $\ell=i$ or $j$. It then holds that
	\begin{equation}\label{eq:ET}
		\mathbb{E}(\mathcal{T})=\frac{1}{2}\sum_{i,j,\ell:\,i\neq j\neq \ell}\int_0^1\mathbb{E}\big[q_\ell\{\bc(t)\}\dot{c}_{(i,j),\ell}(t)\big]\,{\rm d}t\,.
	\end{equation}	
	Notice that $Y_{(i,j),\ell}$ is a function of $\{\mathring{Z}_{i,\ell},\mathring{Z}_{\ell,j},\mathring{Z}_{i,j}\}$. Given $(i,j,\ell)$ such that $i\neq j\neq \ell$, we first consider $\{\mathring{Z}_{i,\ell},\mathring{Z}_{\ell,j},\mathring{Z}_{i,j}\}$ will appear in which $Y_{(i',j'),\ell'}$'s. Recall
	\[
Y_{(i',j'),\ell'}=-\nu_{\ell'}^{-1/2}\bigg\{\bigg(\frac{\mu_{\ell',1}+\mu_{\ell',2}}{2\mu_{\ell',1}\mu_{\ell',2}}\bigg)\mathring{Z}_{i',\ell'}\mathring{Z}_{\ell',j'}\mathring{Z}_{i',j'}-\lambda_{i',\ell'}\mathring{Z}_{i',\ell'}\bigg\}\,.
\]
	Then $\mathring{Z}_{i,\ell}$ will appear in $Y_{(i',j'),\ell'}$ such that either $(i',\ell')=(i,\ell)$, $(\ell',j')=(i,\ell)$ or $(i',j')=(i,\ell)$ holds. Since $\mathring{Z}_{\ell,i}=\mathring{Z}_{i,\ell}$, we know $\mathring{Z}_{i,\ell}$ is also not independent of $Y_{(i',j'),\ell'}$ such that either $(i',\ell')=(\ell,i)$, $(\ell',j')=(\ell,i)$ or $(i',j')=(\ell,i)$ holds. Given $i$ and $\ell$ such that $i\neq\ell$, let
	$
	\mathcal{S}_*(i,\ell)=\{(i',j',\ell'):\{i',j'\}=\{i,\ell\}\}\cup \{(i',j',\ell'):\{j',\ell'\}=\{i,\ell\}\}\cup \{(i',j',\ell'):\{\ell',i'\}=\{i,\ell\}\}
	$. Then we have $\mathring{Z}_{i,\ell}$ is independent of $\{Y_{(i',j'),\ell'}\}_{(i',j',\ell')\notin\mathcal{S}_*(i,\ell)}$. For any $(i,j,\ell)$ such that $i\neq j\neq \ell$, define
	\[
\mathcal{S}(i,j,\ell)=\mathcal{S}_*(i,j)\cup\mathcal{S}_*(j,\ell)\cup\mathcal{S}_*(\ell,i)\,.
\]
	We know $Y_{(i,j),\ell}$ is independent of $\{Y_{(i',j'),\ell'}\}_{(i',j',\ell')\notin\mathcal{S}(i,j,\ell)}$. For any $(i,j,\ell)$ and $(i',j',\ell')$ such that $i\neq j\neq \ell$ and $i'\neq j'\neq\ell'$, let $a_{(i',j'),\ell'}^{(i,j,\ell)}=I\{(i',j',\ell')\in\mathcal{S}(i,j,\ell)\}$. For given $(i,j,\ell)$ such that $i\neq j\neq\ell$, we set $a_{(i',j'),\ell'}^{(i,j,\ell)}=0$ if $\ell'\in\{i',j'\}$. Write
	\[
\ba_{(i',j')}^{(i,j,\ell)}=\{a_{(i',j'),1}^{(i,j,\ell)},\ldots,a_{(i',j'),p}^{(i,j,\ell)}\}^\T\,
\]
and define
\begin{align*}
\bc^{(i,j),\ell}(t)=\{c_1^{(i,j),\ell}(t),\ldots,c_p^{(i,j),\ell}(t)\}^\T:=\sum_{i',j':\,i'\neq j'}\bc_{(i',j')}(t)\circ \ba_{(i',j')}^{(i,j,\ell)}\,,
\end{align*}
	where $\circ$ denotes the Hadamard product. Let
\begin{align*}
\bc^{-(i,j),\ell}(t)=&~\bc(t)-\bc^{(i,j),\ell}(t)\,.
\end{align*}
	We can see that $\bc^{-(i,j),\ell}(t)$ is independent of $\{\mathring{Z}_{i,\ell},\mathring{Z}_{\ell,j},\mathring{Z}_{i,j}\}$.
	 By Taylor expansion,
		\begin{align*}
		\int_0^1\mathbb{E}\big[q_\ell\{\bc(t)\}\dot{c}_{(i,j),\ell}(t)\big]\,{\rm d}t=&\underbrace{\int_0^1\mathbb{E}\big[q_\ell\{\bc^{-(i,j),\ell}(t)\}\dot{c}_{(i,j),\ell}(t)\big]\,{\rm d}t}_{{\rm I}_1(i,j,\ell)}\\
		&+\sum_{k=1}^{p}\underbrace{\int_0^1\mathbb{E}\big[q_{\ell k}\{\bc^{-(i,j),\ell}(t)\}\dot{c}_{(i,j),\ell}(t)c^{(i,j),\ell}_{k}(t)\big]\,{\rm d}t}_{{\rm I}_2(i,j,\ell,k)}\\
		&+\sum_{k,l=1}^{p}\int_0^1\int_0^1(1-\tau)\mathbb{E}\big[q_{\ell kl}\{\bc^{-(i,j),\ell}(t)+\tau \bc^{(i,j),\ell}(t)\}\\
		&~~~~~~~~~~~~~\underbrace{~~~~~~~~~~~~~~~~~~~~~\times\dot{c}_{(i,j),\ell}(t)c^{(i,j),\ell}_k(t)c^{(i,j),\ell}_l(t)\big]\,{\rm d}\tau{\rm d}t}_{{\rm I}_3(i,j,\ell,k,l)}\,.
	\end{align*}
	Together with \eqref{eq:ET}, we have
	\begin{align*}
	2\mathbb{E}(\mathcal{T})=&~\sum_{i,j,\ell:\,i\neq j\neq\ell}{\rm I}_1(i,j,\ell)+\sum_{i,j,\ell:\,i\neq j\neq\ell}\sum_{k=1}^p{\rm I}_2(i,j,\ell,k)\\
 &+\sum_{i,j,\ell:\,i\neq j\neq\ell}\sum_{k:\,k\neq\ell}\sum_{l=1}^p{\rm I}_3(i,j,\ell,k,l)+\sum_{i,j,\ell:\,i\neq j\neq\ell}\sum_{l=1}^p{\rm I}_3(i,j,\ell,\ell,l)\,.
	\end{align*}
	As shown in Sections \ref{se:f1}--\ref{se:f3}, it holds that
	\begin{align}
&~~~~~~~~~~~~~~~~~~~~~\sum_{i,j,\ell:\,i\neq j\neq\ell}{\rm I}_1(i,j,\ell)=0\,,\label{eq:sumI1}\\
&~~\sum_{i,j,\ell:\,i\neq j\neq\ell}\sum_{k=1}^p|{\rm I}_2(i,j,\ell,k)|\lesssim p^{-1/2}\phi^3\log^{7/2}p\,, \label{eq:sumI2}\\
&\sum_{i,j,\ell:\,i\neq j\neq \ell}\sum_{k,l=1}^p|{\rm I}_3(i,j,\ell,k,l)|\lesssim p^{-1/2}\phi^3\log^{7/2}p\,.\label{eq:sumI3}
\end{align}
Then $\sup_{v\in[0,1]}|\mathbb{E}(\mathcal{T})|\lesssim p^{-1/2}\phi^3\log^{7/2}p$. We complete the proof of Lemma \ref{lemma5}. $\hfill\Box$

	\subsection{Proof of \eqref{eq:sumI1}.}\label{se:f1}
	To simplify the notation and without causing much confusion, we write $\bc(t)$, $\bc^{-(i,j),\ell}(t)$, $c_k^{(i,j),\ell}(t)$, $c_{(i,j),\ell}(t)$ and $\dot{c}_{(i,j),\ell}(t)$ as $\bc$, $\bc^{-(i,j),\ell}$, $c_k^{(i,j),\ell}$, $c_{(i,j),\ell}$ and $\dot{c}_{(i,j),\ell}$, respectively. Recall
\begin{align*}
&q_\ell\{\bc^{-(i,j),\ell}\}\dot{c}_{(i,j),\ell}=\frac{q_\ell\{\bc^{-(i,j),\ell}\}}{\sqrt{N}}\bigg[\frac{1}{\sqrt{t}}\{\sqrt{v}{Y}_{(i,j),\ell}+\sqrt{1-v}V_{(i,j),\ell}\}-\frac{1}{\sqrt{1-t}}W_{(i,j),\ell}\bigg]\,.
\end{align*}
	Since $Y_{(i,j),\ell}$ is independent of $\bc^{-(i,j),\ell}$, then
	\begin{align*}
	\mathbb{E}\big[q_\ell\{\bc^{-(i,j),\ell}\}Y_{(i,j),\ell}\big]=0\,.
	\end{align*}
	 Notice that $V_{(i',j'),\ell'}$ with $i'\neq j'\neq\ell'$ included in $\bc^{-(i,j),\ell}$ satisfies $|\{i',j',\ell'\}\cap\{i,j,\ell\}|\leq1$. It follows from \eqref{eq:rebd1}, \eqref{eq:rebd2} and \eqref{eq:rebd3} that ${\rm Cov}\{Y_{(i',j'),\ell'},Y_{(i,j),\ell}\}=0$ for any $i'\neq j'\neq\ell'$ such that $|\{i',j',\ell'\}\cap\{i,j,\ell\}|\leq1$.	Since ${\rm Cov}\{V_{(i',j'),\ell'},V_{(i,j),\ell}\}={\rm Cov}\{Y_{(i',j'),\ell'},Y_{(i,j),\ell}\}$, then
	$
	{\rm Cov}\{V_{(i',j'),\ell'},V_{(i,j),\ell}\}=0$
	 for any $i'\neq j'\neq\ell'$ such that $|\{i',j',\ell'\}\cap\{i,j,\ell\}|\leq1$. Recall that $\{V_{(i,j),\ell}\}_{i,j,\ell:\,i\neq j\neq \ell}$ are normal random variables. Therefore, $V_{(i,j),\ell}$ is independent of $\bc^{-(i,j),\ell}$. Then
	 \begin{align*}
	 \mathbb{E}\big[q_\ell\{\bc^{-(i,j),\ell}\}V_{(i,j),\ell}\big]=0\,.
	 \end{align*}
	  Analogously, we also have  $\mathbb{E}[q_\ell\{\bc^{-(i,j),\ell}\}W_{(i,j),\ell}]=0$. Hence, ${\rm I}_1(i,j,\ell)\equiv0$ for any $i\neq j\neq \ell$. We complete the proof of \eqref{eq:sumI1}. $\hfill\Box$

	\subsection{Proof of \eqref{eq:sumI2}.}\label{se:f2}
	
	To simplify the notation and without causing much confusion, we write $\bc(t)$, $\bc^{-(i,j),\ell}(t)$, $c_k^{(i,j),\ell}(t)$, $c_{(i,j),\ell}(t)$ and $\dot{c}_{(i,j),\ell}(t)$ as $\bc$, $\bc^{-(i,j),\ell}$, $c_k^{(i,j),\ell}$, $c_{(i,j),\ell}$ and $\dot{c}_{(i,j),\ell}$, respectively. Due to $c_k^{(i,j),\ell}=\sum_{i',j':\,i'\neq j'}c_{(i',j'),k}a_{(i',j'),k}^{(i,j,\ell)}$, then
	\begin{equation}\label{eq:I2ijellk}
		{\rm I}_2(i,j,\ell,k)=\sum_{i',j':\,i'\neq j'}\int_0^1\mathbb{E}\big[q_{\ell k}\{\bc^{-(i,j),\ell}\}\dot{c}_{(i,j),\ell}c_{(i',j'),k}a_{(i',j'),k}^{(i,j,\ell)}\big]\,{\rm d}t\,.
	\end{equation}
	Notice that
	\begin{align*}
\sum_{i,j,\ell:\,i\neq j\neq\ell}\sum_{k=1}^p|{\rm I}_2(i,j,\ell,k)|=&\,\sum_{i,j,\ell:\,i\neq j\neq\ell}\sum_{k:\,k\neq i,j,\ell}|{\rm I}_2(i,j,\ell,k)|+\sum_{i,j,\ell:\,i\neq j\neq\ell}|{\rm I}_2(i,j,\ell,i)|\\
&+\sum_{i,j,\ell:\,i\neq j\neq\ell}|{\rm I}_2(i,j,\ell,j)|+\sum_{i,j,\ell:\,i\neq j\neq\ell}|{\rm I}_2(i,j,\ell,\ell)|\,.
\end{align*}
	As we will show in Sections \ref{se:f2.1} and \ref{se:f2.2}, it holds that
	\begin{align}
&~~~~~~~~~~~~~~~~~\sum_{i,j,\ell:\,i\neq j\neq \ell}\sum_{k:\,k\neq i,j,\ell}|{\rm I}_2(i,j,\ell,k)|\lesssim \frac{\phi^3\log^{7/2}p}{p^{1/2}}\,,\label{eq:F.21}\\
&\sum_{i,j,\ell:\,i\neq j\neq\ell}\big\{|{\rm I}_2(i,j,\ell,i)|+|{\rm I}_2(i,j,\ell,j)|+|{\rm I}_2(i,j,\ell,\ell)|\big\}\lesssim \frac{\phi^3\log^{7/2}p}{p^{1/2}}\,.\label{eq:F.22}
\end{align}	
	
	\subsubsection{Case 1: $k\neq i,j,\ell$.}\label{se:f2.1}
	
	Note that $i'\neq j'$ and $k\neq i,j,\ell$. Due to $a_{(i',j'),k}^{(i,j,\ell)}=I\{(i',j',k)\in\mathcal{S}(i,j,\ell)\}$ for $k\notin\{i',j'\}$ and $a_{(i',j'),k}^{(i,j,\ell)}=0$ for $k\in\{i',j'\}$, then $a_{(i',j'),k}^{(i,j,\ell)}=1$ if and only if $\{i',j'\}\subset\{i,j,\ell\}$. It follows from \eqref{eq:I2ijellk} that
	\begin{align}
		&{\rm I}_2(i,j,\ell,k)\notag\\
		&~~~~~~=\int_0^1\mathbb{E}\big[q_{\ell k}\{\bc^{-(i,j),\ell}\}\dot{c}_{(i,j),\ell}c_{(i,\ell),k}\big]\,{\rm d}t+\int_0^1\mathbb{E}\big[q_{\ell k}\{\bc^{-(i,j),\ell}\}\dot{c}_{(i,j),\ell}c_{(\ell,i),k}\big]\,{\rm d}t\notag\\
		&~~~~~~~+\int_0^1\mathbb{E}\big[q_{\ell k}\{\bc^{-(i,j),\ell}\}\dot{c}_{(i,j),\ell}c_{(i,j),k}\big]\,{\rm d}t+\int_0^1\mathbb{E}\big[q_{\ell k}\{\bc^{-(i,j),\ell}\}\dot{c}_{(i,j),\ell}c_{(j,i),k}\big]\,{\rm d}t\label{eq:i2}\\
		&~~~~~~~+\int_0^1\mathbb{E}\big[q_{\ell k}\{\bc^{-(i,j),\ell}\}\dot{c}_{(i,j),\ell}c_{(j,\ell),k}\big]\,{\rm d}t+\int_0^1\mathbb{E}\big[q_{\ell k}\{\bc^{-(i,j),\ell}\}\dot{c}_{(i,j),\ell}c_{(\ell,j),k}\big]\,{\rm d}t\,.\notag
	\end{align}
	Notice that
	\begin{align*}
\dot{c}_{(i,j),\ell}c_{(i,\ell),k}=&~\frac{1}{N}\bigg[\frac{1}{\sqrt{t}}\{\sqrt{v}{Y}_{(i,j),\ell}+\sqrt{1-v}V_{(i,j),\ell}\}-\frac{1}{\sqrt{1-t}}W_{(i,j),\ell}\bigg]\\
&\times \big[\sqrt{t}\{\sqrt{v}Y_{(i,\ell),k}+\sqrt{1-v}V_{(i,\ell),k}\}+\sqrt{1-t}W_{(i,\ell),k}\big]\,.
\end{align*}
As shown in Section \ref{se:f1}, $\{Y_{(i,j),\ell},V_{(i,j),\ell},W_{(i,j),\ell}\}$ is independent of $\bc^{-(i,j),\ell}$. Then
	\begin{equation}\label{eq:term0}
		\begin{split}
			&N\cdot\mathbb{E}\big[q_{\ell k}\{\bc^{-(i,j),\ell}\}\dot{c}_{(i,j),\ell}c_{(i,\ell),k}\big]\\
			&~~~~~~=v\mathbb{E}\big[q_{\ell k}\{\bc^{-(i,j),\ell}\}Y_{(i,j),\ell}Y_{(i,\ell),k}\big]+(1-v)\mathbb{E}\big[q_{\ell k}\{\bc^{-(i,j),\ell}\}V_{(i,j),\ell}V_{(i,\ell),k}\big]\\
			&~~~~~~~~~~-\mathbb{E}\big[q_{\ell k}\{\bc^{-(i,j),\ell}\}W_{(i,j),\ell}W_{(i,\ell),k}\big]\,.
		\end{split}
	\end{equation}
	Recall $Y_{(i,\ell),k}$ is a function of $\{\mathring{Z}_{i,k},\mathring{Z}_{k,\ell},\mathring{Z}_{i,\ell}\}$, and $\bc^{-(i,j),\ell}$ is independent of $\mathring{Z}_{i,\ell}$. We will remove the components in $\bc^{-(i,j),\ell}$ that depend on $\mathring{Z}_{i,k}$ and $\mathring{Z}_{k,\ell}$. As we have shown that $\mathring{Z}_{i,j}$ is independent of $\{Y_{(i',j'),\ell'}\}_{(i',j',\ell')\notin\mathcal{S}_*(i,j)}$ for any given $i$ and $j$, we define
	\begin{align*}
\bc^{(i,j),\ell,(i,\ell),k}=&~\sum_{i',j':\,i'\neq j'}\bc_{(i',j')}\circ \big\{\ba_{(i',j')}^{(i,j,\ell)}+\ba_{(i',j')}^{(i,k)}+\ba_{(i',j')}^{(k,\ell)}-\ba_{(i',j')}^{(i,j,\ell)}\circ\ba_{(i',j')}^{(i,k)}\\
&~~~~~~~~~~~~~~~~~~~-\ba_{(i',j')}^{(i,j,\ell)}\circ\ba_{(i',j')}^{(k,\ell)}-\ba_{(i',j')}^{(i,k)}\circ\ba_{(i',j')}^{(k,\ell)}+\ba_{(i',j')}^{(i,j,\ell)}\circ\ba_{(i',j')}^{(i,k)}\circ\ba_{(i',j')}^{(k,\ell)}\big\}\,,
\end{align*}
	where $\ba_{(i',j')}^{(i,j)}=\{a_{(i',j'),1}^{(i,j)},\ldots,a_{(i',j'),p}^{(i,j)}\}^\T$ with $a_{(i',j'),\ell'}^{(i,j)}=I\{(i',j',\ell')\in\mathcal{S}_*(i,j)\}$ for $\ell'\neq i',j'$ and $a_{(i',j'),\ell'}^{(i,j)}=0$ for $\ell'\in\{i',j'\}$. Therefore, we know $\{Y_{(i,j),\ell},Y_{(i,\ell),k}\}$ is independent of $\bc-\bc^{(i,j),\ell,(i,\ell),k}$. Recall $
	\bc^{(i,j),\ell}=\sum_{i',j':\,i'\neq j'}\bc_{(i',j')}\circ \ba_{(i',j')}^{(i,j,\ell)}$. Then
	\begin{align}
&\bc^{(i,j),\ell,(i,\ell),k}-\bc^{(i,j),\ell}\notag\\
&~~~~~~~=\{\bc_{(i,k)}+\bc_{(k,i)}+\bc_{(\ell,k)}+\bc_{(k,\ell)}\}\circ({\bf 1}-\bfe_i-\bfe_j-\bfe_k-\bfe_\ell)\notag\\
&~~~~~~~~~~~+\sum_{m\neq i,j,\ell,k}\{c_{(k,m),i}+c_{(m,k),i}\}\bfe_i+\sum_{m\neq i,j,\ell,k}\{c_{(k,m),\ell}+c_{(m,k),\ell}\}\bfe_\ell\label{eq:cijlilk-ijl}\\
&~~~~~~~~~~~+\sum_{m\neq i,j,\ell,k}\{c_{(i,m),k}+c_{(m,i),k}+c_{(\ell,m),k}+c_{(m,\ell),k}\}\bfe_k\,,\notag
\end{align}
	where ${\bf 1}$ is a $p$-dimensional vector with all components being $1$, and $\bfe_s$ is a $p$-dimensional vector with the $s$-th component being $1$ and other components being $0$. Let
	\[
\bc^{-(i,j),\ell,(i,\ell),k}=\bc-\bc^{(i,j),\ell,(i,\ell),k}\,.
\]
	Recall $
	\bc^{-(i,j),\ell}=\bc-\bc^{(i,j),\ell}$. Then
	\begin{align*}
	    \bc^{(i,j),\ell,(i,\ell),k-(i,j),\ell}:=\bc^{-(i,j),\ell}-\bc^{-(i,j),\ell,(i,\ell),k}=\bc^{(i,j),\ell,(i,\ell),k}-\bc^{(i,j),\ell}\,.
	\end{align*}
Write $\bc^{(i,j),\ell,(i,\ell),k-(i,j),\ell}=\{c_1^{(i,j),\ell,(i,\ell),k-(i,j),\ell},\ldots,c_p^{(i,j),\ell,(i,\ell),k-(i,j),\ell}\}^\T$. Since $\{Y_{(i,j),\ell},Y_{(i,\ell),k}\}$ is independent of $\bc^{-(i,j),\ell,(i,\ell),k}$, by Taylor expansion, it holds that
\begin{align*}
\mathbb{E}\big[q_{\ell k}\{\bc^{-(i,j),\ell}\}Y_{(i,j),\ell}Y_{(i,\ell),k}\big]=&~\mathbb{E}\big[q_{\ell k}\{\bc^{-(i,j),\ell,(i,\ell),k}\}\big]\mathbb{E}\big\{Y_{(i,j),\ell}Y_{(i,\ell),k}\big\}\\
&+\sum_{m=1}^p\int_0^1\mathbb{E}\big[q_{\ell km}\{\bc^{-(i,j),\ell,(i,\ell),k}+\tau \bc^{(i,j),\ell,(i,\ell),k-(i,j),\ell}\}\\
&~~~~~~~~~~~~~~~~~~~~~~\times Y_{(i,j),\ell}Y_{(i,\ell),k}c_m^{(i,j),\ell,(i,\ell),k-(i,j),\ell}\big]\,{\rm d}\tau\,.
\end{align*}
Notice that for any $V_{(i',j'),\ell'}$ and $W_{(i',j'),\ell'}$ with $i'\neq j'\neq \ell'$ included in $\bc^{-(i,j),\ell,(i,\ell),k}$,
\begin{align*}
{\rm Cov}\{V_{(i,j),\ell},V_{(i',j'),\ell'}\}={\rm Cov}\{W_{(i,j),\ell},W_{(i',j'),\ell'}\}=&\,{\rm Cov}\{Y_{(i,j),\ell},Y_{(i',j'),\ell'}\}=0\,,\\
{\rm Cov}\{V_{(i,\ell),k},V_{(i',j'),\ell'}\}={\rm Cov}\{W_{(i,\ell),k},W_{(i',j'),\ell'}\}=&\,{\rm Cov}\{Y_{(i,\ell),k},Y_{(i',j'),\ell'}\}=0\,.
\end{align*}
Since $V_{(i,j),\ell}, V_{(i,\ell),k}, W_{(i,j),\ell}$ and $W_{(i,\ell),k}$ are normal random variables, then $\{V_{(i,j),\ell}, V_{(i,\ell),k},$ $W_{(i,j),\ell}, W_{(i,\ell),k}\}$ is independent of $\bc^{-(i,j),\ell,(i,\ell),k}$. Analogously, we also have
\begin{align*}
\mathbb{E}\big[q_{\ell k}\{\bc^{-(i,j),\ell}\}V_{(i,j),\ell}V_{(i,\ell),k}\big]=&~\mathbb{E}\big[q_{\ell k}\{\bc^{-(i,j),\ell,(i,\ell),k}\}\big]\mathbb{E}\big\{V_{(i,j),\ell}V_{(i,\ell),k}\big\}\\
&+\sum_{m=1}^p\int_0^1\mathbb{E}\big[q_{\ell km}\{\bc^{-(i,j),\ell,(i,\ell),k}+\tau \bc^{(i,j),\ell,(i,\ell),k-(i,j),\ell}\}\\
&~~~~~~~~~~~~~~~~~~~~~~\times V_{(i,j),\ell}V_{(i,\ell),k}c_m^{(i,j),\ell,(i,\ell),k-(i,j),\ell}\big]\,{\rm d}\tau\,,\\
\mathbb{E}\big[q_{\ell k}\{\bc^{-(i,j),\ell}\}W_{(i,j),\ell}W_{(i,\ell),k}\big]=&~\mathbb{E}\big[q_{\ell k}\{\bc^{-(i,j),\ell,(i,\ell),k}\}\big]\mathbb{E}\big\{W_{(i,j),\ell}W_{(i,\ell),k}\big\}\\
&+\sum_{m=1}^p\int_0^1\mathbb{E}\big[q_{\ell km}\{\bc^{-(i,j),\ell,(i,\ell),k}+\tau \bc^{(i,j),\ell,(i,\ell),k-(i,j),\ell}\}\\
&~~~~~~~~~~~~~~~~~~~~~~\times W_{(i,j),\ell}W_{(i,\ell),k}c_m^{(i,j),\ell,(i,\ell),k-(i,j),\ell}\big]\,{\rm d}\tau\,.
\end{align*}
Recall
$
\mathbb{E}\{Y_{(i,j),\ell}Y_{(i,\ell),k}\}=\mathbb{E}\{V_{(i,j),\ell}V_{(i,\ell),k}\}=\mathbb{E}\{W_{(i,j),\ell}W_{(i,\ell),k}\}$.
 Then \eqref{eq:term0} implies that
	\begin{align}
		&~N\cdot\mathbb{E}\big[q_{\ell k}\{\bc^{-(i,j),\ell}\}\dot{c}_{(i,j),\ell}c_{(i,\ell),k}\big]\notag\\
		=&~\sum_{m=1}^p\int_0^1\mathbb{E}\big[q_{\ell km}\{\bc^{-(i,j),\ell,(i,\ell),k}+\tau \bc^{(i,j),\ell,(i,\ell),k-(i,j),\ell}\}c_m^{(i,j),\ell,(i,\ell),k-(i,j),\ell}\label{eq:term1}\\
		&~~~~~~~~~~~~~~~\times \{vY_{(i,j),\ell}Y_{(i,\ell),k}+(1-v)V_{(i,j),\ell}V_{(i,\ell),k}-W_{(i,j),\ell}W_{(i,\ell),k}\}\big]\,{\rm d}\tau\,.\notag
	\end{align}
	Define
	\begin{equation*}\label{eq:E1}
\mathcal{E}_1=\big\{|Y_{(i,j),\ell}|\vee |V_{(i,j),\ell}|\vee |W_{(i,j),\ell}|\leq B~\textrm{for any}~i\neq j\neq \ell\big\}
\end{equation*}
	for some $B>0$ that will be specified later. Write $\nu=p\gamma^{-2}+\gamma^{-6}$. Then $\nu\asymp\max_{s\in[p]}\nu_s\asymp\min_{s\in[p]}\nu_s$. Recall ${\rm Var}\{V_{(i,j),\ell}\}={\rm Var}\{W_{(i,j),\ell}\}={\rm Var}\{Y_{(i,j),\ell}\}$. As shown in Section \ref{sec:le4case1}, ${\rm Var}\{V_{(i,j),\ell}\}\asymp \nu^{-1}\gamma^{-6}\asymp {\rm Var}\{W_{(i,j),\ell}\}$. Since $V_{(i,j),\ell}$ and $W_{(i,j),\ell}$ are normal random variables with mean zero, then
	\[
\max_{i,j,\ell:\,i\neq j\neq \ell}|V_{(i,j),\ell}|=\nu^{-1/2}\gamma^{-3}\cdot O_{\p}(\log^{1/2}p)=\max_{i,j,\ell:\,i\neq j\neq k}|W_{(i,j),\ell}|\,.
\]
	Notice that $\max_{i,j,\ell:\,i\neq j\neq \ell}|Y_{(i,j),\ell}|\lesssim \nu^{-1/2}\gamma^{-3}$. Recall $\nu=p\gamma^{-2}+\gamma^{-6}$. Then $\nu^{-1/2}\gamma^{-3}\lesssim1$. Selecting $B=C_*\log^{1/2}p$ for a sufficiently large constant $C_*>0$, then $\mathbb{P}(\mathcal{E}_1^c)\lesssim p^{-C}$, where $C>0$ can be sufficiently large if we select a sufficiently large $C_*$.
For given $(i,j,\ell,k)$ such that $i\neq j\neq \ell\neq k$, restricted on $\mathcal{E}_1$, it holds that
\begin{align}
		|\bc^{(i,j),\ell}|_\infty\leq&~\frac{60B}{\sqrt{N}}+\sum_{s_1,s_2:\,s_1\neq s_2\atop s_1,s_2\in\{i,j,\ell\}}\bigg|\frac{1}{\sqrt{N}}\sum_{m\neq i,j,\ell,k}Y_{(s_1,m),s_2}\bigg|\notag\\
		&+\sum_{s_1,s_2:\,s_1\neq s_2\atop s_1,s_2\in\{i,j,\ell\}}\bigg|\frac{1}{\sqrt{N}}\sum_{m\neq i,j,\ell,k}Y_{(m,s_1),s_2}\bigg|\notag\\
		&+\sum_{s_1,s_2:\,s_1\neq s_2\atop s_1,s_2\in\{i,j,\ell\}}\bigg|\frac{1}{\sqrt{N}}\sum_{m\neq i,j,\ell,k}V_{(s_1,m),s_2}\bigg|\label{eq:keybound}\\
		&+\sum_{s_1,s_2:\,s_1\neq s_2\atop s_1,s_2\in\{i,j,\ell\}}\bigg|\frac{1}{\sqrt{N}}\sum_{m\neq i,j,\ell,k}V_{(m,s_1),s_2}\bigg|\notag\\
		&+\sum_{s_1,s_2:\,s_1\neq s_2\atop s_1,s_2\in\{i,j,\ell\}}\bigg|\frac{1}{\sqrt{N}}\sum_{m\neq i,j,\ell,k}W_{(s_1,m),s_2}\bigg|\notag\\
		&+\sum_{s_1,s_2:\,s_1\neq s_2\atop s_1,s_2\in\{i,j,\ell\}}\bigg|\frac{1}{\sqrt{N}}\sum_{m\neq i,j,\ell,k}W_{(m,s_1),s_2}\bigg|\,,\notag\\
		|\bc^{(i,j),\ell,(i,\ell),k-(i,j),\ell}|_\infty\leq&~\frac{24B}{\sqrt{N}}+\sum_{s_1,s_2:\,s_1\neq s_2\atop s_1,s_2\in\{i,\ell,k\}}\bigg|\frac{1}{\sqrt{N}}\sum_{m:\,m\neq i,j,\ell,k}Y_{(s_1,m),s_2}\bigg|\notag\\
		&+\sum_{s_1,s_2:\,s_1\neq s_2\atop s_1,s_2\in\{i,\ell,k\}}\bigg|\frac{1}{\sqrt{N}}\sum_{m:\,m\neq i,j,\ell,k}Y_{(m,s_1),s_2}\bigg|\notag\\
		&+\sum_{s_1,s_2:\,s_1\neq s_2\atop s_1,s_2\in\{i,\ell,k\}}\bigg|\frac{1}{\sqrt{N}}\sum_{m:\,m\neq i,j,\ell,k}V_{(s_1,m),s_2}\bigg|\notag\\
		&+\sum_{s_1,s_2:\,s_1\neq s_2\atop s_1,s_2\in\{i,\ell,k\}}\bigg|\frac{1}{\sqrt{N}}\sum_{m:\,m\neq i,j,\ell,k}V_{(m,s_1),s_2}\bigg|\label{eq:term2}\\
		&+\sum_{s_1,s_2:\,s_1\neq s_2\atop s_1,s_2\in\{i,\ell,k\}}\bigg|\frac{1}{\sqrt{N}}\sum_{m:\,m\neq i,j,\ell,k}W_{(s_1,m),s_2}\bigg|\notag\\
		&+\sum_{s_1,s_2:\,s_1\neq s_2\atop s_1,s_2\in\{i,\ell,k\}}\bigg|\frac{1}{\sqrt{N}}\sum_{m:\,m\neq i,j,\ell,k}W_{(m,s_1),s_2}\bigg|\,,\notag
	\end{align}
	which implies that
	\begin{align}
		&|\bc^{(i,j),\ell}-(\tau-1)\bc^{(i,j),\ell,(i,\ell),k-(i,j),\ell}|_\infty\notag\\
		&~~~~~~~~~\leq24\max_{s_1,s_2:\,s_1\neq s_2\atop s_1,s_2\in\{i,j,\ell,k\}}\bigg\{\bigg|\frac{1}{\sqrt{N}}\sum_{m\neq i,j,\ell,k}Y_{(s_1,m),s_2}\bigg|+\bigg|\frac{1}{\sqrt{N}}\sum_{m\neq i,j,\ell,k}Y_{(m,s_1),s_2}\bigg|\bigg\}\notag\\
		&~~~~~~~~~~~~~+24\max_{s_1,s_2:\,s_1\neq s_2\atop s_1,s_2\in\{i,j,\ell,k\}}\bigg\{\bigg|\frac{1}{\sqrt{N}}\sum_{m\neq i,j,\ell,k}V_{(s_1,m),s_2}\bigg|+\bigg|\frac{1}{\sqrt{N}}\sum_{m\neq i,j,\ell,k}V_{(m,s_1),s_2}\bigg|\bigg\}\label{eq:cha}\\
		&~~~~~~~~~~~~~+24\max_{s_1,s_2:\,s_1\neq s_2\atop s_1,s_2\in\{i,j,\ell,k\}}\bigg\{\bigg|\frac{1}{\sqrt{N}}\sum_{m\neq i,j,\ell,k}W_{(s_1,m),s_2}\bigg|+\bigg|\frac{1}{\sqrt{N}}\sum_{m\neq i,j,\ell,k}W_{(m,s_1),s_2}\bigg|\bigg\}\notag\\
		&~~~~~~~~~~~~~+\frac{84B}{\sqrt{N}}\notag
			\end{align}
	under $\mathcal{E}_1$. Given $s_1,s_2\in\{i,j,\ell,k\}$ such that $s_1\neq s_2$, we have
	\begin{align*}
\frac{1}{\sqrt{N}}\sum_{m\neq i,j,\ell,k}Y_{(s_1,m),s_2}=&-\nu_{s_2}^{-1/2}\bigg(\frac{\mu_{s_2,1}+\mu_{s_2,2}}{2\mu_{s_2,1}\mu_{s_2,2}}\bigg)\frac{1}{\sqrt{N}}\sum_{m\neq i,j,\ell,k}\mathring{Z}_{s_1,s_2}\mathring{Z}_{s_2,m}\mathring{Z}_{s_1,m}\\
&+\frac{p-4}{\sqrt{N}}\nu_{s_2}^{-1/2}\lambda_{s_1,s_2}\mathring{Z}_{s_1,s_2}\,.
\end{align*}
Recall that $\max_{s\in[p],k\in\{1,2\}}\mu_{s,k}\asymp\gamma^3\asymp\min_{s\in[p],k\in\{1,2\}}\mu_{s,k}$, $\max_{s_1,s_2:\,s_1\neq s_2}\lambda_{s_1,s_2}\asymp\gamma^{-1}\asymp\min_{s_1,s_2:\,s_1\neq s_2}\lambda_{s_1,s_2}$, $\max_{s\in[p]}\nu_s\asymp\nu\asymp\min_{s\in[p]}\nu_s$ with $\nu=p\gamma^{-2}+\gamma^{-6}$, and $N=(p-1)(p-2)$. 	Then
	\[
\bigg|\frac{1}{\sqrt{N}}\sum_{m\neq i,j,\ell,k}Y_{(s_1,m),s_2}\bigg|\lesssim \frac{1}{\nu^{1/2}\gamma^{3}} \bigg|\frac{1}{\sqrt{N}}\sum_{m\neq i,j,\ell,k}\mathring{Z}_{s_2,m}\mathring{Z}_{s_1,m}\bigg|+\frac{1}{\nu^{1/2}\gamma}\,.
\]
	Notice that $\nu^{-1/2}\gamma^{-1}\lesssim p^{-1/2}$, $\nu^{-1/2}\gamma^{-3}\lesssim 1$ and $\{\mathring{Z}_{s_1,m}\mathring{Z}_{s_2,m}\}_{m\neq i,j,k,\ell}$ is an independent sequence with mean zero. By Bernstein inequality,
	\[
\mathbb{P}\bigg(\bigg|\frac{1}{\sqrt{N}}\sum_{m\neq i,j,\ell,k}\mathring{Z}_{s_2,m}\mathring{Z}_{s_1,m}\bigg|>u\bigg)\lesssim\exp(-{Cpu^2})
\]
	for any $u=o(1)$, which implies
	\begin{align}
		&\max_{i,j,\ell,k:\,i\neq j\neq \ell\neq k}\max_{s_1,s_2:\,s_1\neq s_2\atop s_1,s_2\in\{i,j,\ell,k\}}\bigg|\frac{1}{\sqrt{N}}\sum_{m\neq i,j,\ell,k}Y_{(s_1,m),s_2}\bigg|\notag\\
		&~~~~~~~~~~~~~~\leq \frac{C}{\nu^{1/2}\gamma}+\frac{1}{\nu^{1/2}\gamma^{3}}\cdot O_{\p}\bigg(\frac{\log^{1/2}p}{p^{1/2}}\bigg)=O_{\p}\bigg(\frac{\log^{1/2}p}{p^{1/2}}\bigg)\,.\label{eq:argu1}
	\end{align}
	Analogously, we also have
	\begin{align*}
\max_{i,j,\ell,k:\,i\neq j\neq \ell\neq k}\max_{s_1,s_2:\,s_1\neq s_2\atop s_1,s_2\in\{i,j,\ell,k\}}\bigg|\frac{1}{\sqrt{N}}\sum_{m\neq i,j,\ell,k}Y_{(m,s_1),s_2}\bigg|=O_{\p}\bigg(\frac{\log^{1/2}p}{p^{1/2}}\bigg)\,.
\end{align*}
Note that $
	{\rm Cov}\{V_{(i_1,j_1),\ell_1},V_{(i_2,j_2),\ell_2}\}={\rm Cov}\{Y_{(i_1,j_1),\ell_1},Y_{(i_2,j_2),\ell_2}\}$.	As shown in Sections \ref{sec:le4case1} and \ref{sec:pfBell1ell22}, we have
	\begin{align*}
\textrm{Cov}\{V_{(s_1,m),s_2},V_{(s_1,m),s_2}\}=&~\frac{1}{\nu_{s_2}}\bigg(\frac{\mu_{s_2,1}+\mu_{s_2,2}}{2\mu_{s_2,1}\mu_{s_2,2}}\bigg)^2{\rm Var}(Z_{s_1,s_2}){\rm Var}(Z_{s_2,m}){\rm Var}(Z_{s_1,m})\\
&+\frac{\lambda_{s_1,s_2}^2}{\nu_{s_2}}{\rm Var}(Z_{s_1,s_2})
\end{align*}
	 and $\textrm{Cov}\{V_{(s_1,m),s_2},V_{(s_1,m'),s_2}\}=\nu_{s_2}^{-1}\lambda_{s_1,s_2}^2{\rm Var}(Z_{s_1,s_2})$ for any $m\neq m'$. Recall $\nu=p\gamma^{-2}+\gamma^{-6}$. Then
	 \begin{align*}
	     {\rm Var}\bigg\{\sum_{m\neq i,j,\ell,k}V_{(s_1,m),s_2}\bigg\}\asymp \frac{p}{\nu\gamma^{2}}\bigg(\frac{1}{\gamma^{4}}+p\bigg)\asymp p\,.
	 \end{align*}
	  Since $\{V_{(s_1,m),s_2}\}_{m\neq i,j,\ell,k}$ are normal random variables with mean zero, then it holds that
	  $N^{-1/2}\sum_{m\neq i,j,\ell,k}V_{(s_1,m),s_2}$ is also a normal random variable with mean zero and variance $N^{-1}{\rm Var}\{\sum_{m\neq i,j,\ell,k}V_{(s_1,m),s_2}\}\asymp p^{-1}$. Therefore,
	\begin{align}
		&\max_{i,j,\ell,k:\,i\neq j\neq \ell\neq k}\max_{s_1,s_2:\,s_1\neq s_2\atop s_1,s_2\in\{i,j,\ell,k\}}\bigg|\frac{1}{\sqrt{N}}\sum_{m\neq i,j,\ell,k}V_{(s_1,m),s_2}\bigg|=O_{\p}\bigg(\frac{\log^{1/2}p}{p^{1/2}}\bigg)\,.\label{eq:argu2}
	\end{align}
	Also, as shown in Section \ref{sec:pfBell1ell23}, ${\rm Cov}\{V_{(m,s_1),s_2},V_{(m',s_1),s_2}\}={\rm Cov}\{Y_{(m,s_1),s_2},Y_{(m',s_1),s_2}\}=0$ for any $m\neq m'$ and $m,m'\notin\{i,j,\ell,k\}$. Then $N^{-1/2}\sum_{m\neq i,j,\ell,k}V_{(m,s_1),s_2}$ is a normal random variable with mean zero and variance $N^{-1}\sum_{m\neq i,j,\ell,k}{\rm Var}\{V_{(m,s_1),s_2}\}\asymp p^{-1}\nu^{-1}\gamma^{-6}$. Notice that $\nu\asymp p\gamma^{-2}+\gamma^{-6}$. Then  $N^{-1}\sum_{m\neq i,j,\ell,k}{\rm Var}\{V_{(m,s_1),s_2}\}\lesssim p^{-1}$, which implies
	\begin{align}
		&\max_{i,j,\ell,k:\,i\neq j\neq \ell\neq k}\max_{s_1,s_2:\,s_1\neq s_2\atop s_1,s_2\in\{i,j,\ell,k\}}\bigg|\frac{1}{\sqrt{N}}\sum_{m\neq i,j,\ell,k}V_{(m,s_1),s_2}\bigg|=O_{\p}\bigg(\frac{\log^{1/2}p}{p^{1/2}}\bigg)\,.\label{eq:argu3}
	\end{align}
	Identically,
	\begin{align*}
&\max_{i,j,\ell,k:\,i\neq j\neq \ell\neq k}\max_{s_1,s_2:\,s_1\neq s_2\atop s_1,s_2\in\{i,j,\ell,k\}}\bigg|\frac{1}{\sqrt{N}}\sum_{m\neq i,j,\ell,k}W_{(s_1,m),s_2}\bigg|\\
&~~~~~~~~~~~=O_{\p}\bigg(\frac{\log^{1/2}p}{p^{1/2}}\bigg)=\max_{i,j,\ell,k:\,i\neq j\neq \ell\neq k}\max_{s_1,s_2:\,s_1\neq s_2\atop s_1,s_2\in\{i,j,\ell,k\}}\bigg|\frac{1}{\sqrt{N}}\sum_{m\neq i,j,\ell,k}W_{(m,s_1),s_2}\bigg|\,.
\end{align*}
	Define
	\begin{align*}
		\mathcal{E}_2(Y)=&~\Bigg\{\max_{s_1,s_2:\,s_1\neq s_2\atop s_1,s_2\in\{i,j,\ell,k\}}\bigg|\frac{1}{\sqrt{N}}\sum_{m\neq i,j,\ell,k}Y_{(s_1,m),s_2}\bigg|\leq\frac{C_{**}\log^{1/2}p}{p^{1/2}}~~\textrm{and}\\
		&~~~\max_{s_1,s_2:\,s_1\neq s_2\atop s_1,s_2\in\{i,j,\ell,k\}}\bigg|\frac{1}{\sqrt{N}}\sum_{m\neq i,j,\ell,k}Y_{(m,s_1),s_2}\bigg|\leq\frac{C_{**}\log^{1/2}p}{p^{1/2}}~~\textrm{for any}~i\neq j\neq\ell\neq k\Bigg\}
	\end{align*}
	for some sufficiently large constant $C_{**}>0$. We can also define $\mathcal{E}_2(V)$ and $\mathcal{E}_2(W)$ in the same manner. Let
	\begin{equation}\label{eq:E2}
		\mathcal{E}_2=\mathcal{E}_2(Y)\cap\mathcal{E}_2(V)\cap\mathcal{E}_2(W)\,.
	\end{equation}
	 Then
	$
	\mathbb{P}(\mathcal{E}_2^c)\lesssim p^{-C}$, where $C$ can be sufficiently large if we select a sufficiently large $C_{**}$. Recall $B=C_*\log^{1/2}p$. Restricted on $\mathcal{E}_1\cap\mathcal{E}_2$, by \eqref{eq:cha}, we have
	\[
\max_{i,j,\ell,k:\,i\neq j\neq\ell\neq k}|\bc^{(i,j),\ell}-(\tau-1)\bc^{(i,j),\ell,(i,\ell),k-(i,j),\ell}|_\infty\leq \frac{84C_*\log^{1/2}p}{N^{1/2}}+\frac{144C_{**}\log^{1/2}p}{p^{1/2}}\,.
\]
	As $p\rightarrow\infty$, if $\phi\ll p^{1/2}(\log p)^{-3/2}$, it holds that
	\[
\frac{84C_*\log^{1/2}p}{N^{1/2}}+\frac{144C_{**}\log^{1/2}p}{p^{1/2}}<\frac{3}{4\beta}
\]
	with $\beta=\phi\log p$, which implies that
	\begin{equation}\label{eq:cha1}
		\max_{i,j,\ell,k:\,i\neq j\neq\ell\neq k}|\bc^{(i,j),\ell}-(\tau-1)\bc^{(i,j),\ell,(i,\ell),k-(i,j),\ell}|_\infty<\frac{3}{4\beta}
	\end{equation}
	under $\mathcal{E}_1\cap\mathcal{E}_2$.
	
	As shown in Lemma A.5 of \cite{CCK_2013}, there exists $U_{\ell km}(\bv)$ such that $|q_{\ell km}(\bv)|\leq U_{\ell km}(\bv)$ for any $\bv\in\mathbb{R}^{p}$, where $\sum_{\ell,k,m=1}^{p}U_{\ell km}(\bv)\lesssim \phi\beta^2$ for any $\bv\in\mathbb{R}^p$. Thus, \eqref{eq:term1} leads to
	\begin{align}\label{eq:Tij1Tij2}
\sum_{\ell\neq i,j}\sum_{k\neq i,j,\ell}\big|N\cdot\mathbb{E}\big[q_{\ell k}\{\bc^{-(i,j),\ell}\}\dot{c}_{(i,j),\ell}c_{(i,\ell),k}\big]\big|\lesssim T_{i,j,1}+T_{i,j,2}
\end{align}
with
\begin{align*}
&T_{i,j,1}=(\log p)\sum_{\ell\neq i,j}\sum_{k\neq i,j,\ell}\sum_{m=1}^p\int_0^1\mathbb{E}\big[ U_{\ell km}\{\bc^{-(i,j),\ell,(i,\ell),k}+\tau \bc^{(i,j),\ell,(i,\ell),k-(i,j),\ell}\}\\
&~~~~~~~~~~~~~~~~~~~~~~~~~~~~~~~~~~~~~~~~~~~~~~~~~~~~~~~~\times |c_m^{(i,j),\ell,(i,\ell),k-(i,j),\ell}|\cdot I(\mathcal{E}_1\cap\mathcal{E}_2)\big]\,{\rm d}\tau\,,\\
&T_{i,j,2}=\sum_{\ell\neq i,j}\sum_{k\neq i,j,\ell}\sum_{m=1}^p\int_0^1\mathbb{E}\big[I(\mathcal{E}_1^c\cup\mathcal{E}_2^c)\cdot U_{\ell km}\{\bc^{-(i,j),\ell,(i,\ell),k}+\tau \bc^{(i,j),\ell,(i,\ell),k-(i,j),\ell}\}\\
&~~~~~~~~~~~~~~~~~~~~~~~~~~~~~~~~~~~~~~~~~~~~~\times \{|Y_{(i,j),\ell}||Y_{(i,\ell),k}|+|V_{(i,j),\ell}||V_{(i,\ell),k}|+|W_{(i,j),\ell}||W_{(i,\ell),k}|\}\\
&~~~~~~~~~~~~~~~~~~~~~~~~~~~~~~~~~~~~~~~~~~~~~\times |c_m^{(i,j),\ell,(i,\ell),k-(i,j),\ell}|\big]\,{\rm d}\tau\,.\notag
\end{align*}
Due to
	$
	   \bc^{-(i,j),\ell,(i,\ell),k}+\tau \bc^{(i,j),\ell,(i,\ell),k-(i,j),\ell}=\bc-\bc^{(i,j),\ell}+(\tau-1)\bc^{(i,j),\ell,(i,\ell),k-(i,j),\ell}$, together with \eqref{eq:cha1}, Lemma A.6 of \cite{CCK_2013} implies that, restricted on $\mathcal{E}_1\cap\mathcal{E}_2$,
	\begin{align*}
	  U_{\ell km}(\bc)\lesssim U_{\ell km}\{\bc^{-(i,j),\ell,(i,\ell),k}+\tau \bc^{(i,j),\ell,(i,\ell),k-(i,j),\ell}\}\lesssim U_{\ell km}(\bc)
	\end{align*}
 for any $\tau\in[0,1]$. Thus,
	\begin{align*}
T_{i,j,1}\lesssim \phi\beta^2(\log p)\mathbb{E}\bigg\{I(\mathcal{E}_1\cap\mathcal{E}_2)\max_{k,\ell:\,k\neq\ell} |\bc^{(i,j),\ell,(i,\ell),k-(i,j),\ell}|_\infty\bigg\}\,.
\end{align*}
	Restricted on $\mathcal{E}_1\cap\mathcal{E}_2$, it follows from \eqref{eq:term2} that
	\[
\max_{k,\ell:\,k\neq\ell} |\bc^{(i,j),\ell,(i,\ell),k-(i,j),\ell}|_\infty\leq \frac{24C_*\log^{1/2}p}{p^{1/2}}+\frac{6C_{**}\log^{1/2}p}{p^{1/2}}\,,
\]
	which implies
	\begin{equation}\label{eq:Tij1}
		T_{i,j,1}\lesssim\frac{\phi\beta^2\log^{3/2}p}{p^{1/2}}\,.
	\end{equation}	
	For $T_{i,j,2}$, due to $U_{\ell km}\{\bc^{-(i,j),\ell,(i,\ell),k}(t)+\tau \bc^{(i,j),\ell,(i,\ell),k-(i,j),\ell}(t)\}\lesssim \phi\beta^2$, by Cauchy-Schwarz inequality, it holds that
	\begin{align*}
		T_{i,j,2}\lesssim&~\phi\beta^2\sum_{\ell\neq i,j}\sum_{k\neq i,j,\ell}\sum_{m=1}^p\mathbb{E}\big[I(\mathcal{E}_1^c\cup\mathcal{E}_2^c)|c_m^{(i,j),\ell,(i,\ell),k-(i,j),\ell}|\\
		&~~~~~~~~~~~~~~~~~~~~~~~~~~~~~~~~\times \{|Y_{(i,j),\ell}||Y_{(i,\ell),k}|+|V_{(i,j),\ell}||V_{(i,\ell),k}|+|W_{(i,j),\ell}||W_{(i,\ell),k}|\}\big]\\
		\lesssim&~p^3\phi\beta^2\mathbb{P}^{1/2}(\mathcal{E}_1^c\cup\mathcal{E}_2^c)\\
		&~~~~~~~\times\max_{\ell,k:\,\ell\neq k}\mathbb{E}^{1/2}\big[|\bc^{(i,j),\ell,(i,\ell),k-(i,j),\ell}|_\infty^2\{|Y_{(i,j),\ell}|^2|Y_{(i,\ell),k}|^2+|V_{(i,j),\ell}|^2|V_{(i,\ell),k}|^2\\
		&~~~~~~~~~~~~~~~~~~~~~~~~~~~~~~~~+|W_{(i,j),\ell}|^2|W_{(i,\ell),k}|^2\}\big]\\
		\lesssim&~p^3\phi\beta^2\mathbb{P}^{1/2}(\mathcal{E}_1^c\cup\mathcal{E}_2^c)\max_{\ell,k:\,\ell\neq k}\mathbb{E}^{1/4}\big\{|\bc^{(i,j),\ell,(i,\ell),k-(i,j),\ell}|_\infty^4\big\}\mathbb{E}^{1/8}\{|Y_{(i,j),\ell}|^8\}\\
		&~~~~~~~~~~~~~~~~~~~~~~~~~~~~~~~~~~~~~~~~~~~~~~~\times \mathbb{E}^{1/8}
		\{|Y_{(i,\ell),k}|^8\}\\
		&+p^3\phi\beta^2\mathbb{P}^{1/2}(\mathcal{E}_1^c\cup\mathcal{E}_2^c)\max_{\ell,k:\,\ell\neq k}\mathbb{E}^{1/4}\big\{|\bc^{(i,j),\ell,(i,\ell),k-(i,j),\ell}|_\infty^4\big\}\mathbb{E}^{1/8}\{|V_{(i,j),\ell}|^8\}\\
		&~~~~~~~~~~~~~~~~~~~~~~~~~~~~~~~~~~~~~~~~~~~~~~~\times\mathbb{E}^{1/8}\{|V_{(i,\ell),k}|^8\}\,,
	\end{align*}
	where the last step is based on the fact $\{V_{(i,j),\ell},V_{(i,\ell),k}\}$ and $\{W_{(i,j),\ell},W_{(i,\ell),k}\}$ are identically distributed. Notice that $\max_{i,j,\ell:\,i\neq j\neq \ell}|Y_{(i,j),\ell}|\lesssim \nu^{-1/2}\gamma^{-3}\lesssim1$ and $V_{(i,j),\ell}$ is a normal distributed random variable with ${\rm Var}\{V_{(i,j),\ell}\}\asymp \nu^{-1}\gamma^{-6}\lesssim1$. Thus
	\[
T_{i,j,2}\lesssim p^3\phi\beta^2\mathbb{P}^{1/2}(\mathcal{E}_1^c\cup\mathcal{E}_2^c)\max_{\ell,k:\,\ell\neq k}\mathbb{E}^{1/4}\big\{|\bc^{(i,j),\ell,(i,\ell),k-(i,j),\ell}|_\infty^4\big\}\,.
\]
By \eqref{eq:cijlilk-ijl}, following the same arguments for \eqref{eq:argu1}, \eqref{eq:argu2} and \eqref{eq:argu3}, we have
	\[
\max_{\ell,k:\,\ell\neq k}\mathbb{E}^{1/4}\big\{|\bc^{(i,j),\ell,(i,\ell),k-(i,j),\ell}|_\infty^4\big\}\lesssim1\,,
\]
	which implies $T_{i,j,2}\lesssim  p^{3}\phi\beta^2\mathbb{P}^{1/2}(\mathcal{E}_1^c\cup\mathcal{E}_2^c)$. Recall $\mathbb{P}(\mathcal{E}_1^c\cup\mathcal{E}_2^c)\lesssim p^{-C}$ and $C$ can be sufficiently large if we select two sufficiently large constants $C_*$ and $C_{**}$ in the definition of $\mathcal{E}_1$ and $\mathcal{E}_2$. Hence, with suitable selection of $(C_*,C_{**})$, we have
	\[
T_{i,j,2}\lesssim\frac{\phi\beta^2\log^{3/2}p}{p^{1/2}}\,.
\]
	Together with \eqref{eq:Tij1}, \eqref{eq:Tij1Tij2} implies that
	\begin{align*}
	\sum_{\ell:\,\ell\neq i,j}\sum_{k:\,k\neq i,j,\ell}\bigg|\int_0^1\mathbb{E}\big[q_{\ell k}\{\bc^{-(i,j),\ell}\}\dot{c}_{(i,j),\ell}c_{(i,\ell),k}\big]\,{\rm d}t\bigg|\lesssim \frac{{\phi\beta^2\log^{3/2}p}}{{p^{5/2}}}
\end{align*}
holds uniformly over $(i,j)$ such that $i\neq j$. Following the same arguments, we can bound the other terms in \eqref{eq:i2}. Then
	\[
	\sum_{\ell:\,\ell\neq i,j}\sum_{k:\,k\neq i,j,\ell}|{\rm I}_2(i,j,\ell,k)|\lesssim \frac{{\phi\beta^2\log^{3/2}p}}{{p^{5/2}}}
		\]
		holds uniformly over $(i,j)$ such that $i\neq j$.
	Since $\beta=\phi\log p$, we then have \eqref{eq:F.21}. $\hfill\Box$
	
	\subsubsection{Case 2: $k=i, j, \ell$.}\label{se:f2.2}
	We first consider the case with $k=i$. Notice that $c_{(i',j'),i}=0$ if $i'=i$ or $j'=i$. By \eqref{eq:I2ijellk}, it holds that
	\begin{align*}
{\rm I}_2(i,j,\ell,i)=&~\int_0^1\mathbb{E}\big[q_{\ell i}\{\bc^{-(i,j),\ell}\}\dot{c}_{(i,j),\ell}c_{(\ell,j),i}\big]\,{\rm d}t+\int_0^1\mathbb{E}\big[q_{\ell i}\{\bc^{-(i,j),\ell}\}\dot{c}_{(i,j),\ell}c_{(j,\ell),i}\big]\,{\rm d}t\\
&+\sum_{m\neq i,j,\ell}\int_0^1\mathbb{E}\big[q_{\ell i}\{\bc^{-(i,j),\ell}\}\dot{c}_{(i,j),\ell}c_{(m,j),i}\big]\,{\rm d}t\\
&+\sum_{m\neq i,j,\ell}\int_0^1\mathbb{E}\big[q_{\ell i}\{\bc^{-(i,j),\ell}\}\dot{c}_{(i,j),\ell}c_{(m,\ell),i}\big]\,{\rm d}t\\
&+\sum_{m\neq i,j,\ell}\int_0^1\mathbb{E}\big[q_{\ell i}\{\bc^{-(i,j),\ell}\}\dot{c}_{(i,j),\ell}c_{(j,m),i}\big]\,{\rm d}t\\
&+\sum_{m\neq i,j,\ell}\int_0^1\mathbb{E}\big[q_{\ell i}\{\bc^{-(i,j),\ell}\}\dot{c}_{(i,j),\ell}c_{(\ell,m),i}\big]\,{\rm d}t\,.
\end{align*}
	Note that
	\begin{align*}
\dot{c}_{(i,j),\ell}=&\,\frac{1}{\sqrt{N}}\bigg[\frac{1}{\sqrt{t}}\{\sqrt{v}{Y}_{(i,j),\ell}+\sqrt{1-v}V_{(i,j),\ell}\}-\frac{1}{\sqrt{1-t}}W_{(i,j),\ell}\bigg]\,,\\
c_{(j,\ell),i}=&\,\frac{1}{\sqrt{N}}\big[\sqrt{t}\big\{\sqrt{v}Y_{(j,\ell),i}+\sqrt{1-v}V_{(j,\ell),i}\big\}+\sqrt{1-t}W_{(j,\ell),i}\big]\,.
\end{align*}
	As shown in Sections \ref{sec:pflemma5} and \ref{se:f1}, $\{Y_{(i,j),\ell}, V_{(i,j),\ell}, W_{(i,j),\ell}\}$ is independent of $\bc^{-(i,j),\ell}$. Following the same arguments in Section \ref{se:f1} to show $V_{(i,j),\ell}$ is independent of $\bc^{-(i,j),\ell}$, we also know $\{V_{(j,\ell),i},W_{(j,\ell),i}\}$ is independent of $\bc^{-(i,j),\ell}$. Notice that $Y_{(j,\ell),i}$ is a function of $(\mathring{Z}_{j,i},\mathring{Z}_{i,\ell},\mathring{Z}_{j,\ell})=(\mathring{Z}_{i,j},\mathring{Z}_{i,\ell},\mathring{Z}_{\ell,j})$ and $Y_{(i,j),\ell}$ is a function of $(\mathring{Z}_{i,j},\mathring{Z}_{i,\ell},\mathring{Z}_{\ell,j})$. Based on the arguments to show $Y_{(i,j),\ell}$ is independent of $\bc^{-(i,j),\ell}$, we know $Y_{(j,\ell),i}$ is also independent of $\bc^{-(i,j),\ell}$. Then 
	\begin{align*}
&N\cdot\mathbb{E}\big[q_{\ell i}\{\bc^{-(i,j),\ell}\}\dot{c}_{(i,j),\ell}c_{(j,\ell),i}\big]\\
&~~~~~~~~~~=v\mathbb{E}\big[q_{\ell i}\{\bc^{-(i,j),\ell}\}\big]\mathbb{E}\{Y_{(i,j),\ell}Y_{(j,\ell),i}\}+(1-v)\mathbb{E}\big[q_{\ell i}\{\bc^{-(i,j),\ell}\}\big]\mathbb{E}\{V_{(i,j),\ell}V_{(j,\ell),i}\}\\
&~~~~~~~~~~~~~~-\mathbb{E}\big[q_{\ell i}\{\bc^{-(i,j),\ell}\}\big]\mathbb{E}\{W_{(i,j),\ell}W_{(j,\ell),i}\}\\
&~~~~~~~~~~=0\,.
\end{align*}
	Analogously, we also have $\mathbb{E}[q_{\ell i}\{\bc^{-(i,j),\ell}\}\dot{c}_{(i,j),\ell}c_{(\ell,j),i}]=0$. Then
	\begin{align}\label{eq:i2new}
{\rm I}_2(i,j,\ell,i)=&~\sum_{m\neq i,j,\ell}\int_0^1\mathbb{E}\big[q_{\ell i}\{\bc^{-(i,j),\ell}\}\dot{c}_{(i,j),\ell}c_{(m,j),i}\big]\,{\rm d}t\notag\\
&+\sum_{m\neq i,j,\ell}\int_0^1\mathbb{E}\big[q_{\ell i}\{\bc^{-(i,j),\ell}\}\dot{c}_{(i,j),\ell}c_{(m,\ell),i}\big]\,{\rm d}t\\
&+\sum_{m\neq i,j,\ell}\int_0^1\mathbb{E}\big[q_{\ell i}\{\bc^{-(i,j),\ell}\}\dot{c}_{(i,j),\ell}c_{(j,m),i}\big]\,{\rm d}t\notag\\
&+\sum_{m\neq i,j,\ell}\int_0^1\mathbb{E}\big[q_{\ell i}\{\bc^{-(i,j),\ell}\}\dot{c}_{(i,j),\ell}c_{(\ell,m),i}\big]\,{\rm d}t\notag\,.
\end{align}

	In the sequel, we only need to bound each term in \eqref{eq:i2new}.
	For $\mathbb{E}[q_{\ell i}\{\bc^{-(i,j),\ell}\}\dot{c}_{(i,j),\ell}c_{(\ell,m),i}]$ with $m\neq i,j,\ell$, it holds that
	\begin{align*}
&N\cdot\mathbb{E}\big[q_{\ell i}\{\bc^{-(i,j),\ell}\}\dot{c}_{(i,j),\ell}c_{(\ell,m),i}\big]\\
&~~~~~~~=v\mathbb{E}\big[q_{\ell i}\{\bc^{-(i,j),\ell}\}Y_{(i,j),\ell}Y_{(\ell,m),i}\big]+(1-v)\mathbb{E}\big[q_{\ell i}\{\bc^{-(i,j),\ell}\}V_{(i,j),\ell}V_{(\ell,m),i}\big]\\
&~~~~~~~~~~~-\mathbb{E}\big[q_{\ell i}\{\bc^{-(i,j),\ell}\}W_{(i,j),\ell}W_{(\ell,m),i}\big]\,.
\end{align*}
	Following the same arguments in Section \ref{se:f1} to show $V_{(i,j),\ell}$ is independent of $\bc^{-(i,j),\ell}$, we also have $\{V_{(\ell,m),i},W_{(\ell,m),i}\}$ is independent of $\bc^{-(i,j),\ell}$. Thus,
	\begin{align}
&N\cdot\mathbb{E}\big[q_{\ell i}\{\bc^{-(i,j),\ell}\}\dot{c}_{(i,j),\ell}c_{(\ell,m),i}\big]\notag\\
&~~~~~~~=v\mathbb{E}\big[q_{\ell i}\{\bc^{-(i,j),\ell}\}Y_{(i,j),\ell}Y_{(\ell,m),i}\big]-v\mathbb{E}\big[q_{\ell i}\{\bc^{-(i,j),\ell}\}\big]\mathbb{E}\{Y_{(i,j),\ell}Y_{(\ell,m),i}\}\,.\label{eq:asyex}
\end{align}
	Notice that $Y_{(\ell,m),i}$ is a function of $\{\mathring{Z}_{\ell,i},\mathring{Z}_{i,m},\mathring{Z}_{\ell,m}\}$, and $c_{(i',j'),\ell'}$'s involving $\mathring{Z}_{i,\ell}$ are not included in $\bc^{-(i,j),\ell}$. Similar to the strategy used in Section \ref{se:f2.1}, we can remove $c_{(i',j'),\ell'}$'s that related to $\{\mathring{Z}_{\ell,m},\mathring{Z}_{i,m}\}$ from $\bc^{-(i,j),\ell}$. Define
	\begin{align*}
\bc^{(i,j),\ell,(\ell,m),i}=&~\sum_{i',j':\,i'\neq j'}\bc_{(i',j')}\circ \big\{\ba_{(i',j')}^{(i,j,\ell)}+\ba_{(i',j')}^{(\ell,m)}+\ba_{(i',j')}^{(i,m)}-\ba_{(i',j')}^{(i,j,\ell)}\circ\ba_{(i',j')}^{(\ell,m)}\\
&~~~~~~~~~~~~~~~-\ba_{(i',j')}^{(i,j,\ell)}\circ\ba_{(i',j')}^{(i,m)}-\ba_{(i',j')}^{(\ell,m)}\circ\ba_{(i',j')}^{(i,m)}+\ba_{(i',j')}^{(i,j,\ell)}\circ\ba_{(i',j')}^{(\ell,m)}\circ\ba_{(i',j')}^{(i,m)}\big\}
\end{align*}
for $\ba_{(i',j')}^{(i,j,\ell)}$ and $\ba_{(i',j')}^{(i,j)}$ specified, respectively, in Sections \ref{sec:pflemma5} and \ref{se:f2.1}.
	Then $\{Y_{(i,j),\ell},Y_{(\ell,m),i}\}$ is independent of $\bc^{-(i,j),\ell,(\ell,m),i}:=\bc-\bc^{(i,j),\ell,(\ell,m),i}$. Let $\bc^{(i,j),\ell,(\ell,m),i-(i,j),\ell}=\bc^{-(i,j),\ell}-\bc^{-(i,j),\ell,(\ell,m),i}$. Then it holds that $\bc^{(i,j),\ell,(\ell,m),i-(i,j),\ell}=\bc^{(i,j),\ell,(\ell,m),i}-\bc^{(i,j),\ell}$. Recall
	$
	  	\bc^{(i,j),\ell}\\=\sum_{i',j':\,i'\neq j'}\bc_{(i',j')}\circ \ba_{(i',j')}^{(i,j,\ell)}$.
	We have
	 \begin{align*}
\bc^{(i,j),\ell,(\ell,m),i-(i,j),\ell}=&\,\{\bc_{(\ell,m)}+\bc_{(m,\ell)}+\bc_{(i,m)}+\bc_{(m,i)}\}\circ({\bf 1}-\bfe_i-\bfe_{j}-\bfe_\ell-\bfe_m)\\
&+\sum_{u\neq i,j,\ell,m}\{c_{(u,m),i}+c_{(m,u),i}\}\bfe_i+\sum_{u\neq i,j,\ell,m}\{c_{(u,m),\ell}+c_{(m,u),\ell}\}\bfe_\ell\\
&+\sum_{u\neq i,j,\ell,m}\{c_{(u,\ell),m}+c_{(\ell,u),m}+c_{(u,i),m}+c_{(i,u),m}\}\bfe_m\,.
\end{align*}
	Write
	$
	\bc^{(i,j),\ell,(\ell,m),i-(i,j),\ell}=\{c_1^{(i,j),\ell,(\ell,m),i-(i,j),\ell},\ldots,c_p^{(i,j),\ell,(\ell,m),i-(i,j),\ell}\}^\T
	$.
	By Taylor expansion, 
	\begin{align*}
\mathbb{E}\big[q_{\ell i}\{\bc^{-(i,j),\ell}\}Y_{(i,j),\ell}Y_{(\ell,m),i}\big]=&~\mathbb{E}\big[q_{\ell i}\{\bc^{-(i,j),\ell,(\ell,m),i}\}\big]\mathbb{E}\{Y_{(i,j),\ell}Y_{(\ell,m),i}\}\\
&+\sum_{s=1}^p\int_0^1\mathbb{E}\big[q_{\ell is}\{\bc^{-(i,j),\ell,(\ell,m),i}+\tau \bc^{(i,j),\ell,(\ell,m),i-(i,j),\ell}\}\\
&~~~~~~~~~~~~~~~~~~~~~~\times Y_{(i,j),\ell}Y_{(\ell,m),i}c_s^{(i,j),\ell,(\ell,m),i-(i,j),\ell}\big]\,{\rm d}\tau\,.
\end{align*}
	Together with \eqref{eq:asyex}, it holds that
			\begin{align*}
		&\big|N\cdot\mathbb{E}\big[q_{\ell i}\{\bc^{-(i,j),\ell}\}\dot{c}_{(i,j),\ell}c_{(\ell,m),i}\big]\big|\\
		&~~~~~~~~~~\leq
		\underbrace{\big|\mathbb{E}\big[q_{\ell i}\{\bc^{-(i,j),\ell,(\ell,m),i}\}\big]-\mathbb{E}\big[q_{\ell i}\{\bc^{-(i,j),\ell}\}\big]\big|\big|\mathbb{E}\{Y_{(i,j),\ell}Y_{(\ell,m),i}\}\big|}_{R_1(i,j,\ell,m)}\\
		&~~~~~~~~~~~~~~+\sum_{s=1}^p\int_0^1\mathbb{E}\big[|q_{\ell is}\{\bc^{-(i,j),\ell,(\ell,m),i}+\tau \bc^{(i,j),\ell,(\ell,m),i-(i,j),\ell}\}|\\
		&~~~~~~~~~~~~~~~~~~\underbrace{~~~~~~~~~~~~~~~~~~~\times|Y_{(i,j),\ell}||Y_{(\ell,m),i}||c_s^{(i,j),\ell,(\ell,m),i-(i,j),\ell}|\big]\,{\rm d}\tau}_{R_2(i,j,\ell,m)}\,.
	\end{align*}
	Note that
	\begin{align*}
	   \bc^{-(i,j),\ell,(\ell,m),i}+\tau \bc^{(i,j),\ell,(\ell,m),i-(i,j),\ell}=\bc-\bc^{(i,j),\ell}+(\tau-1)\bc^{(i,j),\ell,(\ell,m),i-(i,j),\ell}\,.
	\end{align*}
 Following the identical arguments in Section \ref{se:f2.1} for bounding the term on the right-hand side of \eqref{eq:term1}, we have
 \begin{align}\label{eq:R2ijlm}
   \sum_{i:\,i\neq j,m}\sum_{\ell:\,\ell\neq i,j,m}R_2(i,j,\ell,m)\lesssim \frac{\phi\beta^2\log^{3/2}p}{p^{1/2}}\,.
 \end{align}
	By Taylor expansion,
	\begin{align*}
&\mathbb{E}\big[q_{\ell i}\{\bc^{-(i,j),\ell}\}\big]-\mathbb{E}\big[q_{\ell i}\{\bc^{-(i,j),\ell,(\ell,m),i}\}\big]\\
&~~~~~~~~~=\sum_{s=1}^p\int_0^1\mathbb{E}\big[q_{\ell is}\{\bc^{-(i,j),\ell,(\ell,m),i}+\tau \bc^{(i,j),\ell,(\ell,m),i-(i,j),\ell}\}c_s^{(i,j),\ell,(\ell,m),i-(i,j),\ell}\big]\,{\rm d}\tau\,.
\end{align*}
Parellel to \eqref{eq:R2ijlm}, we also have
	\[
\sum_{i:\,i\neq j,m}\sum_{\ell:\,\ell\neq i,j,m}R_1(i,j,\ell,m)\lesssim \frac{\phi\beta^2\log^{3/2}p}{p^{1/2}}\,.
\]
	Therefore,
	\[
\sum_{i:\,i\neq j,m}\sum_{\ell:\,\ell\neq i,j,m}\big|N\cdot\mathbb{E}\big[q_{\ell i}\{\bc^{-(i,j),\ell}\}\dot{c}_{(i,j),\ell}c_{(\ell,m),i}\big]\big|\lesssim \frac{\phi\beta^2\log^{3/2}p}{p^{1/2}}\,,
\]
	which implies
	\[
\sum_{i,j,\ell:\,i\neq j\neq\ell}\bigg|\sum_{m\neq i,j,\ell}\int_0^1\mathbb{E}\big[q_{\ell i}\{\bc^{-(i,j),\ell}\}\dot{c}_{(i,j),\ell}c_{(\ell,m),i}\big]\,{\rm d}t\notag\bigg|\lesssim \frac{\phi\beta^2\log^{3/2}p}{p^{1/2}}\,.
\]
	Analogously, we can obtain the same result for other terms in \eqref{eq:i2new}. Recall $\beta=\phi\log p$. Hence, we have that
	\begin{align*}
\sum_{i,j,\ell:\,i\neq j\neq\ell}|{\rm I}_2(i,j,\ell,i)|\lesssim \frac{\phi^3\log^{7/2}p}{p^{1/2}}\,.
\end{align*}
Identically, we also have
\begin{align*}
\sum_{i,j,\ell:\,i\neq j\neq\ell}|{\rm I}_2(i,j,\ell,j)|\lesssim \frac{\phi^3\log^{7/2}p}{p^{1/2}}~~~\textrm{and}~~~
\sum_{i,j,\ell:\,i\neq j\neq\ell}|{\rm I}_2(i,j,\ell,\ell)|\lesssim \frac{\phi^3\log^{7/2}p}{p^{1/2}}\,.
\end{align*}
We complete the proof of \eqref{eq:F.22}. $\hfill\Box$	
	
	\subsection{Proof of \eqref{eq:sumI3}.}\label{se:f3}
	
	To simplify the notation, we write $\bc(t)$, $\bc^{-(i,j),\ell}(t)$, $c_k^{(i,j),\ell}(t)$, $c_{(i,j),\ell}(t)$ and $\dot{c}_{(i,j),\ell}(t)$ as $\bc$, $\bc^{-(i,j),\ell}$, $c_k^{(i,j),\ell}$, $c_{(i,j),\ell}$ and $\dot{c}_{(i,j),\ell}$, respectively. Define
	\begin{align}\label{eq:mathcalE}
\mathcal{E}=\big\{|Y_{(i,j),\ell}|\vee|V_{(i,j),\ell}|\vee|W_{(i,j),\ell}|\leq p^{1/2}/(4\beta)~\textrm{for any}~i\neq j\neq \ell\big\}\,.
\end{align}
	We then have
	\begin{align*}
		{\rm I}_3(i,j,\ell,k,l)=&~\underbrace{\int_0^1\int_0^1(1-\tau)\mathbb{E}\big[I(\mathcal{E}^c)q_{\ell kl}\{\bc^{-(i,j),\ell}+\tau\bc^{(i,j),\ell}\}\dot{c}_{(i,j),\ell}c^{(i,j),\ell}_kc^{(i,j),\ell}_l\big]\,{\rm d}\tau{\rm d}t}_{{\rm I}_{3,1}(i,j,\ell,k,l)}\\
		&+\underbrace{\int_0^1\int_0^1(1-\tau)\mathbb{E}\big[I(\mathcal{E})q_{\ell kl}\{\bc^{-(i,j),\ell}+\tau\bc^{(i,j),\ell}\}\dot{c}_{(i,j),\ell}c^{(i,j),\ell}_kc^{(i,j),\ell}_l\big]\,{\rm d}\tau{\rm d}t}_{{\rm I}_{3,2}(i,j,\ell,k,l)}\,.
	\end{align*}
	Let $\omega(t)=1/(\sqrt{t}\wedge\sqrt{1-t})$ for any $t\in(0,1)$. Notice that
	\begin{align*}
\dot{c}_{(i,j),\ell}=\frac{1}{\sqrt{N}}\bigg[\frac{1}{\sqrt{t}}\{\sqrt{v}{Y}_{(i,j),\ell}+\sqrt{1-v}V_{(i,j),\ell}\}-\frac{1}{\sqrt{1-t}}W_{(i,j),\ell}\bigg]\,.
\end{align*}
	Then
	\begin{align*}
\max_{i,j,\ell:\,i\neq j\neq\ell}|\dot{c}_{(i,j),\ell}|\lesssim \frac{\omega(t)}{p}\max_{i,j,\ell:\,i\neq j\neq\ell}\big\{|Y_{(i,j),\ell}|\vee|V_{(i,j),\ell}|\vee|W_{(i,j),\ell}|\big\}\,.
\end{align*}
	On the other hand, same as \eqref{eq:keybound}, it holds that
	\begin{align}
		\max_{i,j,\ell:\,i\neq j\neq \ell}|\bc^{(i,j),\ell}|_\infty\leq&~\frac{54}{\sqrt{N}}\max_{i,j,\ell:\,i\neq j\neq \ell}\big\{|Y_{(i,j),\ell}|\vee|V_{(i,j),\ell}|\vee|W_{(i,j),\ell}|\big\}\notag\\
		&+\max_{i,j,\ell:\,i\neq j\neq \ell}\sum_{s_1,s_2:\,s_1\neq s_2\atop s_1,s_2\in\{i,j,\ell\}}\Bigg|\frac{1}{\sqrt{N}}\sum_{m\neq i,j,\ell}Y_{(s_1,m),s_2}\Bigg|\notag\\
		&+\max_{i,j,\ell:\,i\neq j\neq\ell}\sum_{s_1,s_2:\,s_1\neq s_2\atop s_1,s_2\in\{i,j,\ell\}}\Bigg|\frac{1}{\sqrt{N}}\sum_{m\neq i,j,\ell}Y_{(m,s_1),s_2}\Bigg|\notag\\
		&+\max_{i,j,\ell:\,i\neq j\neq \ell}\sum_{s_1,s_2:\,s_1\neq s_2\atop s_1,s_2\in\{i,j,\ell\}}\Bigg|\frac{1}{\sqrt{N}}\sum_{m\neq i,j,\ell}V_{(s_1,m),s_2}\Bigg|\label{eq:boundke}\\
		&+\max_{i,j,\ell:\,i\neq j\neq \ell}\sum_{s_1,s_2:\,s_1\neq s_2\atop s_1,s_2\in\{i,j,\ell\}}\Bigg|\frac{1}{\sqrt{N}}\sum_{m\neq i,j,\ell}V_{(m,s_1),s_2}\Bigg|\notag\\
		&+\max_{i,j,\ell:\,i\neq j\neq \ell}\sum_{s_1,s_2:\,s_1\neq s_2\atop s_1,s_2\in\{i,j,\ell\}}\Bigg|\frac{1}{\sqrt{N}}\sum_{m\neq i,j,\ell}W_{(s_1,m),s_2}\Bigg|\notag\\
		&+\max_{i,j,\ell:\,i\neq j\neq \ell}\sum_{s_1,s_2:\,s_1\neq s_2\atop s_1,s_2\in\{i,j,\ell\}}\Bigg|\frac{1}{\sqrt{N}}\sum_{m\neq i,j,\ell}W_{(m,s_1),s_2}\Bigg|\notag\,.
	\end{align}
	We define
	\begin{align*}
		\mathcal{E}(Y)=&~\Bigg\{\max_{s_1,s_2:\,s_1\neq s_2\atop s_1,s_2\in\{i,j,\ell\}}\bigg|\frac{1}{\sqrt{N}}\sum_{m\neq i,j,\ell}Y_{(s_1,m),s_2}\bigg|\leq\frac{C_{*}\log^{1/2}p}{p^{1/2}}~~\textrm{and}\\
		&~~~\max_{s_1,s_2:\,s_1\neq s_2\atop s_1,s_2\in\{i,j,\ell\}}\bigg|\frac{1}{\sqrt{N}}\sum_{m\neq i,j,\ell}Y_{(m,s_1),s_2}\bigg|\leq\frac{C_{*}\log^{1/2}p}{p^{1/2}}~~\textrm{for any}~i\neq j\neq\ell\Bigg\}
	\end{align*}
	for some sufficiently large constant $C_{*}>0$. We can also define $\mathcal{E}(V)$ and $\mathcal{E}(W)$ in the same manner. Let
	$
	\tilde{\mathcal{E}}=\mathcal{E}(Y)\cap\mathcal{E}(V)\cap\mathcal{E}(W)$.
	Using the same arguments in Section \ref{se:f2.1} to derive the upper bound of $\mathbb{P}(\mathcal{E}_2^c)$ for $\mathcal{E}_2$ specified in \eqref{eq:E2}, it holds that $
	\mathbb{P}(\tilde{\mathcal{E}}^c)\lesssim p^{-C}$, where $C$ can be sufficiently large if we select a sufficiently large $C_{*}$. Hence, restricted on $\tilde{\mathcal{E}}$, we have
	\begin{align*}
\max_{i,j,\ell:\,i\neq j\neq \ell}|\bc^{(i,j),\ell}|_\infty\lesssim \frac{1}{p}\max_{i,j,\ell:\,i\neq j\neq \ell}\big\{|Y_{(i,j),\ell}|\vee|V_{(i,j),\ell}|\vee|W_{(i,j),\ell}|\big\}+\frac{\log^{1/2}p}{p^{1/2}}\,,
\end{align*}
which implies that
\begin{align}
\max_{i,j,\ell:\,i\neq j\neq\ell}|\dot{c}_{(i,j),\ell}||\bc^{(i,j),\ell}|_\infty^2\lesssim&~\frac{\omega(t)}{p^3}\max_{i,j,\ell:\,i\neq j\neq \ell}\big\{|Y_{(i,j),\ell}|^3\vee|V_{(i,j),\ell}|^3\vee|W_{(i,j),\ell}|^3\big\}\notag\\
&+\frac{\omega(t)\log p}{p^2}\max_{i,j,\ell:\,i\neq j\neq \ell}\big\{|Y_{(i,j),\ell}|\vee|V_{(i,j),\ell}|\vee|W_{(i,j),\ell}|\big\}\label{eq:f.41}
\end{align}
	under $\tilde{\mathcal{E}}$.
	As shown in Lemma A.5 of \cite{CCK_2013}, there exists $U_{\ell kl}(\bv)$ such that $|q_{\ell kl}(\bv)|\leq U_{\ell kl}(\bv)$ for any $\bv\in\mathbb{R}^{p}$, where $\sum_{\ell,k,l=1}^{p}U_{\ell kl}(\bv)\lesssim \phi\beta^2$ for any $\bv\in\mathbb{R}^p$. Then
	\begin{align}\label{eq:i31sum}
		\sum_{\ell,k,l=1}^p|{\rm I}_{3,1}(i,j,\ell,k,l)|\lesssim&~\phi\beta^2p^3\int_0^1\mathbb{E}\bigg\{I(\mathcal{E}^c)\max_{i,j,\ell:\,i\neq j\neq\ell}|\dot{c}_{(i,j),\ell}||\bc^{(i,j),\ell}|_\infty^2\bigg\}\,{\rm d}t\notag\\
		=&~\phi\beta^2p^3\int_0^1\mathbb{E}\bigg\{I(\mathcal{E}^c\cap\tilde{\mathcal{E}})\max_{i,j,\ell:\,i\neq j\neq\ell}|\dot{c}_{(i,j),\ell}||\bc^{(i,j),\ell}|_\infty^2\bigg\}\,{\rm d}t\\
		&+\phi\beta^2p^3\int_0^1\mathbb{E}\bigg\{I(\mathcal{E}^c\cap\tilde{\mathcal{E}}^c)\max_{i,j,\ell:\,i\neq j\neq\ell}|\dot{c}_{(i,j),\ell}||\bc^{(i,j),\ell}|_\infty^2\bigg\}\,{\rm d}t\,.\notag
	\end{align}
	Notice that $\max_{i,j,\ell:\,i\neq j\neq \ell}|Y_{(i,j),\ell}|\leq C$, and $V_{(i,j),\ell}$ and $W_{(i,j),\ell}$ are normal random variables. It holds that
	\begin{equation}\label{eq:tailbound}
		\mathbb{P}\bigg[\max_{i,j,\ell:\,i\neq j\neq \ell}\big\{|Y_{(i,j),\ell}|\vee|V_{(i,j),\ell}|\vee|W_{(i,j),\ell}|\big\}>u\bigg]\leq Cp^3\exp(-Cu^2)
	\end{equation}
	for any $u>0$. Thus, for $\mathcal{E}$ defined as \eqref{eq:mathcalE}, we have $\mathbb{P}(\mathcal{E}^c)\lesssim p^3\exp(-C\beta^{-2}p)$. By \eqref{eq:f.41},
	\begin{align*}
	    &\int_0^1\mathbb{E}\bigg\{I(\mathcal{E}^c\cap\tilde{\mathcal{E}})\max_{i,j,\ell:\,i\neq j\neq\ell}|\dot{c}_{(i,j),\ell}||\bc^{(i,j),\ell}|_\infty^2\bigg\}\,{\rm d}t\\ &~~~~~~~~~~~~~~\lesssim\frac{1}{p^{3}}\cdot\mathbb{E}\bigg[I(\mathcal{E}^c)\max_{i,j,\ell:\,i\neq j\neq \ell}\big\{|Y_{(i,j),\ell}|^3\vee|V_{(i,j),\ell}|^3\vee|W_{(i,j),\ell}|^3\big\}\bigg]\\
	&~~~~~~~~~~~~~~~~~~+\frac{\log p}{p^{2}}\cdot\mathbb{E}\bigg[I(\mathcal{E}^c)\max_{i,j,\ell:\,i\neq j\neq \ell}\big\{|Y_{(i,j),\ell}|\vee|V_{(i,j),\ell}|\vee|W_{(i,j),\ell}|\big\}\bigg].
	\end{align*}
	By Cauchy-Schwarz inequality, we have
	\begin{align*}
		&\mathbb{E}\bigg[I(\mathcal{E}^c)\max_{i,j,\ell:\,i\neq j\neq \ell}\big\{|Y_{(i,j),\ell}|\vee|V_{(i,j),\ell}|\vee|W_{(i,j),\ell}|\big\}\bigg]\\
		&~~~~~~~~~~\leq
\mathbb{P}^{1/2}(\mathcal{E}^c)\cdot\mathbb{E}^{1/2}\bigg[\max_{i,j,\ell:\,i\neq j\neq \ell}\big\{|Y_{(i,j),\ell}|^2\vee|V_{(i,j),\ell}|^2\vee|W_{(i,j),\ell}|^2\big\}\bigg]\\
&~~~~~~~~~~\lesssim p^{3}\exp\bigg(-\frac{Cp}{\beta^{2}}\bigg)\,.
\end{align*}
Analogously, we also have
\begin{align*}
\mathbb{E}\bigg[I(\mathcal{E}^c)\max_{i,j,\ell:\,i\neq j\neq \ell}\big\{|Y_{(i,j),\ell}|^3\vee|V_{(i,j),\ell}|^3\vee|W_{(i,j),\ell}|^3\big\}\bigg]\lesssim p^{3}\exp\bigg(-\frac{Cp}{\beta^{2}}\bigg)\,.
\end{align*}
	Hence,
	\begin{align}\label{eq:i311}
		\int_0^1\mathbb{E}\bigg\{I(\mathcal{E}^c\cap\tilde{\mathcal{E}})\max_{i,j,\ell:\,i\neq j\neq\ell}|\dot{c}_{(i,j),\ell}||\bc^{(i,j),\ell}|_\infty^2\bigg\}\,{\rm d}t\lesssim p(\log p)\exp\bigg(-\frac{Cp}{\beta^{2}}\bigg)\,.
	\end{align}
	By Cauchy-Schwarz inequality, it holds that
	\begin{align*}
		&\mathbb{E}\bigg\{I(\mathcal{E}^c\cap\tilde{\mathcal{E}}^c)\max_{i,j,\ell:\,i\neq j\neq\ell}|\dot{c}_{(i,j),\ell}||\bc^{(i,j),\ell}|_\infty^2\bigg\}\\
		&~~~~~~~~~\leq \mathbb{P}^{1/2}(\tilde{\mathcal{E}^c})\cdot\mathbb{E}^{1/2}\bigg\{\max_{i,j,\ell:\,i\neq j\neq\ell}|\dot{c}_{(i,j),\ell}|^2|\bc^{(i,j),\ell}|_\infty^4\bigg\}\\
		&~~~~~~~~~\leq \mathbb{P}^{1/2}(\tilde{\mathcal{E}}^c)\cdot\mathbb{E}^{1/4}\bigg\{\max_{i,j,\ell:\,i\neq j\neq\ell}|\dot{c}_{(i,j),\ell}|^4\bigg\}\cdot\mathbb{E}^{1/4}\bigg\{\max_{i,j,\ell:\,i\neq j\neq \ell}|\bc^{(i,j),\ell}|_\infty^8\bigg\}\,.
	\end{align*}
	Notice that $
	\mathbb{P}(\tilde{\mathcal{E}}^c)\lesssim p^{-C}$, where $C$ can be sufficiently large if we select a sufficiently large $C_{*}$ in the definition of $\tilde{\mathcal{E}}$. Thus, with suitable selection of $C_*$, we have
	\begin{align}\label{eq:key1}
		\int_0^1\mathbb{E}\bigg\{I(\mathcal{E}^c\cap\tilde{\mathcal{E}}^c)\max_{i,j,\ell:\,i\neq j\neq\ell}|\dot{c}_{(i,j),\ell}||\bc^{(i,j),\ell}|_\infty^2\bigg\}\,{\rm d}t\lesssim\frac{1}{p^{11/2}}\,.
	\end{align}
	Together with \eqref{eq:i311} and $\beta=\phi\log p$, \eqref{eq:i31sum} implies that
	\begin{align}\label{eq:i31}
		\sum_{i,j,\ell:\,i\neq j\neq \ell}\sum_{k,l=1}^p|{\rm I}_{3,1}(i,j,\ell,k,l)|\lesssim&~\frac{\phi\beta^2}{p^{1/2}}+\phi\beta^2p^6(\log p)\exp(-Cp\beta^{-2})\notag\\
		=&~\frac{\phi^3\log^2p}{p^{1/2}}+\phi^3p^6\log ^3p\cdot\exp\bigg(-\frac{Cp}{\phi^2\log^2p}\bigg)\,.
	\end{align}
	
	In the sequel, we consider ${\rm I}_{3,2}(i,j,\ell,k,l)$. Due to $|q_{\ell kl}(\bv)|\leq U_{\ell kl}(\bv)$ for any $\bv\in\mathbb{R}^{p}$, by triangle inequality,
	\begin{align}\label{eq:i32}
		|{\rm I}_{3,2}(i,j,\ell,k,l)|\lesssim&~\int_0^1\int_0^1(1-\tau)\mathbb{E}\big[I(\mathcal{E}\cap\tilde{\mathcal{E}})U_{\ell kl}\{\bc^{-(i,j),\ell}+\tau\bc^{(i,j),\ell}\}\notag\\
		&~~~~~~~~~~~~~~~~~~~~~~~~~~~~~~\times|\dot{c}_{(i,j),\ell}||c^{(i,j),\ell}_k||c^{(i,j),\ell}_l|\big]\,{\rm d}\tau{\rm d}t\\
		&+\int_0^1\int_0^1(1-\tau)\mathbb{E}\big[I(\mathcal{E}\cap\tilde{\mathcal{E}}^c)U_{\ell kl}\{\bc^{-(i,j),\ell}+\tau\bc^{(i,j),\ell}\}\notag\\
		&~~~~~~~~~~~~~~~~~~~~~~~~~~~~~~\times|\dot{c}_{(i,j),\ell}||c^{(i,j),\ell}_k||c^{(i,j),\ell}_l|\big]\,{\rm d}\tau{\rm d}t\,.\notag
	\end{align}
	Since $\sum_{\ell,k,l=1}^{p}U_{\ell kl}(\bv)\lesssim \phi\beta^2$ for any $\bv\in\mathbb{R}^p$, then
	\begin{align*}
		&\int_0^1\int_0^1(1-\tau)\mathbb{E}\big[I(\mathcal{E}\cap\tilde{\mathcal{E}}^c)U_{\ell kl}\{\bc^{-(i,j),\ell}+\tau\bc^{(i,j),\ell}\}|\dot{c}_{(i,j),\ell}||c^{(i,j),\ell}_k||c^{(i,j),\ell}_l|\big]\,{\rm d}\tau{\rm d}t\\
		&~~~~~~~~~~~~~~\lesssim \phi\beta^2\int_0^1\mathbb{E}\big[I(\mathcal{E}\cap\tilde{\mathcal{E}}^c)|\dot{c}_{(i,j),\ell}||c^{(i,j),\ell}_k||c^{(i,j),\ell}_l|\big]\,{\rm d}t\\
		&~~~~~~~~~~~~~~\lesssim \phi\beta^2\int_0^1\mathbb{E}\bigg\{I(\mathcal{E}\cap\tilde{\mathcal{E}}^c)\max_{i,j,\ell:\,i\neq j\neq\ell}|\dot{c}_{(i,j),\ell}||\bc^{(i,j),\ell}|_\infty^2\bigg\}\,{\rm d}t\,.
	\end{align*}
	Same as \eqref{eq:key1}, we have
	\begin{align}\label{eq:bd0}
		&\sum_{i,j,\ell:\,i\neq j\neq\ell}\sum_{k,l=1}^p\int_0^1\int_0^1(1-\tau)\mathbb{E}\big[I(\mathcal{E}\cap\tilde{\mathcal{E}}^c)U_{\ell kl}\{\bc^{-(i,j),\ell}+\tau\bc^{(i,j),\ell}\}\notag\\
		&~~~~~~~~~~~~~~~~~~~~~~~~~~~~~~~~~~~~~~~~~~~~~~~~~~\times|\dot{c}_{(i,j),\ell}||c^{(i,j),\ell}_k||c^{(i,j),\ell}_l|\big]\,{\rm d}\tau{\rm d}t\lesssim \frac{\phi^3\log^2p}{p^{1/2}}\,.
	\end{align}
	Restricted on $\mathcal{E}\cap\tilde{\mathcal{E}}$, \eqref{eq:boundke} implies that
	\begin{align*}
\max_{i,j,\ell:\,i\neq j\neq \ell}|\bc^{(i,j),\ell}|_\infty\lesssim \frac{1}{p^{1/2}\beta}+\frac{\log^{1/2}p}{p^{1/2}}\leq\frac{3}{4\beta}
\end{align*}
	for sufficiently large $p$ if $\phi\ll p^{1/2}(\log p)^{-3/2}$. Lemma A.6 of \cite{CCK_2013} implies that, restricted on $\mathcal{E}\cap\tilde{\mathcal{E}}$, $U_{\ell kl}(\bc)\lesssim U_{\ell kl}\{\bc^{-(i,j),\ell}+\tau\bc^{(i,j),\ell}\}\lesssim U_{\ell kl}(\bc)$ for any $\tau\in[0,1]$. Hence,
	\begin{align*}
		&\int_0^1\int_0^1(1-\tau)\mathbb{E}\big[I(\mathcal{E}\cap\tilde{\mathcal{E}})U_{\ell kl}\{\bc^{-(i,j),\ell}+\tau\bc^{(i,j),\ell}\}|\dot{c}_{(i,j),\ell}||c^{(i,j),\ell}_k||c^{(i,j),\ell}_l|\big]\,{\rm d}\tau{\rm d}t\\
		&~~~~~~~~~~~~~~~~\lesssim \int_0^1\underbrace{\mathbb{E}\big\{I(\mathcal{E}\cap\tilde{\mathcal{E}})U_{\ell kl}(\bc)|\dot{c}_{(i,j),\ell}||c^{(i,j),\ell}_k||c^{(i,j),\ell}_l|\big\}}_{R_3(i,j,\ell,k,l)}\,{\rm d}t\,.
	\end{align*}
	Notice that
	\begin{align*}
		c_k^{(i,j),\ell}=\left\{\begin{aligned}
			c_{(i,j),k}+c_{(j,i),k}+c_{(j,\ell),k}+c_{(\ell,j),k}+c_{(\ell,i),k}+c_{(i,\ell),k}\,,~&~\textrm{if}~k\neq i,j,\ell\,, \\
			\sum_{m\neq i,j,\ell}\{c_{(m,j),i}+c_{(j,m),i}\}+\sum_{m\neq i,\ell}\{c_{(m,\ell),i}+c_{(\ell,m),i}\}\,,\,\,&~\textrm{if}~k=i\,,\\
			\sum_{m\neq i,j,\ell}\{c_{(m,i),j}+c_{(i,m),j}\}+\sum_{m\neq j,\ell}\{c_{(m,\ell),j}+c_{(\ell,m),j}\}\,,\,\,&~\textrm{if}~k=j\,,\\
			\sum_{m\neq i,j,\ell}\{c_{(m,j),\ell}+c_{(j,m),\ell}\}+\sum_{m\neq i,\ell}\{c_{(m,i),\ell}+c_{(i,m),\ell}\}\,,\,\,&~\textrm{if}~k=\ell\,.
		\end{aligned} \right.
	\end{align*}
Hence, if $k\in\{i, j, \ell\}$, we have
	\begin{align*}
		|c_k^{(i,j),\ell}|\lesssim&~\sum_{s_1,s_2:\,s_1\neq s_2\atop s_1,s_2\in\{i,j,\ell\}}\Bigg|\frac{1}{\sqrt{N}}\sum_{m\neq i,j,\ell}Y_{(s_1,m),s_2}\Bigg|\notag+\sum_{s_1,s_2:\,s_1\neq s_2\atop s_1,s_2\in\{i,j,\ell\}}\Bigg|\frac{1}{\sqrt{N}}\sum_{m\neq i,j,\ell}Y_{(m,s_1),s_2}\Bigg|\notag\\
		&+\sum_{s_1,s_2:\,s_1\neq s_2\atop s_1,s_2\in\{i,j,\ell\}}\Bigg|\frac{1}{\sqrt{N}}\sum_{m\neq i,j,\ell}V_{(s_1,m),s_2}\Bigg|+\sum_{s_1,s_2:\,s_1\neq s_2\atop s_1,s_2\in\{i,j,\ell\}}\Bigg|\frac{1}{\sqrt{N}}\sum_{m\neq i,j,\ell}V_{(m,s_1),s_2}\Bigg|\notag\\
		&+\sum_{s_1,s_2:\,s_1\neq s_2\atop s_1,s_2\in\{i,j,\ell\}}\Bigg|\frac{1}{\sqrt{N}}\sum_{m\neq i,j,\ell}W_{(s_1,m),s_2}\Bigg|\notag+\sum_{s_1,s_2:\,s_1\neq s_2\atop s_1,s_2\in\{i,j,\ell\}}\Bigg|\frac{1}{\sqrt{N}}\sum_{m\neq i,j,\ell}W_{(m,s_1),s_2}\Bigg|\notag\,.
	\end{align*}
Restricted on $\tilde{\mathcal{E}}$,
	\begin{align*}
\max_{i,j,\ell:\,i\neq j\neq \ell}\max_{k\in\{i,j,\ell\}}|c_k^{(i,j),\ell}|\lesssim\frac{\log^{1/2}p}{p^{1/2}}\,.
\end{align*}	
	
Since $\{W_{(i,j),\ell}\}$ is an independent copy of $\{V_{(i,j),\ell}\}$ and $\sum_{\ell,k,l=1}^{p}U_{\ell kl}(\bv)\lesssim \phi\beta^2$ for any $\bv\in\mathbb{R}^p$, by \eqref{eq:cijt1} and \eqref{eq:dotcijt1}, we have
	\begin{align}\label{eq:bd1}		
	&\sum_{i,j,\ell:\,i\neq j\neq\ell}\sum_{k:\,k\neq i,j,\ell}\sum_{l:\,l\neq i,j,\ell}R_3(i,j,\ell,k,l)\notag\\
		&~~~~~~~~~~~\lesssim\phi\beta^2\sum_{i,j:\,i\neq j}\mathbb{E}\bigg\{\max_{\ell:\,\ell\neq i,j}|\dot{c}_{(i,j),\ell}|\cdot\max_{\ell:\,\ell\neq i,j}\max_{k:\,k\neq i,j,\ell}|c_k^{(i,j),\ell}|^2\bigg\}\notag\\
		&~~~~~~~~~~~\lesssim\frac{\phi\beta^2\omega(t)}{p}\cdot\mathbb{E}\bigg[\max_{i,j,\ell:\,i\neq j\neq\ell}\big\{|Y_{(i,j),\ell}|^3\vee|V_{(i,j),\ell}|^3\vee|W_{(i,j),\ell}|^3\big\}\bigg]\\
		&~~~~~~~~~~~\lesssim\frac{\phi\beta^2\omega(t)}{p}+\frac{\phi\beta^2\omega(t)}{p}\cdot\mathbb{E}\bigg\{\max_{i,j,\ell:\,i\neq j\neq\ell}|V_{(i,j),\ell}|^3\bigg\}\notag\\		
		&~~~~~~~~~~~\lesssim\frac{\phi\beta^2\omega(t)\log^{3/2}p}{p}\,,\notag\\
&\sum_{i,j,\ell:\,i\neq j\neq \ell}\sum_{k:\,k\neq i,j,\ell}\sum_{l\in\{i,j,\ell\}}R_3(i,j,\ell,k,l)\notag\\
&~~~~~~~~~~~\lesssim \frac{\phi\beta^2\log^{1/2}p}{p^{1/2}}\sum_{i,j:\,i\neq j}\mathbb{E}\bigg\{\max_{\ell:\,\ell\neq i,j}|\dot{c}_{(i,j),\ell}|\cdot\max_{\ell:\,\ell\neq i,j}\max_{k:\,k\neq i,j,\ell}|c_k^{(i,j),\ell}|\bigg\}\label{eq:boundadd1}\\
&~~~~~~~~~~~\lesssim \frac{\phi\beta^2\omega(t)\log^{1/2}p}{p^{1/2}}\mathbb{E}\bigg[\max_{i,j,\ell:\,i\neq j\neq \ell}\big\{|Y_{(i,j),\ell}|^2\vee|V_{(i,j),\ell}|^2\vee|W_{(i,j),\ell}|^2\big\}\bigg]\notag\\
&~~~~~~~~~~~\lesssim \frac{\phi\beta^2\omega(t)\log^{3/2}p}{p^{1/2}}\,.\notag
\end{align}

When $k\in\{i,j\}$ and $l\notin\{i,j,\ell\}$, we have
	\begin{align*}
		&\sum_{i,j,\ell:\,i\neq j\neq\ell}\sum_{l:\,l\neq i,j,\ell}R_3(i,j,\ell,k,l)\\
		&~~~~~~~~~\lesssim\phi\beta^2p^{1/2}\log^{1/2}p\cdot\mathbb{E}\bigg\{\max_{i,j,\ell:\,i\neq j\neq \ell}|\dot{c}_{(i,j),\ell}|\cdot\max_{i,j,\ell:\,i\neq j\neq \ell}\max_{l:\,l\neq i,j,\ell}|c^{(i,j),\ell}_l|\bigg\}\\
		&~~~~~~~~~\lesssim\frac{\phi\beta^2\omega(t)\log^{1/2}p}{p^{3/2}}\cdot\mathbb{E}\bigg[\max_{i,j,\ell:\,i\neq j\neq\ell}\big\{|Y_{(i,j),\ell}|^2\vee|V_{(i,j),\ell}|^2\vee|W_{(i,j),\ell}|^2\big\}\bigg]\\
		&~~~~~~~~~\lesssim\frac{\phi\beta^2\omega(t)\log^{3/2}p}{p^{3/2}}\,,
	\end{align*}
which implies
	\begin{align}\label{eq:bd2}
		&\sum_{i,j,\ell:\,i\neq j\neq\ell}\sum_{k\in\{i,j\}}\sum_{l:\,l\neq i,j,\ell}R_3(i,j,\ell,k,l)\lesssim \frac{\phi\beta^2\omega(t)\log^{3/2}p}{p^{3/2}}\,.
	\end{align}
	When $k\in\{i,j\}$ and $l\in\{i,j,\ell\}$, we have
	\begin{align*}
		\sum_{i,j,\ell:\,i\neq j\neq \ell}R_3(i,j,\ell,k,l)\lesssim&~\phi\beta^2\log p\cdot\mathbb{E}\bigg\{\max_{i,j,\ell:\,i\neq j\neq \ell}|\dot{c}_{(i,j),\ell}|\bigg\}\\
		\lesssim&~\frac{\phi\beta^2\omega(t)\log p}{p}\cdot\mathbb{E}\bigg[\max_{i,j,\ell:\,i\neq j\neq\ell}\big\{|Y_{(i,j),\ell}|\vee|V_{(i,j),\ell}|\vee|W_{(i,j),\ell}|\big\}\bigg]\\
		\lesssim&~\frac{\phi\beta^2\omega(t)\log^{3/2}p}{p}\,,
	\end{align*}
	which implies
 \begin{align*}
\sum_{i,j,\ell:\,i\neq j\neq \ell}\sum_{k\in\{i,j\}}\sum_{\ell\in\{i,j,\ell\}}R_3(i,j,\ell,k,l)\lesssim\frac{\phi\beta^2\omega(t)\log^{3/2}p}{p}\,.
 \end{align*}
	Together with \eqref{eq:bd1}, \eqref{eq:boundadd1} and \eqref{eq:bd2}, we have
\begin{align*}
\sum_{i,j,\ell:\,i\neq j\neq \ell}\sum_{k:\,k\neq \ell}\sum_{l=1}^pR_3(i,j,\ell,k,l)\lesssim \frac{\phi\beta^2\omega(t)\log^{3/2}p}{p^{1/2}}\,,
\end{align*}
which implies
	\begin{align}\label{eq:bd3}
		&\sum_{i,j,\ell:\,i\neq j\neq\ell}\sum_{k:\,k\neq\ell}\sum_{l=1}^p\int_0^1\int_0^1(1-\tau)\mathbb{E}\big[I(\mathcal{E}\cap\tilde{\mathcal{E}})U_{\ell kl}\{\bc^{-(i,j),\ell}+\tau\bc^{(i,j),\ell}\}\notag\\
&~~~~~~~~~~~~~~~~~~~~~~~~~~~~~~~~~~~~~~~~~~~\times|\dot{c}_{(i,j),\ell}||c^{(i,j),\ell}_k||c^{(i,j),\ell}_l|\big]\,{\rm d}\tau{\rm d}t\,{\rm d}\tau{\rm d}t\lesssim \frac{\phi^3\log^{7/2}p}{p^{1/2}}\,.
	\end{align}
	Combining \eqref{eq:bd0} and \eqref{eq:bd3}, \eqref{eq:i32} implies that
	\begin{align*}
	    \sum_{i,j,\ell:\,i\neq j\neq\ell}\sum_{k:\,k\neq\ell}\sum_{l=1}^p|{\rm I}_{3,2}(i,j,\ell,k,l)|\lesssim \frac{\phi^3\log^{7/2}p}{p^{1/2}}\,.
	\end{align*}
Analogously, we can also show
\[
\sum_{i,j,\ell:\,i\neq j\neq \ell}\sum_{l=1}^p|{\rm I}_{3,2}(i,j,\ell,\ell,l)|\lesssim\frac{\phi^3\log^{7/2}p}{p^{1/2}}\,.
\]
Together with \eqref{eq:i31}, it holds that
	\begin{align*}
	   \sum_{i,j,\ell:\,i\neq j\neq\ell}\sum_{k,l=1}^p|{\rm I}_{3}(i,j,\ell,k,l)|\lesssim \frac{\phi^3\log^{7/2}p}{{p^{1/2}}}+\phi^3p^6\log ^3p\cdot\exp\bigg(-\frac{Cp}{\phi^{2}\log^2 p}\bigg)\,.
	\end{align*}
	If we select $\phi\ll p^{1/2}\log^{-3/2}p$, we have
	\begin{align*}
	    	\sum_{i,j,\ell:\,i\neq j\neq\ell}\sum_{k,l=1}^p|{\rm I}_{3}(i,j,\ell,k,l)|\lesssim \frac{\phi^3\log^{7/2}p}{p^{1/2}}\,.
	\end{align*}
We complete the proof of \eqref{eq:sumI3}. $\hfill\Box$


\subsection{Proof of Lemma {\rm\ref{la:mu1mul2}}.}\label{sec:pfla:mu1mul2}
Notice that
	\begin{align*}
	\mathbb{E}\{\varphi_{(i,j),0}\}=\frac{\gamma}{1+\exp(\xi+\check{\theta}_i+\check{\theta}_j)}~~~\textrm{and}~~~\mathbb{E}\{\varphi_{(i,j),1}\}=\frac{\gamma}{1+\exp(-\xi-\check{\theta}_i-\check{\theta}_j)}\,.
	\end{align*}
	Recall $\xi=-\omega_1\log p+\xi^+$ and $\check{\theta}_\ell=\omega_2\log p+\check{\theta}_\ell^+$ for all $\ell\in\mathcal{S}$, where $\omega_1\in[0,2)$ and $\omega_2\in[0,1)$ such that $0\leq \omega_1-\omega_2<1$, and $|\xi^+|\vee \max_{\ell\in\mathcal{S}}|\check{\theta}_\ell^+|=o(\log p)$. Write $\chi_{p}=\exp(-|\xi^+|\vee\max_{\ell\in\mathcal{S}}|\check{\theta}_\ell^+|)$.

For $i\in\mathcal{S}$ and $j\in\mathcal{S}^c$, it holds that $
	\xi+\check{\theta}_i+\check{\theta}_j=(\omega_2-\omega_1)\log p+\xi^++\check{\theta}_i^+
	$.  If $\omega_1>\omega_2$, due to $0<\exp(\xi+\check{\theta}_i+\check{\theta}_j)\ll 1$ and $-\xi-\check{\theta}_i-\check{\theta}_j\asymp \log p$ , then $\mathbb{E}\{\varphi_{(i,j),0}\}\gtrsim \gamma$ and $\mathbb{E}\{\varphi_{(i,j),1}\}\gtrsim p^{\omega_2-\omega_1}\chi_{p}^2\gamma$.
If  $\omega_1=\omega_2$, then $ \mathbb{E}\{\varphi_{(i,j),0}\}\gtrsim \chi_{p}^2 \gamma$ and $ \mathbb{E}\{\varphi_{(i,j),1}\}\gtrsim \chi_{p}^2 \gamma$. For $i,j\in\mathcal{S}^c$, it holds that  $\xi+\check{\theta}_i+\check{\theta}_j=-\omega_1\log p+\xi^+$. If $\omega_1=0$, then $ \mathbb{E}\{\varphi_{(i,j),0}\}\gtrsim \chi_{p} \gamma$ and $ \mathbb{E}\{\varphi_{(i,j),1}\}\gtrsim \chi_{p} \gamma$. If $\omega_1>0$, due to $0<\exp(\xi+\check{\theta}_i+\check{\theta}_j)\ll 1$ and $-\xi-\check{\theta}_i-\check{\theta}_j\asymp \log p$, then $\mathbb{E}\{\varphi_{(i,j),0}\}\gtrsim \gamma$ and $\mathbb{E}\{\varphi_{(i,j),1}\}\gtrsim p^{-\omega_1}\chi_{p}\gamma$.
Therefore,  we have
\begin{align}\label{eq:ephiij02}
		&\mathbb{E}\{\varphi_{(i,j),0}\}\gtrsim\left\{\begin{aligned}
	\gamma I(\omega_1>\omega_2)+\chi_{p}^{2}\gamma I(\omega_1=\omega_2)\,,~~~~&\textrm{if}~i\in\mathcal{S},\,j\in\mathcal{S}^c\,,
\\
   \gamma I(\omega_1>0)+\chi_{p}\gamma I(\omega_1=0)\,,~~~~~~&\textrm{if}~i,j\in\mathcal{S}^c\,,
		\end{aligned}\right.\\
 & ~~~~~~~~~~~~\mathbb{E}\{\varphi_{(i,j),1}\}\gtrsim\left\{\begin{aligned}\label{eq:ephiij12}
			p^{\omega_2-\omega_1}\chi_{p}^{2}\gamma \,,~~~~~~&\textrm{if}~i\in\mathcal{S},~ j\in\mathcal{S}^c\,,\\
   p^{-\omega_1}\chi_{p}\gamma \,,~~~~~~~&\textrm{if}~i,j\in\mathcal{S}^c\,.
		\end{aligned}\right.
	\end{align}

\subsubsection{Lower bound of $\min_{\ell\in\mathcal{S}}\mu_{\ell,1}$.}\label{sucase1}

Recall that
\begin{align*}
	\mu_{\ell,1}=&\,\frac{1}{|\mathcal{H}_{\ell}|}\sum_{(i,j)\in\mathcal{H}_{\ell}}\mathbb{E}\{\varphi_{(i,\ell),1}\varphi_{(i,j),0}\varphi_{(\ell,j),1}\}\,,
\end{align*}
where $\mathcal{H}_{\ell}=\{(i,j):i,j\neq\ell~\textrm{such that}~i< j\}$. For any $\ell\in\mathcal{S}$, we have
\begin{align*}
2|\mathcal{H}_{\ell}|\mu_{\ell,1}\geq \sum_{i,j\in\mathcal{S}^c:\,i\neq j}\mathbb{E}\{\varphi_{(i,\ell),1}\varphi_{(i,j),0}\varphi_{(\ell,j),1}\}\,.
\end{align*}
By \eqref{eq:ephiij02} and \eqref{eq:ephiij12}, due to $s\ll p$, it holds that
\begin{align*}
\sum_{i,j\in\mathcal{S}^c:\,i\neq j}\mathbb{E}\{\varphi_{(i,\ell),1}\varphi_{(i,j),0}\varphi_{(\ell,j),1}\}\gtrsim p^{2\omega_2-2\omega_1+2}\chi_{p}^5\gamma^3\,.
\end{align*}
Due to $\mathcal{H}_{\ell}\asymp p^2$,  then
\begin{align*}
\min_{\ell \in\mathcal{S}}\mu_{\ell,1}\gtrsim p^{2\omega_2-2\omega_1}\chi_{p}^5\gamma^3\,.
\end{align*}
We obtain the lower bound of $\min_{\ell\in\mathcal{S}}\mu_{\ell,1}$. $\hfill\Box$

\subsubsection{Lower bound of $\min_{\ell\in\mathcal{S}^c}\mu_{\ell,1}$.}\label{sucase2}
For any $\ell\in\mathcal{S}^c$, we have
\begin{align*}
2|\mathcal{H}_{\ell}|\mu_{\ell,1}
\geq \sum_{i,j\in\mathcal{S}^c:\,i\neq \ell,\,j \neq i,\ell}\mathbb{E}\{\varphi_{(i,\ell),1}\varphi_{(i,j),0}\varphi_{(\ell,j),1}\}\,.
\end{align*}
By \eqref{eq:ephiij02} and \eqref{eq:ephiij12}, due to $s\ll p$, it holds that
\begin{align*}
\sum_{i,j\in\mathcal{S}^c:\,i\neq\ell,\,j\neq \ell,i}\mathbb{E}\{\varphi_{(i,\ell),1}\varphi_{(i,j),0}\varphi_{(\ell,j),1}\}\gtrsim p^{-2\omega_1+2}\chi_{p}^3\gamma^3\,.
\end{align*}
Due to $\mathcal{H}_{\ell}\asymp p^2$, then
\begin{align*}
\min_{\ell \in\mathcal{S}^c}\mu_{\ell,1}\gtrsim p^{-2\omega_1}\chi_{p}^3\gamma^3\,.
\end{align*}
We obtain the lower bound of $\min_{\ell\in\mathcal{S}^c}\mu_{\ell,1}$. $\hfill\Box$

\subsubsection{Lower bound of $\min_{\ell\in\mathcal{S}}\mu_{\ell,2}$.}\label{sucase12}
Recall that
\begin{align*}
	\mu_{\ell,2}=&\,\frac{1}{|\mathcal{H}_{\ell}|}\sum_{(i,j)\in\mathcal{H}_{\ell}}\mathbb{E}\{\varphi_{(i,\ell),0}\varphi_{(i,j),1}\varphi_{(\ell,j),0}\}\,,
\end{align*}
where $\mathcal{H}_{\ell}=\{(i,j):i,j\neq\ell~\textrm{such that}~i< j\}$. For any $\ell\in\mathcal{S}$, we have
\begin{align*}
2|\mathcal{H}_{\ell}|\mu_{\ell,2}
\geq \sum_{i,j\in\mathcal{S}^c:\,i\neq j}\mathbb{E}\{\varphi_{(i,\ell),0}\varphi_{(i,j),1}\varphi_{(\ell,j),0}\}\,.
\end{align*}
By \eqref{eq:ephiij02} and \eqref{eq:ephiij12}, due to $s\ll p$, it holds that
\begin{align*}
\sum_{i,j\in\mathcal{S}^c:\,i\neq j}\mathbb{E}\{\varphi_{(i,\ell),0}\varphi_{(i,j),1}\varphi_{(\ell,j),0}\}\gtrsim p^{-\omega_1+2}\chi_{p}^5\gamma^3\,.
\end{align*}
Due to $\mathcal{H}_{\ell}\asymp p^2$, then
\begin{align*}
\min_{\ell \in\mathcal{S}}\mu_{\ell,2}\gtrsim p^{-\omega_1}\chi_{p}^5\gamma^3\,.
\end{align*}
We obtain the lower bound of $\min_{\ell\in\mathcal{S}}\mu_{\ell,2}$. $\hfill\Box$

\subsubsection{Lower bound of $\min_{\ell\in\mathcal{S}^c}\mu_{\ell,2}$.}\label{sucase112}

For any $\ell\in\mathcal{S}^c$, we have
\begin{align*}
2|\mathcal{H}_{\ell}|\mu_{\ell,2}
	\geq \sum_{i,j\in\mathcal{S}^c:\,i\neq \ell,\,j \neq i,\ell}\mathbb{E}\{\varphi_{(i,\ell),0}\varphi_{(i,j),1}\varphi_{(\ell,j),0}\}\,.
\end{align*}
By \eqref{eq:ephiij02} and \eqref{eq:ephiij12}, due to $s\ll p$, it holds that
\begin{align*}
\sum_{i,j\in\mathcal{S}^c:\,i\neq\ell,\,j\neq \ell,i}\mathbb{E}\{\varphi_{(i,\ell),0}\varphi_{(i,j),1}\varphi_{(\ell,j),0}\}\gtrsim p^{-\omega_1+2}\chi_{p}^3\gamma^3\,.
\end{align*}
Due to $\mathcal{H}_{\ell}\asymp p^2$, then
\begin{align*}
\min_{\ell \in\mathcal{S}^c}\mu_{\ell,2}\gtrsim p^{-\omega_1}\chi_{p}^3\gamma^3\,.
\end{align*}
We obtain the lower bound of $\min_{\ell\in\mathcal{S}^c}\mu_{\ell,2}$.
$\hfill\Box$

\subsection{Proof of Lemma {\rm\ref{la:9}}.}\label{sec:pfl7}

 For any $i,j,\ell\in[p]$, let $\mathring{\psi}_1(i,j;\ell)=\psi_1(i,j;\ell)-\mathbb{E}\{\psi_1(i,j;\ell)\}$ and  $\mathring{\psi}_2(i,j;\ell)=\psi_2(i,j;\ell)-\mathbb{E}\{\psi_2(i,j;\ell)\}$, where $\psi_1(i,j;\ell)=\varphi_{(i,\ell),1}\varphi_{(i,j),0}\varphi_{(\ell,j),1}$ and $\psi_2(i,j;\ell)=\varphi_{(i,\ell),0}\varphi_{(i,j),1}\varphi_{(\ell,j),0}$. Write  $\mathscr{F}_{\ell}=\{Z_{i,\ell},Z_{\ell,j}:(i,j)\in\mathcal{H}_{\ell}\}$. As we have shown in Section \ref{sec:pflem1},
	\begin{align*}
	\hat{\mu}_{\ell,1}-\mu_{\ell,1}
 =&~\underbrace{\frac{1}{|\mathcal{H}_\ell|}\sum_{(i,j)\in\mathcal{H}_\ell}[\psi_1(i,j;\ell)-\mathbb{E}\{\psi_1(i,j;\ell)\,|\,\mathscr{F}_\ell\}]}_{I_{\ell,1,1}}\\
		&+\underbrace{\frac{1}{|\mathcal{H}_\ell|}\sum_{(i,j)\in\mathcal{H}_\ell}[\mathbb{E}\{\psi_1(i,j;\ell)\,|\,\mathscr{F}_\ell\}-\mathbb{E}\{\psi_1(i,j;\ell)\}]}_{I_{\ell,1,2}}\,,\\
 \hat{\mu}_{\ell,2}-\mu_{\ell,2}
 =&~\underbrace{\frac{1}{|\mathcal{H}_\ell|}\sum_{(i,j)\in\mathcal{H}_\ell}[\psi_2(i,j;\ell)-\mathbb{E}\{\psi_2(i,j;\ell)\,|\,\mathscr{F}_\ell\}]}_{I_{\ell,2,1}}\\
		&+\underbrace{\frac{1}{|\mathcal{H}_\ell|}\sum_{(i,j)\in\mathcal{H}_\ell}[\mathbb{E}\{\psi_2(i,j;\ell)\,|\,\mathscr{F}_\ell\}-\mathbb{E}\{\psi_2(i,j;\ell)\}]}_{I_{\ell,2,2}}\,.
	\end{align*}
Identical to the arguments stated in Section \ref{sec:pflem11},  we have
	\begin{align}\label{eq:Il11}
		\max_{k\in\{1,2\}}\max_{\ell\in\mathcal{S}}|I_{\ell,k,1}|=O_{\p}\bigg(\frac{\log^{1/2}s}{p}\bigg)~~\textrm{and}~~\max_{k\in\{1,2\}}\max_{\ell\in\mathcal{S}^c}|I_{\ell,k,1}|=O_{\p}\bigg(\frac{\log^{1/2}p}{p}\bigg)\,.
	\end{align}
Define $\mathring{\varphi}_{(i,j),\tau}=\varphi_{(i,j),\tau}-\mathbb{E}\{\varphi_{(i,j),\tau}\}$. Same as \eqref{eq:Iell12},  it holds that
	\begin{align}
		(p-1)(p-2)I_{\ell,1,2}
		=&~\underbrace{2\sum_{i,j:\,i\neq j,\,i,j\neq \ell}\mathring{\varphi}_{(i,\ell),1}\mathbb{E}\{\varphi_{(\ell,j),1}\}\mathbb{E}\{\varphi_{(i,j),0}\}}_{I_{\ell,1,2}(1)}\label{eq:Iell12S}\\
		&+\underbrace{\sum_{i,j:\,i\neq j,\,i,j\neq \ell}\mathring{\varphi}_{(i,\ell),1}\mathring{\varphi}_{(\ell,j),1}\mathbb{E}\{\varphi_{(i,j),0}\}}_{I_{\ell,1,2}(2)}\,,\notag\\
		(p-1)(p-2)I_{\ell,2,2}
		=&~\underbrace{2\sum_{i,j:\,i\neq j,\,i,j\neq \ell}\mathring{\varphi}_{(i,\ell),0}\mathbb{E}\{\varphi_{(\ell,j),0}\}\mathbb{E}\{\varphi_{(i,j),1}\}}_{I_{\ell,2,2}(1)}\label{eq:Iell22Sc}\\
		&+\underbrace{\sum_{i,j:\,i\neq j,\,i,j\neq \ell}\mathring{\varphi}_{(i,\ell),0}\mathring{\varphi}_{(\ell,j),0}\mathbb{E}\{\varphi_{(i,j),1}\}}_{I_{\ell,2,2}(2)}\,.\notag
	\end{align}
In the sequel, we need the following lemmas whose proofs are given in Sections \ref{sec:pfla:bdvarphi} and \ref{sec:pfla:Ail}, respectively.

\begin{lemma}\label{la:bdvarphi}
Let $(\alpha,\beta)\in\mathcal{M}(\gamma, C_1)$ for some fixed constant $C_1\in(0,0.5)$.
	It holds that
	\begin{align*}
		\mathbb{E}\{\varphi_{(i,j),0}\}\lesssim&\left\{\begin{aligned}
			\gamma I(\omega_1\geq2\omega_2)+p^{\omega_1-2\omega_2}\chi_{p}^{-3}\gamma I(\omega_1<2\omega_2)\,,~~&\textrm{if}~i,j\in\mathcal{S}\,,  \\
			\gamma\,,~~~~~~~~~~~~~~~~~~~~~~~~~~~~~~~~&\textrm{if}~i\in\mathcal{S}^c~\textrm{or}~ j\in\mathcal{S}^c\,,
		\end{aligned}\right.\\
  \mathbb{E}\{\varphi_{(i,j),1}\}\lesssim&\left\{\begin{aligned}
			p^{2\omega_2-\omega_1}\chi_{p}^{-3}\gamma  I(\omega_1>2\omega_2)+\gamma I(\omega_1\leq2\omega_2)\,,~~&\textrm{if}~i,j\in\mathcal{S}\,,  \\
			p^{\omega_2-\omega_1}\chi_{p}^{-2}\gamma I(\omega_1>\omega_2)+\gamma I(\omega_1=\omega_2)\,,~~~~~&\textrm{if}~i\in\mathcal{S}, j\in\mathcal{S}^c\,,\\
   p^{-\omega_1}\chi_{p}^{-1}\gamma I(\omega_1>0)+\gamma I(\omega_1=0)\,,~~~~~~~~~&\textrm{if}~i,j\in\mathcal{S}^c\,.
		\end{aligned}\right.
	\end{align*}
	\end{lemma}
\begin{lemma}\label{la:Ail}  Let $(\alpha,\beta)\in\mathcal{M}(\gamma, C_1)$ for some fixed constant $C_1\in(0,0.5)$. It holds that
	\begin{align*}
 \max_{i,\ell\in\mathcal{S}:\,i\neq \ell}A_{i,\ell}\lesssim&~ \{sp^{\min(2\omega_2-\omega_1,\,\omega_1-2\omega_2)}\chi_{p}^{-1}+p^{1+\omega_2-\omega_1}\}\chi_{p}^{-2}\gamma^2\,,\\
\max_{i\in\mathcal{S},\,\ell \in\mathcal{S}^c}A_{i,\ell}\lesssim&~\{sp^{\min(-\omega_2,\,\omega_2-\omega_1)} \chi_{p}^{-4}+p^{1-\omega_1}\}\chi_{p}^{-1}\gamma^2\,,  \\
 \max_{i\in\mathcal{S}^c,\,\ell \in\mathcal{S}}A_{i,\ell}\lesssim&~\{sp^{\min(2\omega_2-\omega_1,\,0)} \chi_{p}^{-1}+p^{1+\omega_2-\omega_1}\}\chi_{p}^{-2}\gamma^2\,,  \\
 \max_{i, \ell\in\mathcal{S}^c:\, i\neq \ell}A_{i,\ell}\lesssim&~ (sp^{\omega_2-\omega_1} \chi_{p}^{-1}+p^{1-\omega_1})\chi_{p}^{-1}\gamma^2\,.
\end{align*}
\end{lemma}

\subsubsection{Convergence rate of $\max_{\ell\in\mathcal{S}}|\hat{\mu}_{\ell,1}-\mu_{\ell,1}|$.}\label{mu1:S}
 Notice that $\hat{\mu}_{\ell,1}-\mu_{\ell,1}=I_{\ell, 1,1}+I_{\ell, 1,2}$.
Given the convergence rate of $\max_{\ell \in \mathcal{S}}|I_{\ell,1,1}|$ in \eqref{eq:Il11}, in order to establish the convergence rate of $\max_{\ell\in\mathcal{S}}|\hat{\mu}_{\ell,1}-\mu_{\ell,1}|$, we only need to derive the convergence rate of $\max_{\ell \in \mathcal{S}}|I_{\ell,1,2}|$.
For $I_{\ell,1,2}(1)$, we have
	\begin{align*}
		I_{\ell,1,2}(1)&=\sum_{i:\,i\neq \ell}\mathring{\varphi}_{(i,\ell),1}\bigg[2\sum_{j:\,j\neq i,\ell}\mathbb{E}\{\varphi_{(\ell,j),1}\}\mathbb{E}\{\varphi_{(i,j),0}\}\bigg]=\sum_{i:\,i\neq \ell}\mathring{\varphi}_{(i,\ell),1}A_{i,\ell}
		\\&=\sum_{i\in\mathcal{S}:\,i\neq \ell}\mathring{\varphi}_{(i,\ell),1}A_{i,\ell}+\sum_{i\in\mathcal{S}^c}\mathring{\varphi}_{(i,\ell),1}A_{i,\ell}\,.
	\end{align*}
For any $i,\ell\in\mathcal{S}$, write
\begin{align*}
\tilde{A}_{i,\ell,\mathcal{S},\mathcal{S}}=&~\frac{A_{i,\ell}}{\max_{i,\ell\in\mathcal{S}:\,i\neq \ell}A_{i,\ell}}\,.
\end{align*}
It then holds that
\[
\sum_{i\in\mathcal{S}:\,i\neq \ell}\mathring{\varphi}_{(i,\ell),1}A_{i,\ell}=\bigg(\max_{i,\ell\in\mathcal{S}:i\neq\ell}A_{i,\ell}\bigg)\sum_{i\in\mathcal{S}:\,i\neq \ell}\mathring{\varphi}_{(i,\ell),1}\tilde{A}_{i,\ell,\mathcal{S},\mathcal{S}}\,.
\]
Due to $A_{i,\ell}\geq 0$ for any $i,\ell$, we have $\max_{i,\ell\in\mathcal{S}}\tilde{A}_{i,\ell,\mathcal{S},\mathcal{S}}\leq 1$. Since $\{\mathring{\varphi}_{(i,\ell),1}\}_{i\in\mathcal{S}:\,i\neq\ell}$ is an independent sequence, by Bernstein inequality,
\[
\mathbb{P}\bigg\{\bigg|\sum_{i\in\mathcal{S}:\,i\neq \ell}\mathring{\varphi}_{(i,\ell),1}\tilde{A}_{i,\ell,\mathcal{S},\mathcal{S}}\bigg|>u\bigg\}\lesssim \exp(-Cs^{-1}u^2)
\]
for any $0<u\leq O(s)$, which implies
\[
\max_{\ell\in\mathcal{S}}\bigg|\sum_{i\in\mathcal{S}:\,i\neq \ell}\mathring{\varphi}_{(i,\ell),1}\tilde{A}_{i,\ell,\mathcal{S},\mathcal{S}}\bigg|=O_{\p}(s^{1/2}\log^{1/2}s)\,.
\]
Together with Lemma \ref{la:Ail},  we have
\begin{align*}
&\max_{\ell\in\mathcal{S}}\bigg|\sum_{i\in\mathcal{S}:\,i\neq \ell}\mathring{\varphi}_{(i,\ell),1}A_{i,\ell}\bigg|\\
&~~~~~~~=O_{\p}\big[\{sp^{\min(2\omega_2-\omega_1,\,\omega_1-2\omega_2)}\chi_{p}^{-1}+p^{1+\omega_2-\omega_1}\}\chi_{p}^{-2}\gamma^2s^{1/2}\log^{1/2}s\big]\,.
\end{align*}
Analogously, we also have 	
\begin{align*}
&\max_{\ell\in\mathcal{S}}\bigg|\sum_{i\in\mathcal{S}^c}\mathring{\varphi}_{(i,\ell),1}A_{i,\ell}\bigg|\\
&~~~~~~~=O_{\p}\big[\{sp^{\min(2\omega_2-\omega_1,\,0)} \chi_{p}^{-1}+p^{1+\omega_2-\omega_1}\}\chi_{p}^{-2}\gamma^2p^{1/2}\log^{1/2}s\big]\,.
\end{align*}
 Due to $p^{\min(2\omega_2-\omega_1,\,\omega_1-2\omega_2)}\leq p^{\min(2\omega_2-\omega_1,\,0)} $ and $s\ll p$, then
 \begin{align*}
 \max_{\ell\in\mathcal{S}}|I_{\ell,1,2}(1)|=O_{\p}\big[\{sp^{\min(2\omega_2-\omega_1,\,0)} \chi_{p}^{-1}+p^{1+\omega_2-\omega_1}\}\chi_{p}^{-2}\gamma^2p^{1/2}\log^{1/2}s\big]\,.
 \end{align*}
	It follows from Lemma \ref{la:bdvarphi} that $\max_{i,j:\,i\neq j,\,i,j\neq \ell}\mathbb{E}\{\varphi_{(i,j),0}\}\lesssim \gamma$. Using the same arguments in Section \ref{sec:pflem12} for deriving the convergence rate of $\max_{\ell\in[p]}|I_{\ell,1,2}(2)|$ there, we have
 $\max_{\ell\in\mathcal{S}}|I_{\ell,1,2}(2)|=O_{\p}(\gamma p\log s)
$.
	Therefore, by \eqref{eq:Iell12S},
	we have
 \begin{align*}
     \max_{\ell \in \mathcal{S}}|I_{\ell,1,2}|=&~O_{\p}\big[\gamma^2p^{-3/2}\{sp^{\min(2\omega_2-\omega_1,\,0)} \chi_{p}^{-1}+p^{1+\omega_2-\omega_1}\}\chi_{p}^{-2}\log^{1/2}s\big]\\&+O_{\p}(\gamma p^{-1}\log s)\,.
 \end{align*}
 Together with \eqref{eq:Il11}, it holds that
	\begin{align}\label{eq:hatmul1S1}
	\max_{\ell\in\mathcal{S}}|\hat{\mu}_{\ell,1}-\mu_{\ell,1}|=&~O_{\p}\big[\gamma^2p^{-3/2}\{sp^{\min(2\omega_2-\omega_1,\,0)} \chi_{p}^{-1}+p^{1+\omega_2-\omega_1}\}\chi_{p}^{-2}\log^{1/2}s\big]\\&+O_{\p}(\gamma p^{-1}\log s)+O_{\p}(p^{-1}\log^{1/2}s)\notag\,.
	\end{align}
We obtain the convergence rate of $\max_{\ell\in\mathcal{S}}|\hat{\mu}_{\ell,1}-\mu_{\ell,1}|$. $\hfill\Box$

\subsubsection{Convergence rate of $\max_{\ell\in\mathcal{S}}|\hat{\mu}_{\ell,2}-\mu_{\ell,2}|$.}\label{mu2:S}	Notice that $\hat{\mu}_{\ell,2}-\mu_{\ell,2}=I_{\ell, 2,1}+I_{\ell, 2,2}$.
Given the convergence rate of
 $\max_{\ell \in \mathcal{S}}|I_{\ell,2,1}|$  in \eqref{eq:Il11}, in order to establish the convergence rate of $\max_{\ell\in\mathcal{S}}|\hat{\mu}_{\ell,2}-\mu_{\ell,2}|$, we only need to derive the convergence rate of $\max_{\ell \in \mathcal{S}}|I_{\ell,2,2}|$.
For $I_{\ell,2,2}(1)$, we have
	\begin{align*}
		I_{\ell,2,2}(1)
		=\sum_{i\in\mathcal{S}:\,i\neq \ell}\mathring{\varphi}_{(i,\ell),0}A_{\ell, i}+\sum_{i\in\mathcal{S}^c}\mathring{\varphi}_{(i,\ell),0}A_{\ell, i}\,.
	\end{align*}
Following the same arguments in Section \ref{mu1:S} for obtaining the convergence rate of $\max_{\ell\in\mathcal{S}}|I_{\ell,1,2}(1)|$ there, it holds that \begin{align}\label{eq:Il221S}
\max_{\ell\in\mathcal{S}}|I_{\ell,2,2}(1)|=&~O_{\p}\big[\gamma^2s^{1/2}\{sp^{\min(2\omega_2-\omega_1,\,\omega_1-2\omega_2)}\chi_{p}^{-1}+p^{1+\omega_2-\omega_1}\}\chi_{p}^{-2}\log^{1/2}s\big]\notag\\&+O_{\p}\big[\gamma^2p^{1/2}\{sp^{\min(-\omega_2,\,\omega_2-\omega_1)} \chi_{p}^{-4}+p^{1-\omega_1}\}\chi_{p}^{-1}\log^{1/2}s\big]\,.
 \end{align}
 For $I_{\ell,2,2}(2)$, it holds that
	\begin{align*}
		|I_{\ell,2,2}(2)|\leq&~\underbrace{\bigg|\sum_{i,\,j\in\mathcal{S}:\,i\neq j,\,i,j\neq \ell}\mathring{\varphi}_{(i,\ell),0}\mathring{\varphi}_{(\ell,j),0}\mathbb{E}\{\varphi_{(i,j),1}\}\bigg|}_{D_{\ell,1}}\\
  &+\underbrace{\bigg|\sum_{i\in\mathcal{S},\,j\in\mathcal{S}^c:\,i\neq  \ell}\mathring{\varphi}_{(i,\ell),0}\mathring{\varphi}_{(\ell,j),0}\mathbb{E}\{\varphi_{(i,j),1}\}\bigg|}_{D_{\ell,2}}
		\\&+\underbrace{\bigg|\sum_{i\in\mathcal{S}^c,\,j\in\mathcal{S}:\,j\neq \ell}\mathring{\varphi}_{(i,\ell),0}\mathring{\varphi}_{(\ell,j),0}\mathbb{E}\{\varphi_{(i,j),1}\}\bigg|}_{D_{\ell,3}}\\
  &+\underbrace{\bigg|\sum_{i,j\in\mathcal{S}^c:\,i\neq j }\mathring{\varphi}_{(i,\ell),0}\mathring{\varphi}_{(\ell,j),0}\mathbb{E}\{\varphi_{(i,j),1}\}\bigg|}_{D_{\ell,4}}\,.
	\end{align*}
	By  Lemma \ref{la:bdvarphi}, $\max_{i,j\in\mathcal{S}:\,i\neq j,\,i,j\neq \ell}\mathbb{E}\{\varphi_{(i,j),1}\}\lesssim p^{2\omega_2-\omega_1}\chi_{p}^{-3}\gamma  I(\omega_1>2\omega_2)+\gamma I(\omega_1\leq2\omega_2)$ and $\max_{i,j\in\mathcal{S}^c:\,i\neq j}\mathbb{E}\{\varphi_{(i,j),1}\}\lesssim p^{-\omega_1}\chi_p^{-1}I(\omega_1>0)+\gamma I(\omega_1=0)$.
 Applying the same arguments in Section \ref{sec:pflem12} for deriving the convergence rate of $\max_{\ell\in[p]}|I_{\ell,1,2}(2)|$ there, we have
 \begin{align*}
	\max_{\ell\in\mathcal{S}}D_{\ell,1}=&~O_{\p}\big[ \gamma s\{p^{2\omega_2-\omega_1}\chi_{p}^{-3} I(\omega_1>2\omega_2)+ I(\omega_1\leq2\omega_2)\}\log s\big]\,,\\
\max_{\ell\in\mathcal{S}}D_{\ell,4}=&~O_{\p}\big[\gamma p\{p^{-\omega_1}\chi_{p}^{-1}I(\omega_1>0)+I(\omega_1=0)\}\log s\big]\,.
	\end{align*}

Given $\ell\in\mathcal{S}$, define $d_{i,j}^\ell=\mathbb{E}\{\varphi_{(i,j),1}\}$ for any $(i,j)\in\mathcal{S}\times\mathcal{S}^c$ with $i\neq \ell$, and $d_{i,j}^\ell=0$ otherwise. Then
\begin{align}\label{eq:Dl2}
D_{\ell,2}=\bigg|\sum_{i,j=1}^p\mathring{\varphi}_{(i,\ell),0}\mathring{\varphi}_{(\ell,j),0}d_{i,j}^\ell\bigg|=\bigg(\max_{i,j\in[p]}d_{i,j}^\ell\bigg)\bigg|\sum_{i,j=1}^p\mathring{\varphi}_{(i,\ell),0}\mathring{\varphi}_{(\ell,j),0}\tilde{d}_{i,j}^\ell\bigg|
\end{align}
with $\tilde{d}_{i,j}^\ell=d_{i,j}^\ell/\max_{i,j\in[p]}d_{i,j}^\ell\in[0,1]$.
By the decoupling inequalities of \cite{DM_1995} and Theorem 3.3 of \cite{Gineetal_2000}, it holds that
	\begin{align*}
	\max_{\ell\in\mathcal{S}}\mathbb{P}\bigg(\bigg|\sum_{i,j=1}^p\mathring{\varphi}_{(i,\ell),0}\mathring{\varphi}_{(\ell,j),0}\tilde{d}_{i,j}^\ell\bigg|>u\bigg)\lesssim&~\exp(-Cu^{1/2})+\exp(-Cp^{-1/3}u^{2/3})\notag
		\\&~~+\exp(-Cs^{-1/2}p^{-1/2}u)
  +\exp(-Cs^{-1}p^{-1}u^2)
	\end{align*}
	for any $u>0$, which implies
 \begin{align*}
\max_{\ell \in\mathcal{S}}\bigg|\sum_{i,j=1}^p\mathring{\varphi}_{(i,\ell),0}\mathring{\varphi}_{(\ell,j),0}\tilde{d}_{i,j}^\ell\bigg|=O_{\p}(s^{1/2}p^{1/2}\log s)\,.
 \end{align*}
It follows from Lemma \ref{la:bdvarphi} that $0<\max_{i,j\in[p]}d_{i,j}^\ell\lesssim p^{\omega_2-\omega_1}\chi_{p}^{-2}\gamma I(\omega_1>\omega_2)+\gamma I(\omega_1=\omega_2)$. By \eqref{eq:Dl2}, we have
\begin{align*}  \max_{\ell\in\mathcal{S}}D_{\ell,2}=O_{\p}\big[ \gamma s^{1/2} p^{1/2}\{p^{\omega_2-\omega_1}\chi_{p}^{-2}I(\omega_1>\omega_2)+ I(\omega_1=\omega_2)\}\log s\big]\,.
\end{align*}
Analogously, we can also show
\begin{align*}  \max_{\ell\in\mathcal{S}}D_{\ell,3}=O_{\p}\big[ \gamma s^{1/2} p^{1/2}\{p^{\omega_2-\omega_1}\chi_{p}^{-2}I(\omega_1>\omega_2)+ I(\omega_1=\omega_2)\}\log s\big]\,.
\end{align*}
Therefore,
 \begin{align*}\label{eq:Il222S}
     \max_{\ell \in \mathcal{S}}|I_{\ell,2,2}(2)|=&~O_{\p}\big[\gamma s\{p^{2\omega_2-\omega_1}\chi_{p}^{-3}  I(\omega_1>2\omega_2)+ I(\omega_1\leq2\omega_2)\}\log s\big]\notag\\
     &+O_{\p}\big[\gamma s^{1/2} p^{1/2}\{p^{\omega_2-\omega_1}\chi_{p}^{-2} I(\omega_1>\omega_2)+I(\omega_1=\omega_2)\}\log s\big]\\
     &+O_{\p}\big[ \gamma p\{p^{-\omega_1}\chi_{p}^{-1} I(\omega_1>0)+I(\omega_1=0)\}\log s\big]\,.\notag
 \end{align*}
By \eqref{eq:Iell22Sc} and \eqref{eq:Il221S},
 we have
	\begin{align*}
\max_{\ell\in\mathcal{S}}|I_{\ell,2,2}|=&~O_{\p}\big[\gamma^2s^{1/2}p^{-2}\{sp^{\min(2\omega_2-\omega_1,\,\omega_1-2\omega_2)}\chi_{p}^{-1}+p^{1+\omega_2-\omega_1}\}\chi_{p}^{-2}\log^{1/2}s\big]\notag\\&+O_{\p}\big[\gamma^2p^{-3/2}\{sp^{\min(-\omega_2,\,\omega_2-\omega_1)} \chi_{p}^{-4}+p^{1-\omega_1}\}\chi_{p}^{-1}\log^{1/2}s\big]\\&+O_{\p}\big[\gamma s p^{-2}\{p^{2\omega_2-\omega_1}\chi_{p}^{-3} I(\omega_1>2\omega_2)+I(\omega_1\leq2\omega_2)\}\log s\big]\notag\\&+O_{\p}\big[\gamma s^{1/2} p^{-3/2}\{p^{\omega_2-\omega_1}\chi_{p}^{-2} I(\omega_1>\omega_2)+I(\omega_1=\omega_2)\}\log s\big]\\&+O_{\p}\big[\gamma p^{-1}\{p^{-\omega_1}\chi_{p}^{-1}I(\omega_1>0)+ I(\omega_1=0)\}\log s\big]\,.
 \end{align*}
Together with \eqref{eq:Il11}, it holds that
\begin{align*}
\max_{\ell\in\mathcal{S}}|\hat{\mu}_{\ell,2}-\mu_{\ell,2}|=&~O_{\p}\big[\gamma^2s^{1/2}p^{-2}\{sp^{\min(2\omega_2-\omega_1,\,\omega_1-2\omega_2)}\chi_{p}^{-1}+p^{1+\omega_2-\omega_1}\}\chi_{p}^{-2}\log^{1/2}s\big]\\
&+O_{\p}\big[\gamma^2p^{-3/2}\{sp^{\min(-\omega_2,\,\omega_2-\omega_1)} \chi_{p}^{-4}+p^{1-\omega_1}\}\chi_{p}^{-1}\log^{1/2}s\big]\\
&+O_{\p}\big[\gamma s p^{-2}\{p^{2\omega_2-\omega_1}\chi_{p}^{-3}  I(\omega_1>2\omega_2)+I(\omega_1\leq2\omega_2)\}\log s\big]\\
&+O_{\p}\big[\gamma s^{1/2} p^{-3/2}\{p^{\omega_2-\omega_1}\chi_{p}^{-2}I(\omega_1>\omega_2)+I(\omega_1=\omega_2)\}\log s\big]\\
&+O_{\p}\big[\gamma p^{-1}\{p^{-\omega_1}\chi_{p}^{-1}I(\omega_1>0)+I(\omega_1=0)\}\log s\big]\\
&+O_{\p}(p^{-1}\log^{1/2}s)\,.
 \end{align*}
We obtain the convergence rate of $\max_{\ell\in\mathcal{S}}|\hat{\mu}_{\ell,2}-\mu_{\ell,2}|$.
$\hfill\Box$

\subsubsection{ Convergence rate of $\max_{\ell\in\mathcal{S}^c}|\hat{\mu}_{\ell,1}-\mu_{\ell,1}|$.}\label{mu1:Sc}
Notice that $\hat{\mu}_{\ell,1}-\mu_{\ell,1}=I_{\ell, 1,1}+I_{\ell, 1,2}$.
Given the convergence rate of $\max_{\ell \in \mathcal{S}^c}|I_{\ell,1,1}|$ in \eqref{eq:Il11}, in order to establish the convergence rate of $\max_{\ell\in\mathcal{S}^c}|\hat{\mu}_{\ell,1}-\mu_{\ell,1}|$, we only need to derive the convergence rate of $\max_{\ell \in \mathcal{S}^c}|I_{\ell,1,2}|$.
%
For $I_{\ell,1,2}(1)$, we have
	\begin{align*}
		I_{\ell,1,2}(1)
		=\sum_{i\in\mathcal{S}}\mathring{\varphi}_{(i,\ell),1}A_{i,\ell}+\sum_{i\in\mathcal{S}^c:\,i\neq \ell}\mathring{\varphi}_{(i,\ell),1}A_{i,\ell}\,.
	\end{align*}
Following the same arguments in Section \ref{mu1:S} for obtaining the convergence rate of $\max_{\ell\in\mathcal{S}}|I_{\ell,1,2}(1)|$ there, it holds that
\begin{align*}
\max_{\ell\in\mathcal{S}^c}|I_{\ell,1,2}(1)|=&~O_{\p}\big[\gamma^2s^{1/2}\{sp^{\min(-\omega_2,\,\omega_2-\omega_1)} \chi_{p}^{-4}+p^{1-\omega_1}\}\chi_{p}^{-1}\log^{1/2}p\big]\\
&+O_{\p}\big\{\gamma^2p^{1/2}(sp^{\omega_2-\omega_1} \chi_{p}^{-1}+p^{1-\omega_1})\chi_{p}^{-1}\log^{1/2}p\big\}\,.
\end{align*}
	It follows from Lemma \ref{la:bdvarphi} that $\max_{i,j:\,i\neq j,\,i,j\neq \ell}\mathbb{E}\{\varphi_{(i,j),0}\}\lesssim \gamma$. Applying the same arguments in Section \ref{sec:pflem12} for deriving the convergence rate of $\max_{\ell\in[p]}|I_{\ell,1,2}(2)|$ there, we have $\max_{\ell\in\mathcal{S}^c}|I_{\ell,1,2}(2)|=O_{\p}(\gamma p\log p)$.  Therefore, by \eqref{eq:Iell12S},
	we have \begin{align*}
	    \max_{\ell \in \mathcal{S}^c}|I_{\ell,1,2}|=&~O_{\p}\big[\gamma^2s^{1/2}p^{-2}\{sp^{\min(-\omega_2,\,\omega_2-\omega_1)} \chi_{p}^{-4}+p^{1-\omega_1}\}\chi_{p}^{-1}\log^{1/2}p\big]\\
     &+O_{\p}\big\{\gamma^2p^{-3/2}(sp^{\omega_2-\omega_1} \chi_{p}^{-1}+p^{1-\omega_1})\chi_{p}^{-1}\log^{1/2}p\big\}+O_{\p}(\gamma p^{-1}\log p)\,.
	\end{align*}
	Together with \eqref{eq:Il11}, we have
\begin{align*}
	\max_{\ell\in\mathcal{S}^c}|\hat{\mu}_{\ell,1}-\mu_{\ell,1}|=&~O_{\p}\big[\gamma^2s^{1/2}p^{-2}\{sp^{\min(-\omega_2,\,\omega_2-\omega_1)} \chi_{p}^{-4}+p^{1-\omega_1}\}\chi_{p}^{-1}\log^{1/2}p\big]\\
 &+O_{\p}\big\{\gamma^2p^{-3/2}(sp^{\omega_2-\omega_1} \chi_{p}^{-1}+p^{1-\omega_1})\chi_{p}^{-1}\log^{1/2}p\big\}\\
 &+O_{\p}(\gamma p^{-1}\log p)+O_{\p}(p^{-1}\log^{1/2}p)\,.
\end{align*}
We obtain the convergence rate of $\max_{\ell\in\mathcal{S}^c}|\hat{\mu}_{\ell,1}-\mu_{\ell,1}|$. $\hfill\Box$

\subsubsection{Convergence rate of $\max_{\ell\in\mathcal{S}^c}|\hat{\mu}_{\ell,2}-\mu_{\ell,2}|$.}\label{mu2:Sc}
Notice that $\hat{\mu}_{\ell,2}-\mu_{\ell,2}=I_{\ell, 2,1}+I_{\ell, 2,2}$.
Given the convergence rate of
 $\max_{\ell \in \mathcal{S}^c}|I_{\ell,2,1}|$  in \eqref{eq:Il11}, in order to establish the convergence rate of $\max_{\ell\in\mathcal{S}^c}|\hat{\mu}_{\ell,2}-\mu_{\ell,2}|$, we only need to derive the convergence rate of $\max_{\ell \in \mathcal{S}^c}|I_{\ell,2,2}|$.
For $I_{\ell,2,2}(1)$, we have
	\begin{align*}
		I_{\ell,2,2}(1)
  =\sum_{i\in\mathcal{S}}\mathring{\varphi}_{(i,\ell),0}A_{\ell, i}+\sum_{i\in\mathcal{S}^c:\,i\neq \ell}\mathring{\varphi}_{(i,\ell),0}A_{\ell, i}\,.
	\end{align*}
Following the same arguments in Section \ref{mu1:S} for obtaining the convergence rate of $\max_{\ell\in\mathcal{S}}|I_{\ell,1,2}(1)|$ there, it holds that
 \begin{align*}
\max_{\ell\in\mathcal{S}^c}|I_{\ell,2,2}(1)|=&~O_{\p}\big[\gamma^2s^{1/2}\{sp^{\min(2\omega_2-\omega_1,\,0)} \chi_{p}^{-1}+p^{1+\omega_2-\omega_1}\}\chi_{p}^{-2}\log^{1/2}p\big]\\
&+O_{\p}\big\{\gamma^2p^{1/2}(sp^{\omega_2-\omega_1} \chi_{p}^{-1}+p^{1-\omega_1})\chi_{p}^{-1}\log^{1/2}p\big\}\,.
 \end{align*}
For $I_{\ell,2,2}(2)$,   by the same arguments for deriving the convergence rate of $\max_{\ell\in\mathcal{S}}|I_{\ell,2,2}(2)|$ in Section \ref{mu2:S}, it holds that
 \begin{align*}
     \max_{\ell \in \mathcal{S}^c}|I_{\ell,2,2}(2)|=&~O_{\p}\big[\gamma s\{p^{2\omega_2-\omega_1}\chi_{p}^{-3}  I(\omega_1>2\omega_2)+I(\omega_1\leq2\omega_2)\}\log p\big]\\
     &+O_{\p}\big[\gamma s^{1/2} p^{1/2}\{p^{\omega_2-\omega_1}\chi_{p}^{-2}I(\omega_1>\omega_2)+ I(\omega_1=\omega_2)\}\log p\big]\\
     &+O_{\p}\big[\gamma p\{p^{-\omega_1}\chi_{p}^{-1}I(\omega_1>0)+ I(\omega_1=0)\}\log p\big]\,.
 \end{align*}
Therefore, by \eqref{eq:Iell22Sc},
 we have
	\begin{align*}
	    \max_{\ell \in \mathcal{S}^c}|I_{\ell,2,2}|=&~O_{\p}\big[\gamma^2s^{1/2}p^{-2}\{sp^{\min(2\omega_2-\omega_1,\,0)} \chi_{p}^{-1}+p^{1+\omega_2-\omega_1}\}\chi_{p}^{-2}\log^{1/2}p\big]\\
     &+O_{\p}\big\{\gamma^2p^{-3/2}(sp^{\omega_2-\omega_1} \chi_{p}^{-1}+p^{1-\omega_1})\chi_{p}^{-1}\log^{1/2}p\big\}\\
     &+O_{\p}\big[\gamma s p^{-2}\{p^{2\omega_2-\omega_1}\chi_{p}^{-3}  I(\omega_1>2\omega_2)+ I(\omega_1\leq2\omega_2)\}\log p\big]\\
     &+O_{\p}\big[\gamma s^{1/2} p^{-3/2}\{p^{\omega_2-\omega_1}\chi_{p}^{-2}I(\omega_1>\omega_2)+I(\omega_1=\omega_2)\}\log p\big]\\
     &+O_{\p}\big[\gamma p^{-1}\{p^{-\omega_1}\chi_{p}^{-1}I(\omega_1>0)+I(\omega_1=0)\}\log p\big]\,.
	\end{align*}
 Together with \eqref{eq:Il11}, it holds that
	\begin{align*}
	\max_{\ell\in\mathcal{S}^c}|\hat{\mu}_{\ell,2}-\mu_{\ell,2}|=&~O_{\p}\big[\gamma^2s^{1/2}p^{-2}\{sp^{\min(2\omega_2-\omega_1,\,0)} \chi_{p}^{-1}+p^{1+\omega_2-\omega_1}\}\chi_{p}^{-2}\log^{1/2}p\big]\\
 &+O_{\p}\big\{\gamma^2p^{-3/2}(sp^{\omega_2-\omega_1} \chi_{p}^{-1}+p^{1-\omega_1})\chi_{p}^{-1}\log^{1/2}p\big\}\\
 &+O_{\p}\big[\gamma sp^{-2}\{p^{2\omega_2-\omega_1}\chi_{p}^{-3}I(\omega_1>2\omega_2)+ I(\omega_1\leq2\omega_2)\}\log p\big]\\
 &+O_{\p}\big[\gamma s^{1/2} p^{-3/2}\{p^{\omega_2-\omega_1}\chi_{p}^{-2} I(\omega_1>\omega_2)+ I(\omega_1=\omega_2)\}\log p\big]\\
 &+O_{\p}\big[\gamma p^{-1}\{p^{-\omega_1}\chi_{p}^{-1} I(\omega_1>0)+I(\omega_1=0)\}\log p\big]\\
 &+O_{\p}(p^{-1}\log^{1/2}p)\,.
	\end{align*}
We obtain the convergence rate of $\max_{\ell\in\mathcal{S}^c}|\hat{\mu}_{\ell,2}-\mu_{\ell,2}|$. $\hfill\Box$

\subsection{Proof of Lemma {\rm\ref{la:bdvarphi}}.}\label{sec:pfla:bdvarphi}

	Notice that $\gamma=1-\alpha-\beta$ and
	\begin{align*}
	\mathbb{E}\{\varphi_{(i,j),0}\}=&\,\gamma\cdot\mathbb{P}(X_{i,j}=0)=\frac{\gamma}{1+\exp(\xi+\check{\theta}_i+\check{\theta}_j)}\,,\\
 \mathbb{E}\{\varphi_{(i,j),1}\}=&\,\gamma\cdot\mathbb{P}(X_{i,j}=1)=\frac{\gamma}{1+\exp(-\xi-\check{\theta}_i-\check{\theta}_j)}\,.
	\end{align*}
	Recall $\xi=-\omega_1\log p+\xi^+$ and $\check{\theta}_\ell=\omega_2\log p+\check{\theta}_\ell^+$ for all $\ell\in\mathcal{S}$, where $\omega_1\in[0,2)$ and $\omega_2\in[0,1)$ such that $0\leq \omega_1-\omega_2<1$, and $|\xi^+|\vee \max_{\ell\in\mathcal{S}}|\check{\theta}_\ell^+|=o(\log p)$. Write $\chi_{p}=\exp(-|\xi^+|\vee\max_{\ell\in\mathcal{S}}|\check{\theta}_\ell^+|)$.

	For $i,j\in\mathcal{S}$, it holds that
	$
	\xi+\check{\theta}_i+\check{\theta}_j=(2\omega_2-\omega_1)\log p+\xi^++\check{\theta}_i^++\check{\theta}_j^+
	$.
	If $\omega_1>2\omega_2$, due to $\exp(\xi+\check{\theta}_i+\check{\theta}_j)>0$ and $-\xi-\check{\theta}_i-\check{\theta}_j\asymp\log p$, then $ \mathbb{E}\{\varphi_{(i,j),0}\}\leq \gamma$ and $\mathbb{E}\{\varphi_{(i,j),1}\}\lesssim p^{2\omega_2-\omega_1}\chi_{p}^{-3}\gamma$. If $\omega_1=2\omega_2$, due to $\exp(\xi+\check{\theta}_i+\check{\theta}_j)>0$ and $\exp(-\xi-\check{\theta}_i-\check{\theta}_j)>0$, then $ \mathbb{E}\{\varphi_{(i,j),0}\}\leq \gamma$ and $\mathbb{E}\{\varphi_{(i,j),1}\}\leq\gamma$.
 If $\omega_1<2\omega_2$, due to $\xi+\check{\theta}_i+\check{\theta}_j\asymp\log p$ and $\exp(-\xi-\check{\theta}_i-\check{\theta}_j)>0$, then $ \mathbb{E}\{\varphi_{(i,j),0}\}\lesssim p^{-2\omega_2+\omega_1}\chi_{p}^{-3}\gamma$ and $\mathbb{E}\{\varphi_{(i,j),1}\}\leq\gamma$.

For $i\in\mathcal{S}$ and $j\in\mathcal{S}^c$, it holds that $
	\xi+\check{\theta}_i+\check{\theta}_j=(\omega_2-\omega_1)\log p+\xi^++\check{\theta}_i^+
	$. Due to $\omega_1\geq\omega_2$ and $\exp(\xi+\check{\theta}_i+\check{\theta}_j)>0$, then $ \mathbb{E}\{\varphi_{(i,j),0}\}\leq \gamma$. If $\omega_1=\omega_2$, due to $\exp(-\xi-\check{\theta}_i-\check{\theta}_j)>0$, then $\mathbb{E}\{\varphi_{(i,j),1}\}\leq\gamma$. If $\omega_1>\omega_2$, due to $-\xi-\check{\theta}_i-\check{\theta}_j\asymp \log p$, then $\mathbb{E}\{\varphi_{(i,j),1}\}\lesssim p^{\omega_2-\omega_1}\chi_{p}^{-2}\gamma$.

 For $i,j\in\mathcal{S}^c$, it holds that  $\xi+\check{\theta}_i+\check{\theta}_j=-\omega_1\log p+\xi^+$. Due to $\omega_1\geq0$ and $\exp(\xi+\check{\theta}_i+\check{\theta}_j)>0$, then $ \mathbb{E}\{\varphi_{(i,j),0}\}\leq \gamma$. If $\omega_1=0$, due to $\exp(-\xi-\check{\theta}_i-\check{\theta}_j)>0$, then $\mathbb{E}\{\varphi_{(i,j),1}\}\leq \gamma$. If $\omega_1>0$, due to $-\xi-\check{\theta}_i-\check{\theta}_j\asymp \log p$, then $\mathbb{E}\{\varphi_{(i,j),1}\}\leq p^{-\omega_1}\chi_{p}^{-1}\gamma$. $\hfill\Box$

\subsection{Proof of Lemma {\rm\ref{la:Ail}}.}\label{sec:pfla:Ail}
Recall
$
A_{i,\ell}=2\sum_{j:\,j\neq i,\ell}\mathbb{E}\{\varphi_{(\ell,j),1}\}\mathbb{E}\{\varphi_{(i,j),0}\}$.

\subsubsection{Upper bound of  $\max_{i,\ell\in\mathcal{S}:\,i\neq \ell}A_{i,\ell}$.}\label{scase1}

For any $i,\ell \in\mathcal{S}$ and $i\neq \ell$,
\begin{align*}
    A_{i,\ell}=
	2\sum_{j\in\mathcal{S}:\,j\neq i,\ell}\mathbb{E}\{\varphi_{(\ell,j),1}\}\mathbb{E}\{\varphi_{(i,j),0}\}+2\sum_{j\in\mathcal{S}^c}\mathbb{E}\{\varphi_{(\ell,j),1}\}\mathbb{E}\{\varphi_{(i,j),0}\}\,.
\end{align*}
By Lemma \ref{la:bdvarphi}, due to $s\ll p$ and $\chi_p\in(0,1]$, it holds that
\begin{align*}
    \sum_{j\in\mathcal{S}:\,j\neq i,\ell}\mathbb{E}\{\varphi_{(\ell,j),1}\}\mathbb{E}\{\varphi_{(i,j),0}\}\lesssim&~ sp^{\min(2\omega_2-\omega_1,\,\omega_1-2\omega_2)}\chi_{p}^{-3}\gamma^2\,,\\
    \sum_{j\in\mathcal{S}^c}\mathbb{E}\{\varphi_{(\ell,j),1}\}\mathbb{E}\{\varphi_{(i,j),0}\}\lesssim&~ p^{1+\omega_2-\omega_1}\chi_{p}^{-2}\gamma^2\,,
\end{align*}
which implies
$$\max_{i,\ell\in\mathcal{S}:\,i\neq \ell}A_{i,\ell}\lesssim \{sp^{\min(2\omega_2-\omega_1,\,\omega_1-2\omega_2)}\chi_{p}^{-1}+p^{1+\omega_2-\omega_1}\}\chi_{p}^{-2}\gamma^2\,.$$
We obtain the upper bound of $\max_{i,\ell\in\mathcal{S}:\,i\neq \ell}A_{i,\ell}$. $\hfill\Box$

\subsubsection{Upper bound of  $\max_{i\in\mathcal{S},\,\ell \in\mathcal{S}^c}A_{i,\ell}$.}\label{scase2}
For any $i\in\mathcal{S}$ and $\ell \in\mathcal{S}^c$,
\begin{align*}
    \,A_{i,\ell}=
	2\sum_{j\in\mathcal{S}:\,j\neq i}\mathbb{E}\{\varphi_{(\ell,j),1}\}\mathbb{E}\{\varphi_{(i,j),0}\}+2\sum_{j\in\mathcal{S}^c:\,j\neq \ell}\mathbb{E}\{\varphi_{(\ell,j),1}\}\mathbb{E}\{\varphi_{(i,j),0}\}\,.
\end{align*}
By Lemma \ref{la:bdvarphi}, due to $s\ll p$ and $\chi_p\in(0,1]$, it holds that
\begin{align*}
    \sum_{j\in\mathcal{S}:\,j\neq i}\mathbb{E}\{\varphi_{(\ell,j),1}\}\mathbb{E}\{\varphi_{(i,j),0}\}\lesssim&~ sp^{\min(-\omega_2,\,\omega_2-\omega_1)} \chi_{p}^{-5}\gamma^2\,,\\
   \sum_{j\in\mathcal{S}^c:\,j\neq \ell}\mathbb{E}\{\varphi_{(\ell,j),1}\}\mathbb{E}\{\varphi_{(i,j),0}\}\lesssim&~ p^{1-\omega_1}\chi_{p}^{-1}\gamma^2\,,
\end{align*}
which implies
$$\max_{i\in\mathcal{S},\,\ell \in\mathcal{S}^c}A_{i,\ell}\lesssim\{sp^{\min(-\omega_2,\,\omega_2-\omega_1)} \chi_{p}^{-4}+p^{1-\omega_1}\}\chi_{p}^{-1}\gamma^2\,.$$
We obtain the upper bound of $\max_{i\in\mathcal{S},\,\ell \in\mathcal{S}^c}A_{i,\ell}$. $\hfill\Box$

\subsubsection{Upper bound of  $\max_{i\in\mathcal{S}^c,\,\ell \in\mathcal{S}}A_{i,\ell}$.}\label{scase23}
For any $i\in\mathcal{S}^c$ and $\ell \in\mathcal{S}$,
\begin{align*}
    A_{i,\ell}=
	2\sum_{j\in\mathcal{S}:\,j\neq \ell}\mathbb{E}\{\varphi_{(\ell,j),1}\}\mathbb{E}\{\varphi_{(i,j),0}\}+2\sum_{j\in\mathcal{S}^c:\,j\neq i}\mathbb{E}\{\varphi_{(\ell,j),1}\}\mathbb{E}\{\varphi_{(i,j),0}\}\,.
\end{align*}
By Lemma \ref{la:bdvarphi}, due to $s\ll p$ and $\chi_p\in(0,1]$, it holds that
\begin{align*}
    \sum_{j\in\mathcal{S}:\,j\neq \ell}\mathbb{E}\{\varphi_{(\ell,j),1}\}\mathbb{E}\{\varphi_{(i,j),0}\}\lesssim&~ sp^{\min(2\omega_2-\omega_1,\,0)} \chi_{p}^{-3}\gamma^2\,,\\
   \sum_{j\in\mathcal{S}^c:\,j\neq i}\mathbb{E}\{\varphi_{(\ell,j),1}\}\mathbb{E}\{\varphi_{(i,j),0}\}\lesssim &~p^{1+\omega_2-\omega_1}\chi_{p}^{-2}\gamma^2\,,
\end{align*}
which implies
$$\max_{i\in\mathcal{S}^c,\,\ell \in\mathcal{S}}A_{i,\ell}\lesssim\{sp^{\min(2\omega_2-\omega_1,\,0)} \chi_{p}^{-1}+p^{1+\omega_2-\omega_1}\}\chi_{p}^{-2}\gamma^2\,.$$
We obtain the upper bound of $\max_{i\in\mathcal{S}^c,\,\ell \in\mathcal{S}}A_{i,\ell}$. $\hfill\Box$

\subsubsection{Upper bound of  $\max_{i, \ell\in\mathcal{S}^c:\, i\neq \ell}A_{i,\ell}$.}\label{scase4}
For any $i,\ell \in\mathcal{S}^c$ and $i\neq \ell$,
\begin{align*}
    A_{i,\ell}=
	2\sum_{j\in\mathcal{S}}\mathbb{E}\{\varphi_{(\ell,j),1}\}\mathbb{E}\{\varphi_{(i,j),0}\}+2\sum_{j\in\mathcal{S}^c:\,j\neq i,\ell}\mathbb{E}\{\varphi_{(\ell,j),1}\}\mathbb{E}\{\varphi_{(i,j),0}\}\,.
\end{align*}
By Lemma \ref{la:bdvarphi}, due to $s\ll p$ and $\chi_p\in(0,1]$, it holds that
\begin{align*}
    \sum_{j\in\mathcal{S}}\mathbb{E}\{\varphi_{(\ell,j),1}\}\mathbb{E}\{\varphi_{(i,j),0}\}\lesssim&~ sp^{\omega_2-\omega_1} \chi_{p}^{-2}\gamma^2\,,\\
   \sum_{j\in\mathcal{S}^c:\,j\neq i,\ell}\mathbb{E}\{\varphi_{(\ell,j),1}\}\mathbb{E}\{\varphi_{(i,j),0}\}\lesssim&~ p^{1-\omega_1}\chi_{p}^{-1}\gamma^2\,,
\end{align*}
which implies $$\max_{i, \ell\in\mathcal{S}^c:\, i\neq \ell}A_{i,\ell}\lesssim (sp^{\omega_2-\omega_1} \chi_{p}^{-1}+p^{1-\omega_1})\chi_{p}^{-1}\gamma^2\,.$$
We obtain the upper bound of $\max_{i, \ell\in\mathcal{S}^c:\, i\neq \ell}A_{i,\ell}$. $\hfill\Box$


\begin{thebibliography}{11}
		
		
		
		\bibitem[Blocki et~al.(2013)]{Blockietal_2013}
		Blocki, J., Blum, A., Datta, A., and Sheffet, O. (2013).
		Differentially private data analysis of social networks via
		restricted sensitivity. {\sl Proceedings of the 4th Conference on Innovations in Theoretical Computer Science}, 87--96.
			
		\bibitem[Benjamini et~al.(1995)]{Benjamini_1995}
		Benjamini, Y., and Hochberg, Y. (1995). Controlling the false discovery rate: a practical and powerful approach to multiple testing. {\sl Journal of the Royal Statistical Society Series B},  {\bf57}, 289--300.	
		
		\bibitem[Chang, Chen and Wu(2024)]{ChangChenWu_2021}
		Chang, J., Chen, X., and Wu, M. (2024). Central limit theorems for high dimensional dependent data. {\sl Bernoulli}, {\bf30}, 712--742.
		
		\bibitem[Chang et~al. (2024)]{CHKYF_suppl_2024}
        Chang, J., Hu, Q., Kolaczyk, E.D., Yao, Q., and Yi, F. (2024).  Supplement to ``Edge Differentially Private Estimation in the $\beta$-model via Jittering and Method of Moments''.  DOI: 10.1214/[provided by typesetter].

        \bibitem[Chang, Kolaczyk and Yao(2020)]{ChangKolaczykYao_2020}
		Chang, J., Kolaczyk, E. D., and Yao, Q. (2020).  Discussion of `Network cross-validation by edge sampling'. {\sl Biometrika}, {\bf107}, 277--280.
		
		\bibitem[Chang, Kolaczyk and Yao(2022)]{ChangKolaczykYao_2018}
		Chang, J., Kolaczyk, E. D., and Yao, Q. (2022).
		Estimation of subgraph densities in noisy networks.
		{\sl Journal of the American Statistical Association}, {\bf117}, 361--374.
		
		\bibitem[Chang et~al.(2018)]{CQYZ_2018}
		Chang, J., Qiu, Y., Yao, Q., and Zou, T. (2018). Confidence regions for entries of a large precision matrix. {\sl Journal of
			Econometrics}, {\bf206}, 57--82.
		
		\bibitem[Chang et~al.(2017)]{CZZZ_2017}
		Chang, J., Zheng, C., Zhou, W.-X., and Zhou, W. (2017). Simulation-based hypothesis testing of high dimensional means
		under covariance heterogeneity. {\sl Biometrics}, {\bf73}, 1300--1310.
		
		\bibitem[Chatterjee, Diaconis and Sly(2011)]{ChatterjeeDiaconisSly_2011}
		Chatterjee, S., Diaconis, P., and Sly, A. (2011).
		Random graphs with a given degree sequence.
		{\sl The Annals of Applied Probability}, {\bf 21}, 1400--1435.
		
		\bibitem[Chen, Kato and Leng(2021)]{ChenKatoLeng_2021}
		Chen, M., Kato, K., and Leng, C. (2021). Analysis of networks via the sparse $\beta$-model. {\sl Journal of the Royal Statistical Society Series B}, {\bf 83}, 887--910.
		
		\bibitem[Chernozhukov, Chetverikov and Kato(2013)]{CCK_2013}
		Chernozhukov, V., Chetverikov, D., and Kato, K. (2013). Gaussian approximations and multiplier bootstrap for maxima of sums of high-dimensional random vectors. {\sl The Annals of Statistics}, {\bf 41}, 2786--2819.
		
		\harvarditem[Chernozhukov, Chetverikov and Kato]{Chernozhukov, Chetverikov  and   Kato}{2017}{CCK2017}
Chernozhukov, V., Chetverikov, D.,  and  Kato, K.  (2017). Central limit theorems and bootstrap in high  dimensions. {\sl The Annals of Probability}, {\bf 45}, 2309--2352.
		
		
		
		
		\bibitem[de la Pe\~{n}a and Montgomery-Smith(1995)]{DM_1995}
		de la Pe\~{n}a, V., and Montgomery-Smith, S. (1995). Decoupling inequalities for the tail probabilities of multivariate $U$-statistics. {\sl The Annals of Probability}, {\bf23}, 806--816.
		
		\bibitem[Dwork(2006)]{Dwork_2006}
		Dwork, C. (2006). Dierential privacy. {\sl 33rd International Colloquium on Automata, languages and programming}, 1--12.
		
		\bibitem[Dwork et~al.(2006)]{Dworketal_2006}
		Dwork, C., McSherry, F., Nissim, K., and Smith, A. (2006). Calibrating noise to sensitivity in private data analysis.
		{\sl Proceedings of the 3rd Theory of Cryptography Conference}, 265--284.
		
		\bibitem[Fan, Zhang and Yan(2020)]{FanZhangYan_2020}
		Fan, Y., Zhang, H., and Yan, T. (2020). Asymptotic theory for differentially private generalized $\beta $-models with parameters increasing. {\sl Statistics and Its Interface},
		{\bf 13}, 385--398.

		
		\bibitem[Gin\'{e}, Lata{\l}a and Zinn(2000)]{Gineetal_2000}
		
		Gin\'{e}, E., Lata{\l}a, R., and Zinn, J. (2000). Exponential and moment inequalities for $U$-statistics. In {\sl High Dimensional Probability} II (E. Gin\'{e}, D. M. Mason and J. A. Wellner, eds.) 13--38. Birkh\"{a}user, Boston. 
		
		
		\bibitem[Hennig(2007)]{Hennig_2007}
		Hennig, C. (2007). Cluster-wise assessment of cluster stability. {\sl Computational Statistics \& Data Analysis},
		{\bf52}, 258--271.
		
		
		\bibitem[Jiang et~al.(2020)]{JiangEtal_2020}
		Jiang, H., Pei, J., Yu, D., Yu, J., Gong, B., and Cheng, X. (2020). Applications of differential privacy in social network analysis: A survey. {\sl IEEE Transactions on Knowledge and Data Engineering}, {\bf 35}, 108--127.
		
		\bibitem[Karwa, Krivitsky and Slavkovi\'c(2017)]{KarwaKrivitskySlavkovic_2017}
		Karwa, V., Krivitsky, P. N., and Slavkovi\'c, A. (2017).
		Sharing social network data: differentially private estimation of exponential family random
		graph models.
		{\sl Journal of the Royal Statistical Society Series C}, {\bf 66}, 481--500.
		
		\bibitem[Karwa and Slavkovi\'c(2016)]{KarwaSlavkovic_2016}
		Karwa, V., and Slavkovi\'c, A. (2016). Inference using noisy degrees: differentially
		private $\beta$-model and synthetic graphs. {\sl The Annals of Statistics}, {\bf 44}, 87--112.
		
		
		
		\bibitem[Mukherjee, Mukherjee and Sen(2018)]{Mukherjeeetal_2018}
		Mukherjee, R., Mukherjee, S., and Sen, S. (2018). Detection thresholds for the $\beta$-model on sparse graphs. {\sl The Annals of Statistics}, {
  \bf 46}, 1288--1317.
		
		\bibitem[Nissim, Raskhodnikova and Smith (2007)]{NissimRaskhodnikovaSmith_2007}
		Nissim, K., Raskhodnikova, S., and Smith, A. (2007). Smooth sensitivity and sampling in private data analysis. {\sl In Proceedings of the Thirty-ninth Annual ACM symposium on Theory of Computing}, 75--84.
		
		
		\bibitem[Rinaldo, Petrovi{\'c} and Fienberg(2013)]{Rinaldoetal_2013}
		Rinaldo, A., Petrovi{\'c}, S., and Fienberg, S. E. (2013). Maximum likelihood estimation in the $\beta$-model. {\sl The Annals of Statistics}, {\bf 41}, 1085--1110.
		
		\bibitem[Stein and Leng(2021)]{SteinLeng_2021}
	    Stein, S., and Leng, C. (2020). A sparse $\beta $-model with covariates for networks. arXiv preprint arXiv:2010.13604.
	
		
		\bibitem[Wasserman and Zhou(2010)]{WassermanZhou_2010}
		Wasserman, L., and Zhou, S. (2010). A statistical framework for differential privacy.
		{\sl Journal of the American Statistical Association}, {\bf 105}, 375--389.
		
		\bibitem[Yan and Xu(2013)]{yx13}
		Yan, T., and Xu, J. (2013). A central limit theorem in the $\beta$-model
		for undirected random
		graphs with a diverging number of vertices.
		{\sl Biometrika}, {\bf 100}, 519--524.
		
		\bibitem[Zhang et al.(2021)]{Zhangetal_2021}
		Zhang, Y., Wang, Q., Zhang, Y., Yan, T., and Luo, J. (2021). L-2 Regularized maximum likelihood for $\beta $-model in large and sparse networks. arXiv preprint arXiv:2110.11856.
		
	\end{thebibliography}
\end{document}